\documentclass[11pt]{article}
\pdfoutput=1

\usepackage{mathrsfs}
\usepackage{amsmath,amssymb}
\usepackage{bm}
\usepackage{natbib}
\usepackage[usenames]{color}
\usepackage{amsthm}

\usepackage{multirow} 
\usepackage{enumitem}

\usepackage{enumitem}
\usepackage{dsfont}
\usepackage{mathtools}

\usepackage[colorlinks,
linkcolor=red,
anchorcolor=blue,
citecolor=blue
]{hyperref}

\usepackage{mylatexstyle}

\usepackage{setspace}
\usepackage[left=1in, right=1in, top=1in, bottom=1in]{geometry}

\setlist[itemize]{itemsep=0pt, topsep=2pt}
\setlist[enumerate]{itemsep=-2pt, topsep=2pt}

\usepackage{xcolor}

\ifdefined\final
\usepackage[disable]{todonotes}
\else
\usepackage[textsize=tiny]{todonotes}
\fi
\newcommand{\cy}[1]{{ #1}}

\newcommand{\xr}[1]{{ #1}}

\allowdisplaybreaks

\title{Multiple Descent in the Multiple Random Feature Model}

\author
{
Xuran Meng\thanks{Department of Statistics and Actuarial Science, The University of Hong Kong; e-mail: {\tt u3007800@connect.hku.hk}}
	~~~and~~~
Jianfeng Yao\thanks{School of Data Science,
	The Chinese University of Hong Kong (Shenzhen); e-mail: {\tt jeffyao@cuhk.edu.cn}}
~~~and~~~
	Yuan Cao\thanks{Department of Statistics and Actuarial Science, The University of Hong Kong;
  e-mail: {\tt yuancao@hku.hk}}
}

\date{}

\begin{document}

\maketitle

\begin{abstract} 
Recent works have demonstrated a \textit{double descent} phenomenon in over-parameterized learning.
Although this phenomenon has been investigated by recent works, it has not been fully understood in theory. In this paper, \cy{we investigate the multiple descent phenomenon in a class of multi-component prediction models. 
We first consider a ``double random feature model'' (DRFM) concatenating two types of random features}, and study the excess risk achieved by the DRFM in ridge regression. We calculate the precise limit of the excess risk under the high dimensional framework where the training sample size, the dimension of data, and the dimension of random features tend to infinity proportionally. Based on the calculation, we further theoretically demonstrate that the risk curves of DRFMs can exhibit \textit{triple descent}. We then provide a thorough experimental study to verify our theory. At last, we extend our study to the ``multiple random feature model'' (MRFM), and show that MRFMs ensembling $K$ types of random features may exhibit \textit{$(K+1)$-fold descent}. Our analysis points out that risk curves with a specific number of descent generally exist in \cy{learning multi-component prediction models}. 
\end{abstract}

\section{Introduction}
Modern machine learning models such as deep neural networks are usually highly over-parameterized so that they can be trained to exactly fit the training data. Such over-parameterized models have gained immense popularity and achieved state-of-the-art performance in various learning tasks. However, in classical statistical learning theory, over-parameterized models are believed to have high excess risks due to overfitting, and hence their success has not been fully explained in theory. This gap between theory and practice has motivated a number of recent works to study the success of over-parameterized models.





\begin{figure}[t!]
    \centering
    \includegraphics[width=6in]{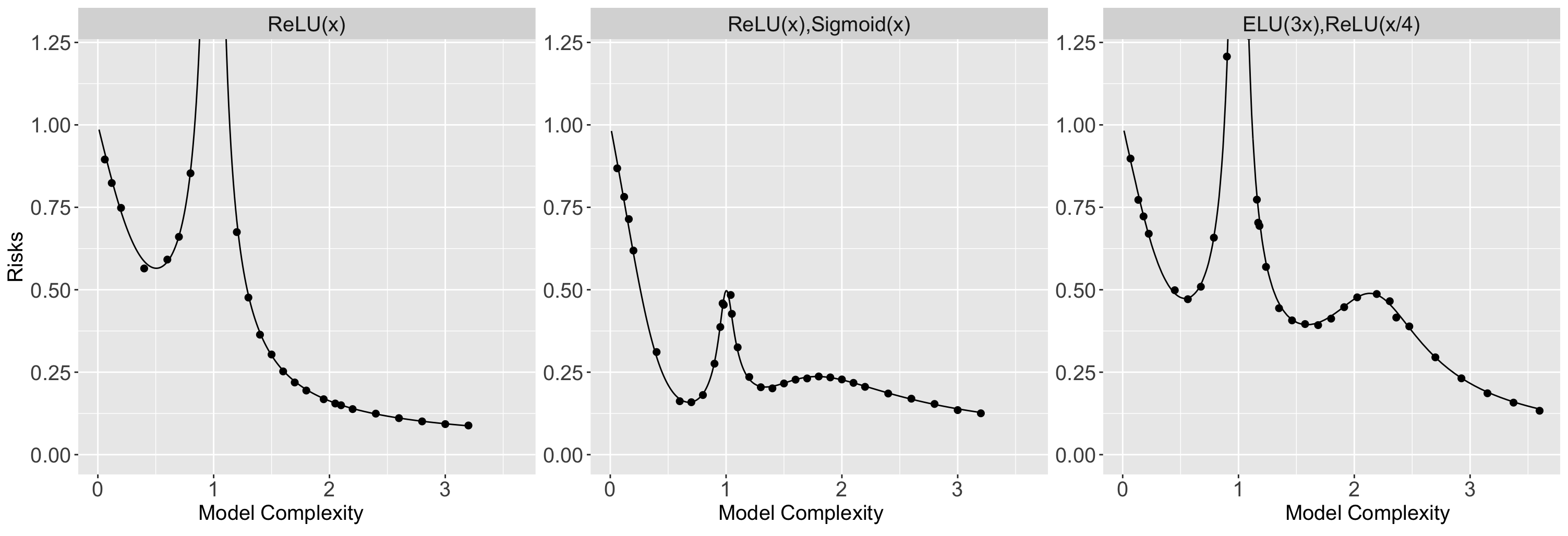}\vspace{-5pt}\\
    \hspace{19pt}(a)\hspace{120pt}(b)\hspace{130pt}(c)
    \vspace{-10pt}
    \caption{Examples of double and triple descent. (a) gives the excess risk of a random feature model with ReLU activation function; (b) shows the excess risk of a double random feature model with ReLU and sigmoid activation functions; (c) shows the excess risk of a double random feature model with ELU and ReLU activation functions. The $x$-axis is the model complexity (number of parameters/sample size) and the $y$-axis is the excess risk. The curve gives our theoretical predictions, and the dots are our numerical results.}
    \label{fig:tripleinIntro}
    \vspace{-10pt}
\end{figure}

\cy{Recent works have pointed out a \textit{double/multiple descent} phenomenon in over-parameterized learning: as the number of parameters in a model increases, the excess risk may increase and decrease multiple times (see Figure~\ref{fig:tripleinIntro} for some examples)}. The double descent phenomenon was first demonstrated experimentally by \citet{belkin2019reconciling} in random feature models, random forests and neural networks, and then studied theoretically by a series of works under different settings. 
Specifically, 
\citet{belkin2020two} 
theoretically demonstrated the double descent shape of the risk curve of the minimum norm predictor in learning \xr{linear  models and Fourier series models}.
\citet{wu2020optimal,mel2021theory,hastie2022surprises} studied the excess risk in linear regression under the setting where the dimension and sample size go to infinity preserving a fixed ratio, and showed that the risk decreases with respect to this ratio in the over-parameterized setting. \citet{mei2022generalization,liao2020random} further studied double descent in random feature models when the sample size, data dimension and the number of random features have fixed ratios \xr{and \citet{adlam2022random} extended the model by adding bias terms.}
\xr{\citet{deng2022model} studied double descent under logistic model. \citet{Enami2020Generalization} studied the asymptotic generalization error of generalized linear models.}
\cy{Several recent works have also studied other learning settings under which the risk curves exhibit triple descent or multiple descent.} 
Specifically, \citet{liang2020multiple} gave an upper bound on the risk of the minimum-norm interpolants in a reproducing kernel Hilbert space and showed that it has a multiple descent shape with infinitely many peaks. \citet{chen2021multiple} showed that with different and well-designed data distributions in linear regression, 
the risk curve can have an arbitrary number of peaks at arbitrary locations as the data dimension increases. \citet{mel2021theory,li2021minimum} showed that the risk curve of linear regression can exhibit multiple descent when learning anisotropic data. 
\citet{pennington2020Triple} demonstrated triple descent for a specific random feature model associated with an over-parameterized two-layer neural network in the so-called \cy{``neural tangent kernel'' \citep{jacot2018neural} regime.} \cy{\citep{misiakiewicz2022spectrum,xiao2022precise} showed that the risk curve of certain kernel predictors can exhibit multiple descent concerning the sample size and data dimension. }

While recent works have provided valuable insights, the double, triple and multiple descent phenomena have not been fully understood in theory. \cy{Specifically, we note that various modern learning methods utilize multi-component predictors of a general form
\begin{align}\label{eq:multi-component}
    f(\xb) = f_1(\xb) + f_2(\xb) + \cdots + f_K(\xb),
\end{align}
where $f_1(\xb),\ldots, f_K(\xb)$ are individual prediction models. Such a multi-component formulation covers different learning methods. For example, ensemble methods \citep{hansen1990neural,dietterich2000ensemble,krogh1994neural} can naturally be formulated as \eqref{eq:multi-component}; two-layer neural networks utilizing feature concatenation is also a summation of multiple components defined by different features; two-layer ResNet \citep{he2016deep} models can be formulated as \eqref{eq:multi-component} by treating the feedforward part and the skip-connection part of the model as two components; a class of semi-parametric methods consider a parametric component and a non-parametric component in the model \citep{zhao2016partially,chernozhukov2018double}; \xr{Neural network models with the exact form of \eqref{eq:multi-component} can also be applied to solve partial differential equations \citep{Liu2020multi}.}}

\cy{ 
In this work, we aim to study the double/multiple descent phenomenon in learning multi-component predictors. We theoretically demonstrate that
}
\begin{center}
\emph{
\cy{There exists a learning problem, such that for any $K\in \mathbb{N}_+$, there exists a $K$-component prediction model whose risk curve exhibits $(K+1)$-fold descent.}
}    
\end{center}
\cy{The learning problem mentioned in the claim above is the same learning problem where recent works have demonstrated double descent for random feature models \citet{mei2022generalization}, and is also essentially the same learning problem (with slight modification) studied in \citet{hastie2022surprises} analyzing double descent in linear regression. Therefore, demonstrating this claim provides new insights into how complicated prediction model structures can affect the risk curve.}




\cy{
This paper aims to study the double/multiple descent phenomena in learning multi-component predictors of the form \eqref{eq:multi-component} through the lens of random feature models. Specifically, }we introduce double and multiple random feature models (DRFMs and MRFMs), \xr{which ensemble} two or more types of random features defined by different nonlinear activation functions. Under the setting where the training sample size, the dimension of data, and the dimension of random features tend to infinity proportionally, we establish an asymptotic limit of the excess risk achieved by DRFMs and MRFMs\xr{, and} demonstrate that the risk curve of a DRFM can exhibit triple descent: an example for the DRFM with \xr{different} activation functions is given in Figure~\ref{fig:tripleinIntro}. More generally, we also show that the risk curve of an MRFM with $K$ types of random features can exhibit $(K+1)$-fold descent. 

We summarize the contributions of this paper as follows.
\begin{enumerate}[leftmargin = *]
    \item \cy{Our first contribution is to demonstrate the existence of multiple descent in learning certain multi-component predictors. Specifically, we  demonstrate that DRFMs  may exhibit triple descent, and then extend the analysis to MRFMs and show that MRFMs consisting of $K$ types of random features may have a risk curve with $(K+1)$-fold descent. To the best of our knowledge, such multiple descent risk curves with a specific number of peaks have not been well understood in random feature models or other multi-component learning models, and therefore we believe that DRFMs and MRFMs can serve as important examples in the literature of multiple descent.}
    \item \cy{We provide a natural and intuitive explanation of multiple descent in DRFMs and MRFMs. For example, for DRFMs, we point out that the existence of triple descent risk curves is predictable by considering the two extreme cases: (i) the DRFM uses two random features of the same type and scale, and (ii) one type of random feature in the DRFM has a very small scale and is thus negligible. This scale difference refers to a large gap in magnitude between the two random features, such as the activation pair $(\sigma_1(x),c_0\sigma_2(x))$ where the constant $c_0$ is small. We point out that these two cases both lead to double descent but with different peak locations. Therefore, for DRFMs where the scale difference between the two parts of random features is neither too big nor too small, we can expect triple descent to appear. Following this intuition, we successfully anticipate multiple descent in various simulations, and correctly predict the number of peaks in the risk curves and the locations of all peaks in the risk curves.}
    

    \item \cy{We also establish comprehensive theoretical results to back up our intuitive explanation. We calculate the precise limit of the excess risk achieved by DRFMs and MRFMs. This is an extension of the study of \citet{mei2022generalization} which analyzed the vanilla random feature model with a single activation function.
    We also establish a novel type of theoretical proof of multiple descent which is based on the comparison between excess risk values at different over-parameterization levels. }

\end{enumerate}

\xr{Our calculation of the theoretical limit of the excess risks of DRFMs and MRFMs follow the blueprint of \citet{mei2022generalization} and expand upon it by constructing new linear pencil matrices and giving new calculation for the related Stieltjes transforms. In essence, \citet{mei2022generalization} accomplished the following:
\begin{enumerate}
\item introduced a risk function decomposition and proved convergence in $L_1$ distance; \label{item:decompose}
\item\label{item:pencil} used a linear pencil matrix and its partial derivatives of logarithmic potential to express the decomposed terms;
\item\label{item:asymptotic} provided an asymptotic approximation of the logarithmic potential and proved that partial derivatives are also approximated in $L_1$ distance;
\item\label{item:calculation} calculated theoretical values from the asymptotic approximation.
\end{enumerate}
Our theoretical analysis of the excess risk  follow the decomposition method in \ref{item:decompose}, but due to the increased complexity of our model, we develop several new technical lemmas to overcome   this higher complexity; such examples include  Proposition~\ref{prop:decomp-expectation}  and Lemma~\ref{lemma:decomp2} in the Appendix.} \xr{Moreover, the increased complexity in the main terms of the decomposition necessitates a more complex linear pencil matrix, as defined in Definition~\ref{def:linear pencil} and \ref{def:linear pencil-m}. Although the construction of the linear pencil matrices is inspired by Item \ref{item:pencil} above, the higher complexity of our model results in a more intricate  construction and  a more complex calculation of the  related Stieltjes transforms and their logarithmic potentials than in \ref{item:asymptotic}. Specifically, Proposition~\ref{prop:implicit1} provides the calculation of the Stieltjes transforms in DRFM and serves as the inspiration for the calculation of the Stieltjes transforms in MRFM. In MRFM, we utilize mathematical induction to complete this calculation.}

\xr{Besides the calculation of the theoretical limits of excess risks, this paper also presents a novel theory in the demonstration of multiple descent (given in Propositions~\ref{prop:property_lbd=0} and \ref{prop:property_psi}). Instead of directly investigating the theoretical limits, our approach focuses on taking limits within specific parameter ranges to observe the resulting behavior. Specifically, we employ the following steps:
\begin{enumerate}
\item We give a fixed ratio between the number of training parameters and the sample size.
\item Within this ratio, we set the regularization parameter $\lambda$ to approach zero, which allows us to approximate the implicit $\bnu$-system introduced later. We then replace the approximate solution with the theoretical limits.
\item To assess the impact of scale differences, we let one of the activation function scales tend towards zero, and examine the resulting theoretical limits.
\end{enumerate}
By employing this method, we successfully utilize the $\varepsilon-\delta$ language to accurately depict the presence of two peaks and determine their precise locations.
}


The remaining of the paper is organized as follows. We first give some additional references and  notations below. Section~\ref{sec:Problemsetting} introduces the problem settings. 
Section~\ref{sec:mainresult2} establishes the theoretical limits of the excess risks of double random feature models. Section~\ref{sec:simulation} gives theoretical analyses and simulations to demonstrate triple descent in some DRFMs. Section~\ref{sec:generalcase} extends the results to multiple random feature models and gives numerical simulations to demonstrate multiple descent.
Finally, Section~\ref{sec:Conclusion} concludes the paper and discusses some related questions for future investigation. Proofs of the main results and some additional experiments are presented in the appendix. 


\subsection{Additional related works} 
Besides the works we previously discussed, a series of recent works have also studied the double and triple descent phenomena.  \citet{montanari2020interpolation} considered a two-layer neural network in the neural tangent regime, showed an interpolation phase transition, and gave a characterization of the generalization error which decreases with the number of training parameters. \citet{adlam2020understanding} developed a novel bias-variance decomposition, and utilized the decomposition to show double descent in random feature regression. \citet{d2020double} developed a quantitative theory for the double descent phenomenon in the lazy learning regime of two-layer neural networks, and showed that overfitting is beneficial when the noise level in the data is low. \citet{geiger2020scaling} utilized the intuition of double descent to show that the smallest generalization error can sometimes be achieved by the ensemble of several neural networks of intermediate sizes. 
\citet{nakkiran2020optimal,patil2022mitigating} studied how an appropriately chosen parameters or suitable cross validation procedure can mitigate multiple descent in the prediction models. \citet{d2020triple} investigated the parameter-wise double descent and sample-wise triple descent phenomena in random feature regression. \citet{deng2021model} showed double descent phenomenon in logistic regression.


Our paper is also closely related to the recent studies of the ``benign overfitting'' phenomenon. \citet{tsigler2020benign} showed that for certain regression problems, the risk achieved by the minimum norm linear interpolator can be asymptotically optimal. \citet{bartlett2020benign} further extended the results in \citet{tsigler2020benign} to the setting of linear ridge regression. \citet{chatterji2021finite} studied the risk of the maximum margin linear classifier in learning sub-Gaussian mixtures with additional label-flipping noises. \citet{cao2021risk} established matching upper and lower bounds of the risk achieved by the maximum margin linear classifier. \citet{frei2022benign} showed that fully-connected two-layer networks trained to achieve a zero training error can still achieve an asymptotically optimal test
error. \citet{cao2022benign} studied signal learning and noise memorization during the training of a two-layer convolutional neural network and revealed a phase transition between benign and harmful overfitting. Note that most studies along this line of research focus on the setting where the number of parameters $N$ is much larger than the sample size $n$ (e.g., $N = \Omega(n^2)$). In comparison, our work considers the setting where $N$ and $n$ go to infinity in comparable magnitudes, and studies how the excess risk changes with respect to their ratio.

\subsection{Notations}
We use lower case letters to denote scalars, and use 
bold face letters to denote vectors and matrices.
For functions $f,g$ and a probability measure $\nu$, we denote $\la f,g\ra_{\nu}=\int f(\xb)g(\xb)\nu(d\xb)$. 
 The $\ell_2$-norm of a vector $\vb$ is $\lVert\vb \rVert_2$.  
 For a matrix $\Ab$, we use $\lVert\Ab \rVert_{\star}$, $\lVert\Ab\rVert_{\max}$, $\lVert\Ab\rVert_{\op}$ and $\lVert\Ab\rVert_{F}$ to denote its nuclear norm, maximum norm, operator norm, and Frobinuous norm, respectively, and use $\tr(\Ab)$ to denote its trace. 
 A sub-matrix of $\Ab$ with row indices in $I$ and column indices in $J$ is denoted by $\Ab_{I,J}$, and $\tr_{I}(\Ab)=\tr(\Ab_{I,I})$ is the trace of the square sub-matrix with indices in $I$. 
\if TF
 $\tr_{\cA}(\Ab)=\sum_{i\in\cA}\Ab(i,i)$. $\vb^\T $ and $\Ab^\T $ denote the transpose. $\Ab_{\cA_1,\cA_2}$ denote its sub-matrix with row in $\cA_1$ and column in $\cA_2$. For example, $\Ab=\begin{bmatrix}
1&2&3\\4&5&6\\7&8&9
\end{bmatrix}$, $\Ab_{[1:2],[2:3]}=\begin{bmatrix}
2&3\\5&6
\end{bmatrix}$ is the sub-matrix with row 1,2 and column 2, 3.
\fi 

The sets of natural, real and complex numbers are denoted by $\NN$,  $\RR$ and $\mathbb{C}$, respectively. 
For $z\in \bbC$, we use $\Re(z)$ and $\Im(z)$ to denote its real and imaginary part. $\bbC_+=\{z\in\bbC:\Im(z)>0\}$ denotes the upper half complex plane with positive imaginary part. Let $\rmi=\sqrt{-1}$ be the imaginary unit. 
The unit sphere of $\RR^d$ is denoted by  $\SSS^{d-1}=\{\xb\in\RR^d:\lVert\xb\rVert_2=1\}$ and $c\cdot\SSS^{d-1}$ denotes the sphere with radius  $c>0$. The set of integers from  $n_1$ to $n_2$ is denoted by $[n_1:n_2]=\{n_1,\ldots,n_2\}$ and 
$[n]=[1:n] =\{1,\ldots,n\}$.  Moreover, $\1_q\in\RR^q$ denotes $q$-dimensional all-one vectors.

We use the standard asymptotic notations $\Theta_d(\cdot)$, $O_d(\cdot)$, $o_d(\cdot)$ and $\Omega_d(\cdot)$, where the subscript $d$ emphasizes the asymptotic variable. We write $X_1(d)=O_{\PP}(X_2(d))$ if for any
$\varepsilon>0$, there exists $C>0$ such that $\PP(|X_1(d)/X_2(d)|>C)\leq \varepsilon$ for all $d$. Similarly, we denote $X_1(d)=o_{\PP}(X_2(d))$ if $\{X_1(d)/X_2(d)\}_d$ converges to $0$ in probability.



\section{The double random feature model}
\label{sec:Problemsetting}



We consider regression problems where, for a data pair $(\xb, y)$, the goal is to predict the scalar response $y$ using the input vector $\xb\in\RR^d$. We analyze the prediction performance of a {\em double random feature model}, or DRFM, constructed as follows.
The random features are based on two nonlinear activation functions $\sigma_1,\sigma_2$ and $N$ random feature parameter vectors 
$\btheta_i\sim \Unif(\sqrt{d}\cdot\SSS^{d-1})$, $i\in[N]$. We let $a_i\in\RR$, $i\in[N]$ be the linear combination coefficients of the random features, and denote $\bTheta = [\btheta_1,\ldots,\btheta_N]^\top \in \RR^{N\times d}$, $\ab=[a_1,\ldots, a_N]^\top \in \RR^N$. 
Then a DRFM predictor has the form
\begin{equation}\label{eq:DRFM}
   \widehat y = f(\xb;\ab,\bTheta)=\sum_{i=1}^{N_1}a_i\sigma_1\big(\la\btheta_i,\xb\ra/\sqrt{d}\big)+\sum_{i=N_1+1}^{N}a_i\sigma_2\big(\la\btheta_i,\xb\ra/\sqrt{d}\big).
\end{equation}
In \eqref{eq:DRFM}, the first $N_1$ units use the activation function $\sigma_1$ and the first part of the random feature parameters $\bTheta_1 = [\btheta_1,\ldots,\btheta_{N_1}]^\T$,
while the remaining $N_2=N-N_1$ units use the second activation function $\sigma_2$ and the second part of the random feature parameters $\bTheta_2 = [\btheta_{N_1 + 1},\ldots,\btheta_{N}]^\T$.
Note that the coefficients $a_1,\ldots,a_N$ are the trainable parameters, while $\btheta_1,\ldots,\btheta_N$ are randomly generated parameters to define the random features.

Note that in our definition of $f(\xb;\ab,\bTheta)$,  we have introduced the factor $1/ \sqrt{d}$ inside the activation functions $\sigma_j(\cdot)$.
This normalization facilitates our analysis using random matrix theory. Note also that the random feature parameters $\btheta_i$ are imposed to both have a fixed length $\sqrt d$, but the setting covers a more general situation where the parameters can have different lengths, say $c_1\sqrt d$ and $c_2\sqrt d$, respectively. Indeed, if $\|\btheta_i \|_2=c_j\sqrt{d} $, we can introduce $\tilde{\sigma}_{j}(z) = \sigma_{j}(c_j z)$ so that $ \sigma_{j}\big(\la \btheta_i,\xb\ra/\sqrt{d}\big)=\tilde{\sigma}_{j}\big(\la{\bm{\tau}}_i,\xb\ra/\sqrt{d}\big)$ where $\btau_j =\btheta_j/c_j $ has length $\sqrt d$.

To go further, we specify the data we aim to learn with double random feature models. We assume the data are generated from a distribution defined as follows. 
\begin{definition}[Data generation model]\label{def:data_model}
The distribution of the data pair $(\xb, y)$ is given as follows: 
\begin{enumerate}[leftmargin = *]
    \item The input vector  $\xb$ follows the uniform distribution on the sphere $\sqrt{d}\cdot\SSS^{d-1}$ of raidus   $ \sqrt{d}$.
    \item The output is $y =  \la \bbeta_{1,d}, \xb \ra+F_0+\varepsilon$, 
    where  
    $\bbeta_{1,d}\in\RR^d$,  $F_0\in\RR$, and $\varepsilon$ is 
    a  noise independent of $\xb$. We assume that  $\EE(\varepsilon)=0$, $\EE(\varepsilon^2)=\tau^2$, and $\EE(\varepsilon^4)<+\infty$.
\end{enumerate}
The parameters of the data generation  model are  $\bbeta_d=[F_0,\bbeta_{1,d}^\T ]^\T $ and we hereafter denote by $\cD(\bbeta_d)$  the probability distribution of the pair  $(\xb, y)$. 
\end{definition}
This data generation model is standard in recent literature on double descent. Similar settings have been studied in a number of recent works \citep{hamsici2007spherical,marinucci2011random,di2014nonparametric,mei2022generalization}.

Given a training data set $S = \{(\xb_i,y_i)\}_{i=1}^n$ consisting of $n$ independent samples from the data generation model in Definition~\ref{def:data_model}, we denote the data matrix by $\Xb=[\xb_1,...,\xb_n]^\T \in\RR^{n\times d}$, the label vector by $\yb=[y_1,...,y_n]^\T $ and the noise vector by $\bvarepsilon=[\varepsilon_1,\ldots,\varepsilon_n]^\T$.  
Then we fit a DRFM predictor $f(\cdot;\ab,\bTheta)$ based on the training data set $S$ via the principle of ridge regression. Specifically, we learn the coefficient vector $\ab$ by minimizing the $\ell_2$-regularized square loss:
\begin{equation}\label{eq:def_hata}
\hat{\ab}=\argmin_{\ab}\left\{\frac1n\sum_{i=1}^n\Big(y_i-f(\xb_i;\ab,\bTheta)\Big)^2+\frac{d}{n}\lambda\lVert\ab\rVert_2^2 \right\}, 
\end{equation}
where $\lambda>0$ is the regularization parameter. We here use the factor $d / n$ in the regularization term to simplify our analysis. Removing the factor   does not affect the results in this paper, because we consider the setting where $d / n$ has a positive limit. This fact will be formally clarified in Section~\ref{sec:mainresult2}.

The excess risk of the predictor $f(\cdot; \hat\ab, \bTheta)$ can be written as
\begin{align}
    R_d(\Xb,\bTheta,\lambda,\bbeta_d,\bvarepsilon) 
    = \EE_{\xb\sim\Unif(\sqrt{d}\cdot\SSS^{d-1}) } \big[ F_0+\xb^\T \bbeta_{1,d} - f(\xb; \hat\ab, \bTheta)\big]^2.\label{eq:def_Rd}
\end{align}
This notation of the excess risk specifically highlights the dependency of the risk on $\Xb,\bTheta,\lambda,\bbeta_d,\bvarepsilon$. Note that we do not take average over the randomness of the training data $\Xb$, the noise vector $\bvarepsilon$ or the random features $\bTheta$, but aim to show the convergence of the risk towards a fixed value as $d,N,n\rightarrow \infty$ in an appropriate manner.

\section{\xr{Excess risks of} double random feature models}
\label{sec:mainresult2}
In this section we present our main results on the excess risks of DRFMs. We first give a definition. 
\begin{definition}
\label{def:someconstant}
The spherical moments of the activation functions $\sigma_j$ ($j=1,2$) are 
 \begin{equation*}
        \mu_{j,0}\triangleq \EE\{\sigma_j(G)\},\quad \mu_{j,1}\triangleq\EE\{G\sigma_j(G)\},\quad
        \mu_{j,2}^2\triangleq{\EE\{\sigma_j(G)^2\}-\mu_{j,0}^2-\mu_{j,1}^2}, 
 \end{equation*}
 where $ G\sim\rmN(0,1)$ is standard normal. We collect the six constants $\mu_{j,0},\mu_{j,1},\mu_{j,2}^2$, $j=1,2$ in a vector $\bmu$. 
\end{definition}

\xr{In Definition~\ref{def:someconstant}, the first index $j$ points out the corresponding activation function, and the second index $k$ links to the specific spherical moment.} We now introduce the main assumptions in this paper. 
\begin{assumption}
\label{assump1} The nonlinear activation functions  $\sigma_{j}:\RR\rightarrow\RR$ $(j=1,2)$ are weakly differentiable, with weak derivative $\sigma_{j}'$. Moreover, for some constants $0< C_0,C_1<+\infty$, $|\sigma_j(u)|\vee |\sigma_j'(u)|\leq C_0e^{C_1|u|}$,  $u\in\RR$. 
\end{assumption}
It is easy to see that commonly used activation functions such as ReLU, sigmoid, and hyperbolic tangent functions all satisfy Assumption~\ref{assump1}. Therefore this is a mild assumption. 
\begin{assumption}
\label{assump2}
The data dimension $d$, random feature dimensions $N_1, N_2$, and sample size $n$ are such that $d\rightarrow\infty$, $N_1=N_1(d)\rightarrow\infty$, $N_2=N_2(d)\rightarrow\infty$, $n=n(d)\rightarrow \infty$. Moreover,  
when $d\rightarrow\infty$,  the following limits exist:
\begin{equation*}
    \lim_{d\rightarrow+\infty}N_1/d=\psi_1>0,\quad
    \lim_{d\rightarrow+\infty}N_2/d=\psi_2>0,\quad
    \lim_{d\rightarrow+\infty}n/d=\psi_3>0.\quad
\end{equation*} 
\end{assumption}

Assumption~\ref{assump2} defines the asymptotic framework for our analysis where $N_1, N_2, n, d$ go to infinity proportionally to each other. 
We let $\psi=\psi_1+\psi_2$ and  $\bpsi=[\psi_1,\psi_2,\psi_3]$.

\begin{assumption}
\label{assump3}
Let $F_{1,d}=\lVert\bbeta_{1,d}\rVert_2$. Then 
$\lim\limits_{d\rightarrow+\infty}  F_{1,d}=F_1>0$. Moreover, if $F_0\neq0$, then $\mu_{1,0}^2 + \mu_{2,0}^2 > 0$.
\end{assumption}
The condition $F_1>0$ fixes the asymptotic scale of $\bbeta_{1,d}$.  
The second condition means that when $F_0=\EE(y)\neq 0$, 
we need either $\mu_{1,0}^2>0$ or  $ \mu_{2,0}^2 > 0$ so that the predictor $f(\xb; \hat\ab, \bTheta)$ can approximate the response $y$ well when $d\rightarrow\infty$.

The statement of the main results needs some further preparation. For any $\xi\in\bbC_+$, we consider the following system of equations for the unknowns $\nu_1,\nu_2,\nu_3$:
\begin{align}
\label{eq:implicit2}
\left\{
\begin{aligned}
&\nu_{1}\cdot\bigg(-\xi-\mu_{1,2}^2\nu_{3}-\frac{\mu_{1,1}^2\nu_{3}}{1-\mu_{2,1}^2\nu_{2}\nu_{3}-\mu_{1,1}^2\nu_{1}\nu_{3}} \bigg) = \psi_1,\\
&\nu_{2}\cdot\bigg(-\xi-\mu_{2,2}^2\nu_{3}-\frac{\mu_{2,1}^2\nu_{3}}{1-\mu_{1,1}^2\nu_{1}\nu_{3}-\mu_{2,1}^2\nu_{2}\nu_{3}}\bigg)= \psi_2,\\
&\nu_{3}\cdot\bigg(-\xi-\mu_{1,2}^2\nu_{1}-\mu_{2,2}^2\nu_{2}-\frac{\mu_{1,1}^2\nu_{1}+\mu_{2,1}^2\nu_{2}}{1-\mu_{1,1}^2\nu_{1}\nu_{3}-\mu_{2,1}^2\nu_{2}\nu_{3}}\bigg) = \psi_3.
\end{aligned}
\right.
\end{align}
\xr{This system will be hereafter referred as the  $\bnu$\textit{-system}}. For different values of $\xi\in\bbC_+$, the solutions of the above system can be viewed as functions of $\xi$. We let $\bnu(\xi)= [\nu_1,\nu_2,\nu_3]^\T(\xi):~\bbC_+\rightarrow\bbC_+^3$ be the analytic function defined on $\bbC_+$ satisfying 
(i) for any $\xi\in\bbC_+$, $\bnu(\xi)$ is a solution to $\bnu$\textit{-system}~\eqref{eq:implicit2}, (ii) there exists a sufficiently large constant $\xi_0$, such that $| \nu_j(\xi) | \leq 2\psi_j / \xi_0$, for all $\xi$ with $\Im(\xi) \geq \xi_0$ and $j=1,2,3$. 
It can be shown that such a function $\bnu$ exists and is unique, and therefore our definition of $\bnu$ is valid.  
The details are given in Proposition~\ref{prop:existence_uniqueness_nu}. We hereafter denote $\bnu=\bnu(\xi,\bmu)$ to emphasize the dependence in $\bmu$.

\begin{definition}[Auxiliary matrices]\label{def:mainresultsdefinition} 
Define $\xi^*=\sqrt{\lambda}\cdot\rmi$, and
$$\nu^*_j \triangleq  \nu_j(\xi^*;\bmu), \quad j=1,2,3.
$$
\xr{Here, $\nu_j$ is the solution of $\bnu$\textit{-system}~\eqref{eq:implicit2}}. Moreover, let $M_N\triangleq \nu^*_1\mu_{1,1}^2+\nu^*_2\mu_{2,1}^2$, $M_D\triangleq \nu^*_3M_N-1$, and define the  matrices 
\if UT
\begin{equation*}
\begin{split}
      \Hb\triangleq \begin{bmatrix}
    -\frac{\nu^{*2}_3\mu_{1,1}^4}{M_D^2}+\frac{\psi_1}{\nu^{*2}_1}&-\frac{\nu^{*2}_3\mu_{1,1}^2\mu_{2,1}^2}{M_D^2}&-\frac{\mu_{1,1}^2}{M_D^2}-\mu_{1,2}^2\\
    -\frac{\nu^{*2}_3\mu_{1,1}^2\mu_{2,1}^2}{M_D^2}&-\frac{\nu^{*2}_3\mu_{2,1}^4}{M_D^2}+\frac{\psi_2}{\nu^{*2}_2}&-\frac{\mu_{2,1}^2}{M_D^2}-\mu_{2,2}^2\\
    -\frac{\mu_{1,1}^2}{M_D^2}-\mu_{1,2}^2&-\frac{\mu_{2,1}^2}{M_D^2}-\mu_{2,2}^2&-\frac{M_N^2}{M_D^2}+\frac{\psi_3}{\nu^{*2}_3}
    \end{bmatrix}, 
    \quad 
    \Vb\triangleq \begin{bmatrix}
    \mu_{1,2}^2&0&\frac{\mu_{1,1}^2}{M_D^2}&\frac{\nu^{*2}_3\mu_{1,1}^2}{M_D^2}\\
    \mu_{2,2}^2&0&\frac{\mu_{2,1}^2}{M_D^2}&\frac{\nu^{*2}_3\mu_{2,1}^2}{M_D^2}\\
    0&1&\frac{M_N^2}{M_D^2}&\frac{1}{M_D^2}
    \end{bmatrix},  
\end{split}
\end{equation*}
\fi
\begin{equation*}
\begin{split}
      \Hb\triangleq 
      \begin{bmatrix}
        -\frac{\nu^{*2}_3\mu_{1,1}^4}{M_D^2}+\frac{\psi_1}{\nu^{*2}_1}&-\frac{\nu^{*2}_3\mu_{1,1}^2\mu_{2,1}^2}{M_D^2}&\quad -\frac{\mu_{1,1}^2}{M_D^2}-\mu_{1,2}^2\\
          * & -\frac{\nu^{*2}_3\mu_{2,1}^4}{M_D^2}+\frac{\psi_2}{\nu^{*2}_2}&\quad -\frac{\mu_{2,1}^2}{M_D^2}-\mu_{2,2}^2\\
         * &  * & \quad -\frac{M_N^2}{M_D^2}+\frac{\psi_3}{\nu^{*2}_3}
    \end{bmatrix}, 
    \quad 
    \Vb\triangleq \begin{bmatrix}
    \mu_{1,2}^2&~0~&~\frac{\mu_{1,1}^2}{M_D^2}&~~\frac{\nu^{*2}_3\mu_{1,1}^2}{M_D^2}\\
    \mu_{2,2}^2&0&\frac{\mu_{2,1}^2}{M_D^2}&\frac{\nu^{*2}_3\mu_{2,1}^2}{M_D^2}\\
    0&1&\frac{M_N^2}{M_D^2}&\frac{1}{M_D^2}
    \end{bmatrix},  
\end{split}
\end{equation*}
($\Hb $ is symmetric). Finally, let $\Lb\triangleq \Vb^\T \Hb^{-1}\Vb$. 
\end{definition}

\xr{See Proposition~\ref{prop:connectGd} for the reason of selecting $\xi=\sqrt{\lambda}\cdot \rmi$}. We are now in the position to state our main theorem which establishes the 
theoretical risk curve for the double random feature model.  
\begin{theorem}
\label{thm:mainthm}
Let the data matrix $\Xb$, noise vector $\bvarepsilon$, and the DRFM model $f(\cdot;\ab,\bTheta)$ with random feature parameter matrix $\bTheta$ be defined as in Section~\ref{sec:Problemsetting}. 
Moreover, let $M_D$ and $\Lb$ be defined in Definition~\ref{def:mainresultsdefinition}. Then under Assumptions~\ref{assump1},~\ref{assump2} and \ref{assump3}, for any regularization parameter $\lambda>0$, the asymptotic excess risk  $R_d(\Xb,\bTheta,\lambda,\bbeta_d,\bvarepsilon)$ of the DRFM defined in \eqref{eq:def_Rd} satisfies
\begin{equation*}
    \EE_{\Xb,\bTheta,\bvarepsilon}\big|R_d(\Xb,\bTheta,\lambda,\bbeta_d,\bvarepsilon)-\cR(\lambda,\bpsi,\bmu,F_1,\tau)\big|=o_d(1),
\end{equation*}
where 
\begin{equation}
\label{eq:TheoremR}
\begin{split}
        \cR(\lambda,\bpsi,\bmu,F_1,\tau)=&F_1^2\bigg(\frac{1}{M_D^2}+\Lb_{3,4}+\Lb_{1,4}\bigg)+\tau^2\big(\Lb_{2,3}+\Lb_{1,2}\big),
\end{split}
\end{equation}
\xr{and $\Lb_{i,j}$ are the elements in the matrix $\Lb$ which is defined in Definition~\ref{def:mainresultsdefinition}.}
\end{theorem}
The proof of Theorem~\ref{thm:mainthm} is given in Appendix~\ref{sec:proofmainresult}. \xr{In Theorem~\ref{thm:mainthm}, the regularization parameter $\lambda$ is treated as a constant that does not depend on $d,n,p$. Note that the first three terms in \eqref{eq:TheoremR} correspond to the estimation bias, and the last two terms are the variance terms.} It can be checked that the values in $\bnu^*=[\nu_1^*,\nu_2^*,\nu_3^*]^\T$ are all purely imaginary numbers in $\bbC_+$. As the matrices $\Hb$ and $\Vb$ only depend on \xr{$\nu_j^{*2}$ (which are all negative)}, their elements are real-valued, so do the elements of the matrix $\Lb$. 
Moreover, given $\nu_{\sfc}^*$, $\sfc = 1,2,3$, the terms $\Lb_{3,4},\Lb_{1,4},\Lb_{2,3},\Lb_{1,2}$ in \eqref{eq:TheoremR} all have closed form solutions. Due to the complexity of the solutions, we defer the calculation to Appendix~\ref{sec:proofmainresult}.

\begin{remark}\label{rmk:vanillaRFM}
By inspecting the expressions of the matrices $\Hb$, $\Vb$ and $\Lb$, we see that the dependence of the asymptotic excess risk \eqref{eq:TheoremR} on the activation functions is expressed through their spherical moments $\mu_{j,1}$ and $\mu_{j,2}$, $j=1,2$. In particular, if we let $\mu_{1,1}=\mu_{2,1}$ and $\mu_{1,2}=\mu_{2,2}$, we are led to the case of a single activation function, and the asymptotic excess risk  \eqref{eq:TheoremR} coincides with the one found in \citet{mei2022generalization} 
for vanilla random feature models.
\end{remark}

\begin{remark}
Theorem~\ref{thm:mainthm} shows that the excess risk converges to $\cR(\lambda,\bpsi,\bmu,F_1,\tau)$ in $L_1$ distance, which is a  type of strong convergence. It directly implies convergence in probability: for any $\rho,\delta > 0$, there exists $d_0\in \NN$ such that for all $d \geq d_0$, $\PP\big( \big|R_d(\Xb,\bTheta,\lambda,\bbeta_d,\bvarepsilon)-\cR(\lambda,\bpsi,\bmu,F_1,\tau)\big| \leq \rho\big) \geq  1 - \delta$.
\end{remark}

\section{The phenomenon of triple descent in DRFMs}
\label{sec:simulation}
In this section we establish theoretical results showing the existence of DRFMs with triple descent risk curves and use simulations to verify our results.

\cy{Before we propose the detailed results, we first explain our intuition by considering the two extreme cases below:  

\begin{itemize}[leftmargin = *]
    \item Case 1 (no scale difference):
    As discussed in Remark~\ref{rmk:vanillaRFM}, 
    if the two activation functions have identical spherical moments, that is, $\mu_{1,1} = \mu_{2,1}$ and $\mu_{1,2} = \mu_{2,2}$, 
    the risk curve should be identical to that of a vanilla (single) random feature model. Hence, according to the study of vanilla random feature models in \citet{mei2022generalization}, the risk curve commonly has a double descent shape, with the peak at the interpolation threshold $(N_1 + N_2) / n = 1$.
    
    \item Case 2 (large scale difference): If one of the two types of random features is too small in scale compared to the other, then we can expect that this small-scale part of random features is almost negligible. For example, under the extreme case that $N_1 = N_2$ and $\sigma_2(\cdot)\equiv0 $, the second type of random features can never contribute to the learned predictor, and this case also reduces to a vanilla random feature model. Therefore we can expect the risk curve to reach the peak at $N_1 / n = 1$, that is, $(N_1 + N_2) / n = 2$. 
\end{itemize}
We can see that the two extreme cases above both lead to double descent. However, in the first case, the peak is at $(N_1 + N_2)/ n = 1$, while in the second case, the peak is at  $(N_1 + N_2)/ n = 2$. When the scales of the two parts of random features are neither too similar nor too different, we can expect the risk curve to exhibit certain characteristics from both extreme cases, possibly having two peaks at $(N_1 + N_2)/ n = 1$ and $(N_1 + N_2)/ n = 2$ respectively -- this is exactly triple descent. This motivates us to conjecture that triple descent can occur when the two parts of random features have appropriate scale differences.}

\subsection{Triple descent: theoretical results}
 The asymptotic excess risk function $\cR(\lambda,\bpsi,\bmu,F_1,\tau)$
established in Theorem~\ref{thm:mainthm} can imply the existence of triple descent in double random feature models. Note that this risk function depends on several parameters including the smoothing parameter $\lambda$, the number of features in the model and some spherical moments of the involved activation functions. Here we focus on the ``ridgeless regression'' setting where $\lambda \rightarrow 0$, and we aim to construct specific configurations of $\bmu$ such that for any fixed values of $F_1$ and $\tau$, the risk function exhibits (at least) triple descent as $(\psi_1 + \psi_2) / \psi_3$ increases.

For convenience, we use in this section the shorthand $\cR:=\cR(\lambda,\bpsi,\bmu,F_1,\tau)$. The following proposition demonstrates triple descent by considering the asymptotic regime where $\lambda \to 0 $ and ${\mu_{2,1},\mu_{2,2}\to0}$:
the former points to a limiting ridgeless regression model and the latter signifies the scale differences between two activation functions by shrinking  the second activation to 0.


 

\begin{proposition}[$\lambda\to0$]
\label{prop:property_lbd=0}
Consider the same assumptions as in Theorem~\ref{thm:mainthm} and the asymptotic excess risk function $\cR:=\cR(\lambda,\bpsi,\bmu,F_1,\tau)$.   For fixed $0<\psi_1,\psi_2,\psi_3<+\infty$ we have:
\begin{enumerate}[leftmargin = *]
    \item When  $(\psi_1+\psi_2)/\psi_3=c_1<1$,
    ~$\displaystyle     \lim\limits_{\lambda\to0}\cR <+\infty;
    $
    \item When  $(\psi_1+\psi_2)/\psi_3=1$,  
    ~$\displaystyle 
        \lim\limits_{\lambda\to0}\cR =+\infty;
    $
    \item When   $1<(\psi_1+\psi_2)/\psi_3=c_2<1+\psi_2/\psi_1$, 
    ~$\displaystyle 
    \varliminf\limits_{\mu_{2,1},\mu_{2,2}\to0}\lim\limits_{\lambda\to0}\cR <+\infty;
    $
    \item When  $(\psi_1+\psi_2)/\psi_3=1+\psi_2/\psi_1$,   
    ~$\displaystyle 
        \lim\limits_{\mu_{2,1},\mu_{2,2}\to0}\lim\limits_{\lambda\to0}\cR =+\infty.
    $
\end{enumerate}
\end{proposition}

The proof of Proposition~\ref{prop:property_lbd=0} is given in Appendix~\ref{sec:appendix_property}. \cy{This proposition theoretically demonstrate the existence of triple descent for certain DRFMs. Note that the risk function $\cR$ depends on $\psi_1,\psi_2,\psi_3$. To simplify and standardize the setting, we specifically consider the case where $\psi_3$ and the ratio $\psi_1/\psi_2$ are both fixed. In this case, the change of model complexity has a single degree of freedom, which can be characterized by $c : = (\psi_1+\psi_2)/\psi_3$. Now we can investigate the curve of the risk function with respect to $c$ and see if the shape exhibits triple descent.}

\cy{To see how Proposition~\ref{prop:property_lbd=0} demonstrates triple descent, we pick two fixed ``reference points'' $0<c_1<1$ and $1< c_2 < 1 + \psi_2/\psi_1$ (recall that we are considering the setting where $\psi_2/\psi_1$ is fixed.) By the third and fourth conclusions above, we can choose a large enough constant $M_1>0$ (Not related to $\psi_1,\psi_2$ and $\psi_3$), for which there exist $\mu_{2,1}$ and $\mu_{2,2}$ such that 
\begin{align*}
    \lim\limits_{\lambda\to0}\cR >M_1~\text{when }(\psi_1+\psi_2)/\psi_3=1+\psi_2/\psi_1,~ \text{and} ~ \lim\limits_{\lambda\to0}\cR <M_1~\text{when }(\psi_1+\psi_2)/\psi_3=c_2.
\end{align*}
For these chosen spectral moments  $\mu_{2,1}$ and $\mu_{2,2}$ and by  the first and second conclusions of the proposition,
one can find a large constant $M_2>M_1$ such that 
\begin{align*}
    \lim\limits_{\lambda\to0}\cR >M_2~\text{when }(\psi_1+\psi_2)/\psi_3=1,\quad\text{and}\quad \lim\limits_{\lambda\to0}\cR <M_2~\text{when }(\psi_1+\psi_2)/\psi_3=c_1.
\end{align*}
Recall that $\psi_1\sim N_1/d$,   $\psi_2\sim N_2/d$ and $\psi_3\sim n/d$ 
in the limits. It is customary to consider the asymptotic excess risk function  $\cR $ with respect to the ``model complexity parameter'' $c=\lim_{d\to+\infty}(N_1+N_2)/n=(\psi_1+\psi_2)/\psi_3$.}  \xr{Based on this analysis, we are able to find constants $0 < M_1 < M_2$ and $\mu_{2,1}$, $\mu_{2,2}$ that do not depend on $\psi_1,\psi_2,\psi_3$, so that the following four results hold:}
\begin{enumerate}[leftmargin = *]
    \item \cy{$c=c_1$, $\lim\limits_{\lambda\to0}\cR <M_2$;}
    \item $c=1$, $\lim\limits_{\lambda\to0}\cR >M_2$;
    \item \cy{$c=c_2$, $\lim\limits_{\lambda\to0}\cR <M_1$;}
    \item $c=1+\psi_2/\psi_1$, $\lim\limits_{\lambda\to0}\cR >M_1$.
\end{enumerate}
\xr{These four cases above correspond to four situations with different model complexities: each case is for a specific value of $c = (\psi_1 + \psi_2)/\psi_3= \lim_{d\to+\infty} (N_1 + N_2)/n$.}
The next proposition shows that the risk function has a finite limit  when the model complexity parameter $c$ tends to infinity, or in other words, in the infinitely over-parameterized regime.

\begin{proposition}[$\psi_1,\psi_2\to+\infty$]
\label{prop:property_psi}
Consider the same assumptions as in Theorem~\ref{thm:mainthm} and the asymptotic excess risk function $\cR:=\cR(\lambda,\bpsi,\bmu,F_1,\tau)$ with non-degenerate activation functions.    
For fixed $\psi_3$ and $r_1,r_2>0$, we have
\begin{align*}
    \lim\limits_{\substack{ \psi_1,\psi_2 \to+\infty \\ \psi_1/r_1=\psi_2/r_2 } }\cR =
     \frac{F_1^2\psi_3+\tau^2\chi_0^2} 
    {(\chi_0+1)^2\psi_3-\chi_0^2}, 
\end{align*}
where
\begin{align*}
    &\chi_0=\frac{(r_1\mu_{1,1}^{2}+r_2\mu_{2,1}^2)\chi_1}
   {2\sum\limits_{i,j=1}^2r_ir_j\mu_{i,1}^2\mu_{j,2}^2}, \\
    &\chi_1=(\psi_3-1)\sum\limits_{i=1}^2r_i\mu_{i,1}^2-\sum\limits_{i=1}^2r_i\mu_{i,2}^2+\sqrt{\Big((\psi_3-1)\sum\limits_{i=1}^2r_i\mu_{i,1}^2-\sum\limits_{i=1}^2r_i\mu_{i,2}^2\Big)^2+4\psi_3\sum\limits_{i,j=1}^2r_ir_j\mu_{i,1}^2\mu_{j,2}^2}.
\end{align*}
\end{proposition}

\begin{figure}[t!]
    \centering
    \includegraphics[width=4in]{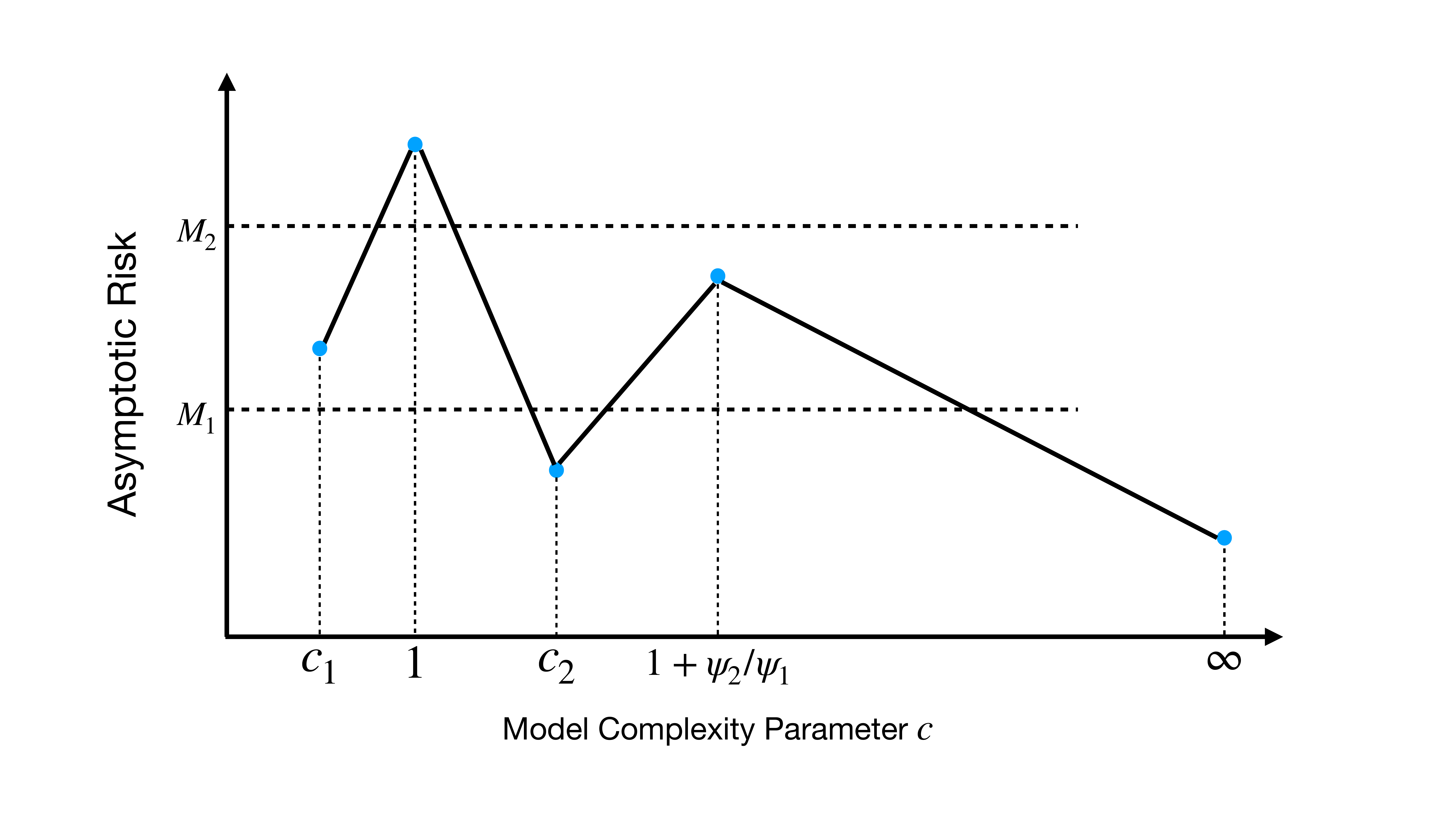}\\
    \vspace{-10pt}
    \caption{Existence of  triple descent in a double random feature model: the four points $c_1$ to $1+\psi_2/\psi_1$ for the model complexity parameter \xr{$c=(\psi_1+\psi_2)/\psi_3$} are found  in  Proposition~\ref{prop:property_lbd=0} and  the last point depicts  the limit found in  Proposition~\ref{prop:property_psi} when $c\to\infty$. }
    \label{fig:example}
    \vspace{-10pt}
\end{figure}

\xr{The proof of Proposition~\ref{prop:property_psi} can be found in Appendix~\ref{sec:appendix_property}. By combining the derived limiting risk value from Proposition~\ref{prop:property_psi}, where $c$ tends to infinity, with the summary provided after Proposition~\ref{prop:property_lbd=0}, we can observe that the asymptotic excess risk function $\cR$ exhibits (at least) triple descent with the chosen parameter values. This behavior occurs as the model complexity parameter $c$ increases from 0 to $c_1,1,c_2,1+\psi_2/\psi_1 $, and eventually tends to infinity. A visual representation of this phenomenon can be seen in Figure~\ref{fig:example}. Furthermore, when the model complexity is $c<1$, the asymptotic risk takes the form of a $\sfU$ shape, which aligns with classical theory.}

\xr{
\begin{remark}
\label{remark:tripledescent}
We can also consider the case where $\psi_1,\psi_2$ goes to zero, and this case corresponds to the setting where random feature model is almost reduced to a constant predictor. In this case, it is easy to show that $\lim\limits_{ \psi_1,\psi_2 \to0}\cR =F_1^2$. In classical statistical theory, the first descent occurs here when the model complexity gradually increases: as the predictor becomes more complicated than a constant predictor, the asymptotic risk will first decrease below $F_1^2$.
\end{remark}
}

\subsection{Triple descent: empirical evidence}
\label{subsec:tripleinDRFM}
In this subsection, we 
empirically demonstrate the triple descent phenomenon in double random feature models.
The simulation design is as follows.
\begin{itemize}[leftmargin = *]
    \item  Training data  $\{(\xb_i, y_i)\}_{i=1}^n$ are generated independently following   Definition~\ref{def:data_model} with $\tau = 0.1$: each $\xb_i$ is uniformly generated from the sphere $\sqrt{d}\cdot \SSS^{d-1}$, and the corresponding response is given as $y_i=\la\bbeta_1,\xb_i\ra+F_0+\varepsilon_i$, 
    where $\bbeta_1$ is a randomly chosen unit vector;
    \item $F_0=0.2$, $\lambda=10^{-5}$; 
    \item Training sample size $n=1000$, data dimension $d=300$  and $N_1 = N_2$ varying   from $0 $ to $1.6 n$.
\end{itemize}
As we gradually increase the dimensions of random features $N_1=N_2$ from $0 $ to $1.6 n$, the model complexity parameter $c(d)=(N_1+N_2)/n$ varies from 0 to 3.2.  
The  empirical and finite-horizon values for the limiting
 excess risk $\cR(\lambda,\bpsi,\bmu,F_1,\tau)$  in  Theorem~\ref{thm:mainthm}  are obtained  on  a test data set of size $700$ and  averaged from 30 independent replications.
 
\begin{figure}[t!]
    \centering
    \includegraphics[width=0.95\textwidth]{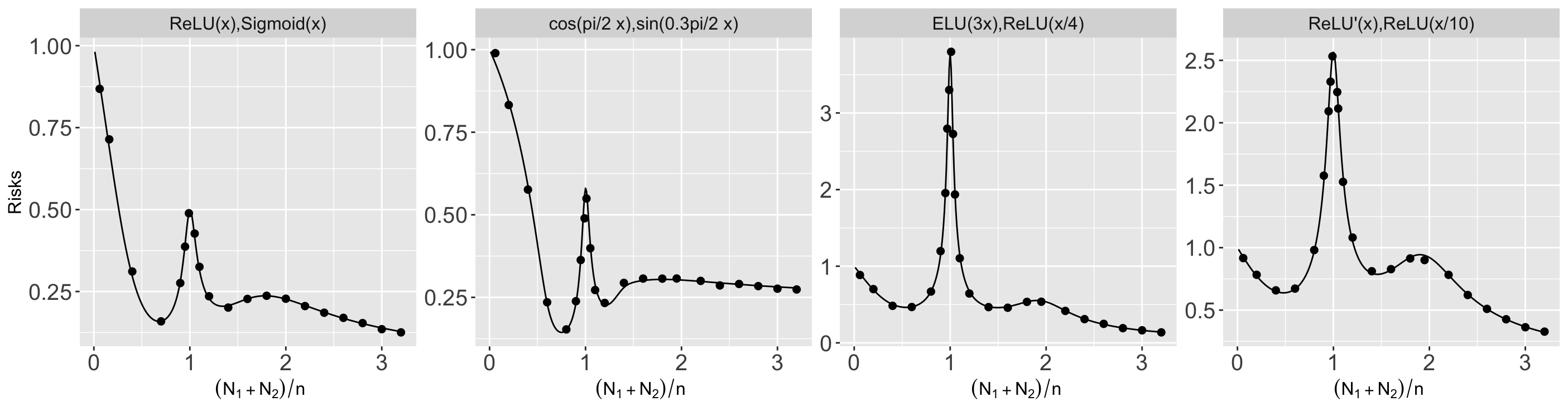}\vspace{-5pt}
    
   \hspace{12pt}(a)\hspace{86pt}\hspace{12pt}(b)\hspace{80pt}\hspace{10pt}(c)\hspace{80pt}\hspace{15pt}(d)
    \vspace{-10pt}
    \caption{Triple descent in double random feature models with different activation functions. 
    The plots show both the asymptotic excess risks (curves) and empirical excess risks (dots).
    From (a) to (d), the activation functions are $\big(\ReLU(x),\Sigmoid(x)\big)$, $\big(\cos(\frac{\pi}{2}x),\sin(\frac{0.3\pi}{2}x)\big)$, $\big(\ELU(3x),\ReLU(x/4)\big)$ and $\big(\ReLU'(x),\ReLU(x/10)\big)$.}
    \label{fig:exist}
    \vspace{-10pt}
\end{figure}
The results are given in Figure~\ref{fig:exist}. 
In this figure (and all other figures of this section),  the values of the asymptotic risk  $\cR(\lambda,\bpsi,\bmu,F_1,\tau)$ are shown as continuous curves while empirical  risk values are plotted using black dots. We consider activation functions     $\ReLU(x) = x_+$, 
$\ReLU'(x) = \ind\{ x > 0\}$, $\Sigmoid(x) = 1 / (1 + e^{-x})$, $\ELU(x) = x_+ - (1-e^x)_-$, 
as well as trigonometric functions $\cos(x)$ and $\sin(x)$. We slightly scale the activation functions to show clearer shapes of triple descent: the four plots in  Figure~\ref{fig:exist} represent  DRFMs with activation pairs $\big(\ReLU(x),\Sigmoid(x)\big)$, $\big(\cos(\frac{\pi}{2}x),\sin(\frac{0.3\pi}{2}x)\big)$, $\big(\ELU(3x),\ReLU(x/4)\big)$ and $\big(\ReLU'(x),\ReLU(x/10)\big)$, respectively.

Clearly, the empirical risk values well match their theoretical counterparts in all the examined settings, which empirically validates the asymptotic risks established  in Theorem~\ref{thm:mainthm}. 
More importantly, these risk curves all exhibit triple descent
as predicted by Propositions~\ref{prop:property_lbd=0} and  
\ref{prop:property_psi} (see also Figure~\ref{fig:example}),  where the four critical constants have the following  values under the present experimental design: 
\[  c_1 <1,  \quad c_2=1, \quad 1<c_3<2, \quad c_4=2.
\]

\subsection{Impact of scale difference on triple descent} 
\label{subsec:explaintriple}  
\begin{figure}[t!]
    \centering
    \includegraphics[width=0.95\textwidth]{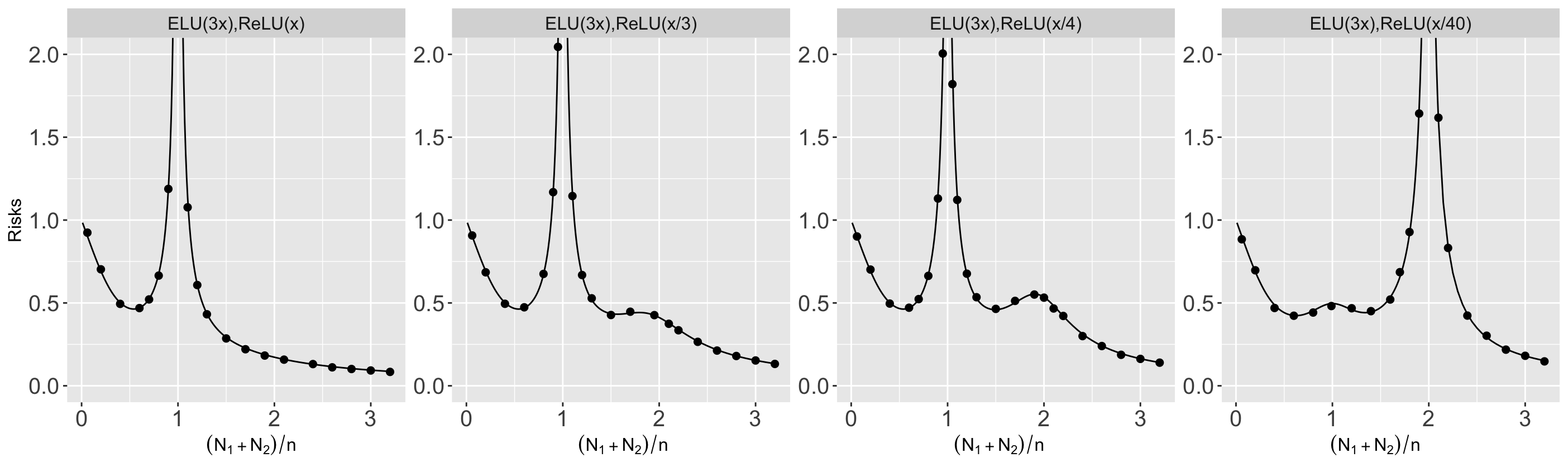}\vspace{-5pt}
    
    \hspace{10pt}(a)\hspace{86pt}\hspace{12pt}(b)\hspace{80pt}\hspace{12pt}(c)\hspace{80pt}\hspace{15pt}(d)
    \vspace{-10pt}
    \caption{Risk curves of DRFMs with scaled $\ReLU$ and $\ELU$ activation functions. The plots show both the asymptotic excess risks (curves) and empirical excess risks (dots).
     From (a) to (d), the activation functions are $(\ELU(3x),\ReLU(x))$, $(\ELU(3x),\ReLU(x/3))$, $(\ELU(3x),\ReLU(x/4))$ and $(\ELU(3x),\ReLU(x/40))$ respectively.}
    \label{fig:activations1}
    \vspace{-10pt}
\end{figure}

\begin{figure}[t!]
    \centering
    \includegraphics[width=0.95\textwidth]{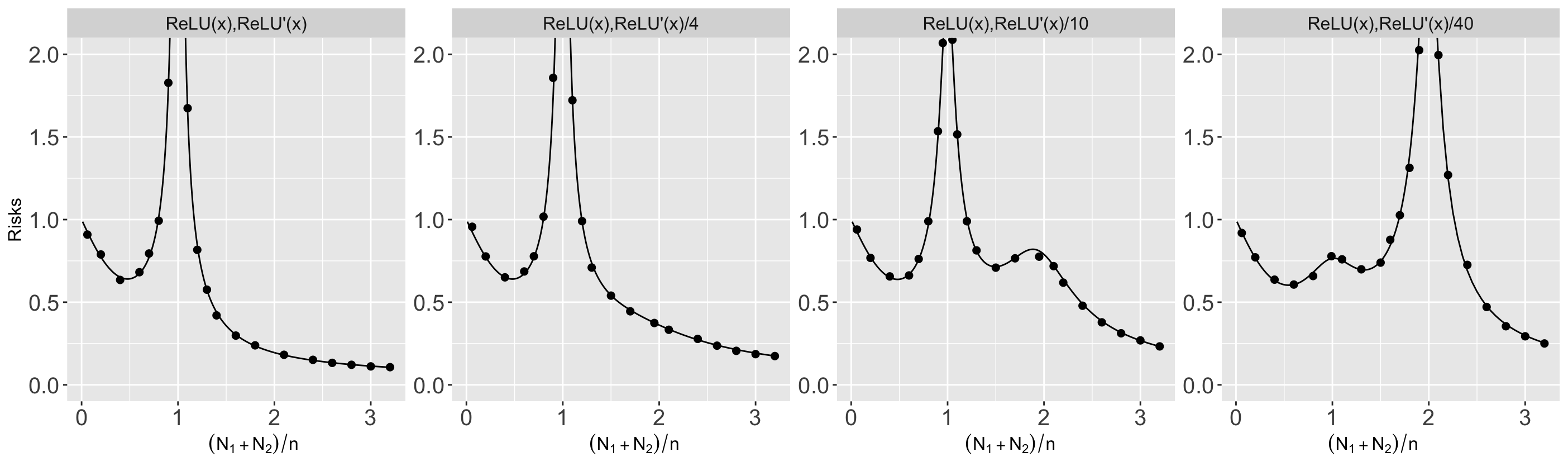}\vspace{-5pt}
    
    \hspace{10pt}(a)\hspace{86pt}\hspace{10pt}(b)\hspace{80pt}\hspace{15pt}(c)\hspace{80pt}\hspace{15pt}(d)
    \vspace{-10pt}
    \caption{Risk curves of DRFMs with scaled $\ReLU$ and $\ReLU'$ activation functions. The plots show both the asymptotic excess risks (curves) and empirical excess risks (dots). From (a) to (d), the activation functions are $(\ReLU(x),\ReLU'(x))$, $(\ReLU(x),\ReLU'(x)/4)$, $(\ReLU(x),\ReLU'(x)/10)$ and $(\ReLU(x),\ReLU'(x)/40)$ respectively.}
    \label{fig:activations2}
    \vspace{-10pt}
\end{figure}

As demonstrated in Propositions~\ref{prop:property_lbd=0} and  
\ref{prop:property_psi}, when the magnitude of a random feature is of a smaller order than the other feature,  triple descent appears in a DRFM. In this section, we  use our theoretical predictions as well as simulations to verify this result. The experiment setups are the same as the experiments in Section~\ref{subsec:tripleinDRFM}, except that here we use different pairs of activation functions. For two activation functions $\sigma_1, \sigma_2$, we gradually decrease the scale of $\sigma_2$ by using activation pairs $(\sigma_1(x),c_0\sigma_2(x))$ with a smaller and smaller factor $c_0$. Results for activation pairs $(\ELU,\ReLU)$ and $(\ReLU,\ReLU')$ are reported in Figure~\ref{fig:activations1} and Figure~\ref{fig:activations2}, respectively. Clearly, in both figures, the empirical errors (dots) well match their theoretical counterparts (curves). Moreover, In Figure~\ref{fig:activations1} (a) and Figure~\ref{fig:activations2} (a), we present the result when we appropriately balance the two activation functions such that the two parts of the random features have similar scales, and the resulting risk curves exhibit double descent with a peak at $(N_1 + N_2 ) / n = 1$. As the scale of the second random feature decreases, the risk curves transit from double descent curves to triple descent curves in Figure~\ref{fig:activations1} (b), (c) and Figure~\ref{fig:activations2} (b), (c). Finally, in  Figure~\ref{fig:activations1} (d) and Figure~\ref{fig:activations2} (d) when the scale differences are large, the risk curves have a large peak near $c = 2$ but only a very small peak near $c = 1$. Clearly, these results perfectly match Proposition~\ref{prop:property_lbd=0}, and thus backs up the triple descent phenomena in DRFMs. 

\subsection{Impact of the ratio between random feature dimensions} 
Our previous experiments are all under the setting where $N_1 = N_2$, which corresponds to the case where the two parts of the random features have the same dimensions. In fact,  we can study more general settings where $N_1$ and $N_2$ hold a ratio other than $1$. Specifically, suppose that $\sigma_1$ has larger scale compared to $\sigma_2$. Then based on Proposition~\ref{prop:property_lbd=0}, it is clear that the first peak should be around $c= 1$, while the second peak should be around $c=1+\psi_2/\psi_1$. 

\begin{figure}[t!]
    \centering
    \includegraphics[width=0.95\textwidth]{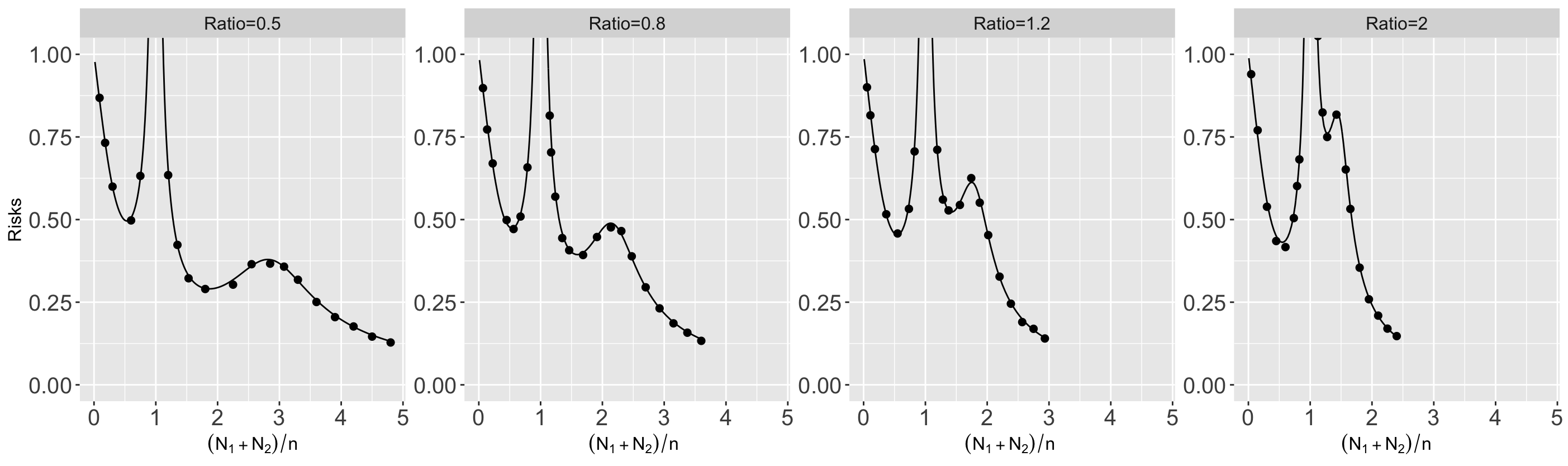}\vspace{-5pt}
    
    \hspace{13pt}(a)\hspace{86pt}\hspace{10pt}(b)\hspace{80pt}\hspace{15pt}(c)\hspace{80pt}\hspace{15pt}(d)
    \vspace{-10pt}
    \caption{Risk curves of DRFMs with different ratios between random feature dimensions. The plots show both the asymptotic excess risks (curves) and empirical excess risks (dots). From (a) to (d), the ratios $N_1/N_2$ are 0.5, 0.8, 1.2 and 2, respectively. The activation functions are chosen as $\sigma_1(x)=\ELU(3x)$ and $\sigma_2(x)=\ReLU(x/4)$.}
    \label{fig:ratio}
\end{figure}
 
We now consider the same experiment setup as in Section~\ref{subsec:tripleinDRFM}, except that here we focus on the activation pair $\big(\ELU(3x),\ReLU(x/4)\big)$, and no longer require $N_1 =N_2$. Instead, we consider the ratios $N_1/N_2\in\{0.5, 0.8, 1.2,2 \}$ and plot the corresponding risk curves. Note that the coordinates in the first part of random features are about 10 times those in the second part (in magnitude), and the second peak in the risk curve is expected to be around the position $1 + (N_1 / N_2)^{-1}$. 




The simulation results are reported in Figure~\ref{fig:ratio}.  It can be seen that the second peaks in Figure~\ref{fig:ratio} (a), (b), (c), (d) are around $c=1 + (N_1 / N_2)^{-1} = 3,9/4,11/6,3/2$, respectively. This further verifies Proposition~\ref{prop:property_lbd=0}, and shows how one can design double random feature models with specific peak locations. 
We have also studied other key factors affecting the risk curve, such as the regularization parameter and the signal-to-noise ratio. Details of experimental results are reported in Appendix~\ref{sec:otherfactor}.



%


\subsection{Further discussion}
\label{sec:NTK}
\xr{Due to the complexity of the theoretical expressions involving almost 10 variables, it is difficult to provide a precise characterization such as under what conditions is the 2nd descent is lower than the 1st descent, or under what conditions is the peak of the second descent lower than the bottom of the 1st descent. While we have found empirically that the second peak tends to appear when the scale of $\sigma_2$ is small enough, it is hard to make a formal statement on  the general conditions that guarantee this fact. We would also like to note that Appendix~\ref{sec:otherfactor} provides some
analysis on the effects of SNR and regularization on the multiple descent phenomenon. Specifically, we have observed that SNR affects the trend of the risks in the under-parameterized regime ($(N_1+N_2)/n<1$) and the highly over-parameterized regime ($(N_1+N_2)/n>2$, while $\lambda$ affects the existence of the peak. Additionally,  benign overfitting tends to occur when the SNR is
high, while optimal regularization can help mitigate the multiple descent, as has been shown in
previous literature \citep{nakkiran2020optimal,mei2022generalization}.}

\section{The multiple random feature model}
\label{sec:generalcase}
In the previous sections, we have studied double random feature models   based on two activation functions. In this section, we extend our results to the case with  
$K$ activation functions ($K \in \NN_+$).  


Suppose that for $j\in [K]$, there are $N_j$ random feature units using activation function $\sigma_j$. Then we let $N = N_1 + \cdots + N_K$ be the total dimension of the random features. Moreover, we define the index set of the random feature units using the activation function $\sigma_j$ as
$$
\cN_{\sfc} = \Bigg\{ i\in[N] : 1 + \sum_{r=1}^{\sfc - 1} N_r \leq i \leq \sum_{r=1}^{\sfc} N_r \Bigg\},~\sfc\in[K].
$$
Let $\btheta_i\sim \Unif(\sqrt{d}\cdot\SSS^{d-1})$, $i\in[N]$ be the random feature parameter vectors and $a_i\in\RR$, $i\in[N]$ be the linear combination coefficients of the random features. Then we denote $\bTheta = [\btheta_1,\ldots,\btheta_N]^\top \in \RR^{N\times d}$, $\ab=[a_1,\ldots, a_N]^\top \in \RR^N$. 
A \textit{multiple random feature model} (MRFM) predictor is defined as
$$
f(\xb;\ab,\bTheta)=\sum_{\sfc=1}^K\sum_{i\in\cN_{\sfc}}a_i\sigma_\sfc\big(\la\btheta_i,\xb\ra/\sqrt{d}\big).
$$
We also denote by $\bTheta_\sfc = [\btheta_{\cN_\sfc}]^\T\in\RR^{N_{\sfc}\times d}$ the collection of the random feature parameter vectors using the activation function $\sigma_{\sfc}$.
We learn the same data model in Definition~\ref{def:data_model} by fitting a training data set $S = \{(\xb_i,y_i)\}_{i=1}^n$ with the function $f(\xb;\ab,\bTheta)$ using ridge regression. 
Similar to Section~\ref{sec:Problemsetting}, we learn the coefficient vector $\ab$ by minimizing the $\ell_2$-regularized square loss function:
\begin{equation*}
\hat{\ab}=\argmin_{\ab}\left\{\frac1n\sum_{j=1}^n\Big(y_j-f(\xb_j;\ab,\bTheta)\Big)^2+\frac{d}{n}\lambda\lVert\ab\rVert_2^2 \right\}. 
\end{equation*}
The excess risk is denoted by  $R_d(\Xb,\bTheta,\lambda,\bbeta_d,\bvarepsilon)$  highlighting its dependence on $\Xb,\bTheta,\lambda,\bbeta_d$ and $\bvarepsilon$:
\begin{align}
    R_d(\Xb,\bTheta,\lambda,\bbeta_d,\bvarepsilon) 
    = \EE_{\xb\sim\Unif(\sqrt{d}\cdot\SSS^{d-1}) } \big[ F_0+\xb^\T \bbeta_{1,d} - f(\xb; \hat\ab, \bTheta)\big]^2.\label{eq:def_Rd-m}
\end{align}

\subsection{\xr{Excess risks of}  MRFMs}
The definitions and assumptions  below are similar to  those   previously used for DRFMs  in Section~\ref{sec:mainresult2}. 
\begin{definition}
For $\sfc=1,2,\ldots,K$ and $G\sim\rmN(0,1)$, define
\begin{align*}
    \mu_{\sfc,0}\triangleq\EE\sigma_{\sfc}(G),\qquad\mu_{\sfc,1}\triangleq\EE G\sigma_{\sfc}(G),\qquad\mu_{\sfc,2}^2\triangleq{\EE\{ \sigma_{\sfc}^2(G)\}-\mu_{\sfc,0}^2-\mu_{\sfc,1}^2}.
\end{align*}
 These spherical moments are collected into a vector $\bmu$. 
\end{definition}
 
\begin{assumption}
\label{assump1-m}
Let $\sigma_{\sfc}~:~\RR\rightarrow\RR$ $(\sfc=1,2,\ldots,K)$ be weakly differentiable, with weak derivative $\sigma_{\sfc}'$. Assume $|\sigma_\sfc(u)|\vee |\sigma_\sfc'(u)|\leq C_0e^{C_1|u|}$ for some constants $C_0,C_1<+\infty$. 
\end{assumption}

\begin{assumption}
\label{assump2-m}
We consider sequences of parameters $N_1,N_2,\ldots,N_K,n,d$ that go to infinity proportionally to each other. Without loss of generality, let the sequences be indexed by $d$, and assume  for $\sfc=1,\ldots,K$, the following limits exist:
\begin{equation*}
    \lim_{d\rightarrow+\infty}N_\sfc/d=\psi_\sfc\in(0,\infty),\qquad
    \lim_{d\rightarrow+\infty}n/d=\psi_{K+1}\in(0,\infty).
\end{equation*}
These limits are collected into the vector   $\bpsi=[\psi_1,\ldots,\psi_{K},\psi_{K+1}]$. 
\end{assumption}

\begin{assumption}
\label{assump3-m}
Let $F_{1,d}=\lVert\bbeta_{1,d}\rVert_2$. Then 
$\lim\limits_{d\rightarrow+\infty}  F_{1,d}=F_1>0$. Moreover, if $F_0\neq0$, then $\sum\limits_{\sfc=1}^{K}\mu_{\sfc,0}^2 > 0$. 
\end{assumption}


All these assumptions are natural, and parallel Assumptions \ref{assump1}-\ref{assump3} in Section~\ref{sec:mainresult2}, respectively. The presentation of the results for the MRFM also relies on a system of self-consistent equations  as follows. 
For $\xi\in\bbC_+$, consider the following system of equations with unknown functions  $(\nu_1,\ldots,\nu_{K+1})$: $\bbC_+\rightarrow\bbC_+^{K+1}$ (as functions of the complex variable $\xi$):
\begin{align}
\label{eq:implicit2-m}
\left\{
\begin{aligned}
&\nu_\sfc\cdot \bigg(-\xi-\mu_{\sfc,2}^2\nu_{K+1}-\frac{\mu_{\sfc,1}^2\nu_{K+1}}{1-\sum_{\sfc=1}^{K}\mu_{\sfc,1}^2\nu_{\sfc}\nu_{K+1}} \bigg) = \psi_{\sfc},\quad\sfc=1,\ldots,K\\
&\nu_{K+1}\cdot \Bigg(-\xi-\sum_{\sfc=1}^{K}\mu_{\sfc,2}^2\nu_{\sfc}-\frac{\sum_{\sfc=1}^{K}\mu_{\sfc,1}^2\nu_{\sfc}}{1-\sum_{\sfc=1}^{K}\mu_{\sfc,1}^2\nu_{\sfc}\nu_{K+1}}\Bigg) = \psi_{K+1}.
\end{aligned}
\right.
\end{align}
We let $\bnu= [\nu_1,\ldots,\nu_{K+1}]^\T:~\bbC_+\rightarrow\bbC_+^{K+1}$ be the analytic function defined on $\bbC_+$ satisfying 
(i) for any $\xi\in\bbC_+$, $\bnu(\xi)$ is a solution to $\bnu$\textit{-system}~\eqref{eq:implicit2-m}, (ii) there exists a sufficiently large constant $\xi_0$, such that $| \nu_j(\xi) | \leq 2\psi_j / \xi_0$ for all $\xi$ with $\Im(\xi) \geq \xi_0$ and $j\in [K]$. 
It can be shown that such a function $\bnu$ exists and is unique, and therefore our definition of $\bnu$ is valid. 
The full justification is given in Proposition~\ref{prop:existence_uniqueness_nu-m} in the appendix. We also denote $\bnu=\bnu(\xi,\bmu)$ to emphasize the dependence on $\bmu$.

\begin{definition}[Auxiliary matrices]
\label{def:mainresultsdefinition-m}
Define $\xi^*=\sqrt{\lambda}\cdot{\rmi}$,
$$\bnu^*=[\nu_1^*,\ldots,\nu_{K+1}^*]^\T=\big[\nu_1,\ldots,\nu_{K+1}\big]^\T (\xi^*;\bmu)$$
\xr{where $\nu_j$ is the solution of $\bnu$\textit{-system}~\eqref{eq:implicit2-m}}, and let
\[ M_N=\sum\limits_{\sfc=1}^{K}\mu_{\sfc,1}^2\nu_{\sfc}^*,\quad  M_D=\nu_{K+1}^*M_N-1.
\]
We then let $\Hb\in\RR^{(K+1)\times(K+1)}$ be a real symmetric matrix whose $(i,j)$-th entry ($i\le j$) is
 \[
\Hb_{i,j}=\left\{  
\begin{aligned}
&-\frac{\nu_{K+1}^{*2}\mu_{i,1}^4}{M_D^2}+\frac{\psi_{i}}{\nu_i^{*2}},\quad &1\leq i=j\leq K,\\
&-\frac{\nu_{K+1}^{*2}\mu_{i,1}^2\mu_{j,1}^2}{M_D^2},\quad &1\leq i<j\leq K,\\
&-\frac{\mu_{i,1}^2}{M_D^2}-\mu_{i,2}^2  ,\quad & 1\leq i\leq K,~j=K+1,\\
& -\frac{M_N^2}{M_D^2}+\frac{\psi_{K+1}}{\nu_{K+1}^{*2}} ,\quad& i=j=K+1.
\end{aligned}\right. 
\]
Moreover, define $\Vb=[\vb_1,\vb_2,\vb_3,\vb_4]\in\RR^{(K+1)\times4}$, where 
\begin{align*}
\vb_1&=\Big[\mu_{1,2}^2,\mu_{2,2}^2,\ldots,\mu_{K,2}^2,0\Big]^\T,\quad\vb_2=\Big[0,\ldots,0,1\Big]^\T,\\
\vb_3&=\Big[\frac{\mu_{1,1}^2}{M_D^2},\ldots,\frac{\mu_{K,1}^2}{M_D^2},\frac{M_N^2}{M_D^2}\Big]^\T,\quad\vb_4=\Big[\nu_{K+1}^{*2}\frac{\mu_{1,1}^2}{M_D^2},\ldots,\nu_{K+1}^{*2}\frac{\mu_{K,1}^2}{M_D^2},\frac{1}{M_D^2}\Big]^\T.
\end{align*}
Finally, let  
$
    \Lb=\Vb^\T \Hb^{-1}\Vb\in\RR^{4\times4}  
$. 
\end{definition}
It is clear that the above definitions are consistent with Definition~\ref{def:mainresultsdefinition} for the case of $K =2$. Based on these definitions, the asymptotic limit of the excess risk can be expressed as function of the elements of the matrix $\Lb$. Our main result for MRFMs is given in the following theorem.
\begin{theorem}
\label{thm:mainthm-m}
Let the data matrix $\Xb$ and the noise vector $\bvarepsilon$ be generated as in Definition~\ref{def:data_model}.  
Then under Assumptions~\ref{assump1-m},~\ref{assump2-m} and \ref{assump3-m}, for any regularization parameter $\lambda>0$, the asymptotic excess risk $R_d(\Xb,\bTheta,\lambda,\bbeta_d,\bvarepsilon)$ of the MRFM defined in \eqref{eq:def_Rd-m} satisfies 
\begin{equation*}
    \begin{split}
        \EE_{\Xb,\bTheta,\bvarepsilon}\big|R_d(\Xb,\bTheta,\lambda,\bbeta_d,\bvarepsilon)-\cR(\lambda,\bpsi,\bmu,F_1,\tau)\big|=o_d(1),
    \end{split}
\end{equation*}
where, with  $M_D$ and the matrix $\Lb$ defined in Definition~\ref{def:mainresultsdefinition-m}, 
\begin{equation}
\label{eq:TheoremR-m}
\begin{split}
        \cR(\lambda,\bpsi,\bmu,F_1,\tau)=&F_1^2\bigg(\frac{1}{M_D^2}+\Lb_{3,4}+\Lb_{1,4}\bigg)+\tau^2\big(\Lb_{2,3}+\Lb_{1,2}\big).
\end{split}
\end{equation}
\xr{Here, $\Lb_{i,j}$ are the elements in the matrix $\Lb$ which is defined in Definition~\ref{def:mainresultsdefinition-m}}.
\end{theorem}
Theorem~\ref{thm:mainthm-m} is proved  in Appendix~\ref{sec:appendixthm-m}. The asymptotic excess risk for the MRFM given in Equation~\eqref{eq:TheoremR-m} is similar to \eqref{eq:TheoremR} for the DRFM. It is clear that Theorem~\ref{thm:mainthm-m} covers Theorem~\ref{thm:mainthm} and the results in \citet{mei2022generalization} as special cases with $K = 2$ and $K = 1$, respectively. 

\subsection{Multiple descent in MRFMs}
\label{subsec:verifyinmulti}
We now  demonstrate the existence of multiple descent in MRFMs.
The experimental setting is similar to the previous experiments reported in Section~\ref{sec:simulation}. We set $d=300$, $n=1000$, and $\lambda=10^{-4}$. 
In simulation, the training data $\{(\xb_i, y_i)\}_{i=1}^n$ are generated independently according to Definition~\ref{def:data_model}: each $\xb_i$ is uniformly generated from the sphere $\sqrt{d}\cdot \SSS^{d-1}$, and the corresponding response is given as $y_i=\la\bbeta_1,\xb_i\ra+F_0+\varepsilon_i$, 
where $\bbeta_1$ is a randomly chosen unit vector, $F_0=0.2$ and $\tau=0.1$. We estimate the excess risks of the MRFMs with a test data set of size $700$, and take average over $30$ independent runs. We consider two MRFMs with $K=3$ and $K=4$,  respectively. For the case $K = 3$, we consider three activation functions $\sigma_1(x)=\ReLU(9x)$, $\sigma_2(x)=\ReLU(x)$ and $\sigma_3(x)=\ReLU(0.1x)$, and set the ratios between dimensions of random features as $N_1=N_2=N_3/3$. 
For the case $K = 4$, we use four activation functions $\sigma_1(x)=\ReLU(80x)$, $\sigma_2(x)=\ReLU(9x)$, $\sigma_3(x)=\ReLU(x)$ and $\sigma_4(x)=\ReLU(0.1x)$, and keep the ratios $N_1=N_2=N_3=N_4/3$.

\begin{figure}[t!]
    \centering
    \includegraphics[width=5in]{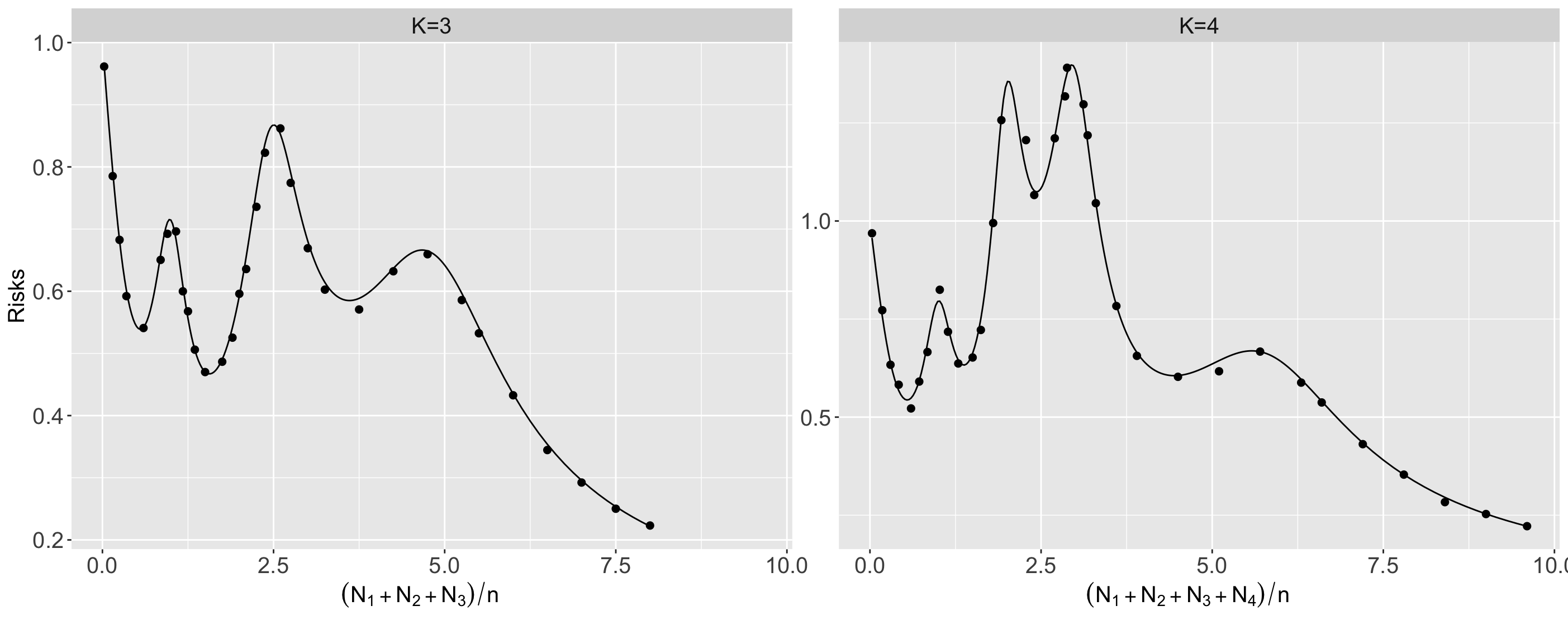}\vspace{-5pt}
    
   \hspace{10pt}(a)\hspace{178pt}(b)
   \vspace{-10pt}
    \caption{Multiple descent in multiple random feature models. (a) gives the risk curve for the MRFM with three activation functions, which exhibits quadruple descent; (b) shows the risk curve for the MRFM with four activation functions, which exhibits quintuple descent.
    }
    \label{fig:multi1}
    \vspace{-10pt}
\end{figure}

The results are given in Figure~\ref{fig:multi1}. We can see that the simulation results (dots) well match the theoretically derived risks (curves), which validates our results in Theorem~\ref{thm:mainthm-m}. Moreover, Figure~\ref{fig:multi1} (a) (where we use three different activation functions) shows quadruple descent, while Figure~\ref{fig:multi1} (b) (where we use four different activation functions) shows quintuple descent. With these observations, we believe an MRFM using $K$ activation functions may exhibit $(K+1)$-fold descent.


Following a similar analysis as in Section~\ref{sec:simulation}, we can also study the locations of each peak in the risk curves as follows. 
First consider the experiment with $K = 3$. Clearly, the first peak always locates around $(N_1 + N_2 + N_3)/ n = 1$. Regarding the second peak, note that the scales of the activation functions are set in descending order. Under this case, the first two types of random features will mainly contribute to the predictor and the third type of random features is negligible, therefore we have $(N_1 + N_2)/ n = 1$ around the second peak. Since $N_1=N_2=N_3/3$, we have $N_1 = N_2 = n / 2$ and $N_3 = 3n/2$. Hence we conclude that the second peak should be around $(N_1 + N_2 + N_3)/ n = 2.5$. Similarly, regarding the third peak, we have $N_1/n=1$, which indicates that the peak locates around $(N_1 + N_2 + N_3)/ n = 5$. These predicted locations clearly match the results shown in Figure~\ref{fig:multi1} (a). For the case $K = 4$, with a similar argument, we can expect that the four peaks are located around $1,2,3,6$, respectively. This also matches the result in Figure~\ref{fig:multi1} (b). 

\section{Conclusion}
\label{sec:Conclusion}
This paper considers the learning of double random feature models and multiple random feature models. We give the explicit formulas for the asymptotic excess risks achieved by DRFMs and MRFMs. 
These theoretical results are further well confirmed by empirical simulations in various settings. We provide an explanation of the triple descent and multiple descent phenomena based on the scale difference between activation functions, and discuss how the ratio between random feature dimensions control the location of the second peaks in the risk curves. By showing that MRFMs with $K$ types of random features may exhibit $(K+1)$-fold descent, we demonstrate that risk curves with a specific number of descent generally exist in random feature based regression.

An immediate future work direction is to study ridgeless regression where $\lambda = 0$. Moreover, our result can help future studies on the advantages and disadvantages of overfitting by quantitatively comparing the risks achieved by over-parameterized/under-parameterized models with different regularization levels. Extending our findings to deep learning is another important future direction.


\bibliography{ref}
\bibliographystyle{ims}

\appendix

\section{Proof of Theorem~\ref{thm:mainthm}}
\label{sec:proofmainresult}

 The proof  is presented in the following four steps. 
\begin{enumerate}[leftmargin = *]
    \item  We  first develop a decomposition of the  risk and find an asymptotic approximation whose main terms are expressed as traces of several random matrices,
    see Proposition~\ref{prop:decomp-expectation}; 
    \item  We then  create a new random matrix called the  linear pencil matrix, which includes all the fundamental random matrices involved in the asymptotic approximation found in the first step,  so that the needed traces are all functions of the limiting spectrum of the linear pencil matrix, see Proposition~\ref{prop:connectGd}; 
    \item Next,  we find the key limiting spectral functions of the linear pencil matrix including its Stieltjes transform and logarithmic potential, and show that the needed traces converge to some specific partial derivatives of the limiting logarithmic potential, see Propositions~\ref{prop:implicit1} and  \ref{prop:substitutiongd}.
    \item The last step collects the results of the previous three steps and establishes the limit of the excess risk (with respect to  the $L_1$ distance).
\end{enumerate}
 The four steps are given  in the following subsections, respectively.  A few technical lemmas and propositions used in these steps are stated without proofs; these proofs are deferred to appendix later. Before proceeding further, we remind the reader  the following notations:  $\Xb=[\xb_1,...,\xb_n]^\T \in\RR^{n\times d}$ with $(\xb_i)_{i\in[n]}\sim \Unif(\sqrt{d}\cdot\SSS^{d-1})$, $\yb=[y_1,...,y_n]^\T $,  $\bTheta=[\bTheta_1^\T ,\bTheta_2^\T ]^\T =[\btheta_1,...,\btheta_N]^\T \in\RR^{N\times d}$  with $(\btheta_i)_{i\in[N]}\sim\Unif(\sqrt{d}\cdot\SSS^{d-1})$. 
Some new  notations are given in the following definition. 
\begin{definition}
\label{def:notation1}
Define
\begin{equation*}
    \begin{split}
       &\Zb _j=\sigma_1\left(\Xb\bTheta_j^\T /\sqrt d\right)/\sqrt d\in\RR^{n\times N_j},~ j=1,2,\quad 
       \Zb =\left( \Zb _1, \Zb _2\right)\in\RR^{n\times N}, \\
      & \bUpsilon=( \Zb ^\T  \Zb +\lambda\Ib_N)^{-1}; \quad  \bsigma(\xb)=\big(\sigma_1(\xb^\T \bTheta_1^\T /\sqrt{d}),\sigma_2(\xb^\T \bTheta_2^\T /\sqrt{d})\big)^\T \in\RR^N;\\
      &\Mb_1=\diag\big(\mu_{1,1}\Ib_{N_1},\mu_{2,1}\Ib_{N_2}\big),~\quad  \Mb_2=\diag\big(\mu_{1,2}\Ib_{N_1},\mu_{2,2}\Ib_{N_2}\big).  
    \end{split}
\end{equation*}
Furthermore, for any matrix $\Wb\in\RR^{N\times N}$, we define a bracket $  [\Wb]_{\Zb}\triangleq \Zb\bUpsilon{\Wb}\bUpsilon  \Zb ^\T$.
\end{definition}

\subsection{Step 1:  bias-variance decomposition of the excess risk} 
\label{step1}

By the definition of $\hat{\ab}$ in \eqref{eq:def_hata}, we have
\begin{align}\label{eq:def_hata1}
    \hat{\ab}=&\argmin_{\ab}\left\{\frac1n\sum_{j=1}^n\Big(y_j-f(\xb_j;\ab,\bTheta)\Big)^2+\frac{d}{n}\lambda\lVert\ab\rVert_2^2 \right\} =\frac{1}{\sqrt{d}}\bUpsilon \Zb ^\T \yb.
\end{align}
The excess risk is then of the form
\begin{align*}
    R_d(\Xb,\bTheta,\lambda,\bbeta_d,\bvarepsilon)=&\EE_{\xb}\big[\xb^\T \bbeta_{1,d}+F_0-\hat{\ab}^\T \bsigma(\xb)\big]^2.
\end{align*}
The goal of Theorem~\ref{thm:mainthm} is to calculate this risk.  
One of the major challenges in this calculation is the nonlinearities of the activation functions. To overcome this challenge, we introduce a decomposition of the risk in the proposition below. We remind readers that $F_{1,d}=\lVert\bbeta_{1,d}\rVert_2$.



\begin{proposition}
\label{prop:decomp-expectation}
For any $\lambda > 0$, let
\begin{equation*}
    \begin{split}
        \overline{R}_d(\Xb,\bTheta,\lambda,F_{1,d},\tau)=F_{1,d}^2-\frac{2F_{1,d}^2}{d}\tr \bigg(\Mb_1\frac{\bTheta\Xb^\T }{d}\Zb\bUpsilon\bigg)+\frac{F_{1,d}^2}{d}\tr\bigg( \big[\tilde{\Ub}\big]_\Zb \frac{\Xb\Xb^\T }{d}\bigg)+\frac{\tau^2}{d}\tr(\big[\tilde{\Ub}\big]_\Zb ),
    \end{split}
\end{equation*}
where $\tilde{\Ub}=\Mb_1 \bTheta\bTheta^\T \Mb_1 / d + \Mb_2\Mb_2$. Then under the same conditions as Theorem~\ref{thm:mainthm},
\begin{align*}
    \begin{split}
        \EE_{\Xb,\bTheta,\bvarepsilon}&\Big|R_d(\Xb,\bTheta,\lambda,\bbeta_d,\bvarepsilon)-\overline{R}_d(\Xb,\bTheta,\lambda,F_{1,d},\tau)\Big|=o_d(1).
    \end{split}
\end{align*}
\end{proposition}
The proof of Proposition~\ref{prop:decomp-expectation} is given in Section~\ref{sec:appendixDecomposition}. It presents the bias-variance decomposition as the sum of four terms: \xr{the first three terms with $F_{1,d}^2$ together} give the bias in the asymptotic excess risk, while the last term \xr{with $\tau^2$} is the variance.

\subsection{Step 2: approximation of the risk decomposition via a linear pencil matrix} \label{step2}
 The approximating function  $\overline{R}_d(\Xb,\bTheta,\lambda,F_{1,d},\tau)$ found in 
 Proposition~\ref{prop:decomp-expectation} depends on three traces of certain random matrices. In this step, we calculate these traces via a special random matrix, namely the linear pencil matrix defined as follows.  
 \begin{definition}
\label{def:linear pencil}
(1) Let 
\[\mathcal{Q}:=\{\qb=[q_1,q_2,q_3,q_4,q_5]\in\RR_+^5:~ q_4,q_5\leq(1+q_1)/2,~ \| \qb \|_2 \leq 1\}.
\]
Depending on $\qb\in\mathcal{Q}$ and $\bmu$, 
the linear pencil matrix $\Ab(\qb,\bmu)$ is 
\begin{equation*}
\Ab(\qb,\bmu)=\begin{bmatrix}
q_2\mu_{1,2}^2\Ib_{N_1}+q_4\mu_{1,1}^2\frac{\bTheta_1\bTheta_1^\T }{d}&q_4\mu_{1,1}\mu_{2,1}\frac{\bTheta_1\bTheta_2^\T }{d}& \Zb_1^\T +q_1\tilde{\Zb}_1^\T \\
q_4\mu_{1,1}\mu_{2,1}\frac{\bTheta_2\bTheta_1^\T }{d}&q_2\mu_{2,2}^2\Ib_{N_2}+q_4\mu_{2,1}^2\frac{\bTheta_2\bTheta_2^\T }{d}& \Zb_2^\T +q_1\tilde{\Zb}_2^\T \\
 \Zb_1+q_1\tilde{\Zb}_1& \Zb_2+q_1\tilde{ \Zb}_2&q_3\Ib_n+q_5\frac{\Xb\Xb^\T }{d}
\end{bmatrix}\in\RR^{\PO\times \PO},
\end{equation*}
where $\PO=N+n$, and  $\tilde{\Zb}_j=\frac{\mu_{j,1}}{d}\Xb\bTheta_j^\T $ for $j=1,2$.

(2)  The Stieltjes transform of the empirical eigenvalue  distribution of $\Ab=\Ab(\qb,\bmu)$ (up to the factor $P/d$)  is
\begin{equation*}
    \begin{split}
    M_d(\xi;\qb,\bmu)=\frac{1}{d}\tr\big[(\Ab-\xi\Ib_{\PO})^{-1}\big], \quad \xi\in\bbC_+, 
    \end{split}
\end{equation*}
and its logarithmic potential is  
\begin{equation*}
\begin{split}
    G_d(\xi;\qb,\bmu)=\frac1d\log\det (\Ab-\xi\Ib_{\PO})=\frac1d  \sum_{i=1}^{\PO}\Log(\lambda_i(\Ab)-\xi),\quad \xi\in\bbC_+.
\end{split}
\end{equation*}
Here $\lambda_1(\Ab)\ge \cdots\ge \lambda_P(\Ab)$ are   the  eigenvalues of $\Ab$, and
$\Log(z):=\Log(|z|)+\rmi\arg(z)$, for $z\in\bbC$, $-\pi<\arg(z)\leq\pi$
is the principal value of a complex logarithmic function. 
\end{definition}
We assume that $\qb\in\cQ$ throughout the paper. 
The three traces in the definition of $\overline{R}_d(\Xb,\bTheta,\lambda,F_{1,d},\tau)$ in Proposition~\ref{prop:decomp-expectation} are now expressed as partial derivatives of the logarithmic potential $G_d$ as shown in the proposition below. 
\begin{proposition}
\label{prop:connectGd}
Let  $\tilde{\Ub}$ be defined in Proposition~\ref{prop:decomp-expectation}. Then we have  
\begin{equation*}
    \begin{split}
        \frac1d\tr \bigg(\Mb_1\frac{\bTheta\Xb^\T }{d}\Zb\bUpsilon\bigg)&=\frac{1}{2}\partial_{q_1} G_d(\xi^*;\qb,\bmu) |_{\qb = \mathbf{0}},\\
        \frac{1}{d}\tr\bigg( \big[\tilde{\Ub}\big]_\Zb \frac{\Xb\Xb^\T }{d}\bigg)&=-\partial^2_{q_4,q_5}G_d(\xi^*;\qb,\bmu)|_{\qb = \mathbf{0}}-\partial^2_{q_5,q_2}G_d(\xi^*;\qb,\bmu)|_{\qb = \mathbf{0}},\\
        \frac{1}{d}\tr(\big[\tilde{\Ub}\big]_\Zb )&=-\partial^2_{q_4,q_3}G_d(\xi^*;\qb,\bmu)|_{\qb = \mathbf{0}}-\partial^2_{q_2,q_3}G_d(\xi^*;\qb,\bmu)|_{\qb = \mathbf{0}}.
    \end{split}
\end{equation*}
\end{proposition}
We remind readers that $\xi^*=\sqrt{\lambda}\cdot\rmi$. The proof of Proposition~\ref{prop:connectGd} is given in Section~\ref{sec:appendixconnectedGd}.


\subsection{Step 3: key limiting spectral functions of the linear pencil matrix}
\label{step3}
Proposition \ref{prop:connectGd} shows that the excess risk can be calculated based on  $G_d(\xi^*;\qb,\bmu)$. Moreover, by Definition~\ref{def:linear pencil}, we have $\frac{d}{d\xi}G_d(\xi;\qb,\bmu)=-M_d(\xi;\qb,\bmu)$, which shows that $G_d(\xi;\qb,\bmu)$ is related to $M_d(\xi;\qb,\bmu)$. Therefore, we study the Stieltjes transform $ M_d(\xi;\qb,\bmu)$ and calculate its limit as $d,n,N \rightarrow \infty$. To do so,
we define the following system of equations.

\begin{definition}
\label{def:implicit1}
For $\xi\in\bbC_+$, define a function $\sfFb(\cdot;\xi,\qb,\bmu)$ from  $\bbC^3$ to  $\bbC^3$ by 
\begin{equation*}
    \mb=[m_1,m_2,m_3] \longmapsto \sfFb(\mb;\xi,\qb,\bmu)  
    = \begin{bmatrix}
    \psi_1\Big\{-\xi+q_2\mu_{1,2}^2-\mu_{1,2}^2m_{3}+\frac{H_{1}}{H_{D}}\Big\}^{-1}\\
    \psi_2\Big\{-\xi+q_2\mu_{2,2}^2-\mu_{2,2}^2m_{3}+\frac{H_{2}}{H_{D}}\Big\}^{-1}\\
    \psi_3\Big\{-\xi+q_3-\mu_{1,2}^2m_{1}-\mu_{2,2}^2m_{2}+\frac{H_{3}}{H_{D}}\Big\}^{-1}
    \end{bmatrix},
\end{equation*}
where
$$\begin{aligned}
H_{1}=&\mu_{1,1}^2q_4(1+m_{3}q_5)-\mu_{1,1}^2(1+q_1)^2m_{3},\\
H_{2}=&\mu_{2,1}^2q_4(1+m_{3}q_5)-\mu_{2,1}^2(1+q_1)^2m_{3},\\
H_{3}=&q_5(1+\mu_{1,1}^2m_{1}q_4+\mu_{2,1}^2m_{2}q_4)-\mu_{2,1}^2(1+q_1)^2m_{2}-\mu_{1,1}^2(1+q_1)^2m_{1},\\
H_{D}=&(1+\mu_{1,1}^2m_{1}q_4+\mu_{2,1}^2m_{2}q_4)(1+m_{3}q_5)-\mu_{2,1}^2(1+q_1)^2m_{2}m_{3}-\mu_{1,1}^2(1+q_1)^2m_{1}m_{3}.
\end{aligned}$$
We write the three coordinates of $\sfFb$ as $\sfFb(\mb;\xi,\qb,\bmu) = [\sfF_1 , \sfF_2 , \sfF_3 ]^\T (\mb;\xi,\qb,\bmu)$.
\end{definition}

We give in Section~\ref{sec:appendixuniquesolution} some properties of the function $\sfFb$. 
In particular, we show that there exists a constant $\xi_0 > 0$, such that for all $\xi$ with $\Im(\xi) > \xi_0$  and $\qb\in\cQ$, $\sfFb(\cdot;\xi, \qb,\bmu)$ has a unique fixed point $\mb(\xi;\qb,\bmu)=[m_1,m_2,m_3]^\T(\xi;\qb,\bmu)$ satisfying $| m_j(\xi) |\leq2\psi_j/\xi_0$
for $j=1,2,3$. 
Note that this fixed point result only defines  $\mb(\xi;\qb,\bmu)$ on $\{\xi: \Im(\xi) > \xi_0\}$.
To extend its definition to $\bbC_+$, we aim to show that $\mb$ is an analytic function on $\{\xi: \Im(\xi) > \xi_0\}$, and its analytic continuation to $\bbC_+$ is still a fixed point of $\sfFb(\cdot;\xi, \qb,\bmu)$, i.e.,
\begin{equation} 
\label{eq:m123}
\mb(\xi;\qb,\bmu)\equiv\sfFb[\mb(\xi;\qb,\bmu);\xi,\qb,\bmu] 
\end{equation}
for all $\xi\in \bbC_+$. More importantly, 
by using 
random matrix theory, we also aim to show that the limiting spectral distribution (LSD) of the matrix $\Ab$ exists and its Stieltjes transform is   
\begin{equation*}
m(\xi;\qb,\bmu)= \sum\limits_{i=1}^3m_i(\xi;\qb,\bmu).
\end{equation*} 
These results are formally given in the following proposition.  





\begin{proposition}
\label{prop:implicit1}
Under Assumptions~\ref{assump1} and \ref{assump2},
$\mb(\xi;\qb,\bmu)$ is analytic on $\{\xi:\Im(\xi)>\xi_0\}$, and has a unique analytic continuation to $\bbC_+$. Moreover, this analytic continuation (still denoted as $\mb(\xi;\qb,\bmu)$) satisfies the following properties:
\begin{enumerate}[leftmargin = *]
    \item $\mb(\xi;\qb,\bmu) \in \bbC_+^3$ for all $\xi\in \bbC_+$.
    \item $\mb(\xi,\qb,\bmu)\equiv\sfFb[\mb(\xi,\qb,\bmu);\xi,\qb,\bmu]$ for all $\xi\in \bbC_+$.
    \item Let $M_d(\xi;\qb,\bmu)$ be defined in Definition~\ref{def:linear pencil}. Then for any compact set $\Omega\subset\bbC_+$,  \begin{equation*}
    \lim_{d\rightarrow+\infty}\EE\bigg[\sup_{\xi\in\Omega}\big|M_d(\xi;\qb,\bmu)-m(\xi;\qb,\bmu)\big|\bigg]=0.
\end{equation*} 
\end{enumerate}
\end{proposition}


The proof of Proposition~\ref{prop:implicit1} is given in Section~\ref{sec:appendixpropimplicit1}. It shows that $M_d(\xi;\qb,\bmu)$ has a deterministic limit equal to  $m(\xi;\qb,\bmu)$. This result, together with the connection between $M_d(\xi;\qb,\bmu)$ and the logarithmic potential $G_d(\xi;\qb,\bmu)$ in Definition~\ref{def:linear pencil}, further indicates that 
$G_d$ may also have a deterministic limit, and its deterministic limit can possibly be expressed as a function of $\mb(\xi;\qb,\bmu)$. In fact, this  limit is found to be 
\begin{equation}
\label{def:gd}
    g(\xi;\qb,\bmu)\triangleq L(\xi,m_1(\xi;\qb,\bmu),m_2(\xi;\qb,\bmu),m_3(\xi;\qb,\bmu);\qb,\bmu),
\end{equation}
where
\begin{equation}\label{eq:L-function}
    \begin{split}
        & L(\xi,z_1,z_2,z_3;\qb,\bmu)  \triangleq \\
        &~ \Log\big[(1+\mu_{1,1}^2z_1q_4+\mu_{2,1}^2z_2q_4)(1+z_3q_5)-\mu_{1,1}^2(1+q_1)^2z_1z_3-\mu_{2,1}^2(1+q_1)^2z_2z_3\big]\\
&~ -\mu_{1,2}^2z_1z_3-\mu_{2,2}^2z_2z_3+q_2\mu_{1,2}^2z_1+q_2\mu_{2,2}^2z_2+q_3z_3-\xi(z_1+z_2+z_3)\\
&~ -\psi_1\Log(z_1/\psi_1)-\psi_2\Log(z_2/\psi_2)-\psi_3\Log(z_3/\psi_3)-\psi_1-\psi_2-\psi_3.
    \end{split}
\end{equation}
The following proposition formally shows that $g(\xi;\qb,\bmu)$ and its partial derivatives are the deterministic limit of the $G_d$ and the partial derivatives of $G_d$, respectively. 
\begin{proposition}
\label{prop:substitutiongd}
Let $G_d(\xi;\qb,\bmu)$ be defined in Definition~\ref{def:linear pencil} and  $g(\xi;\qb,\bmu)$  defined in \eqref{def:gd}. Then 
for any fixed $\xi\in\bbC_+$, $\qb\in\mathcal{Q}$ and $u\in\RR_+$,
$$\begin{aligned}
&\lim_{d\rightarrow+\infty}\EE[|G_d(\xi;\qb,\bmu)-g(\xi;\qb,\bmu)|]=0,\\
&\lim_{d\rightarrow+\infty}\EE[\lVert \nabla_{\qb}G_d(\rmi u;\qb,\bmu)|_{\qb=\0}-\nabla_{\qb}g(\rmi u;\qb,\bmu)|_{\qb=\0}\rVert_2]=0,\\
&\lim_{d\rightarrow+\infty}\EE[\lVert \nabla_{\qb}^2G_d(\rmi u;\qb,\bmu)|_{\qb=\0}-\nabla_{\qb}^2g(\rmi u;\qb,\bmu)|_{\qb=\0}\rVert_{\op}]=0.\end{aligned}$$
\end{proposition}
Proposition~\ref{prop:substitutiongd} is proved in Section~\ref{sec:appendixsubstitution}.



\subsection{Step 4: completion of the proof} 
\label{step4}

According to Propositions~\ref{prop:decomp-expectation}, \ref{prop:connectGd}, and \ref{prop:substitutiongd}, the key terms in the excess risk can be calculated as the partial derivatives of the function $g(\xi;\qb,\bmu)$ at $\qb = \0$. However, $g(\xi;\qb,\bmu)$ is   based on $\mb(\xi;\qb,\bmu)$, and the calculation of the partial derivatives of $g(\xi;\qb,\bmu)$ is non-trivial: $\mb(\xi;\qb,\bmu)$ is originally defined on $\{\xi: \Im(\xi) > \xi_0\}$ as the fixed point of $\sfFb$, and its definition is then extended to $\bbC_+$ in Proposition~\ref{prop:implicit1}. 
To finalize the proof, we first present the following proposition relating $\mb(\xi;\qb,\bmu)$ to the function $\nu(\xi;\bmu)$ defined in Section~\ref{sec:mainresult2}. 





\begin{proposition}\label{prop:existence_uniqueness_nu}
There exists a unique analytic function $\bnu= [\nu_1,\nu_2,\nu_3]^\T: \bbC_+ \rightarrow \bbC_+^3$ such that:
\begin{enumerate}[leftmargin = *]
    \item For any $\xi\in\bbC_+$, $\bnu(\xi;\bmu)$ is a solution to $\bnu$\textit{-system}~\eqref{eq:implicit2}.
    \item There exists $\xi_0>0$, such that $| \nu_j(\xi;\bmu) | \leq 2\psi_j / \xi_0$, for all $\xi$ with $\Im(\xi) \geq \xi_0$ and $j=1,2,3$. Moreover, it holds that $ \bnu(\xi;\bmu) = \mb(\xi;\0,\bmu) $ for all $\xi\in\bbC_+$.
    \item $\bnu^*=\bnu(\sqrt{\lambda}\cdot\rmi;\bmu)$ in Definition~\ref{def:mainresultsdefinition} satisfies \xr{$\nu_j^*=b_j\cdot\rmi$ with $b_j>0$ for all $j=1,2,3$}. 
\end{enumerate}

\end{proposition}

The proof of Proposition~\ref{prop:existence_uniqueness_nu} is given in Section~\ref{sec:appendix_nu}. The proposition thus justifies the definition of $\bnu(\xi;\bmu)$ in Section~\ref{sec:mainresult2} by demonstrating its existence and uniqueness. Moreover, it also relates $\bnu(\xi;\bmu)$ to the function $\mb(\xi;\qb,\bmu)$ introduced in step 3 of the proof. With this result, we can finalize the proof of
Theorem~\ref{thm:mainthm} as follows.


\begin{proof}[Proof of Theorem~\ref{thm:mainthm}]
Let
\begin{align}
     \cR(\lambda,\bpsi,\bmu,F_1,\tau) &=
        F_1^2\cdot \big[ 1-\partial_{q_1} g(\xi^*;\qb,\bmu)  -\partial^2_{q_4,q_5}g(\xi^*;\qb,\bmu)  -\partial^2_{q_2,q_5}g(\xi^*;\qb,\bmu) \big ] \big|_{\qb = \mathbf{0}} \nonumber\\
        &\quad   -\tau^2\cdot \big[\partial^2_{q_3,q_4}g(\xi^*;\qb,\bmu) + \partial^2_{q_2,q_3}g(\xi^*;\qb,\bmu) \big ] \big|_{\qb = \mathbf{0}},\label{eq:finalgd}
\end{align}
where $g$ is defined in \eqref{def:gd}, and $\xi^*=\sqrt{\lambda}\cdot\rmi$ is given in Definition~\ref{def:notation1}. 
Then by Propositions~\ref{prop:decomp-expectation}, \ref{prop:connectGd} and \ref{prop:substitutiongd}, we have 
\begin{align*}
    \begin{split}
        \EE_{\Xb,\bTheta,\bvarepsilon}&\Big|R_d(\Xb,\bTheta,\lambda,\bbeta_d,\bvarepsilon)-\cR(\lambda,\bpsi,\bmu,F_1,\tau)\Big|=o_d(1).
    \end{split}
\end{align*}
Therefore to complete the proof, it suffices to calculate the partial derivative terms of $g(\xi^*;\qb,\bmu)$ at $\qb = \mathbf{0}$. For this calculation, we first note that by the definition of $L(\xi,\zb;\qb,\bmu)$ in \eqref{eq:L-function} and the definition of $\mb$ in \eqref{eq:m123} as the fixed point of $\sfFb(\cdot;\xi,\qb,\bmu)$, we have that
\begin{align}\label{eq:L_property}
    \nabla_{\zb}L(\xi,\zb;\qb,\bmu)|_{\zb=\mb} \equiv\0.
\end{align}
Readers can refer to Lemma~\ref{lemma:q_derivative_m} and its proof for the detailed derivation of \eqref{eq:L_property}. 
Let $ \mb^*(\qb,\bmu) = [ m_1(\xi^*;\qb,\bmu),m_2(\xi^*;\qb,\bmu),m_3(\xi^*;\qb,\bmu) ]^\T $. Then by Proposition~\ref{prop:existence_uniqueness_nu}, we have $\bnu^*=\mb^*(\0,\bmu)$. Therefore,
\begin{align}
    \partial_{q_1} g(\xi^*;\qb,\bmu)\big|_{\qb=\0} &= \partial_{q_1} \big[  L(\xi^*, \mb^*(\qb,\bmu) ;\qb,\bmu) \big] \big|_{\qb=\0} \nonumber\\
    & = \big[ \big\la  \nabla_{\zb} L(\xi^*,\zb;\qb,\bmu) |_{ \zb = \mb^*} , \partial_{q_1} \mb^* \big\ra + \partial_{q_1} L(\xi^*,\zb;\qb,\bmu) |_{\zb = \mb^*} \big] \big|_{\qb = \0}\nonumber\\
    & = 0 + \partial_{q_1} L(\xi^*,\zb;\qb,\bmu) |_{\qb = \0,\zb = \bnu^*}=\frac{2\nu_3^*M_N}{M_D},\label{eq:gradient_g1}
\end{align}
where the first equality is by the definition of $g$, the second equality follows by the chain rule, the third equality follows by \eqref{eq:L_property}, and the last equality is by direct calculation and the definition that $M_N=\nu_1^*\mu_{1,1}^2+\nu_2^*\mu_{2,1}^2$, $M_D=\nu_3^*M_N-1$.

For the second order derivatives, let $q_i$ $q_j$ be the $i^{\text{th}}$ and $j^{\text{th}}$ element in $\qb$ for $i,j=2,3,4,5$. Then by \eqref{eq:L_property}, with similar calculation as \eqref{eq:gradient_g1}, we have
\begin{align}\label{eq:derivative_g2}
        \frac{\partial^2g(\xi^*;\qb,\bmu)}{\partial q_i\partial q_j}=&\frac{\partial^2L(\xi^*,\zb;\qb,\bmu)}{\partial q_i\partial q_j}\bigg|_{\zb=\mb^*}+\bigg\la\nabla_{\zb}\bigg[\frac{\partial L(\xi^*,\zb;\qb,\bmu)}{\partial q_i}\bigg]\bigg|_{\zb=\mb^*},\frac{\partial \mb^*}{\partial q_i}\bigg\ra.
\end{align}
Moreover, by \eqref{eq:L_property} and the formula for implicit differentiation, we have
\begin{align}\label{eq:m_derivative}
    \frac{\partial \mb^*}{\partial q_i}=-\big[\big(\nabla_{\zb}^2L(\xi^*,\zb;\qb,\bmu)\big)\big|_{\zb=\mb^*}\big]^{-1}\frac{\partial\big[\nabla_\zb L(\xi^*,\zb;\qb,\bmu)\big]}{\partial q_i}\bigg|_{\zb=\mb^*}.
\end{align}
In addition, we let $\ub=[q_2,q_3,q_4,q_5,z_1,z_2,z_3]^\T$, and define the symmetric matrix  
\begin{eqnarray}\label{eq:Wmatrix_def}
        \Wb & = & \Wb(\bnu^*,\bmu)=\nabla_{\ub}^2 L(\xi,\zb;\qb,\bmu)|_{\zb=\bnu^*,\qb=0} \nonumber \\
        &  = & \begin{bmatrix}
         ~~0~~&~~0~~&~~0~~&~~0~~&\mu_{1,2}^2&\mu_{2,2}^2&0\\
          * &0&0&0&0&0&1\\
         * & * &-\frac{M_N^2}{M_D^2}&-\frac{\nu_3^2M_N^2}{M_D^2}&\frac{\mu_{1,1}^2}{M_D^2}&\frac{\mu_{2,1}^2}{M_D^2}&\frac{M_N^2}{M_D^2}\\
        * & * & * &-\frac{\nu_3^2}{M_D^2}&\frac{\nu_3^2\mu_{1,1}^2}{M_D^2}&\frac{\nu_3^2\mu_{2,1}^2}{M_D^2}&\frac{1}{M_D^2}\\
        *  & * &  * & * &- \frac{\nu_3^2\mu_{1,1}^4}{M_D^2}+\frac{\psi_1}{\nu_1^2}&-\frac{\nu_3^2\mu_{1,1}^2\mu_{2,1}^2}{M_D^2}&-\frac{\mu_{1,1}^2}{M_D^2}-\mu_{1,2}^2\\
        * & * &  * &  * &  * &-\frac{\nu_3^2\mu_{2,1}^4}{M_D^2}+\frac{\psi_2}{\nu_2^2}~&-\frac{\mu_{2,1}^2}{M_D^2}-\mu_{2,2}^2~\\
        * & * &  * & * &  * &  * &-\frac{M_N^2}{M_D^2}+\frac{\psi_3}{\nu_3^2}
        \end{bmatrix}.
\end{eqnarray}
\if UT
\begin{equation}
    \begin{split}
        \Wb(\mb,\bmu)=\begin{bmatrix}
        0&0&0&0&\mu_{1,2}^2&\mu_{2,2}^2&0\\
        0&0&0&0&0&0&1\\
        0&0&-\frac{M_N^2}{M_D^2}&-\frac{m_3^2M_N^2}{M_D^2}&\frac{\mu_{1,1}^2}{M_D^2}&\frac{\mu_{2,1}^2}{M_D^2}&\frac{M_N^2}{M_D^2}\\
        0&0&-\frac{m_3^2M_N^2}{M_D^2}&-\frac{m_3^2}{M_D^2}&\frac{m_3^2\mu_{1,1}^2}{M_D^2}&\frac{m_3^2\mu_{2,1}^2}{M_D^2}&\frac{1}{M_D^2}\\
        \mu_{1,2}^2&0&\frac{\mu_{1,1}^2}{M_D^2}&\frac{m_3^2\mu_{1,1}^2}{M_D^2}&-\frac{m_3^2\mu_{1,1}^4}{M_D^2}+\frac{\psi_1}{m_1^2}&-\frac{m_3^2\mu_{1,1}^2\mu_{2,1}^2}{M_D^2}&-\frac{\mu_{1,1}^2}{M_D^2}-\mu_{1,2}^2\\
        \mu_{2,2}^2&0&\frac{\mu_{2,1}^2}{M_D^2}&\frac{m_3^2\mu_{2,1}^2}{M_D^2}&-\frac{m_3^2\mu_{1,1}^2\mu_{2,1}^2}{M_D^2}&-\frac{m_3^2\mu_{2,1}^4}{M_D^2}+\frac{\psi_2}{m_2^2}&-\frac{\mu_{2,1}^2}{M_D^2}-\mu_{2,2}^2\\
        0&1&\frac{M_N^2}{M_D^2}&\frac{1}{M_D^2}&-\frac{\mu_{1,1}^2}{M_D^2}-\mu_{1,2}^2&-\frac{\mu_{2,1}^2}{M_D^2}-\mu_{2,2}^2&-\frac{M_N^2}{M_D^2}+\frac{\psi_3}{m_3^2}
        \end{bmatrix}.
    \end{split}
\end{equation}
\fi
Then by \eqref{eq:derivative_g2}, \eqref{eq:m_derivative} and \eqref{eq:Wmatrix_def}, we have
\begin{align}
    &\frac{\partial^2 g(\xi^*;\qb,\bmu)}{\partial q_2 \partial q_5}\Big|_{\qb=0}=\Wb_{1,4}-\Wb_{1,[5:7]}\Big(\Wb_{[5:7],[5:7]}\Big)^{-1}\Wb_{[5:7],4},\label{eq:2nd_derivative_g1}\\
        &\frac{\partial^2 g(\xi^*;\qb,\bmu)}{\partial q_3\partial q_4}\Big|_{\qb=0}=\Wb_{2,3}-\Wb_{2,[5:7]}\Big(\Wb_{[5:7],[5:7]}\Big)^{-1}\Wb_{[5:7],3},\label{eq:2nd_derivative_g2}\\
        &\frac{\partial^2 g(\xi^*;\qb,\bmu)}{\partial q_2\partial q_3}\Big|_{\qb=0}=\Wb_{1,2}-\Wb_{1,[5:7]}\Big(\Wb_{[5:7],[5:7]}\Big)^{-1}\Wb_{[5:7],2},\label{eq:2nd_derivative_g3}\\
        &\frac{\partial^2 g(\xi^*;\qb,\bmu)}{\partial q_4\partial q_5}\Big|_{\qb=0}=\Wb_{3,4}-\Wb_{3,[5:7]}\Big(\Wb_{[5:7],[5:7]}\Big)^{-1}\Wb_{[5:7],4}.\label{eq:2nd_derivative_g4}
\end{align}
Now the terms on the right hand side above can be directly calculated: (recalling $\Vb, \Hb$ given in Definition~\ref{def:mainresultsdefinition}, \xr{and $\nu_j^*$ is the solution of $\bnu$\textit{-system}~\eqref{eq:implicit2} given $\xi=\sqrt{\lambda}\cdot \rmi$}) we have
\begin{align*}
   & \Wb_{1,4}=\Wb_{2,3}=\Wb_{1,2}=0,   \quad  \Wb_{3,4}=-\frac{\nu_3^{*2}M_N^2}{M_D^2}, \\
   &  \Wb_{[5:7],[1:4]}=\Wb_{[1:4],[5:7]}^\T = \Vb,  \quad \text{and}  \quad \Wb_{[5:7],[5:7]} = \Hb. 
\end{align*}
Plugging \eqref{eq:gradient_g1} and \eqref{eq:2nd_derivative_g1}-\eqref{eq:2nd_derivative_g4} into \eqref{eq:finalgd} completes the proof of
Theorem~\ref{thm:mainthm}.
\end{proof}
Finally, recall $\Lb= \Vb^\T \Hb^{-1}\Vb$, we give the closed form expression for the terms $\Lb_{1,4}, \Lb_{2,3}, \Lb_{1,2}, \Lb_{3,4}$ in Theorem~\ref{thm:mainthm}. Let  $[\nu_1^*,\nu_2^*,\nu_3^*]$, $M_N$ and $M_D$ be defined in Definition~\ref{def:mainresultsdefinition}, and
\begin{equation}
\label{eq:explicit_S}
\begin{split}
    S=&\nu_3^{*4}\Big(\nu_2^{*2}M_N^2\mu_{2,1}^4\psi_1+\nu_1^{*2}M_N^2\mu_{1,1}^4\psi_2+\nu_1^{*2}\nu_2^{*2}M_D^2\big(\mu_{1,2}^2\mu_{2,1}^2-\mu_{1,1}^2\mu_{2,2}^2\big)^2\Big)\\
      &-\nu_3^{*2}\nu_2^{*2}\psi_1\big(2M_D^2\mu_{2,1}^2\mu_{2,2}^2+M_D^4\mu_{2,2}^4+\mu_{2,1}^4(1+M_D^2\psi_3)\big)\\
      &-\nu_3^{*2}\nu_1^{*2}\psi_2\big(2M_D^2\mu_{1,1}^2\mu_{1,2}^2+M_D^4\mu_{1,2}^4+\mu_{1,1}^4(1+M_D^2\psi_3)\big)\\
      & -\nu_3^{*2}\psi_1\psi_2M_D^2M_N^2+M_D^4\psi_1\psi_2\psi_3.
      \end{split}
\end{equation}
Then by direct calculation, the terms $\Lb_{1,4}, \Lb_{2,3}, \Lb_{1,2}, \Lb_{3,4}$ satisfy the following equations:
\begin{equation}
\label{eq:explicit_L}
\begin{split}
    \frac{S\cdot\Lb_{1,4}}{\nu_3^{*2}}=&-\nu_3^{*2}M_N^2\Big(\nu_2^{*2}\mu_{2,1}^2\mu_{2,2}^2\psi_1+\nu_1^{*2}\mu_{1,1}^2\mu_{1,2}^2\psi_2\Big)\\
    &+\nu_1^{*2}\mu_{1,2}^2\psi_2\Big(M_D^2\mu_{1,2}^2+\mu_{1,1}^2\big(1+M_D^2\psi_3\big)\Big)\\
      &+\nu_2^{*2}\mu_{2,2}^2\psi_1\Big(M_D^2\mu_{2,2}^2+\mu_{2,1}^2\big(1+M_D^2\psi_3\big)\Big),\\
    \frac{S\cdot\Lb_{2,3}}{\nu_3^{*2}}=&\nu_2^{*2}\mu_{2,1}^2\big(\mu_{2,1}^2+M_D^2\mu_{2,2}^2\big)\psi_1+\nu_1^{*2}\mu_{1,1}^2\big(\mu_{1,1}^2+M_D^2\mu_{1,2}^2\big)\psi_2\\
      & -\nu_3^{*2}M_N^2\big(\nu_2^{*2}\mu_{2,1}^4\psi_1+\nu_1^{*2}\mu_{1,1}^4\psi_2\big)  +M_D^2M_N^2\psi_1\psi_2,\\
    \frac{S\cdot\Lb_{1,2}}{\nu_3^{*2}}=&M_D^2\Big(\nu_2^{*2}\mu_{2,2}^2\big(\mu_{2,1}^2+M_D^2\mu_{2,2}^2\big)\psi_1+\nu_1^{*2}\mu_{1,2}^2\big(\mu_{1,1}^2+M_D^2\mu_{1,2}^2\big)\psi_2\\
      &  -\nu_1^{*2}\nu_2^{*2}\nu_3^{*2}\big(\mu_{1,2}^2\mu_{2,1}^2-\mu_{1,1}^2\mu_{2,2}^2\big)^2\Big),\\
    \frac{M_D^2S\cdot\Lb_{3,4}}{\nu_3^{*2}}=&\nu_3^{*2}\Big(\nu_2^{*2}M_N^2\mu_{2,1}^2\big(M_D^2\mu_{2,2}^2-\mu_{2,1}^2\big)\psi_1+\nu_1^{*2}M_N^2\mu_{1,1}^2\big(M_D^2\mu_{1,2}^2-\mu_{1,1}^2\big)\psi_2\Big)\\
    & +\psi_1\psi_2M_D^2M_N^2   -\nu_1^{*2}\nu_2^{*2}\nu_3^{*2}M_D^2\big(\mu_{1,2}^2\mu_{2,1}^2-\mu_{1,1}^2\mu_{2,2}^2\big)^2\\
        &+\nu_2^{*2}\mu_{2,1}^2\psi_1\big(M_D^2\mu_{2,2}^2+\mu_{2,1}^2+M_D^2\mu_{2,1}^2\psi_3\big)\\
        &+\nu_1^{*2}\mu_{1,1}^2\psi_2\big(M_D^2\mu_{1,2}^2+\mu_{1,1}^2+M_D^2\mu_{1,1}^2\psi_3\big).
      \end{split}
\end{equation}
Clearly, the equations above give explicit calculations of $\Lb_{1,4}, \Lb_{2,3}, \Lb_{1,2}, \Lb_{3,4}$ given the solution $[\nu_1^*,\nu_2^*,\nu_3^*]$ of the self consistent system $\bnu$\textit{-system}~\eqref{eq:implicit2}. \xr{ Readers may keep in mind that $\nu_j^{*2}$ is negative since $\nu_j^*$ is purely imaginary.} 
\subsection{Discussion on the proof of Theorem~\ref{thm:mainthm}}
\xr{ 
In this section, we briefly discuss the proof of Theorem~\ref{thm:mainthm} and highlight the novel challenges we encountered in our own proof and setting compared to \cite{mei2022generalization}. We compare the differences in the various steps of the proof to better understand the unique aspects of our extension.}

\begin{enumerate}[leftmargin = *]
    \item  \xr{The first step in our proof is to directly calculate the excess risk according to its definition, and identify key terms which require further analysis. To do so, we perform a bias-variance decomposition of the risk and find an asymptotic approximation whose main terms are expressed as traces of several random matrices, as detailed in Proposition~\ref{prop:decomp-expectation}. Compared to \cite{mei2022generalization}, our analysis on the DRFM addresses the impact of different activation functions. As shown in Lemma~\ref{lemma:decomp2}, the  terms become more complex for DRFMs, and it requires additional treatment and careful justification to prove the specific negligible terms. Furthermore, the decomposition in Proposition~\ref{prop:decomp-expectation} includes additional diagonal matrices $\Mb_1$ and $\Mb_2$, whereas in \cite{mei2022generalization}, it is only a scalar. To address this difference, we have extended several technical lemmas to accommodate the inclusion of $\Mb_1$ and $\Mb_2$. For further details, please refer to Section~\ref{sec:appendixDecomposition}}
    \item  \xr{The second step involves the construction of a new random matrix, known as the linear pencil matrix. While a similar technique is also used in \cite{mei2022generalization}, 
    the linear pencil matrix is more intricate in our setting, see Definition~\ref{def:linear pencil}, due to the greater complexity of the terms involving $\Mb_1$ and $\Mb_2$. 
    Specifically, the linear pencil matrix in our case is a 3 by 3 block matrix with a more complicated structure than the matrix proposed in \cite{mei2022generalization}.   
}
    \item \xr{The third step is a standard procedure that involves identifying the critical limiting spectral functions of the linear pencil matrix, including its Stieltjes transform and logarithmic potential. However, in our case, these calculations   differ from the reference and require additional investigation  due to the increased complexity of the linear pencil matrix and the more intricate formula of the related implicit equations.    Furthermore, these new calculations of the Stieltjes transform and logarithmic potential provide inspiration for the study of  the  multiple random feature models and we find  a  mathematical induction method to complete the study.}
\end{enumerate}
\section{Proof of Propositions~\ref{prop:property_lbd=0} and \ref{prop:property_psi}}
\label{sec:appendix_property}
In this section we present the detailed proofs of Propositions~\ref{prop:property_lbd=0} and \ref{prop:property_psi}.
We denote by $\bnu^*=\bnu(\sqrt{\lambda}\cdot\rmi;\bmu)=\mb(\sqrt{\lambda}\cdot\rmi;\0,\bmu)$, Proposition~\ref{prop:existence_uniqueness_nu} shows that the three numbers $\nu_j^*$, $j=1,2,3$, are all purely imaginary with positive imaginary parts, that is,  $\nu_j^* = \rmi\nu_j$ where  $\nu_j>0$. 
Moreover by $\bnu$\textit{-system}~\eqref{eq:implicit2}, we also have the following self-consistent equations:
    \begin{align}
\label{eq:implicit_trans}
    \left\{\begin{aligned}
    &\sqrt{\lambda}\nu_1+\mu_{1,2}^2\nu_1\nu_3+\frac{\mu_{1,1}^2\nu_1\nu_3}{1+\mu_{1,1}^2\nu_1\nu_3+\mu_{2,1}^2\nu_2\nu_3}=\psi_1,\\
    &\sqrt{\lambda}\nu_2+\mu_{2,2}^2\nu_2\nu_3+\frac{\mu_{2,1}^2\nu_2\nu_3}{1+\mu_{1,1}^2\nu_1\nu_3+\mu_{2,1}^2\nu_2\nu_3}=\psi_2,\\
    &\sqrt{\lambda}\nu_3+\mu_{1,2}^2\nu_1\nu_3+\mu_{2,2}^2\nu_2\nu_3+\frac{\mu_{1,1}^2\nu_1\nu_3+\mu_{2,1}^2\nu_2\nu_3}{1+\mu_{1,1}^2\nu_1\nu_3+\mu_{2,1}^2\nu_2\nu_3}=\psi_3.
    \end{aligned}\right.
\end{align}
The system \eqref{eq:implicit_trans} can be further rewritten    as 
\begin{align}
\label{eq:reformulate}
    \left\{\begin{aligned}
    &\lambda\nu_1\nu_3=\bigg(\psi_1-\mu_{1,2}^2\nu_1\nu_3-\frac{\mu_{1,1}^2\nu_1\nu_3}{1+\mu_{1,1}^2\nu_1\nu_3+\mu_{2,1}^2\nu_2\nu_3}\bigg)\\
    &\qquad\qquad\cdot\bigg(\psi_3-\mu_{1,2}^2\nu_1\nu_3-\mu_{2,2}^2\nu_2\nu_3-\frac{\mu_{1,1}^2\nu_1\nu_3+\mu_{2,1}^2\nu_2\nu_3}{1+\mu_{1,1}^2\nu_1\nu_3+\mu_{2,1}^2\nu_2\nu_3}\bigg),\\
    &\lambda\nu_2\nu_3=\bigg(\psi_2-\mu_{2,2}^2\nu_2\nu_3-\frac{\mu_{2,1}^2\nu_2\nu_3}{1+\mu_{1,1}^2\nu_1\nu_3+\mu_{2,1}^2\nu_2\nu_3}\bigg)\\
    &\qquad\qquad\cdot\bigg(\psi_3-\mu_{1,2}^2\nu_1\nu_3-\mu_{2,2}^2\nu_2\nu_3-\frac{\mu_{1,1}^2\nu_1\nu_3+\mu_{2,1}^2\nu_2\nu_3}{1+\mu_{1,1}^2\nu_1\nu_3+\mu_{2,1}^2\nu_2\nu_3}\bigg),\\
    &\sqrt{\lambda}(\nu_1+\nu_2-\nu_3)=\psi_1+\psi_2-\psi_3.
    \end{aligned}\right.
\end{align}
Our proofs of Propositions~\ref{prop:property_lbd=0} and \ref{prop:property_psi} mainly study the asymptotic properties of $\nu_1\nu_3$ and $\nu_2\nu_3$ based on \eqref{eq:reformulate}. Specifically, we define
\begin{align*}
    \chi_1(\bmu)=\lim_{\lambda\to0}\nu_1\nu_3,\quad \chi_2(\bmu)=\lim_{\lambda\to0}\nu_2\nu_3.
\end{align*}
Note that the existence of these limits with values in  $[0,+\infty)\cup\{+\infty\}$ 
is guaranteed by the property of Stieltjes transform, and the limit value $\chi_1(\bmu)$, $\chi_2(\bmu)$ are related with the moment vector $\bmu$. In the following proof, we drop the argument $\bmu$ in $\chi_1$, $\chi_2$ for simplicity.

\subsection{Proof of Proposition~\ref{prop:property_lbd=0}}
  We first prove the second and fourth conclusions of Proposition~\ref{prop:property_lbd=0} where the excess risk tends to infinity, and then we prove its first and third conclusions. Readers may keep in mind that when we let $\lambda\to0$, the moment vector $\bmu$ is fixed. 
  

\begin{proof}[\textbf{Second conclusion}]
 If $\psi_3=\psi_1+\psi_2$, then by \eqref{eq:reformulate} we have $\nu_1+\nu_2=\nu_3$. We first use a proof by contradiction to show that $\chi_1= \lim_{\lambda\to0}\nu_1\nu_3 >0$. It is obvious by definition that $\chi_1\geq 0$. If $\chi_1=0$, then from the first equation in  \eqref{eq:reformulate} we have
\begin{align*}
    0=\lim_{\lambda\to0}\lambda\nu_1\nu_3&=\lim_{\lambda\to0}\bigg(\psi_1-\mu_{1,2}^2\nu_1\nu_3-\frac{\mu_{1,1}^2\nu_1\nu_3}{1+\mu_{1,1}^2\nu_1\nu_3+\mu_{2,1}^2\nu_2\nu_3}\bigg)\\
    &\qquad\cdot\bigg(\psi_3-\mu_{1,2}^2\nu_1\nu_3-\mu_{2,2}^2\nu_2\nu_3-\frac{\mu_{1,1}^2\nu_1\nu_3+\mu_{2,1}^2\nu_2\nu_3}{1+\mu_{1,1}^2\nu_1\nu_3+\mu_{2,1}^2\nu_2\nu_3}\bigg),\\
    &= \psi_1\cdot \lim_{\lambda\to0}(\psi_1+\sqrt{\lambda}\nu_2)\geq\psi_1^2.
\end{align*}
This is impossible and hence we have $\chi_1>0$. Moreover, if $\chi_1=+\infty$, we have
\begin{align*}
    0=\lim_{\lambda\to0}&\lambda=\lim_{\lambda\to0}\bigg(\psi_1/(\nu_1\nu_3)-\mu_{1,2}^2-\frac{\mu_{1,1}^2}{1+\mu_{1,1}^2\nu_1\nu_3+\mu_{2,1}^2\nu_2\nu_3}\bigg)\\
    &\qquad\qquad\cdot\bigg(\psi_3-\mu_{1,2}^2\nu_1\nu_3-\mu_{2,2}^2\nu_2\nu_3-\frac{\mu_{1,1}^2\nu_1\nu_3+\mu_{2,1}^2\nu_2\nu_3}{1+\mu_{1,1}^2\nu_1\nu_3+\mu_{2,1}^2\nu_2\nu_3}\bigg)\gg 0,
\end{align*}
which is also a contradiction. Therefore
$0< \chi_1 <\infty$. Similarly we conclude that $0<\chi_2<\infty$.

Furthermore, the relation $\nu_1+\nu_2=\nu_3$ implies that $\nu_1,\nu_2<\nu_3$. Then we have  $\lim\limits_{\lambda\to0}\nu_1,\nu_2<+\infty$ and $\sqrt{\lambda}\nu_1,\sqrt{\lambda}\nu_2\to0$ when $\lambda\to0$.  Therefore \eqref{eq:implicit_trans} gives us the following equations when $\lambda\to0$:
\begin{align}
\label{eq:implicit_lbd=0}
    \left\{
    \begin{aligned}
    &\mu_{1,2}^2\chi_1+\frac{\mu_{1,1}^2\chi_1}{1+\mu_{1,1}^2\chi_1+\mu_{2,1}^2\chi_2}=\psi_1,\\
    &\mu_{2,2}^2\chi_2+\frac{\mu_{2,1}^2\chi_2}{1+\mu_{1,1}^2\chi_1+\mu_{2,1}^2\chi_2}=\psi_2.
    \end{aligned}
    \right.
\end{align}
By  \eqref{eq:implicit_lbd=0}, we can express $\psi_1$, $\psi_2$ and $\psi_3=\psi_1+\psi_2$ by $\chi_1$ and $\chi_2$. Moreover, note that when $\lambda\to0$, 
\begin{align}
\label{eq:Prop4replace}
 \nu_1^*\nu_3^*=-\chi_1,  ~ \nu_2^*\nu_3^*=-\chi_2,~ M_N\nu_3^*=-\mu_{1,1}^2\chi_1-\mu_{2,1}^2\chi_2,~ M_D=-\mu_{1,1}^2\chi_1-\mu_{2,1}^2\chi_2-1.
\end{align}
Therefore $S$ and $\Lb_{i,j}$ in \eqref{eq:explicit_S} and \eqref{eq:explicit_L} can also be expressed by $\chi_1$ and $\chi_2$ when $\lambda\to0$. With direct algebraic calculations, we obtain 
\begin{align*}
    \lim\limits_{\lambda\to0}S=0, \quad \lim\limits_{\lambda\to0}S\cdot (\Lb_{3,4}+\Lb_{1,4})\neq0,\quad
    \lim\limits_{\lambda\to0}S\cdot (\Lb_{2,3}+\Lb_{1,2})\neq0.
\end{align*}
This implies  that $\Lb_{3,4}+\Lb_{1,4}\to\infty$ and $\Lb_{2,3}+\Lb_{1,2}\to\infty$ when $\lambda\to0$. Since $\cR\geq0$, we have 
\begin{align*}
  \lim\limits_{\lambda\to0}\cR(\lambda,\bpsi,\bmu,F_1,\tau)=\lim\limits_{\lambda\to0}F_1^2\bigg(
 \frac{1}{M_D^2}+\Lb_{3,4}+\Lb_{1,4}\bigg)+\tau^2(\Lb_{2,3}+\Lb_{1,2})=+\infty.
\end{align*}
This gives the second conclusion in Proposition~\ref{prop:property_lbd=0}.
\end{proof}
\begin{proof}[\textbf{Fourth conclusion}]
If $(\psi_1+\psi_2)/\psi_3=1+\psi_2/\psi_1$, then $\psi_1=\psi_3$,  and  \eqref{eq:reformulate} gives $\sqrt{\lambda}(\nu_1+\nu_2-\nu_3)=\psi_2$. By substitution of  $\sqrt{\lambda}(\nu_1+\nu_2-\nu_3)=\psi_2$ into the second equation in \eqref{eq:implicit_trans} we obtain
\begin{align*}
    \sqrt{\lambda}\nu_3+\mu_{2,2}^2\nu_2\nu_3+\frac{\mu_{2,1}^2\nu_2\nu_3}{1+\mu_{1,1}^2\nu_1\nu_3+\mu_{2,1}^2\nu_2\nu_3}=\sqrt{\lambda}\nu_1.
\end{align*}
Thus  $\nu_3<\nu_1$. Moreover, if $\chi_1=\lim_{\lambda\to0}\nu_1\nu_3=+\infty$, then the first equation in  \eqref{eq:reformulate} indicates that 
\begin{align*}
    0=\lim\limits_{\lambda\to0}&\lambda=\lim\limits_{\lambda\to0}\bigg(\psi_1/(\nu_1\nu_3)-\mu_{1,2}^2-\frac{\mu_{1,1}^2}{1+\mu_{1,1}^2\nu_1\nu_3+\mu_{2,1}^2\nu_2\nu_3}\bigg)\\
    &\qquad\qquad\cdot\bigg(\psi_3-\mu_{1,2}^2\nu_1\nu_3-\mu_{2,2}^2\nu_2\nu_3-\frac{\mu_{1,1}^2\nu_1\nu_3+\mu_{2,1}^2\nu_2\nu_3}{1+\mu_{1,1}^2\nu_1\nu_3+\mu_{2,1}^2\nu_2\nu_3}\bigg) \gg 0,
\end{align*}
which is impossible.
Therefore $\chi_1<+\infty$. Similarly, the second equation in \eqref{eq:reformulate} gives   $\chi_2=\lim_{\lambda\to0}\nu_2\nu_3<+\infty$. Here, $\chi_1,\chi_2<+\infty$ is obtained under a given moment vector $\bmu$. Combined $\chi_1<+\infty$ with $\nu_3<\nu_1$, we get $\sqrt{\lambda}\nu_3\to0$ as $\lambda\to0$. Therefore the third equation in \eqref{eq:implicit_trans} gives us 
\begin{align}
\label{eq:implicit1_lbd=0}
    \psi_3=\mu_{1,2}^2\chi_1+\mu_{2,2}^2\chi_2+\frac{\mu_{1,1}^2\chi_1+\mu_{2,1}^2\chi_2}{1+\mu_{1,1}^2\chi_1+\mu_{2,1}^2\chi_2}.
\end{align}
We remind the readers that we aim at proving 
\begin{align}\displaystyle 
        \lim\limits_{\mu_{2,1},\mu_{2,2}\to0}\lim\limits_{\lambda\to0}\cR(\lambda,\bpsi,\bmu,F_1,\tau) =+\infty.\label{eq:4thconcolutiongoal}
\end{align}
To show this, we rely on the following claim (recall that $\chi_2$ depends on $\mu_{2,1},\mu_{2,2}$):
\begin{equation}
    \lim\limits_{\mu_{2,1},\mu_{2,2}\to0}\mu_{2,2}^2\chi_2+\mu_{2,1}^2\chi_2=0.
    \label{eq:H-claim}
\end{equation}
In the following, we first explain how \eqref{eq:H-claim} can be used to show \eqref{eq:4thconcolutiongoal}, then give the proof of \eqref{eq:H-claim}. 

\noindent\textit{Proof of \eqref{eq:4thconcolutiongoal} based on \eqref{eq:H-claim}.} By \eqref{eq:implicit1_lbd=0} and \eqref{eq:H-claim}, we have
\begin{align}\label{eq:psi3equationFourthconclusion}
    \psi_3=\lim\limits_{\mu_{2,1},\mu_{2,2}\to0}\mu_{1,2}^2\chi_1+\frac{\mu_{1,1}^2\chi_1}{1+\mu_{1,1}^2\chi_1}.
\end{align}
Recall that in Theorem~\ref{thm:mainthm}, $\cR(\lambda,\bpsi,\bmu,F_1,\tau)$ is defined based on the quantities $S$ and $\Lb_{i,j}$, $i,j=1,\ldots,4$. The analytical expressions of these quantities are given in \eqref{eq:explicit_S} and \eqref{eq:explicit_L} respectively. 
We replace the terms $\bpsi$ and $\bnu^*$ in \eqref{eq:explicit_S} and \eqref{eq:explicit_L} with terms consisting of $\chi_1$ and $\chi_2$ by using equations \eqref{eq:Prop4replace}, \eqref{eq:H-claim}, \eqref{eq:psi3equationFourthconclusion}, and then get that
\begin{align*}
    \lim\limits_{\mu_{2,1},\mu_{2,2}\to0}\lim\limits_{\lambda\to0}S=0,~\lim\limits_{\mu_{2,1},\mu_{2,2}\to0}\lim\limits_{\lambda\to0}S\cdot (\Lb_{3,4}+\Lb_{1,4})\neq0,
    ~\lim\limits_{\mu_{2,1},\mu_{2,2}\to0}\lim\limits_{\lambda\to0}S\cdot(\Lb_{2,3}+\Lb_{1,2})\neq0.
\end{align*}
Therefore the limits   $\Lb_{3,4}+\Lb_{1,4}=\infty$ and $\Lb_{2,3}+\Lb_{1,2}=\infty$ when $\lambda\to0$ and $\mu_{2,1},\mu_{2,2}\to0$. Since $\cR>0$, we have 
\begin{align*}
  &\lim\limits_{\mu_{2,1},\mu_{2,2}\to0}\lim\limits_{\lambda\to0}\cR(\lambda,\bpsi,\bmu,F_1,\tau)\\
  &\qquad=\lim\limits_{\mu_{2,1},\mu_{2,2}\to0}\lim\limits_{\lambda\to0}F_1^2\bigg(
 \frac{1}{M_D^2}+\Lb_{3,4}+\Lb_{1,4}\bigg)+\tau^2(\Lb_{2,3}+\Lb_{1,2})=+\infty.
\end{align*}

\noindent\textit{Proof of \eqref{eq:H-claim}.} We first show that $\lim\limits_{\lambda\to0}\nu_3=0$. From the analysis above, we have $\lim\limits_{\lambda\to0}\nu_3<+\infty$ due to $\nu_3<\nu_1$ and $\chi_1=\lim\limits_{\lambda\to0}\nu_1\nu_3<+\infty$. If $\lim\limits_{\lambda\to0}\nu_3>0$, then combined with $\lim\limits_{\lambda\to0}\nu_1\nu_3<+\infty$ we have $\sqrt{\lambda}\nu_1$, $\sqrt{\lambda}\nu_2\to0$, the first and second equations in \eqref{eq:implicit_trans} give us\begin{align*}
    \left\{
    \begin{aligned}
    &\mu_{1,2}^2\chi_1+\frac{\mu_{1,1}^2\chi_1}{1+\mu_{1,1}^2\chi_1+\mu_{2,1}^2\chi_2}=\psi_1,\\
    &\mu_{2,2}^2\chi_2+\frac{\mu_{2,1}^2\chi_2}{1+\mu_{1,1}^2\chi_1+\mu_{2,1}^2\chi_2}=\psi_2.
    \end{aligned}
    \right.
\end{align*}
Combining the equations above with \eqref{eq:implicit1_lbd=0}, we have $\psi_1+\psi_2=\psi_3$, which is a contradiction to the condition $\psi_1=\psi_3$. Therefore $\lim\limits_{\lambda\to0}\nu_3=0$.

Combining the limit $\lim\limits_{\lambda\to0}\nu_3=0$ with  \eqref{eq:reformulate} yields that   $\lim\limits_{\lambda\to0}\sqrt{\lambda}(\nu_1+\nu_2)=\psi_2$. \eqref{eq:implicit_trans} further indicates the existence of $\lim\limits_{\lambda\to0}\sqrt{\lambda}\nu_1$ and $\lim\limits_{\lambda\to0}\sqrt{\lambda}\nu_2$ respectively due to the existence of $\chi_1$ and $\chi_2$ (The existence can also be guaranteed by the property of Stieltjes transform).
Next we show that $\lim\limits_{\lambda\to0}\sqrt{\lambda}\nu_1,\lim\limits_{\lambda\to0}\sqrt{\lambda}\nu_2>0$. We use a proof by contradiction:
\begin{itemize}[leftmargin = *]
    \item If $\lim\limits_{\lambda\to0}\sqrt{\lambda}\nu_1=0$, then it holds that $\lim\limits_{\lambda\to0}\sqrt{\lambda}\nu_2=\psi_2$,  then we conclude that $\nu_2\gg\nu_1$, and  $\lim\limits_{\lambda\to0}\nu_1\nu_3>0$, $\lim\limits_{\lambda\to0}\nu_2\nu_3=0$ from \eqref{eq:implicit_trans}. This is a contradiction because  
    $\lim\limits_{\lambda\to0}\nu_1\nu_3>0$ and $\lim\limits_{\lambda\to0}\nu_2\nu_3=0$ indicate $\nu_1\gg\nu_2$.
    \item  If $\lim\limits_{\lambda\to0}\sqrt{\lambda}\nu_2=0$, we have $\lim\limits_{\lambda\to0}\sqrt{\lambda}\nu_1=\psi_2$, then the second equation in \eqref{eq:implicit_trans} indicates that  $\lim\limits_{\lambda\to0}\nu_2\nu_3>0$. Moreover, $\lim\limits_{\lambda\to0}\sqrt{\lambda}\nu_2=0$  and $\lim\limits_{\lambda\to0}\sqrt{\lambda}\nu_1=\psi_2$
indicate that $\nu_1\gg\nu_2$, therefore $\lim\limits_{\lambda\to0}\nu_1\nu_3=+\infty$, which contradicts to the conclusion $\chi_1=\lim\limits_{\lambda\to0}\nu_1\nu_3<+\infty$ above.
\end{itemize}
From the analysis above we prove that $\nu_1$ and  $\nu_2$ have the same order when $\lambda\to0$. If $\chi_1=\lim\limits_{\lambda\to0}\nu_1\nu_3=0$, then $\chi_2=\lim\limits_{\lambda\to0}\nu_2\nu_3=0$. The first and second equations in \eqref{eq:implicit_trans} give us $\lim\limits_{\lambda\to0}\sqrt{\lambda}(\nu_1+\nu_2)\to\psi_1+\psi_2$ which contradicts  the third equation in \eqref{eq:reformulate} which  indicates that $\lim\limits_{\lambda\to0}\sqrt{\lambda}(\nu_1+\nu_2)\to\psi_2$.  Therefore we have $\chi_1,\chi_2>0$. Here, we utilize the fact $\lim\limits_{\lambda\to0}\nu_3=0$. Finally  we have 
\begin{align*}
    \nu_1=\Theta\bigg(\frac{1}{\sqrt{\lambda}}\bigg),~ \nu_2=\Theta\bigg(\frac{1}{\sqrt{\lambda}}\bigg),~
    \nu_3=\Theta(\sqrt\lambda).
\end{align*}
Then we can assume that 
\begin{align*}
    \lim\limits_{\lambda\to0}\sqrt{\lambda}\nu_1=\psi_1-n_1,~
    \lim\limits_{\lambda\to0}\sqrt{\lambda}\nu_2=\psi_2-n_2,~
    \lim\limits_{\lambda\to0}\nu_3/\sqrt{\lambda}=k,
\end{align*}
where $0\leq n_1<\psi_1$, $0\leq n_2<\min(\psi_1,\psi_2)$, $k>0$ and $n_1,n_2,k$ satisfy
\begin{align}
\label{eq:implicitprop_psi1=psi3}
\left\{\begin{aligned}
    &\mu_{1,2}^2(\psi_1-n_1)k+\frac{\mu_{1,1}^2(\psi_1-n_1)k}{1+\mu_{1,1}^2(\psi_1-n_1)k+\mu_{2,1}^2(\psi_2-n_2)k}=n_1,\\
    &\mu_{2,2}^2(\psi_2-n_2)k+\frac{\mu_{2,1}^2(\psi_2-n_2)k}{1+\mu_{1,1}^2(\psi_1-n_1)k+\mu_{2,1}^2(\psi_2-n_2)k}=n_2,\\
    &n_1+n_2=\psi_3=\psi_1.
    \end{aligned}\right.
\end{align}
It is easy to see that $\chi_1=(\psi_1-n_1)\cdot k$ and $\chi_2=(\psi_2-n_2)\cdot k$. We can also show that $\lim\limits_{\mu_{2,1},\mu_{2,2}\to0} n_2 = 0$. 
Indeed if $\limsup\limits_{\mu_{2,1},\mu_{2,2}\to0} n_2 > 0$, the second equation in \eqref{eq:implicitprop_psi1=psi3} gives $k\to+\infty$. However, $n_2\cdot k\to+\infty$ leads to a contradiction to the first equation in \eqref{eq:implicitprop_psi1=psi3}. 
 Next, using the second equation in \eqref{eq:implicitprop_psi1=psi3}, we have $\lim\limits_{\mu_{2,1},\mu_{2,2}\to0}\mu_{2,2}^2\chi_2=0$.  
 
As for $\lim\limits_{\mu_{2,1},\mu_{2,2}\to0}\mu_{2,1}^2\chi_2$, note that $\chi_1$ is bounded by the inequality $\chi_1<\psi_3/\mu_{1,2}^2$ due to the first equation in \eqref{eq:implicitprop_psi1=psi3}, thus we have 
\[
0=\lim\limits_{\mu_{2,1},\mu_{2,2}\to0}\frac{\mu_{2,1}^2(\psi_2-n_2)k}{1+\mu_{1,1}^2(\psi_1-n_1)k+\mu_{2,1}^2(\psi_2-n_2)k}= \lim\limits_{\mu_{2,1},\mu_{2,2}\to0}\frac{\mu_{2,1}^2\chi_2}{1+\mu_{1,1}^2\chi_1+\mu_{2,1}^2\chi_2},
\]
which indicates that  $\lim\limits_{\mu_{2,1},\mu_{2,2}\to0}\mu_{2,1}^2\chi_2=0$. Hence the claim \eqref{eq:H-claim} holds and the proof is complete. 
\end{proof}
\begin{proof}[\textbf{Third conclusion}]
Let $r=1-(c_2-1)\psi_3/\psi_2$, then $\psi_3=\psi_1+r\psi_2$ with $0<r<1$. An  analysis similar  to the previous case of $\psi_3=\psi_1$ leads to  $\nu_3<\nu_1+\nu_2$, and $\nu_1=\Theta(\frac{1}{\sqrt{\lambda}}),\nu_2=\Theta(\frac{1}{\sqrt{\lambda}})$, and $\nu_3=\Theta(\sqrt\lambda)$. 
We still assume that  
\begin{align*}
    \lim\limits_{\lambda\to0}\sqrt{\lambda}\nu_1=\psi_1-n_1,~
    \lim\limits_{\lambda\to0}\sqrt{\lambda}\nu_2=\psi_2-n_2,~
    \lim\limits_{\lambda\to0}\nu_3/\sqrt{\lambda}=k,
\end{align*}
where $0\leq n_1<\psi_1$, $r\psi_2\leq n_2<\psi_2$, $k>0$ and $n_1,n_2,k$ satisfy
\begin{align}
\label{eq:implicitprop_psi1+rpsi2=psi3}
\left\{\begin{aligned}
    &\mu_{1,2}^2(\psi_1-n_1)k+\frac{\mu_{1,1}^2(\psi_1-n_1)k}{1+\mu_{1,1}^2(\psi_1-n_1)k+\mu_{2,1}^2(\psi_2-n_2)k}=n_1,\\
    &\mu_{2,2}^2(\psi_2-n_2)k+\frac{\mu_{2,1}^2(\psi_2-n_2)k}{1+\mu_{1,1}^2(\psi_1-n_1)k+\mu_{2,1}^2(\psi_2-n_2)k}=n_2,\\
    &n_1+n_2=\psi_3=\psi_1+r\psi_2.
    \end{aligned}\right.
\end{align}
It is easy to see  that $\chi_1=(\psi_1-n_1)\cdot k$ and $\chi_2=(\psi_2-n_2)\cdot k$. Let ${\mu_{2,1},\mu_{2,2}\to0}$ and  note that $n_2\geq r\psi_2$.  We must have $k\to+\infty$ by  the second equation in \eqref{eq:implicitprop_psi1+rpsi2=psi3}. Therefore the first equation in \eqref{eq:implicitprop_psi1+rpsi2=psi3} indicates that $n_1=\psi_1$ and $n_2=r\psi_2$ as ${\mu_{2,1},\mu_{2,2}\to0}$. Now it is easy to prove the third conclusion in Proposition~\ref{prop:property_lbd=0} if we further assume that $\mu_{2,1}/\mu_{2,2}\to0$ due to 
\begin{align*}
    \varliminf\limits_{\mu_{2,1},\mu_{2,2}\to0}\lim\limits_{\lambda\to0}\cR(\lambda,\bpsi,\bmu,F_1,\tau)\leq\lim\limits_{\substack{\mu_{2,1},\mu_{2,2}\to0\\\mu_{2,1}/\mu_{2,2}\to0}}\cR(\lambda,\bpsi,\bmu,F_1,\tau).
\end{align*}
Define $\bar{\chi}_1=\lim\limits_{\substack{\mu_{2,1},\mu_{2,2}\to0\\\mu_{2,1}/\mu_{2,2}\to0}}\chi_1$. Then we have
\begin{align*}
    &\lim\limits_{\substack{\mu_{2,1},\mu_{2,2}\to0\\\mu_{2,1}/\mu_{2,2}\to0}}\mu_{2,2}^2\chi_2=r\psi_2,~\mu_{1,2}^2\bar{\chi}_1+\frac{\mu_{1,1}^2\bar{\chi}_1}{1+\mu_{1,1}^2\bar{\chi}_1}=\psi_1.
\end{align*}
Combining  the expression of $S$ and $\Lb_{i,j}$ in \eqref{eq:explicit_S} and $\eqref{eq:explicit_L}$ gives us
   \begin{align*}
    &\lim\limits_{\substack{\mu_{2,1},\mu_{2,2}\to0\\\mu_{2,1}/\mu_{2,2}\to0}}\lim\limits_{\lambda\to0}S=(1-r)r\bar{\chi}_1(\psi_2+\mu_{1,1}^2\psi_2\bar{\chi}_1)^2(\mu_{1,2}^2+\mu_{1,1}^4\mu_{1,2}^2\bar{\chi}_1^2+\mu_{1,1}^2(1+2\mu_{1,2}^2\bar{\chi}_1))>0,\\
    &\lim\limits_{\substack{\mu_{2,1},\mu_{2,2}\to0\\\mu_{2,1}/\mu_{2,2}\to0}}\lim\limits_{\lambda\to0}|S\cdot (\Lb_{3,4}M_D^2+\Lb_{1,4})|<+\infty,
    ~\lim\limits_{\substack{\mu_{2,1},\mu_{2,2}\to0\\\mu_{2,1}/\mu_{2,2}\to0}}\lim\limits_{\lambda\to0}|S\cdot(\Lb_{2,3}+\Lb_{1,2})|<+\infty.
\end{align*}
Therefore
$
    \varliminf\limits_{\mu_{2,1},\mu_{2,2}\to0}\lim\limits_{\lambda\to0}\cR(\lambda,\bpsi,\bmu,F_1,\tau)\leq\lim\limits_{\substack{\mu_{2,1},\mu_{2,2}\to0\\\mu_{2,1}/\mu_{2,2}\to0}}\cR(\lambda,\bpsi,\bmu,F_1,\tau)<+\infty.
$
This completes the proof of the third conclusion in Proposition~\ref{prop:property_lbd=0}.
\end{proof}
\begin{proof}[\textbf{First conclusion}]
Let $r=1+(1-c_1)\psi_3/\psi_2$, then we have $\psi_3=\psi_1+r\psi_2$ with $r>1$.  Similarly to the previous arguments, we obtain   $\nu_1\nu_3=\bTheta_{\lambda}(1)$ and $\nu_2\nu_3=\bTheta_{\lambda}(1)$. Also note that $\nu_1+\nu_2<\nu_3$, therefore it holds that $\sqrt{\lambda}\nu_1\to0$ and $\sqrt{\lambda}\nu_1\to0$.
Recall that we defined  $\chi_1=\lim_{\lambda\to0}\nu_1\nu_3$ and $\chi_2=\lim_{\lambda\to0}\nu_2\nu_3$, and the system \eqref{eq:implicit_lbd=0} still holds in the current  case. Substituting  \eqref{eq:implicit_lbd=0} into \eqref{eq:explicit_S} and \eqref{eq:explicit_L}, and after some simple calculation,  we obtain 
\begin{align*}
    \lim\limits_{\lambda\to0}S>(r-1)\mu_{1,1}^4\mu_{2,1}^2(1+\mu_{2,2}^2\chi_2)\chi_1^2\chi_2^2>0,
~\lim\limits_{\lambda\to0}S\cdot (\Lb_{3,4}M_D^2+\Lb_{1,4})<+\infty,
    \lim\limits_{\lambda\to0}S\cdot (\Lb_{2,3}+\Lb_{1,2})<+\infty.
\end{align*}
Therefore when $\psi_3=\psi_1+r\psi_2$, $r>1$, 
$
    \lim\limits_{\lambda\to0}\cR(\lambda,\bpsi,\bmu,F_1,\tau)<+\infty.
$
This completes the proof of the first conclusion in Proposition~\ref{prop:property_lbd=0}.
\end{proof}
\subsection{Proof of Proposition~\ref{prop:property_psi}}
For this proposition, we let  $\psi_0=\psi_1/r_1=\psi_2/r_2\to+\infty$.
By the system~\eqref{eq:implicit_trans} we have
\begin{align}
\label{eq:3inH1}
    \sqrt{\lambda}\nu_3=\psi_3-\mu_{1,2}^2\nu_1\nu_3-\mu_{2,2}^2\nu_2\nu_3-\frac{\mu_{1,1}^2\nu_1\nu_3+\mu_{2,1}^2\nu_2\nu_3}{1+\mu_{1,1}^2\nu_1\nu_3+\mu_{2,1}^2\nu_2\nu_3}.
\end{align}
Therefore $\nu_3<\psi_3/\sqrt{\lambda}$ with fixed $\psi_3$. Then from the first and second equations in \eqref{eq:implicit_trans} we easily  get that $\lim\limits_{\psi_0\to+\infty}\nu_1,\lim\limits_{\psi_0\to+\infty}\nu_2=+\infty$. 
If $\varlimsup\limits_{\psi_0\to+\infty}\nu_3>0$,  further from \eqref{eq:3inH1} we will get
\begin{align*}
    \varlimsup\limits_{\psi_0\to+\infty}\sqrt{\lambda}\bigg(\nu_3+\mu_{1,2}^2\nu_1\nu_3+\mu_{2,2}^2\nu_2\nu_3+\frac{\mu_{1,1}^2\nu_1\nu_3+\mu_{2,1}^2\nu_2\nu_3}{1+\mu_{1,1}^2\nu_1\nu_3+\mu_{2,1}^2\nu_2\nu_3}\bigg)=\psi_3.
\end{align*}
This is a contradiction because 
the left hand side of the equation above tends to infinity  while  the right hand  side is fixed. Therefore we have $\lim\limits_{\psi_0\to+\infty}\nu_3=0$. Combining this result with
\begin{align*}
    &\sqrt{\lambda}\nu_1+\mu_{1,2}^2\nu_1\nu_3+\frac{\mu_{1,1}^2\nu_1\nu_3}{1+\mu_{1,1}^2\nu_1\nu_3+\mu_{2,1}^2\nu_2\nu_3}=\psi_1,\\
    &\sqrt{\lambda}\nu_2+\mu_{2,2}^2\nu_2\nu_3+\frac{\mu_{2,1}^2\nu_2\nu_3}{1+\mu_{1,1}^2\nu_1\nu_3+\mu_{2,1}^2\nu_2\nu_3}=\psi_2,
\end{align*}
 we conclude that
\begin{align*}
    \lim\limits_{\psi_0\to+\infty}\nu_1/\psi_0=r_1/\sqrt{\lambda},~\lim\limits_{\psi_0\to+\infty}\nu_2/\psi_0=r_2/\sqrt{\lambda}.
\end{align*}
We further define 
\begin{align*}
    \varlimsup\limits_{\psi_0\to+\infty}\nu_3\psi_0=\overline{\chi},~\varliminf\limits_{\psi_0\to+\infty}\nu_3\psi_0=\underline{\chi}.
\end{align*}
Then we have 
\begin{align*}
    \varlimsup\limits_{\psi_0\to+\infty}\nu_1\nu_3=r_1\overline{\chi},~
    \varlimsup\limits_{\psi_0\to+\infty}\nu_2\nu_3=r_2\overline{\chi},~
    \varliminf\limits_{\psi_0\to+\infty}\nu_1\nu_3=r_1\underline{\chi},~
    \varliminf\limits_{\psi_0\to+\infty}\nu_2\nu_3=r_2\underline{\chi}.
\end{align*}
Taking the superior and inferior limit when $\psi_0\to+\infty$ in the third equation of \eqref{eq:implicit_trans}, we have
\begin{align*}
\left\{\begin{aligned}
    &\psi_3=\mu_{1,2}^2r_1\overline{\chi}+\mu_{2,2}^2r_2\overline{\chi}+\frac{\mu_{1,1}^2r_1\overline{\chi}+\mu_{2,1}^2r_2\overline{\chi}}{1+\mu_{1,1}^2r_1\overline{\chi}+\mu_{2,1}^2r_2\overline{\chi}},\\
    &\psi_3=\mu_{1,2}^2r_1\underline{\chi}+\mu_{2,2}^2r_2\underline{\chi}+\frac{\mu_{1,1}^2r_1\underline{\chi}+\mu_{2,1}^2r_2\underline{\chi}}{1+\mu_{1,1}^2r_1\underline{\chi}+\mu_{2,1}^2r_2\underline{\chi}}.
    \end{aligned}\right.
\end{align*}
Therefore $\overline{\chi}$ and $\underline{\chi}$ are both the solution of the equation
\begin{align}\label{eq:quadratic_equation_of_chi}
    \psi_3(1+\mu_{1,1}^2r_1x+\mu_{2,1}^2r_2x)=(\mu_{1,2}^2r_1x+\mu_{2,2}^2r_2x)(1+\mu_{1,1}^2r_1x+\mu_{2,1}^2r_2x)+\mu_{1,1}^2r_1x+\mu_{2,1}^2r_2x.
\end{align}
Note that $\overline{\chi}$ and $\underline{\chi}$ are both positive, and the equation above only has one positive root. Therefore we conclude that $\overline{\chi}=\underline{\chi}$, and we can write $\chi:=\overline{\chi}=\underline{\chi}$.
By calculating the positive root of \eqref{eq:quadratic_equation_of_chi}, we easily see that $(r_1\mu_{1,1}^2+r_2\mu_{2,1}^2)\chi=\chi_0$ where $\chi_0$ is defined in Proposition~\ref{prop:property_psi}. Plugging the limits $\nu_1\nu_3\to r_1\chi$ and $\nu_2\nu_3\to r_2\chi$ into $M_D$ and $M_N$ in \eqref{eq:explicit_S} and \eqref{eq:explicit_L}, we obtain $M_D\to -\chi_0-1$, $\nu_3^*M_N\to -\chi_0$ when  $\psi_0\to+\infty$. Direct algebraic calculation then gives 
\begin{align*}
    \Lb_{2,3}\to\frac{\chi_0^2}{(\chi_0+1)^2\psi_3-\chi_0^2},\quad~\Lb_{3,4}\to\frac{\chi_0^2}{(\chi_0+1)^4\psi_3-\chi_0^2(\chi_0+1)^2},\quad~\Lb_{1,2},\Lb_{1,4}\to0
\end{align*}
when $\psi_0\to+\infty$. Then we have
\begin{align*}
    \lim\limits_{\psi_0\to\infty}\cR(\lambda,\bpsi,\bmu,F_1,\tau)&=\lim\limits_{\psi_0\to\infty}F_1^2\bigg(
 \frac{1}{M_D^2}+\Lb_{3,4}+\Lb_{1,4}\bigg)+\tau^2(\Lb_{2,3}+\Lb_{1,2})\\
 &=F_1^2\bigg(\frac{1}{(\chi_0+1)^2}+ \frac{\chi_0^2}{(\chi_0+1)^4\psi_3-\chi_0^2(\chi_0+1)^2}\bigg)+\tau^2\bigg(\frac{\chi_0^2}{(\chi_0+1)^2\psi_3-\chi_0^2}\bigg)\\
 &=\frac{F_1^2\psi_3+\tau^2\chi_0^2} 
    {(\chi_0+1)^2\psi_3-\chi_0^2}.
\end{align*}
This  proves Proposition~\ref{prop:property_psi}.
\section{Proof of Theorem~\ref{thm:mainthm-m}}
\label{sec:appendixthm-m}
\xr{Here, we provide the proof of Theorem~\ref{thm:mainthm-m} for the MRFM. The proof for MRFM bears significant resemblance to the previous proof of Theorem~\ref{thm:mainthm}. In this section, we offer a brief overview of the proof for MRFM, highlighting the key distinctions between these two theorems. Here we will focus on several key steps in the proof that are significantly different from the proof of DRFMs.}



\subsection{Step 1: bias-variance decomposition of the excess risk}
\label{step1-m} 
We first give some notations as follows.  
\begin{definition}
\label{def:notaion-m}
Define
\begin{align*}
    &\Zb _\sfc=\sigma_\sfc\left(\Xb\bTheta_\sfc^\T /\sqrt d\right)/\sqrt d\in\RR^{n\times N_\sfc},~\Zb=\big[ \Zb _1,\ldots, \Zb _K\big],\\
    & \bsigma(\xb)=\big[\sigma_1(\xb^\T \bTheta_1^\T /\sqrt{d}),\ldots,\sigma_K(\xb^\T \bTheta_K^\T /\sqrt{d})\big]^\T \in\RR^N,~\bUpsilon=( \Zb ^\T  \Zb +\lambda\Ib_N)^{-1},\\
 &\Vb_0(F_0)=\EE_{\xb}\big[\bsigma(\xb) F_0\big]\in\RR^{N\times 1},~\Vb(\bbeta_{1,d})=\EE_{\xb}\big[\bsigma(\xb)\xb^\T \bbeta_{1,d}\big]\in\RR^{N\times 1}, \Ub=\EE_{\xb}\big[\bsigma(\xb)\bsigma(\xb)^\T \big]\in\RR^{N\times N}. \qquad\qquad
\end{align*}
\end{definition}
Clearly, these notations are consistent with Definition~\ref{def:notation1} and Proposition~\ref{prop:decomp-expectation}. Based on these notations with direct calculation, we can express the excess risk $R_d(\Xb,\bTheta,\lambda,\bbeta_d,\bvarepsilon)$ as
\begin{align}
\label{eq:Rdecomposition1-m}
  R_d(\Xb,\bTheta,\lambda,\bbeta_d,\bvarepsilon)=F_0^2+F_{1,d}^2-2\yb^\T \Zb\bUpsilon[\Vb(\bbeta_{1,d})+\Vb_0(F_0)]/\sqrt{d}+\yb^\T \big[\Ub\big]_\Zb \yb/d.
\end{align}
To continue the calculation, we consider the Gegenbauer decompositions of the activation functions. Suppose that the Gegenbauer decompositions of  $\sigma_{\sfc}(\cdot)$, $\sfc=1,\ldots,K,$ are 
\begin{align*}
  \sigma_{\sfc}(x)=\sum_{k=0}^{+\infty}\lambda_{d,k}\big(\sigma_{\sfc}\big)B(d,k)\cdot Q_k^{(d)}(\sqrt{d}\cdot x),\quad\sfc=1,\ldots,K,  
\end{align*}
where $\lambda_{d,k}(\sigma_{\sfc})$ are the decomposition coefficients, $Q_k^{(d)}$, $k\in \NN$ are the Gegenbauer polynomials, and $B(d,0)=1$,
$B(d,k)=k^{-1} (2k+d-2)\binom{k+d-3}{k-1}$ for $k\geq1$. Let   
\begin{align}
  &\bLambda_{d,k}=\diag\big(\lambda_{d,k}(\sigma_1)\Ib_{N_1},\ldots,\lambda_{d,k}(\sigma_K)\Ib_{N_K}\big),\quad k\in\NN=\{0,1,...\}, \label{eq:def_lambda_k-m}\\
&\Mb_1=\diag\big(\mu_{1,1}\Ib_{N_1},\ldots,\mu_{K,1}\Ib_{N_K}\big),\quad \Mb_2=\diag\big(\mu_{1,2}\Ib_{N_1},\ldots,\mu_{K,2}\Ib_{N_K}\big).\label{eq:def_M_k-m}
\end{align}
Now we present Proposition~\ref{prop:decomp-expectation-m} below, which is the counterpart of Proposition~\ref{prop:decomp-expectation}.
\begin{proposition}
\label{prop:decomp-expectation-m}
For any given $\lambda$, let
\begin{equation*}
    \begin{split}
        \overline{R}_d(\Xb,\bTheta,\lambda,F_{1,d},\tau)=F_{1,d}^2-\frac{2F_{1,d}^2}{d}\tr \Mb_1\frac{\bTheta\Xb^\T }{d}\Zb\bUpsilon+\frac{F_{1,d}^2}{d}\tr\Big( \big[\tilde{\Ub}\big]_\Zb \frac{\Xb\Xb^\T }{d}\Big)+\frac{\tau^2}{d}\tr(\big[\tilde{\Ub}\big]_\Zb ),
    \end{split}
\end{equation*}
where $\tilde{\Ub}=\Mb_1 \bTheta\bTheta^\T \Mb_1 / d + \Mb_2\Mb_2$. Then under the same conditions as Theorem~\ref{thm:mainthm-m},
\begin{align*}
    \begin{split}
        \EE_{\Xb,\bTheta,\bvarepsilon}&\Big|R_d(\Xb,\bTheta,\lambda,\bbeta_d,\bvarepsilon)-\overline{R}_d(\Xb,\bTheta,\lambda,F_{1,d},\tau)\Big|=o_d(1).
    \end{split}
\end{align*}
\end{proposition}

The proof for Proposition~\ref{prop:decomp-expectation-m} is exactly the same as the proof for Proposition~\ref{prop:decomp-expectation}, except the definitions of $\bLambda_{d,k}$, $\Mb_1$ and $\Mb_2$ are changed. We therefore omit the proof details.


\subsection{Step 2: approximation of the risk decomposition via a linear pencil matrix} \label{step2-m}
The approximating function $ \overline{R}_d(\Xb,\bTheta,\lambda,F_{1,d},\tau)$ established in  Proposition~\ref{prop:decomp-expectation-m} again depends on traces of several random matrices. These traces are next evaluated using a new linear pencil matrix, which is a bit more involved compared with the linear pencil matrix for DRFMs. 


\begin{definition}
\label{def:linear pencil-m}
(1) Let $\mathcal{Q}:=\{\qb=[q_1,q_2,q_3,q_4,q_5]\in\RR_+^5: q_4,q_5\leq(1+q_1)/2, \| \qb \|_2 \leq 1\}$. 
Depending on $\qb\in\cQ$ and $\bmu$,
the linear pencil matrix $\Ab(\qb,\bmu)\in\RR^{\PO\times \PO}$ ($\PO=N+n$) is 
\begin{align*}
    \Ab(\qb,\bmu)=&\begin{bmatrix}
    q_2\Mb_2\Mb_2+q_4\Mb_1\frac{\bTheta\bTheta^\T}{d}\Mb_1&\Zb^\T+q_1\widetilde{\Zb}^\T\\
    \Zb+q_1\widetilde{\Zb}&q_3\Ib_n+q_5\frac{\Xb\Xb^\T}{d}
    \end{bmatrix}\\
    =&\begin{bmatrix}
q_2\mu_{1,2}^2\Ib_{N_1}+q_4\mu_{1,1}^2\frac{\bTheta_1\bTheta_1^\T }{d}&\cdots&q_4\mu_{1,1}\mu_{K,1}\frac{\bTheta_1\bTheta_K^\T }{d}& \Zb_1^\T +q_1\tilde{\Zb}_1^\T \\\vdots&\ddots&\vdots&\vdots\\
q_4\mu_{K,1}\mu_{1,1}\frac{\bTheta_K\bTheta_1^\T }{d}&\cdots&q_2\mu_{K,2}^2\Ib_{N_K}+q_4\mu_{K,1}^2\frac{\bTheta_K\bTheta_K^\T }{d}& \Zb_K^\T +q_1\tilde{\Zb}_K^\T \\
 \Zb_1+q_1\tilde{\Zb}_1&\cdots& \Zb_K+q_1\tilde{ \Zb}_K&q_3\Ib_n+q_5\frac{\Xb\Xb^\T }{d}
\end{bmatrix},
\end{align*}
where $\tilde{\Zb}_{\sfc}=\frac{\mu_{\sfc,1}}{d}\Xb\bTheta_{\sfc}^\T $, $\sfc=1,\ldots,K+1$. 

(2)  The Stieltjes transform of the empirical eigenvalue distribution of $\Ab$ (up to a $P/d$ factor)  is
\begin{equation*}
    \begin{split}
    M_d(\xi;\qb,\bmu)=\frac{1}{d}\tr\big[(\Ab-\xi\Ib_{\PO})^{-1}\big], \quad \xi\in\bbC_+, 
    \end{split}
\end{equation*}
and its logarithmic potential is  
\begin{equation*}
\begin{split}
    G_d(\xi;\qb,\bmu)=\frac1d\log\det(\Ab-\xi\Ib_{\PO})=\frac1d  \sum_{i=1}^{\PO}\Log(\lambda_i(\Ab)-\xi),\quad \xi\in\bbC_+.
\end{split}
\end{equation*}
Here $\{\lambda_i(\Ab)\}_{i\in[\PO]}$ are the  eigenvalues of $\Ab$ in decreasing order, and
$\Log(z):=\Log(|z|)+\rmi\arg(z)$, for $z\in\bbC$, $-\pi<\arg(z)\leq\pi$
is the principal value of a complex logarithmic function. 
\end{definition}

The three traces appearing in the definition of $\overline{R}_d(\Xb,\bTheta,\lambda,F_{1,d},\tau)$ in Proposition~\ref{prop:decomp-expectation-m} are now expressed as partial derivatives of the logarithmic potential $G_d$ as shown in the proposition below.
\begin{proposition}
\label{prop:connectGd-m}
Let $\xi^*$ be defined in Definition~\ref{def:notaion-m} and $\tilde{\Ub}$ be defined in Proposition~\ref{prop:decomp-expectation-m}. Then
\begin{equation*}
    \begin{split}
        \frac1d\tr \Mb_1\frac{\bTheta\Xb^\T }{d}\Zb\bUpsilon&=\frac{1}{2}\partial_{q_1} G_d(\xi^*;\qb,\bmu) |_{\qb = \mathbf{0}},\\
        \frac{1}{d}\tr( \big[\tilde{\Ub}\big]_\Zb \frac{\Xb\Xb^\T }{d})&=-\partial^2_{q_4,q_5}G_d(\xi^*;\qb,\bmu)|_{\qb = \mathbf{0}}-\partial^2_{q_2,q_5}G_d(\xi^*;\qb,\bmu)|_{\qb = \mathbf{0}},\\
        \frac{1}{d}\tr(\big[\tilde{\Ub}\big]_\Zb )&=-\partial^2_{q_3,q_4}G_d(\xi^*;\qb,\bmu)|_{\qb = \mathbf{0}}-\partial^2_{q_2,q_3}G_d(\xi^*;\qb,\bmu)|_{\qb = \mathbf{0}}.
    \end{split}
\end{equation*}
\end{proposition}
The proof for Proposition~\ref{prop:connectGd-m} is the same as the proof for Proposition~\ref{prop:connectGd}.  We omit the details for simplicity.

\subsection{Step 3: key limiting spectral functions of the linear pencil matrix}
\label{step3-m}
Proposition \ref{prop:connectGd-m} shows that the excess risk depends on the limiting spectral properties of the linear pencil matrix $\Ab$. Therefore we study the Stieltjes transform $ M_d(\xi;\qb,\bmu)$ of the empirical eigenvalue distribution of $\Ab$ and calculate its limit as $d,n,N \rightarrow \infty$. We first give the following definition. 
\begin{definition}
\label{def:implicit1-m}
Define $\sfFb(\cdot;\xi,\qb,\bmu)= [\sfF_1(\cdot;\xi,\qb,\bmu),\cdots,
    \sfF_{K+1}(\cdot;\xi,\qb,\bmu)]^\top:\bbC^{K+1}\to\bbC^{K+1}$ as
\begin{equation*}
    \begin{split}
      &\sfF_{\sfc}(\mb;\xi,\qb,\bmu)=\psi_\sfc\bigg\{-\xi+q_2\mu_{\sfc,2}^2-\mu_{\sfc,2}^2m_{K+1}+\frac{H_{\sfc}}{H_{D}}\bigg\}^{-1},\quad\sfc=1,\ldots,K,\\
&\sfF_{K+1}(\mb;\xi,\qb,\bmu)=\psi_{K+1}\bigg\{-\xi+q_3-\sum\limits_{\sfc=1}^K\mu_{\sfc,2}^2m_{\sfc}+\frac{H_{{K+1}}}{H_{D}}\bigg\}^{-1},
    \end{split}
\end{equation*}
where $\xi\in\bbC_+$,  $\mb=[m_1,\ldots,m_{K+1}]\in\bbC^{K+1}$, and
$$\begin{aligned}
&H_{\sfc}=\mu_{\sfc,1}^2q_4(1+m_{K+1}q_5)-\mu_{\sfc,1}^2(1+q_1)^2m_{K+1},\quad\sfc=1,\ldots,K,\\
&H_{{K+1}}=q_5\Big(1+\sum\limits_{\sfc=1}^K\mu_{\sfc,1}^2m_{\sfc}q_4\Big)-(1+q_1)^2\sum\limits_{\sfc=1}^K\mu_{\sfc,1}^2m_{\sfc},\\
&H_{D}=\Big(1+\sum\limits_{\sfc=1}^K\mu_{\sfc,1}^2m_{\sfc}q_4\Big)(1+m_{K+1}q_5)-(1+q_1)^2\sum\limits_{\sfc=1}^K\mu_{\sfc,1}^2m_{\sfc}m_{K+1}.
\end{aligned}$$
\end{definition}
Note that the function $\sfFb(\mb;\xi,\qb,\bmu)$ in Definition~\ref{def:implicit1-m} above is not related to $d$.
Lemma~\ref{lemma:uniquesolution-m} below ensures the existence and uniqueness of the fixed point of  $\sfFb(\mb;\xi,\qb,\bmu)$ for $\xi\in \{\xi\in\bbC: \Im(\xi)>\xi_0 \}$ with some sufficiently large constant $\xi_0$. 

\begin{lemma}
\label{lemma:uniquesolution-m}
For $\sfFb(\mb;\xi,\qb,\bmu)$ in Definition~\ref{def:implicit1-m},
there exists $\xi_0>0$ such that, for any $\xi\in\bbC_+ $ with $\Im(\xi)>\xi_0$, the equation~$\mb=\sfFb(\mb;\xi,\qb,\bmu)$ admits a unique solution in $\DD(2\psi_1/\xi_0)\times\ldots\times\DD(2\psi_{K+1}/\xi_0)$. 
\end{lemma}
The proof of Lemma~\ref{lemma:uniquesolution-m} is given in Section~\ref{sec:appendixuniquesolution-m}. Define the fixed point of $\sfFb(\mb;\xi,\qb,\bmu)$ as the function of $\xi$ on $\{\xi:\Im(\xi)>\xi_0\}$: 
\begin{align}\mb(\xi;\qb,\bmu)=\begin{bmatrix}
m_1(\xi;\qb,\bmu)\\
\vdots \\
m_{K+1}(\xi;\qb,\bmu)
\end{bmatrix}\label{eq:m_def-m}\end{align}
The following proposition shows that $\mb$ is an analytic function on $\{\xi: \Im(\xi) > \xi_0\}$, and its analytic continuation to $\bbC_+$ is still a fixed point of $\sfFb(\cdot;\xi, \qb,\bmu)$.

\begin{proposition}
\label{prop:implicit1-m}
Under Assumptions~\ref{assump1-m} and \ref{assump2-m},
$\mb(\xi;\qb,\bmu)$ is analytic on $\{\xi:\Im(\xi)>\xi_0\}$, and has a unique analytic continuation to $\bbC_+$. Moreover, this analytic continuation (still denoted as $\mb(\xi;\qb,\bmu)$) satisfies the following properties:
\begin{enumerate}[leftmargin = *]
    \item $\mb(\xi;\qb,\bmu) \in \bbC_+^{K+1}$ for all $\xi\in \bbC_+$.
    \item $\mb(\xi,\qb,\bmu)\equiv\sfFb[\mb(\xi,\qb,\bmu);\xi,\qb,\bmu]$ for all $\xi\in \bbC_+$.
    \item Let $M_d(\xi;\qb,\bmu)$ be defined in Definition~\ref{def:linear pencil-m}. Then for any compact set $\Omega\subset\bbC_+$,  \begin{equation*}
    \lim_{d\rightarrow+\infty}\EE\bigg[\sup_{\xi\in\Omega}\big|M_d(\xi;\qb,\bmu)-\sum\limits_{\sfc=1}^{K+1}m_{\sfc}(\xi;\qb,\bmu)\big|\bigg]=0.
\end{equation*} 
\end{enumerate}
\end{proposition}
The proof of Proposition~\ref{prop:implicit1-m} is given in Section~\ref{subsec:appendixpropimplicit1-m}. 
The study of the limiting spectral distribution also leads to a deterministic limit for the logarithmic potential $G_d$.  This limit logarithmic potential is found to be 
\begin{equation}
\label{eq:gd-m}
    g(\xi;\qb,\bmu)\triangleq L(\xi,m_1(\xi;\qb,\bmu),\ldots,m_{K+1}(\xi;\qb,\bmu);\qb,\bmu),
\end{equation}
where the function $L$ is 
\begin{equation}\label{eq:L-function-m}
    \begin{split}
        & L(\xi,z_1,\ldots,z_{K+1};\qb,\bmu)  \triangleq \\
        &\quad \Log\bigg[\Big(1+q_4\sum\limits_{\sfc=1}^K\mu_{\sfc,1}^2z_\sfc\Big)(1+z_{K+1}q_5)-\sum\limits_{\sfc=1}^K\mu_{\sfc,1}^2(1+q_1)^2z_\sfc z_{K+1}\bigg]-\sum\limits_{\sfc=1}^K\mu_{\sfc,2}^2z_\sfc z_{K+1}\\
&\qquad +q_2\sum\limits_{\sfc=1}^K\mu_{\sfc,2}^2z_\sfc+q_3z_{K+1} -\sum\limits_{\sfc=1}^{K+1}\psi_\sfc\Log(z_\sfc/\psi_\sfc)-\xi\bigg(\sum\limits_{\sfc=1}^{K+1}z_\sfc\bigg)-\sum\limits_{\sfc=1}^{K+1}\psi_\sfc.
    \end{split}
\end{equation}
This convergence, together with those of the  partial derivatives of our interest, are  formally established in the following proposition.

\begin{proposition}
\label{prop:substitutiongd-m}
Let $G_d(\xi;\qb,\bmu)$ be defined in Definition~\ref{def:linear pencil-m}, and $g(\xi;\qb,\bmu)$ be defined in equation~\eqref{eq:gd-m}. For any fixed $\qb\in\cQ$, $\xi\in\bbC_+$ and $u\in\RR_+$,
\begin{align*}
&\lim_{d\rightarrow+\infty}\EE[|G_d(\xi;\qb,\bmu)-g(\xi;\qb,\bmu)|]=0,\\
&\lim_{d\rightarrow+\infty}\EE[\lVert \nabla_{\qb}G_d(\rmi u;\qb,\bmu)|_{\qb=\0}-\nabla_{\qb}g(\rmi u;\qb,\bmu)|_{\qb=\0}\rVert_2]=0,\\
&\lim_{d\rightarrow+\infty}\EE[\lVert \nabla_{\qb}^2G_d(\rmi u;\qb,\bmu)|_{\qb=\0}-\nabla_{\qb}^2g(\rmi u;\qb,\bmu)|_{\qb=\0}\rVert_{\op}]=0.\end{align*}
\end{proposition}
The proof of Proposition~\ref{prop:substitutiongd-m} utilizes the key observation that
$ \nabla_{\zb}L(\xi,\zb;\qb,\bmu)|_{\zb=\mb}\equiv\0$. 
We omit the details here since it is similar to the proof of Proposition~\ref{prop:substitutiongd}.

\subsection{Step 4: complete the proof} 
\label{step4-m}
Similar to the previous proof of Theorem~\ref{thm:mainthm}, we give the following proposition to ensure the existence and uniqueness of $\bnu$ defined in Section~\ref{sec:generalcase}.
\begin{proposition}\label{prop:existence_uniqueness_nu-m}
There exists a unique analytic function $\bnu= [\nu_1,\ldots,\nu_{K+1}]^\T: \bbC_+ \rightarrow \bbC_+^{K+1}$ such that:
\begin{enumerate}[leftmargin = *]
    \item For any $\xi\in\bbC_+$, $\bnu(\xi;\bmu)$ is a solution to $\bnu$\textit{-system}~\eqref{eq:implicit2-m}.
    \item There exists $\xi_0>0$, such that $| \nu_j(\xi;\bmu) | \leq 2\psi_j / \xi_0$, for all $\xi$ with $\Im(\xi) \geq \xi_0$ and $j=1,\ldots,K+1$.
\end{enumerate}
Moreover, it holds that $ \bnu(\xi;\bmu) = \mb(\xi;\0,\bmu) $ for all $\xi\in\bbC_+$.
\end{proposition}
\begin{proof}[Proof of Proposition~\ref{prop:existence_uniqueness_nu-m}]
By Proposition~\ref{prop:implicit1-m}, the existence is directly verified as   $\mb(\xi;\0,\bmu)$ is a solution. For the uniqueness of $\bnu$,
note that $\bnu(\xi;\bmu)$ and $\mb(\xi;\0,\bmu)$ are analytic. By Lemma~\ref{lemma:uniquesolution-m},  they are identical on $\{\xi:\Im(\xi)>\xi_0\}$ with some sufficiently large $\xi_0$. The uniqueness of $\bnu$ thus results  from the uniqueness of the analytic continuation.
\end{proof}
Proposition~\ref{prop:existence_uniqueness_nu-m} justifies the definition of $\bnu(\xi;\bmu)$ in Section~\ref{sec:generalcase} by demonstrating its existence and uniqueness. Moreover, it also relates $\bnu(\xi;\bmu)$ to the function $\mb(\xi;\qb,\bmu)$ introduced in step 3 of the proof. With this result, we can finalize the proof of
Theorem~\ref{thm:mainthm-m} as follows. 

\begin{proof}[Proof of Theorem~\ref{thm:mainthm-m}]
Let
\begin{equation}\label{eq:finalgd-m}
    \begin{split}
        \cR(\lambda,\bpsi,\bmu,F_1,\tau)
        & =
        F_1^2\cdot \big[ 1-\partial_{q_1} g(\xi^*;\qb,\bmu)  -\partial^2_{q_4,q_5}g(\xi^*;\qb,\bmu)  -\partial^2_{q_2,q_5}g(\xi^*;\qb,\bmu) \big ] \big|_{\qb = \mathbf{0}}\\
        &\quad  -\tau^2\cdot \big[\partial^2_{q_3,q_4}g(\xi^*;\qb,\bmu) + \partial^2_{q_2,q_3}g(\xi^*;\qb,\bmu) \big ] \big|_{\qb = \mathbf{0}},
    \end{split}
\end{equation}
where $g$ is defined in \eqref{eq:gd-m}. 
Then by Propositions~\ref{prop:decomp-expectation-m}, \ref{prop:connectGd-m} and \ref{prop:substitutiongd-m}, we have 
\begin{align*}
    \begin{split}
        \EE_{\Xb,\bTheta,\bvarepsilon}&\Big|R_d(\Xb,\bTheta,\lambda,\bbeta_d,\bvarepsilon)-\cR(\lambda,\bpsi,\bmu,F_1,\tau)\Big|=o_d(1).
    \end{split}
\end{align*}
Recall equations~\eqref{eq:gd-m} and \eqref{eq:L-function-m}, for any $\xi\in\bbC_+$ we have
$$
\nabla_{\zb}L(\xi,\zb;\qb,\bmu)|_{\zb=\mb}=\0.
$$
Here $\zb=[z_1,\ldots,z_{K+1}]^\T$. Then from the formula for implicit differentiation, we have
\begin{align}\label{eq:g_derivative1-m}
\partial_{q_1} g(\xi^*;\qb,\bmu)|_{\qb=0}=\partial_{q_1} L(\xi^*,\zb;\qb,\bmu)|_{\zb=\bnu^*,\qb=0}=\frac{2\nu_{K+1}^*M_N}{M_D}.
\end{align}
We remind readers  that $M_N=\sum\limits_{\sfc=1}^K\nu_\sfc^*\mu_{\sfc,1}^2$, $M_D=\nu_{K+1}^*M_N-1$ and $\bnu^*=\mb(\xi^*;\0,\bmu)$.
Denote $\ub=(q_2,q_3,q_4,q_5,\zb)$, and construct the matrix $\Wb(\bnu^*,\bmu)=\nabla_{\ub}^2 L(\xi^*,\zb;\qb,\bmu)|_{\zb=\bnu^*,\qb=0}$. 
Note that \eqref{eq:derivative_g2} and \eqref{eq:m_derivative} in our proof of DRFMs still hold for the case of MRFMs. Therefore we have (to simplify the writing, we drop the arguments in the matrix $\Wb$):
{\small
\begin{align}
        &\frac{\partial^2 g(\xi^*;\qb,\bmu)}{\partial q_2 \partial q_5}\Big|_{\qb=0}=\Wb_{1,4}-\Wb_{1,[5:(K+5)]}\Big(\Wb_{[5:(K+5)],[5:(K+5)]}\Big)^{-1}\Wb_{[5:(K+5)],4},\label{eq:g_derivative2-m}\\
        &\frac{\partial^2 g(\xi^*;\qb,\bmu)}{\partial q_3\partial q_4}\Big|_{\qb=0}=\Wb_{2,3}-\Wb_{2,[5:(K+5)]}\Big(\Wb_{[5:(K+5)],[5:(K+5)]}\Big)^{-1}\Wb_{[5:(K+5)],3},\label{eq:g_derivative3-m}\\
        &\frac{\partial^2 g(\xi^*;\qb,\bmu)}{\partial q_2\partial q_3}\Big|_{\qb=0}=\Wb_{1,2}-\Wb_{1,[5:(K+5)]}\Big(\Wb_{[5:(K+5)],[5:(K+5)]}\Big)^{-1}\Wb_{[5:(K+5)],2},\label{eq:g_derivative4-m}\\
        &\frac{\partial^2 g(\xi^*;\qb,\bmu)}{\partial q_4\partial q_5}\Big|_{\qb=0}=\Wb_{3,4}-\Wb_{3,[5:(K+5)]}\Big(\Wb_{[5:(K+5)],[5:(K+5)]}\Big)^{-1}\Wb_{[5:(K+5)],4}.\label{eq:g_derivative5-m}
\end{align}
}
Similar to the case of DRFMs, we have
\begin{align*}
   & \Wb_{1,4}=\Wb_{2,3}=\Wb_{1,2}=0,   \quad  \Wb_{3,4}=-\frac{\nu_{K+1}^{*2}M_N^2}{M_D^2}, \\
   &  \Vb=\Wb_{[5:(K+5)],[1:4]}=\Wb_{[1:4],[5:(K+5)]}^\T,  \quad \text{and}  \quad \Hb=\Big(\Wb_{[5:(K+5)],[5:(K+5)]}\Big). 
\end{align*}
Plugging \eqref{eq:g_derivative1-m} and \eqref{eq:g_derivative2-m}-\eqref{eq:g_derivative5-m} into \eqref{eq:finalgd-m} proves Theorem~\ref{thm:mainthm-m}.
\end{proof}


\section{Other key factors affecting the risk curve}
\label{sec:otherfactor}
Here we investigate several other factors that affect the shape of the risk curve. By studying how these factors affect the risk, we aim to provide a clearer understanding of Proposition~\ref{prop:property_lbd=0}, Proposition~\ref{prop:property_psi} and the triple descent phenomena. Our analysis also shows how we can design DRFMs to achieve a specific risk curve shape. Unlike  \citet{chen2021multiple} which requires designing a specific data distribution, our study shows that various risk curves can be achieved by different random feature models on a fixed data distribution.





\noindent\textbf{The regularization parameter $\lambda$. }We investigate how the regularization parameter  $\lambda$ affect the shape of the risk curve. We again use the same experiment setup as in Section~\ref{subsec:tripleinDRFM}, expect that we focus on activation functions $\ELU(3x)$ and $\ReLU(x/4)$, and calculate the risk curves w.r.t. different regularization parameters $\lambda = 10^{-1},10^{-2}, 10^{-3}$ and $10^{-4}$.

\begin{figure}[t!]
    \centering
    \includegraphics[width=0.95\textwidth]{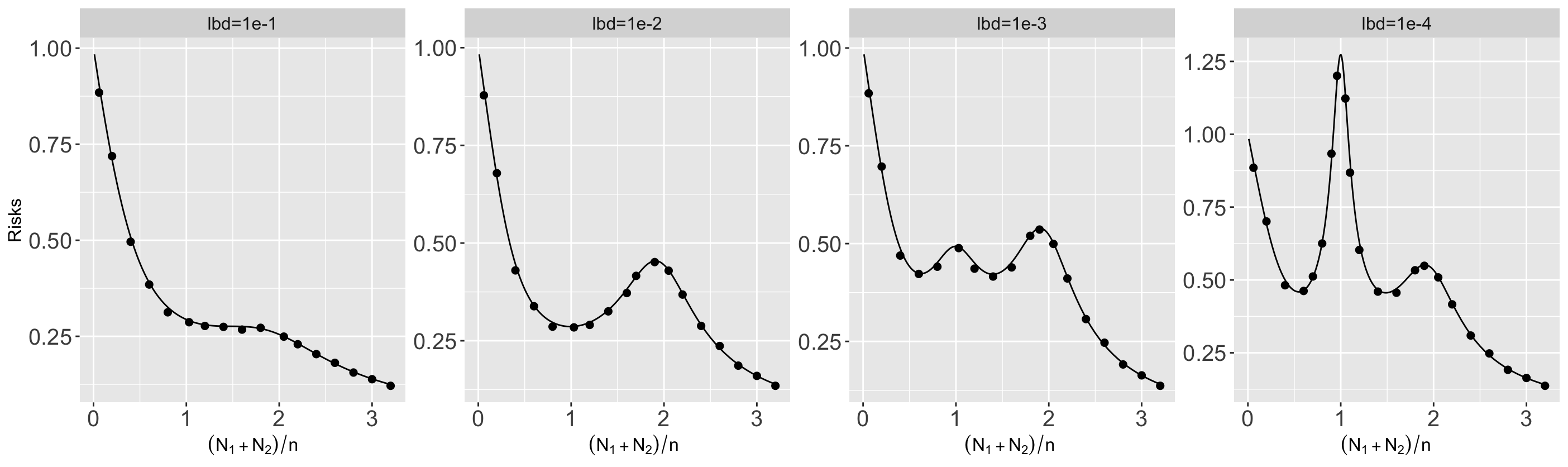}\vspace{-5pt}
    
    \hspace{12pt}(a)\hspace{86pt}\hspace{10pt}(b)\hspace{80pt}\hspace{15pt}(c)\hspace{80pt}\hspace{15pt}(d)
    \vspace{-10pt}
    \caption{Risk curves of DRFMs trained with different regularization parameters. The plots show both the asymptotic excess risks (curves) and empirical excess risks (dots). From (a) to (d), we set $\lambda = 10^{-1},10^{-2}, 10^{-3}$ and $10^{-4}$, respectively. The activation functions are chosen as $\sigma_1(x)=\ELU(3x)$ and $\sigma_2(x)=\ReLU(x/4)$ in all these experiments.}
    \label{fig:lbd}
    \vspace{-10pt}
\end{figure}

The results are given in Figure~\ref{fig:lbd}.  Note that Proposition~\ref{prop:property_lbd=0} holds under the condition $\lambda$ tends to 0. When the regularization parameter $\lambda$ is large, the risk  decreases with the model complexity parameter $c\sim(N_1 + N_2 )/ n$. As $\lambda$ decreases, the peak at $c = 2$ first appears, and then the peak at $c = 1$ also appears when $\lambda = 10^{-3}$. Finally when $\lambda = 10^{-4}$, the risk around $c = 1$ becomes very high. From these experiments, we can conclude that (i) Double/triple descent happens particularly when there is no regularization or when the regularization is very weak. (ii) the risk value of the first peak around $c = 1$ is more sensitive to $\lambda$ then that of the second peak.

\smallskip\noindent\textbf{Signal-to-noise ratio.} We also study how the signal-to-noise ratio (SNR) in the data, which we define as $\| \bbeta_1 \|_2 /\tau$, affects the shape of the risk curve. We again use the same experimental setup as in Section~\ref{subsec:tripleinDRFM}, except that (i) we focus on activation functions $(\ELU(3x)$ and $\ReLU(x/4))$, and (ii) we perform experiments with different values of $\| \bbeta_1 \|_2 = F_1$ and $\tau$.



\begin{figure}[t!]
    \centering
    \includegraphics[width=0.95\textwidth]{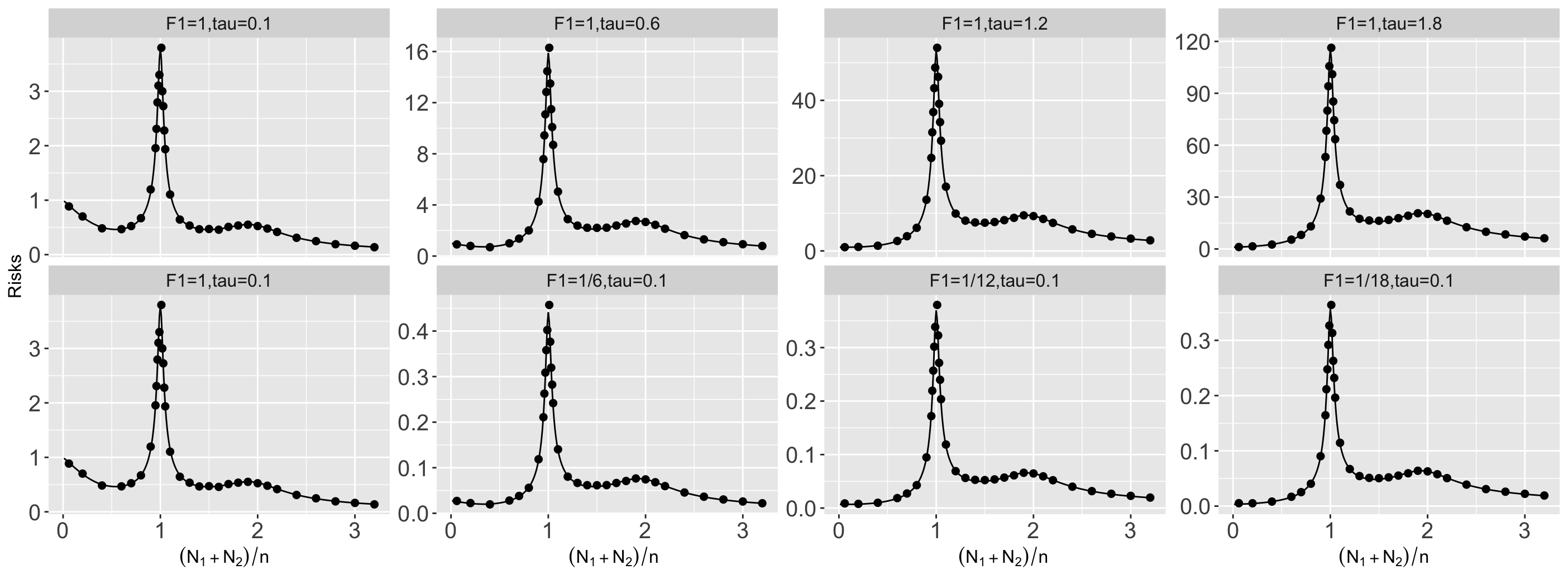}\vspace{-5pt}
    
    \hspace{6pt}(a)\hspace{86pt}\hspace{10pt}(b)\hspace{80pt}\hspace{18pt}(c)\hspace{80pt}\hspace{20pt}(d)
    \vspace{-10pt}
    \caption{Risk curves of DRFMs under different SNRs. The plots show both the asymptotic excess risks (curves) and empirical excess risks (dots).
    In the top row, we set $\lVert\bbeta_1\rVert_2=1$ and $\tau= 0.1, 0.6, 1.2$ and $1.8$ (from (a) to (d)).
    In the bottom row, we set $\tau=0.1$ and $\lVert\bbeta_1\rVert_2= 1, 1/6, 1/12$ and $1/18$ (from (a) to (d)). 
    The parameter values are chosen such that the two figures in each column have the same SNR.}
    \label{fig:SNR}
    \vspace{-10pt}
\end{figure}

The results are given in Figure~\ref{fig:SNR}.  We first see that the risk curves in each column have the same shapes. This matches our theoretical result that the risk has the form $ R = \tau^2( a \cdot \mathrm{SNR} + b) $ for some positive functions $a,b$ depending on the other parameters. Moreover, the SNR has a particularly high impact on the trend of the risks in the under-parameterized regime ($(N_1+N_2)/n<1$) and the highly over-parameterized regime ($(N_1+N_2)/n>2$, shown in Proposition~\ref{prop:property_psi}). Specifically, in column (a) when the SNR is large, we can see that the lowest risk is achieved in the highly over-parameterized regime; on the other hand, in columns (c) and (d) when the SNR is relatively small, the lowest risk is achieved in the under-parameterized regime.

\newpage

$ $

\newpage
The following sections in the appendix  are technical details, and  we briefly summarize the structure below:
\begin{itemize}
    \item In Section~\ref{sec:appendixDecomposition}, we give the proof of Proposition~\ref{prop:decomp-expectation} which introduces the decomposition of the asymptotic excess risk.
    \item In Section~\ref{sec:appendixconnectedGd}, we present the proof of Proposition~\ref{prop:connectGd} by showing how the key terms in the risk decomposition are related to the logarithmic potential of the linear pencil matrix.
    \item In Section~\ref{sec:appendixuniquesolution}, we establish basic properties of the fixed point equation  \eqref{eq:m123} to justify the definition of $\mb(\xi;\qb,\bmu)$ below Definition~\ref{def:implicit1}.
    \item In Section~\ref{sec:appendixpropimplicit1}, we provide the proof of Proposition~\ref{prop:implicit1}, which extends the definition of $m(\xi;\qb,\bmu)$ to $\bbC_{+}$ and shows that it is the asymptotic limit of $M_d(\xi;\qb,\bmu)$ as $d\rightarrow \infty$.
    \item In Section~\ref{sec:appendixsubstitution}, we give the proof of Proposition~\ref{prop:substitutiongd} by relating the logarithmic potential $G_d(\xi;\qb,\bmu)$ to $g(\xi;\qb,\bmu)$ in \eqref{def:gd}.
    \item In Section~\ref{sec:appendix_nu}, we prove Proposition~\ref{prop:existence_uniqueness_nu} to justify the definition of $\bnu(\xi)$ as the unique solution to the system \eqref{eq:implicit2}.
    \item In Section~\ref{sec:addionalproofs-m}, we display the proof of the lemmas and propositions given in Appendix~\ref{sec:appendixthm-m}.
\end{itemize}

\section{Proof of Proposition~\ref{prop:decomp-expectation}}
\label{sec:appendixDecomposition}
Proposition~\ref{prop:decomp-expectation} gives a decomposition of the risk $R_d(\Xb,\bTheta,\lambda,\bbeta_d,\bvarepsilon)$. To prove this decomposition, we first introduce some additional notations and preliminary lemmas. 
\begin{definition}
\label{def:notation2}
Define
\begin{equation*}
    \begin{split}
      \Vb_0(F_0)=F_0\EE_{\xb}\big[\bsigma(\xb) \big] \in\RR^{N\times 1},~ \Vb(\bbeta_{1,d})=\EE_{\xb}\big[\bsigma(\xb)\xb^\T \bbeta_{1,d}\big]\in\RR^{N\times 1},~ \Ub=\EE_{\xb}\big[\bsigma(\xb)\bsigma(\xb)^\T \big]\in\RR^{N\times N}, 
    \end{split}
\end{equation*}
where $\xb$ is a random vector uniformed distributed on the sphere $\sqrt{d}\cdot\SSS^{d-1}$ and $\bsigma(\xb)$ is defined in Definition~\ref{def:notation1}. 
\end{definition}
Note that by the definition of $\bsigma(\xb)$ in Definition~\ref{def:notation1},  $\bsigma(\xb)$ also depends on the random feature parameter matrix $\bTheta$. Therefore, $\Vb_0(F_0)$, $\Vb(\bbeta_{1,d})$ and $\Ub$ also depends on $\bTheta$. 
Now with these notations, and by the definition of $\hat\ab$ in \eqref{eq:def_hata1}, we can rewrite the risk as follows:
\begin{align}
\label{eq:Rdecomposition1}
  R_d(\Xb,\bTheta,\lambda,\bbeta_d,\bvarepsilon)=&\EE_{\xb}\big[\xb^\T \bbeta_{1,d}+F_0-\hat{\ab}^\T \bsigma(\xb)\big]^2\nonumber\\
=&F_0^2+F_{1,d}^2-2\yb^\T \Zb\bUpsilon[\Vb(\bbeta_{1,d})+\Vb_0(F_0)]/\sqrt{d}+\yb^\T \big[\Ub\big]_\Zb \yb/d.
\end{align}
Therefore, to prove Proposition~\ref{prop:decomp-expectation}, it suffices to further decompose the terms $\Ub,\Vb(\bbeta_{1,d})$ and $\Vb_0(F_0)$. 
To handle these terms, we consider the Gegenbauer decomposition \citep{hua1963harmonic} of the nonlinear activation functions. For $\sfc=1,2$, let $\lambda_{d,k}(\sigma_{\sfc})$ be the coefficients of the Gegenbauer decomposition of $\sigma_{\sfc}$, i.e.,
\[\sigma_{\sfc}(x)=\sum_{k=0}^{+\infty}\lambda_{d,k}\big(\sigma_{\sfc}\big)B(d,k)\cdot Q_k^{(d)}(\sqrt{d}\cdot x),\]
where $B(d,0)=1$,
$B(d,k)=k^{-1}(2k+d-2)\binom{k+d-3}{k-1}$ with $k\geq1$, and $Q_k^{(d)}$, $k\in \NN$ are the Gegenbauer polynomials forms an orthogonal basis on $L^2([-d,d],\tau_d)$. $\tau_d$ is the distribution of $\la\xb_1,\xb_2\ra$ where $\xb_1,\xb_2\sim\sqrt{d}\cdot\Unif(\SSS^{d-1})$.
Then  define 
\begin{align}
  &\bLambda_{d,k}=\diag\big(\lambda_{d,k}(\sigma_1)\Ib_{N_1},\lambda_{d,k}(\sigma_2)\Ib_{N_2}\big),\quad k\in\NN=\{0,1,...\}. \label{eq:def_lambda_k}
\end{align}
The following lemma decomposes the three terms in Definition~\ref{def:notation2}.
\begin{lemma}[\xr{ Decomposition on nonlinear activation functions}]
\label{lemma:decomp1}
With  $\Mb_1$ and $\Mb_2$  in Definition~\ref{def:notation1}, and $\bLambda_{d,k}$ in equation~\eqref{eq:def_lambda_k}, we have 
\begin{align*}
&\Vb_0(F_0)=F_0\bLambda_{d,0}\1_N,\\
&\Vb(\bbeta_{1,d})=\bLambda_{d,1}\bTheta\bbeta_{1,d} = \bigg(\frac{\Mb_{1}+\bDelta'}{\sqrt d}\bigg)\bTheta\bbeta_{1,d} ,\\
&\Ub = \bLambda_{d,0}\1_N\1_N^\T \bLambda_{d,0}+\Mb_1\frac{\bTheta\bTheta^\T }{d}\Mb_1+ \Mb_2\Mb_2  + \bDelta.
\end{align*}
where the remainder matrices $\bDelta, \bDelta'$ satisfy $\EE\lVert\bDelta\rVert_{\op}^2\vee \EE\lVert\bDelta'\rVert_{\op}^2=o_d(1)$.
\end{lemma}
Lemma~\ref{lemma:decomp1} is proved  in Section~\ref{sec:appendixlemma1}. 
Plugging the decompositions in Lemma~\ref{lemma:decomp1} into \eqref{eq:Rdecomposition1} will then give a decomposition of the risk consisting of multiple terms. 
The next lemma establishes useful  moment estimations for some of the terms in \eqref{eq:Rdecomposition1}, which helps us get rid of the negligible terms in the decomposition.

\begin{lemma} [\xr{Negligible terms}]
\label{lemma:decomp2}
For any fixed $k\in\NN\backslash\{0\}$, let ${\bGamma}_1\in\RR^{N\times N}$ and ${\bGamma}_2\in\RR^{n\times n}$ 
be symmetric random matrices with $\big[\EE\lVert{{\bGamma}_j }\rVert_{\op}^{k}\big]^{1/k}=O_d(1)$, $j = 1,2$. Define
\begin{equation*}
    \begin{split}
     \cB=&\frac{1}{d}   \1_n^\T \big[{\bGamma}_1\big]_\Zb \1_n,\\
     \cC=&1-\frac{2}{\sqrt d}\tr(\bLambda_{d,0}\1_N\1_n^\T  \Zb \bUpsilon)+\frac{1}{d} \1_n^\T \big[\bLambda_{d,0}\1_N\1_N^\T \bLambda_{d,0}\big]_\Zb \1_n ,\\
     \cD=&\frac{1}{d}\tr(\big[\bLambda_{d,0}\1_N\1_N^\T \bLambda_{d,0}\big]_\Zb \bGamma_2   ),
    \end{split}
\end{equation*}
where $\bLambda_{d,0}$ is defined in equation~\eqref{eq:def_lambda_k}. Then if $\sum_{\sfc}\mu_{\sfc,0}^2>0$, for any fixed $\lambda>0$, there exists a constant $C>0$ such that  
$$
\big(\EE |\cB|^k\big)^{1/k}  \vee \EE |\cC| \vee \big(\EE |\cD|^k\big)^{1/k}  =O_d\Big(d^{-1} e^{C \sqrt{\log d}} \Big)=o_d(1).  
$$
If $\sum_{\sfc}\mu_{\sfc,0}^2=0$, it still holds that $\big(\EE |\cD|^k\big)^{1/k}=o_d(1)$.
\end{lemma}

\if UT
\begin{lemma} 
\label{lemma:decomp3}
Let ${\bGamma}\in\RR^{n\times n}$ any symmetric random matrix with $\big[\EE\lVert{\bGamma }\rVert_{\op}^{2k}\big]^{1/2k}=O_d(1)$ for any fixed $k\in\NN\backslash\{0\}$. Define
$$\begin{aligned}
\cD=&\frac{1}{d}\tr(\big[\bLambda_{d,0}\1_N\1_N^\T \bLambda_{d,0}\big]_\Zb \bGamma   ),
\end{aligned}$$
with $\1_N\in\RR^N$ whose elements are all $1$. $\Zb$ and $\bUpsilon$ are defined in Definition~\ref{def:notation1}, $\bLambda_{d,0}$ is defined in equation~\eqref{eq:def_lambda_k}. Then for any fixed $\lambda>0$, there exists constant $C>0$, s.t
$$\begin{aligned}
\big(\EE& |\cD|^k\big)^{1/k}=O_d\Big(\frac{\exp (C\sqrt{\log d})}{d}\Big)=o_d(1).
\end{aligned}$$
\end{lemma}
\fi 
The proof of the lemma is given in Section~\ref{sec:appendixdecomp2}. To further decompose and calculate the risk, we also need to study the impact of fixed vector $\bbeta_{1,d}$ on the risk. To do so, we aim to show that the risk only depends on $F_{1,d}$ ($=\lVert\bbeta_{1,d}\rVert_2$) due to rotation invariance of the learning problem.  The result is given in the following lemma.

\begin{lemma}[\xr{Preliminary for proof of   Proposition~\ref{prop:decomp-expectation}}]
\label{lemma:replace} 
Suppose $\tilde{\bbeta}_{1,d}\sim\Unif( F_{1,d}\cdot\SSS^{d-1})$ is independent of $(\Xb,\bTheta,\bvarepsilon)$, and denote $\tilde{\bbeta}_d=[F_0,\tilde{\bbeta}_{1,d}^\T ]^\T $. Then for any fixed $\bbeta_{1,d}$, under the assumptions of Proposition~\ref{prop:decomp-expectation}, we have
\begin{align*}
        &\EE_{\Xb,\bTheta,\bvarepsilon} \big| R_d(\Xb,\bTheta,\lambda,\bbeta_d,\bvarepsilon) - \overline{R}_d(\Xb,\bTheta,\lambda,F_{1,d},\tau)\big|\\
        &\qquad\qquad\qquad\qquad\qquad\qquad=\EE_{\Xb,\bTheta,\bvarepsilon,\tilde{\bbeta}_{d}} \big|R_d(\Xb,\bTheta,\lambda,\tilde{\bbeta}_d,\bvarepsilon)-\overline{R}_d(\Xb,\bTheta,\lambda,F_{1,d},\tau)\big|,\\
        &\EE_{\Xb,\bTheta}\big[\Var_{\tilde{\bbeta}_{d},\bvarepsilon}(R_d(\Xb,\bTheta,\lambda,\tilde{\bbeta}_d,\bvarepsilon))\big]=o_d(1).
\end{align*}

\end{lemma}

The proof of the lemma is given in Section~\ref{sec:appendixreplace}. Based on the above lemmas, we are ready to present the proof of Proposition~\ref{prop:decomp-expectation} as follows.

\begin{proof}[Proof of Proposition~\ref{prop:decomp-expectation}] Let  $\tilde{\bbeta}_d=[F_0,\tilde{\bbeta}_{1,d}^\T ]^\T $ with $\tilde{\bbeta}_{1,d}\sim\Unif( F_{1,d}\cdot\SSS^{d-1})$. Then we have
\begin{align*}
    &\EE_{\Xb,\bTheta,\bvarepsilon} \big| R_d(\Xb,\bTheta,\lambda,\bbeta_d,\bvarepsilon) - \overline{R}_d(\Xb,\bTheta,\lambda,F_{1,d},\tau)\big|\\
    &\qquad\qquad =\EE_{\Xb,\bTheta,\bvarepsilon,\tilde{\bbeta}_{d}} \big|R_d(\Xb,\bTheta,\lambda,\tilde{\bbeta}_d,\bvarepsilon)-\overline{R}_d(\Xb,\bTheta,\lambda,F_{1,d},\tau)\big|\\
    &\qquad\qquad \leq\EE_{\Xb,\bTheta,\bvarepsilon,\tilde{\bbeta}_{d}} \big|R_d(\Xb,\bTheta,\lambda,\tilde{\bbeta}_d,\bvarepsilon)-\EE_{\bvarepsilon,\tilde{\bbeta}_{d}}R_d(\Xb,\bTheta,\lambda,\tilde{\bbeta}_d,\bvarepsilon)\big|\\
    &\qquad\qquad \quad+\EE_{\Xb,\bTheta}\big|\EE_{\bvarepsilon,\tilde{\bbeta}_{d}}R_d(\Xb,\bTheta,\lambda,\tilde{\bbeta}_d,\bvarepsilon)-\overline{R}_d(\Xb,\bTheta,\lambda,F_{1,d},\tau)\big|\\
    &\qquad\qquad \leq \EE_{\Xb,\bTheta} \Big[ \sqrt{ \Var_{\tilde{\bbeta}_{d},\bvarepsilon}(R_d(\Xb,\bTheta,\lambda,\tilde{\bbeta}_d,\bvarepsilon)) } \Big]\\
    &\qquad\qquad \quad+\EE_{\Xb,\bTheta}\big|\EE_{\bvarepsilon,\tilde{\bbeta}_{d}}R_d(\Xb,\bTheta,\lambda,\tilde{\bbeta}_d,\bvarepsilon)-\overline{R}_d(\Xb,\bTheta,\lambda,F_{1,d},\tau)\big|\\
    &\qquad\qquad \leq \sqrt{ \EE_{\Xb,\bTheta} \big[ \Var_{\tilde{\bbeta}_{d},\bvarepsilon}(R_d(\Xb,\bTheta,\lambda,\tilde{\bbeta}_d,\bvarepsilon)) \big] }\\
    &\qquad\qquad \quad+\EE_{\Xb,\bTheta}\big|\EE_{\bvarepsilon,\tilde{\bbeta}_{d}}R_d(\Xb,\bTheta,\lambda,\tilde{\bbeta}_d,\bvarepsilon)-\overline{R}_d(\Xb,\bTheta,\lambda,F_{1,d},\tau)\big|\\
    &\qquad\qquad = o_d(1) + \EE_{\Xb,\bTheta}\big|\EE_{\bvarepsilon,\tilde{\bbeta}_{d}}R_d(\Xb,\bTheta,\lambda,\tilde{\bbeta}_d,\bvarepsilon)-\overline{R}_d(\Xb,\bTheta,\lambda,F_{1,d},\tau)\big|,
\end{align*}
where the first equality follows by the first equation in Lemma~\ref{lemma:replace}, the first inequality follows by triangle inequality, the second and third inequalities are by Jensen's inequality, and the last equality follows by Lemma~\ref{lemma:replace} again.
Therefore, to prove the proposition, it suffices to show that
\begin{align*}
\EE_{\Xb,\bTheta}\big|\EE_{\bvarepsilon,\tilde{\bbeta}_{d}}R_d(\Xb,\bTheta,\lambda,\tilde{\bbeta}_d,\bvarepsilon)-\overline{R}_d(\Xb,\bTheta,\lambda,F_{1,d},\tau)\big|=o_d(1).
\end{align*}
Similar to \eqref{eq:Rdecomposition1}, we have
\begin{align}
    R_d(\Xb,\bTheta,\lambda,\tilde\bbeta_{d},\bvarepsilon)=F_0^2+F_{1,d}^2-\frac{2\tilde\yb^\T \Zb\bUpsilon(\Vb(\tilde\bbeta_{1,d})+\Vb_0(F_0))}{\sqrt{d}}+\frac{\tilde\yb^\T \big[\Ub\big]_\Zb \tilde\yb}{d},\label{eq:Rd_def_tilde}
\end{align}
where $\Tilde{\yb}=\1_nF_0+\Xb\tilde\bbeta_{1,d}+\bvarepsilon$. 
From Lemma~\ref{lemma:decomp1}, we further have
\begin{align}\Vb_0(F_0)F_0=F_0^2\bLambda_{d,0}\1_N , \quad \EE_{\tilde{\bbeta}_{1,d}}\big(\Vb(\tilde{\bbeta}_{1,d})\tilde{\bbeta}_{1,d}^\T \big)=F_{1,d}^2\bigg(\frac{\Mb_{1}+\bDelta'}{\sqrt d}\bigg)\frac{\bTheta}{d},\label{eq:decomp1_proof}
\end{align}
and 
\begin{align}
 \Ub = \bLambda_{d,0}\1_N\1_N^\T \bLambda_{d,0}+\Mb_1\frac{\bTheta\bTheta^\T }{d}\Mb_1+ \Mb_2\Mb_2+\bDelta.\label{eq:decomp2_proof}
\end{align}
By \eqref{eq:Rd_def_tilde}, \eqref{eq:decomp1_proof}, \eqref{eq:decomp2_proof} and the definition of $\overline{R}_d(\Xb,\bTheta,\lambda,F_{1,d},\tau)$, we obtain the following equation with direct calculation:
\begin{align*}
&\EE_{\tilde\bbeta_{d},\bvarepsilon}R_d(\Xb,\bTheta,\lambda,\tilde\bbeta_{d},\bvarepsilon)-\overline{R}_d(\Xb,\bTheta,\lambda,F_{1,d},\tau)=\\
    &\qquad \underbrace{F_{0}^2-\frac{2 F_{0}^2}{\sqrt d}\tr\big(\bLambda_{d,0}\1_N\1_n^\T  \Zb \bUpsilon\big)+\frac{ F_{0}^2}{d}\tr\big( [\bLambda_{d,0}\1_N\1_N^\T \bLambda_{d,0}]_\Zb \1_n\1_n^\T \big)}_{I_1}\\
&\qquad -\underbrace{ \frac{2 F_{1,d}^2}{ d}\tr\bigg(\bDelta'\frac{\bTheta \Xb^\T }{d} \Zb \bUpsilon\bigg)}_{I_2}+\underbrace{\frac{ F_{0}^2}{d}
\tr\bigg( \bigg[\Mb_1\frac{\bTheta\bTheta^\T }{d}\Mb_1+\Mb_2\Mb_2\bigg]_\Zb  \1_n\1_n^\T \bigg)}_{I_3}\\
&\qquad +\underbrace{ \frac{ F_{1,d}^2}{d}\tr\bigg( [\bLambda_{d,0}\1_N\1_N^\T \bLambda_{d,0}]_\Zb  \frac{\Xb\Xb^\T }{d}\bigg)}_{I_4}+\underbrace{\frac{\tau^2}{d}
\tr\big(\big[  \bLambda_{d,0}\1_N\1_N^\T \bLambda_{d,0}   \big]_\Zb\big)}_{I_5}\\
&\qquad +\underbrace{\frac{ F_{1,d}^2}{d}  \tr\bigg( [\bDelta]_\Zb \frac{\Xb\Xb^\T }{d}\bigg)}_{I_6}+
   \underbrace{\frac{ F_{1,d}^2}{d}  \tr [\bDelta]_\Zb}_{I_7}
    +\underbrace{\frac{ F_{0}^2}{d} \tr\big( [\bDelta]_\Zb  \1_n\1_n^\T \big)}_{I_8}. 
\end{align*}
We now show that all the terms $I_1,\ldots,I_8$ on the right hand side above are negligible terms. We note that by definition, $\lVert\Zb\bUpsilon\lVert_{\op} = \lVert\Zb (\Zb^\T\Zb + \lambda \Ib)^{-1}\lVert_{\op} \leq 1/ (2\sqrt{\lambda}) $ is deterministically bounded. Therefore we have
\begin{align*}
    \EE|I_2|\leq2F_{1,d}^2\cdot\EE\bigg\|{ \bigg(\bDelta'\frac{\bTheta \Xb^\T }{d} \Zb \bUpsilon\bigg)}\bigg\|_{\op}\leq O_d\Big(\frac{1}{2\sqrt{\lambda}}\Big)\cdot\big(\EE\lVert\bDelta'\rVert_{\op}^2\big)^{\frac{1}{2}}\cdot\Big(\EE\Big\lVert\frac{\bTheta \Xb^\T }{d}\Big\rVert_{\op}^2\Big)^{\frac{1}{2}}=o_d(1),
\end{align*}
where the last equality follows by $\EE\|\bDelta'\|_{\op}^2 = o_d(1)$ in Lemma~\ref{lemma:decomp1}. Moreover, by definition, we have
\begin{align*}
    \| [\bDelta]_{\Zb} \|_{\op} = \| \Zb \bUpsilon  \bDelta (\Zb \bUpsilon )^\T  \|_{\op} \leq \frac{1}{4\lambda} \| \bDelta \|_{\op}.
\end{align*}
Therefore, by Lemma~\ref{lemma:decomp1} that $\EE\|\bDelta\|_{\op}^2 = o_d(1)$, we have
\begin{align*}
    &\EE|I_6| \leq F_{1,d}^2 \cdot \EE \bigg\| [\bDelta]_\Zb \frac{\Xb\Xb^\T }{d} \bigg\|_{\op}  \leq F_{1,d}^2 \cdot \EE \bigg[ \| [\bDelta]_{\Zb} \|_{\op} \cdot  \bigg\|  \frac{\Xb\Xb^\T }{d} \bigg\|_{\op} \bigg] = O_d(\EE\|\bDelta\|_{\op}^2) = o_d(1), \\
    &\EE|I_7| \leq F_{1,d}^2 \cdot \EE \| [\bDelta]_\Zb \|_{\op} \leq \frac{F_{1,d}^2}{4\lambda} \cdot  \EE \| \bDelta \|_{\op} = o_d(1).
\end{align*}
For the remaining terms, we discuss them according to the value of $F_0$. When $F_0=0$, it is clear that $I_1=I_3=0$. Note that under this situation, the condition $\sum_{\sfc}\mu_{\sfc,0}^2>0$ in Lemma~\ref{lemma:decomp2}  may not hold. If $\sum_{\sfc}\mu_{\sfc,0}^2=0$, from Lemma~\ref{lemma:decomp2}, it still holds that
\begin{align*}
    \EE |I_4|=o_d(1),\qquad
    \EE |I_5|=o_d(1).
\end{align*}
Therefore, when   $F_0=0$, Proposition~\ref{prop:decomp-expectation} holds.

When $F_0\neq0$, $\sum_{\sfc}\mu_{\sfc,0}^2>0$ holds from Assumption~\ref{assump3}, the result for $\cC$ in Lemma~\ref{lemma:decomp2} gives the bound for $I_1$, the result for $\cB$ in Lemma~\ref{lemma:decomp2} gives the bounds on $I_3$ and $I_8$, and the result for $\cD$ in Lemma~\ref{lemma:decomp2} gives the bounds on $I_4$ and $I_5$. 
Therefore we have
\begin{align*}
    \EE_{\Xb,\bTheta}\big|\EE_{\tilde\bbeta_{d},\bvarepsilon}R_d(\Xb,\bTheta,\lambda,\tilde\bbeta_{d},\bvarepsilon)-\overline{R}_d(\Xb,\bTheta,\lambda,F_{1,d},\tau)\big|=o_d(1),
\end{align*}
which proves Proposition~\ref{prop:decomp-expectation}.
\end{proof}

\subsection{Proof of Lemma~\ref{lemma:decomp1}}
\label{sec:appendixlemma1}
The proof of Lemma~\ref{lemma:decomp1} is mainly based on the decomposition of the nonlinear activation function.
We first present several classical  lemmas about Gegenbauer polynomials and their relation to Hermite polynomials. 
The following lemma can be found in \citet{mei2022generalization} (see Lemma~9.4 and its proof in the reference). 
\begin{lemma}
\label{lemma:meisongdecomp1}
Let $Q_k^{(d)}(\cdot)$, $k\in\NN$ be the Gegenbauer polynomials. The following properties hold:
\begin{enumerate}[leftmargin = *]
    \item For $\vb_1$,$\vb_2\in\sqrt{d}\cdot\SSS^{d-1}$, suppose $\xb\sim\Unif(\sqrt{d}\cdot\SSS^{d-1})$, then for $k,l \in \NN$,  $$\EE_{\xb}\big[Q_k^{(d)}(\vb_1^\T\xb)Q_l^{(d)}(\xb^\T \vb_2 )\big]\\
        =\frac{\delta_{kl}}{B(d,k)}\cdot Q_k^{(d)}(\vb_1^\T\vb_2 ),$$
    where $\delta_{kl} = 1$ if $k = l$ and  $\delta_{kl} = 0$ if $k \neq l$. 
    \item For $\bTheta_1$ and $\bTheta_2$ defined in  Section~\ref{sec:Problemsetting}, $Q_k^{(d)}(\cdot)$ the point wise function on matrices, the following equality holds: $$\begin{aligned}
&\EE\Big[\sup_{k\geq2}\big\lVert Q_k^{(d)}(\bTheta_{\sfc}\bTheta_{\sfc}^\T )-\Ib_{N_{\sfc}}\big\rVert_{\op}^2\Big]=o_d(1),\quad \sfc=1,2,\\
&\EE\Big[\sup_{k\geq2}\big\lVert Q_k^{(d)}(\bTheta_1\bTheta_2^\T )\big\rVert_{\op}^2\Big]=o_d(1).
\end{aligned}$$
\end{enumerate}
\end{lemma}
The next lemma gives the connection between the coefficients in Hermite polynomials $H_k$ and the coefficients in Gegenbauer polynomials $Q_k^{(d)}$.
\begin{lemma}[\xr{Gegenbauer decomposition}]
\label{lemma:connectioninGegenandHermite} 
Let $Q_k^{(d)}(\cdot)$, $H_k(\cdot)$, $k\in\NN$ be the Gegenbauer and Hermite polynomials respectively. For $\sfc=1,2$, suppose that $\sigma_{\sfc}(x)$ has Gegenbauer decomposition
\[
\sigma_{\sfc}(x)=\sum_{k=0}^{+\infty}\lambda_{d,k}\big(\sigma_{\sfc}\big)B(d,k)\cdot Q_k^{(d)}(\sqrt{d}\cdot x)
\] 
and Hermite polynomial decomposition 
\[
\sigma_{\sfc}(x)=\sum_{k=0}^{+\infty}\alpha_{k}\big(\sigma_{\sfc}\big)/k!\cdot H_k(x).
\]
Then for each $k\in \NN$,  $\lambda_{d,k}^2(\sigma_{\sfc})B(d,k)k!\rightarrow\alpha_{k}^2(\sigma_{\sfc})$ as $d\rightarrow+\infty$. 
\end{lemma}
The proof of Lemma~\ref{lemma:connectioninGegenandHermite} can be found in Appendix A.3 in \citet{mei2022generalization}. Note that the orthogonality of the standard Hermite polynomials ($H_1(x)=x$) implies that for $G\sim N(0,1) $,
\begin{align*}
    \EE[ H_k(G) H_l(G)] = \delta_{kl} \cdot k!.
\end{align*}
Based on this property, let $\alpha_{k}(\sigma_{\sfc})$ be defined in Lemma~\ref{lemma:connectioninGegenandHermite}.  
Then for $\sfc = 1,2$, we have
\begin{align*}
    \alpha_{k}(\sigma_{\sfc})=\mu_{\sfc,k}, ~k = 0,1, \qquad
    \mu_{\sfc,2}^2=\sum_{k\geq2}\frac{\alpha_{k}^2(\sigma_{\sfc})}{k!}, 
\end{align*}
where  
the constants $\mu_{\sfc,k}$ and $\mu_{j,2}$ are defined in Definition~\ref{def:someconstant}. 
Therefore, by Lemma~\ref{lemma:connectioninGegenandHermite}, we further have
 \begin{align}\label{eq:Gegen_mu_property}
    \sum_{k\geq2}\lambda_{d,k}^2(\sigma_{\sfc})B(d,k)\rightarrow\mu_{\sfc,2}^2.
\end{align}
Recall that $\bsigma(\xb)=\big(\sigma_1(\xb^\T \bTheta_1^\T /\sqrt{d}),\sigma_2(\xb^\T \bTheta_2^\T /\sqrt{d})\big)^\T$. Moreover, note that the zeroth order Gegenbauer polynomial $Q_0^{d}(x)=1$. Therefore by Lemma~\ref{lemma:meisongdecomp1} and the Gegenbauer decomposition of $\bsigma_j$ in Lemma~\ref{lemma:connectioninGegenandHermite}, we have
\begin{align*}
    \Vb_0(F_0)&=F_0\EE_{\xb}[\bsigma(\xb)\cdot Q_0^{d}(\xb^\T\1_d)]=\frac{F_0}{B(d,0)}\cdot  \bLambda_{d,0} \cdot Q_0^{d}(\bTheta\1_d)\cdot B(d,0)=F_0\bLambda_{d,0}\1_N.
\end{align*}
Here,  the equality holds from the fact that  $Q_0^{d}(\xb^\T\1_d)=1$ and $Q_0^{d}(\bTheta\1_d)=\1_N$. Similarly, $Q_1^{d}(x)=x/d$ holds. Again from Lemma~\ref{lemma:meisongdecomp1} and Lemma~\ref{lemma:connectioninGegenandHermite}, we have
\begin{align*}
\Vb(\bbeta_{1,d})&=\EE_{\xb}\bsigma(\xb)\xb^\T\bbeta_{1,d}=d\cdot\EE_{\xb}\bsigma(\xb)Q_1^{d}(\xb^\T\bbeta_{1,d})=\frac{d}{B(d,1)}\cdot  \bLambda_{d,1} \cdot Q_1^{d}(\bTheta\bbeta_{1,d})\cdot B(d,1)\\
&=\bLambda_{d,1}\bTheta\bbeta_{1,d} = \bigg(\frac{\Mb_{1}+\bDelta'}{\sqrt d}\bigg)\bTheta\bbeta_{1,d}.
\end{align*}
Here, $\bDelta'=\sqrt{d}\cdot\bLambda_{d,1}-\Mb_1$. From Lemma~\ref{lemma:connectioninGegenandHermite}, set $k=1$ and we have $\sqrt{d}\lambda_{d,1}(\sigma_{\sfc})\rightarrow\mu_{\sfc,1}$. Thus $\Delta'$ satisfies $\EE\lVert\bDelta'\rVert_{\op}^2=o_d(1)$.  As for   $\Ub=\EE_{\xb}\big[\bsigma(\xb)\bsigma(\xb)^\T \big]$, $\Ub$ could be divided into the following block matrix:
\begin{equation*}
    \Ub=\begin{bmatrix}
    \Ub_{1,1}&\quad \Ub_{1,2}\\
    \Ub_{2,1}&\quad \Ub_{2,2}
    \end{bmatrix},
\end{equation*}
where
\begin{equation*}
     \Ub_{i,j} =\EE_{\xb}[\sigma_i( \bTheta_i\xb /\sqrt{d})\sigma_j(\xb^\T \bTheta_j^\T /\sqrt{d})],
        \quad i,j=1,2.
\end{equation*}
Now by Lemma~\ref{lemma:meisongdecomp1}, we have
\begin{align}
    \Ub_{i,j}= \sum_{k=0}^{+\infty}\lambda_{d,k}(\sigma_i)\lambda_{d,k}(\sigma_j)B(d,k)Q_k^{(d)}(\bTheta_i\bTheta_j^\T ),
    \quad i,j=1,2.
    \label{eq:Uij_def}
\end{align}
\if UT
where
\begin{equation*}
    \begin{split}
        &\Ub_{1,1}=\EE_{\xb}[\sigma_1( \bTheta_1\xb /\sqrt{d})\sigma_1(\xb^\T \bTheta_1^\T /\sqrt{d})],\\
        &\Ub_{1,2}=\Ub_{2,1}^\T=\EE_{\xb}[\sigma_1( \bTheta_1\xb /\sqrt{d})\sigma_2(\xb^\T \bTheta_2^\T /\sqrt{d})],\\
        &\Ub_{2,2}=\EE_{\xb}[\sigma_2( \bTheta_2\xb /\sqrt{d})\sigma_2(\xb^\T \bTheta_2^\T /\sqrt{d})].
    \end{split}
\end{equation*}
Now by Lemma~\ref{lemma:meisongdecomp1}, we have
\begin{align}
    \Ub_{1,1} =&\sum_{k,l=0}^{+\infty}\lambda_{d,k}(\sigma_1)\lambda_{d,l}(\sigma_1)B(d,k)B(d,l)\EE_\xb\big[Q_k^{(d)}(\bTheta_1\xb)Q_l^{(d)}(\xb^\T \bTheta_1^\T )\big] \nonumber\\
        =&\sum_{k=0}^{+\infty}\lambda_{d,k}^2(\sigma_1)B(d,k)Q_k^{(d)}(\bTheta_1\bTheta_1^\T ). \label{eq:U11_def}
\end{align}
Similarly, we also have
\begin{align}
    & \Ub_{1,2}=\Ub_{2,1}^\T = \sum_{k=0}^{+\infty}\lambda_{d,k}(\sigma_1)\lambda_{d,k}(\sigma_2)B(d,k)Q_k^{(d)}(\bTheta_1\bTheta_2^\T ), \label{eq:U1221_def} \\
    & \Ub_{2,2} = \sum_{k=0}^{+\infty}\lambda_{d,k}^2(\sigma_2)B(d,k)Q_k^{(d)}(\bTheta_2\bTheta_2^\T ).\label{eq:U22_def}
\end{align}
\fi
Note that $Q_0^{(d)}(x)=1$, $Q_1^{(d)}(x)=x/d$, so   that the first two terms in  the decomposition \eqref{eq:Uij_def}  have a simple form. 
For $k \geq 2$, we approximate the terms using the approximation given in the second item of   Lemma~\ref{lemma:meisongdecomp1}.
\if UT
$$\begin{aligned}
\EE\Big[\sup_{k\geq2}\big\lVert Q_k^{(d)}(\bTheta_{\sfc}\bTheta_{\sfc}^\T )-\Ib_{N_{\sfc}}\big\rVert_{\op}^2\Big]=o_d(1),\quad
\EE\Big[\sup_{k\geq2}\big\lVert Q_k^{(d)}(\bTheta_1\bTheta_2^\T )\big\rVert_{\op}^2\Big]=o_d(1).
\end{aligned}$$ 
\fi 
Consider first $\Ub_{1,1}$. We have
\begin{align}
    \Ub_{1,1} 
        & = \lambda_{d,0}^2(\sigma_1)\cdot \1_{N_1}\1_{N_1}^\T +  \lambda_{d,1}^2(\sigma_1) \cdot B(d,1) \cdot \frac{\bTheta_1\bTheta_1^\T }{d} + \sum_{k = 2}^{+\infty}\lambda_{d,k}^2(\sigma_1)\cdot B(d,k) \cdot Q_k^{(d)}(\bTheta_1\bTheta_1^\T )\nonumber\\
        & = \lambda_{d,0}^2(\sigma_1)\cdot \1_{N_1}\1_{N_1}^\T +  \lambda_{d,1}^2(\sigma_1) \cdot B(d,1) \cdot \frac{\bTheta_1\bTheta_1^\T }{d} + \sum_{k = 2}^{+\infty}\lambda_{d,k}^2(\sigma_1) \cdot B(d,k) \cdot \Ib_{N_1} \nonumber \\
        &\quad  + \sum_{k = 2}^{+\infty}\lambda_{d,k}^2(\sigma_1) \cdot B(d,k) \cdot \big[ Q_k^{(d)}(\bTheta_1\bTheta_1^\T ) - \Ib_{N_1}\big],\label{eq:U11_decomp_eq1}
\end{align}
where  we have used the fact that  $\sum\limits_{k=2}^{+\infty}\lambda_{d,k}^2(\sigma_{1})B(d,k)<+\infty$ for sufficiently large $d$, which is implied by  
\eqref{eq:Gegen_mu_property}. Moreover, by Lemma~\ref{lemma:meisongdecomp1}, the convergence of this series 
  also implies that 
\begin{align}\label{eq:U11_decom_eq2}
    \EE\Bigg\lVert\sum_{k = 2}^{+\infty}\lambda_{d,k}^2(\sigma_1) \cdot B(d,k) \cdot \big[ Q_k^{(d)}(\bTheta_1\bTheta_1^\T ) - \Ib_{N_1}\big]\Bigg\rVert_{\op}^2 = o_d(1).
\end{align}
Therefore by \eqref{eq:U11_decomp_eq1} and \eqref{eq:U11_decom_eq2}, we have
\begin{align}\label{eq:U11_decomp_eq3}
    \EE&\Bigg\lVert\Ub_{1,1}-\lambda_{d,0}^2(\sigma_1)\cdot \1_{N_1}\1_{N_1}^\T - \lambda_{d,1}^2(\sigma_1) \cdot B(d,1) \cdot \frac{\bTheta_1\bTheta_1^\T }{d} - \sum_{k = 2}^{+\infty}\lambda_{d,k}^2(\sigma_1) \cdot B(d,k) \cdot \Ib_{N_1}\Bigg\rVert_{\op}^2=o_d(1).
\end{align}
Now by Lemma~\ref{lemma:connectioninGegenandHermite} and equations \eqref{eq:Gegen_mu_property}, \eqref{eq:U11_decomp_eq3}, we have
\begin{align*}
    \EE&\bigg\lVert\Ub_{1,1}-\lambda_{d,0}^2(\sigma_1)\cdot \1_{N_1}\1_{N_1}^\T - \mu_{1,1}^2 \cdot \frac{\bTheta_1\bTheta_1^\T }{d} - \mu_{1,2}^2 \cdot \Ib_{N_1}\bigg\rVert_{\op}^2=o_d(1).
\end{align*}
This establishes  the approximation for $\Ub_{1,1}$.

For the other sub-matrices $\Ub_{1,2}$, $\Ub_{2,1}$ and $\Ub_{2,2}$, the derivations 
are  exactly the same, and we obtain the following results:
\begin{align*}
        \EE&\bigg\lVert\Ub_{1,2}-\lambda_{d,0}(\sigma_1)\lambda_{d,0}(\sigma_2)\1_{N_1}\1_{N_2}^\T -\mu_{1,1}\mu_{2,1}\cdot\frac{\bTheta_1\bTheta_2^\T }{d}\bigg\rVert_{\op}^2=o_d(1),\\
        \EE&\bigg\lVert\Ub_{2,1}-\lambda_{d,0}(\sigma_1)\lambda_{d,0}(\sigma_2)\1_{N_1}\1_{N_2}^\T -\mu_{1,1}\mu_{2,1}\cdot\frac{\bTheta_2\bTheta_1^\T }{d}\bigg\rVert_{\op}^2=o_d(1),\\
        \EE&\bigg\lVert\Ub_{2,2}-\lambda_{d,0}^2(\sigma_2)\1_{N_2}\1_{N_2}^\T -\mu_{2,1}^2\frac{\bTheta_2\bTheta_2^\T }{d}-\mu_{2,2}^2\cdot\Ib_{N_2}\bigg\rVert_{\op}^2=o_d(1).
\end{align*}
Note that the collection of the approximations for the four blocks
$\Ub_{1,1}, \Ub_{1,2}, \Ub_{2,1}$ and  $\Ub_{1,1}$ gives the matrix 
\[\bLambda_{d,0}\1_N\1_N^\T \bLambda_{d,0}+\Mb_1\frac{\bTheta\bTheta^\T }{d}\Mb_1+\Mb_2\Mb_2,
\] 
so finally we have  
\begin{align*}&\EE\bigg\lVert\Ub-\bLambda_{d,0}\1_N\1_N^\T \bLambda_{d,0}-\Mb_1\frac{\bTheta\bTheta^\T }{d}\Mb_1-\Mb_2\Mb_2\bigg\rVert_{\op}^2 =  o_d(1). 
\end{align*}
\if UT
{\small
\begin{align*}&\EE\bigg\lVert\Ub-\bLambda_{d,0}\1_N\1_N^\T \bLambda_{d,0}-\Mb_1\frac{\bTheta\bTheta^\T }{d}\Mb_1-\Mb_2\Mb_2\bigg\rVert_{\op}^2\\
=&\EE\Bigg\lVert\begin{bmatrix}
\Ub_{1,1}-\lambda_{d,0}^2(\sigma_1)\cdot \1_{N_1}\1_{N_1}^\T - \mu_{1,1}^2 \cdot \frac{\bTheta_1\bTheta_1^\T }{d} - \mu_{1,2}^2 \cdot \Ib_{N_1}&\Ub_{1,2}-\lambda_{d,0}(\sigma_1)\lambda_{d,0}(\sigma_2)\1_{N_1}\1_{N_2}^\T -\mu_{1,1}\mu_{2,1}\cdot\frac{\bTheta_1\bTheta_2^\T }{d}\\
\Ub_{2,1}-\lambda_{d,0}(\sigma_1)\lambda_{d,0}(\sigma_2)\1_{N_1}\1_{N_2}^\T -\mu_{1,1}\mu_{2,1}\cdot\frac{\bTheta_2\bTheta_1^\T }{d}&\Ub_{2,2}-\lambda_{d,0}^2(\sigma_2)\1_{N_2}\1_{N_2}^\T -\mu_{2,1}^2\frac{\bTheta_2\bTheta_2^\T }{d}-\mu_{2,2}^2\cdot\Ib_{N_2}
\end{bmatrix}\Bigg\rVert_{\op}^2\\=&o_d(1).\end{align*}}
\fi 
The proof of Lemma~\ref{lemma:decomp1} is complete.

\subsection{Proof of Lemma~\ref{lemma:decomp2}}
\label{sec:appendixdecomp2}
We first prove that $\big(\EE |\cD|^k\big)^{1/k}=o_d(1)$ if $\sum_{\sfc}\mu_{\sfc,0}^2=0$. Note that the rank-1 matrix $\Ab$ satisfies $|\tr(\Ab)|=\|\Ab\|_{\op}$. Moreover, $\sum_{\sfc}\mu_{\sfc,0}^2=0$ implies $\|\Lambda_{d,0}\|_{\op}=o_d(1)$.  We have $\big(\EE |\cD|^k\big)^{1/k}=O_d(\|\bLambda_{d,0}\frac{\1_N\1_N^\T}{d}\bLambda_{d,0}\|_{\op})\cdot(\EE\|\bGamma_2\|_{\op}^k)^{1/k}=o_d(1)\cdot O_d(1)=o_d(1)$.

In the following proof of Lemma~\ref{lemma:decomp2}, we have the condition $\sum_{\sfc}\mu_{\sfc,0}^2>0$. We separate the proof into two parts, {estimating  $\cB$ and $\cC$},  and { $\cD$}, respectively. 
\subsubsection{Estimation for $\cB$ and $\cC$}
\label{subsubsection:BandC}
Let
\begin{align*}
L_1=\frac{1}{\sqrt d}\tr(\bLambda_{d,0}\1_N\1_n^\T  \Zb \bUpsilon),\qquad L_2({\bGamma})=\frac{1}{d}\tr( \big[{\bGamma}\big]_\Zb \1_n\1_n^\T )=\frac{1}{d}\tr( \Zb \bUpsilon{\bGamma}\bUpsilon  \Zb ^\T \1_n\1_n^\T ),
\end{align*}
where  $\bGamma\in\RR^{N\times N}$ is a  symmetric matrix.
Then  we have 
 \[
  \cB=L_2(\bGamma),
  \quad 
  \cC=1-2L_1+L_2(\bLambda_{d,0}\1_N\1_N^\T\bLambda_{d,0}).
\] 
Define further the following terms:
$$
\begin{array}{l@{\quad}l@{\quad}l@{\quad}l}
&K_{11}=\Tb_1^\T \Eb_0^{-1}\Tb_1,&K_{12}=\Tb_1^\T \Eb_0^{-1}\Tb_2,&K_{22}=\Tb_2^\T \Eb_0^{-1}\Tb_2,\\
&G_{11}=\Tb_1^\T \Eb_0^{-1}{\bGamma}\Eb_0^{-1}\Tb_1,&G_{12}=\Tb_1^\T \Eb_0^{-1}{\bGamma}\Eb_0^{-1}\Tb_2,&G_{22}=\Tb_2^\T \Eb_0^{-1}{\bGamma}\Eb_0^{-1}\Tb_2,
\end{array}
$$
where
$$
\begin{array}{l@{\qquad}l}
\Jb= \Zb-\1_n\1_N^\T \bLambda_{d,0}/\sqrt{d}, &\Eb_0=\Jb^\T \Jb+\lambda\Ib_N,\\
\Tb_1=\psi_3^{1/2}\bLambda_{d,0}\1_N,
&\Tb_2=\frac{1}{\sqrt{n}}\Jb^\T \1_n.
\end{array}
$$
We denote $\psi_3=n/d$ for notation simplification. The proof is organized in two steps: 
\begin{enumerate}[leftmargin = *]
    \item Express $\cB$ and $\cC$ in function of $K_{ij}$ and $G_{ij}$, $i,j\in\{1,2\}$.
    \item Estimate the order of $K_{ij}$ and $G_{ij}$, and show that $\EE|\cB|$ and $\EE|\cC|$ are both $o_d(1)$.
\end{enumerate}
Denote $\Fb_1=[\Tb_1,\Tb_1,\Tb_2] \in \RR^{N\times 3}$, $\Fb_2=[\Tb_1,\Tb_2,\Tb_1] \in \RR^{N\times 3}$, it is easy to see 
\begin{align}
    \bUpsilon=&\big(\Zb^\T\Zb+\lambda\Ib_N\big)^{-1}=\big((\Jb+\1_n\1_N^\T \bLambda_{d,0})^\T(\Jb+\1_n\1_N^\T \bLambda_{d,0})+\lambda\Ib_N\big)^{-1}\nonumber\\
    =&\big(\psi_3\bLambda_{d,0}\1_N\1_N^\T
+\psi_3^{1/2}\bLambda_{d,0}\1_N\Tb_2^\T +\psi_3^{1/2}\Tb_2\1_N^\T \bLambda_{d,0}+\Jb^\T \Jb+\lambda\Ib_N\big)^{-1}\nonumber\\
=&\big(\Eb_0+\Fb_1\Fb_2^\T \big)^{-1}.\nonumber
\end{align}
For $L_1$, replacing $\Zb$ by $\Jb+\1_n\1_N^\T \bLambda_{d,0}/\sqrt{d}$, we have
\begin{align}
L_1=&\tr\big[(\psi_3\bLambda_{d,0}\1_N\1_N^\T \bLambda_{d,0}+\psi_3^{1/2}\bLambda_{d,0}\1_N\Tb_2^\T )\cdot\big(\Eb_0+\Fb_1\Fb_2^\T \big)^{-1}\big]\nonumber\\
=&\tr\big[(\Tb_1\Tb_1^\T +\Tb_1\Tb_2^\T )\cdot\big(\Eb_0+\Fb_1\Fb_2^\T \big)^{-1}\big].\label{eq:L1_cal1}
\end{align}
By the Sherman-Morrison-Woodbury formula, 
\begin{align}
\bUpsilon=\big(\Eb_0+\Fb_1\Fb_2^\T \big)^{-1}=\Eb_0^{-1}-\Eb_0^{-1}\Fb_1(\Ib_3+\Fb_2^\T \Eb_0^{-1}\Fb_1)^{-1}\Fb_2^{T}\Eb_0^{-1}.\label{eq:woodberryfomula}
\end{align}
Plugging \eqref{eq:woodberryfomula} into \eqref{eq:L1_cal1}, we have   
$$\begin{aligned}
L_1=~&(\Tb_1^\T \Eb_0^{-1}\Tb_1-\Tb_1^\T \Eb_0^{-1}\Fb_1(\Ib_3+\Fb_2^\T \Eb_0^{-1}\Fb_1)^{-1}\Fb_2^\T \Eb_0^{-1}\Tb_1)\\
&+(\Tb_2^\T \Eb_0^{-1}\Tb_1-\Tb_2^\T \Eb_0^{-1}\Fb_1(\Ib_3+\Fb_2^\T \Eb_0^{-1}\Fb_1)^{-1}\Fb_2^\T \Eb_0^{-1}\Tb_1)\\
=~&(K_{11}-[K_{11},K_{11},K_{12}](\Ib_3+\bK)^{-1}[K_{11},K_{12},K_{11}]^\T )\\
&+(K_{12}-[K_{12},K_{12},K_{22}](\Ib_3+\bK)^{-1}[K_{11},K_{12},K_{11}]^\T )\\
=~&[K_{11},K_{11},K_{12}](\Ib_3+\bK)^{-1}[1,0,0]^\T +[K_{12},K_{12},K_{22}](\Ib_3+\bK)^{-1}[1,0,0]^\T ,
\end{aligned}$$
where
\[ \bK=\Fb_2^\T \Eb_0^{-1}\Fb_1=\begin{bmatrix}
K_{11}&K_{11}&K_{12}\\
K_{12}&K_{12}&K_{22}\\
K_{11}&K_{11}&K_{12}
\end{bmatrix}.
\]
Thus by simple calculation, 
\begin{align}\label{eq:L1_calculation}
    L_1=1-\frac{K_{12}+1}{K_{11}(1-K_{22})+(K_{12}+1)^2}.
\end{align}
As for $L_2({\bGamma})$, we have 
\begin{align*}
\Zb ^\T \1_n\1_n^\T\Zb/d=&(\Jb+\1_n\1_N^\T \bLambda_{d,0}/\sqrt{d})^\T\1_n\1_n^\T(\Jb+\1_n\1_N^\T \bLambda_{d,0}/\sqrt{d})/d\\ =&\psi_3(\psi_3^{1/2}\bLambda_{d,0}\1_N+\frac{1}{\sqrt{n}}\Jb^\T \1_n)(\psi_3^{1/2}\bLambda_{d,0}\1_N+\frac{1}{\sqrt{n}}\Jb^\T \1_n)^\T\\=&\psi_3(\Tb_1+\Tb_2)(\Tb_1+\Tb_2)^\T .
\end{align*}
Then after similar calculation by  \eqref{eq:woodberryfomula}.
\begin{align}
\cB=L_2({\bGamma})=&\frac{1}{d}\tr\big(\Zb ^\T \1_n\1_n^\T \Zb \bUpsilon{\bGamma}\bUpsilon\big)=\tr\big(\psi_3(\Tb_1+\Tb_2)(\Tb_1+\Tb_2)^\T \bUpsilon{\bGamma}\bUpsilon\big)\nonumber\\
=&\psi_3(\Tb_1+\Tb_2)^\T \big(\Eb_0+\Fb_1\Fb_2^\T \big)^{-1}{\bGamma}\big(\Eb_0+\Fb_1\Fb_2^\T \big)^{-1}(\Tb_1+\Tb_2)\nonumber\\
=&\psi_3(\Tb_1+\Tb_2)^\T\big(\Eb_0^{-1}-\Eb_0^{-1}\Fb_1(\Ib_3+\Fb_2^\T \Eb_0^{-1}\Fb_1)^{-1}\Fb_2^{T}\Eb_0^{-1}\big)\nonumber\\&\quad\cdot{\bGamma}\big(\Eb_0^{-1}-\Eb_0^{-1}\Fb_1(\Ib_3+\Fb_2^\T \Eb_0^{-1}\Fb_1)^{-1}\Fb_2^{T}\Eb_0^{-1}\big)(\Tb_1+\Tb_2)\nonumber  \\
=&\psi_3\frac{G_{11}(1-K_{22})^2+G_{22}(K_{12}+1)^2+2G_{12}(K_{12}+1)(1-K_{22})}{(K_{11}(1-K_{22})+(K_{12}+1)^2)^2}.\label{eq:Bcalcu}
\end{align}
When $\bGamma=\bLambda_{d,0}\1_N\1_N^\T \bLambda_{d,0} $, the $G_{11},G_{12}$ and $G_{22}$ above can be given as
\begin{equation*}
    \begin{split}
        G_{11}=K_{11}^2/\psi_3,\qquad
        G_{12}=K_{11}K_{12}/\psi_3,\qquad
        G_{22}=K_{12}^2/\psi_3.
    \end{split}
\end{equation*}
Then by \eqref{eq:L1_calculation}, \eqref{eq:Bcalcu},
we have
\begin{equation}\label{eq:Ccalcu}
    \begin{split}
        \cC =1-2L_1+L_2(\bLambda_{d,0}\1_N\1_N^\T \bLambda_{d,0})=\frac{(K_{12}+1)^2}{(K_{11}(1-K_{22})+(K_{12}+1)^2)^2}.
    \end{split}
\end{equation}

We next estimate the order for $K_{11}$, $K_{12}$, $K_{22}$, $G_{11}$, $G_{12}$ and $G_{22}$ respectively.
By the inequality 
$\big\lVert
\big[
\Ab~~\Bb
\big]\big\rVert_{\op}\leq \lVert \Ab\rVert_{\op}+\lVert \Bb\rVert_{\op} 
$
for any matrices $\Ab$ and $\Bb$, we have
\begin{align}
    \lVert \Jb\rVert_{\op}&\leq \lVert \Zb_1-\lambda_{d,0}(\sigma_1)\1_n\1_{N_1}^\T /\sqrt{d}\rVert_{\op}+\lVert \Zb_2-\lambda_{d,0}(\sigma_2)\1_n\1_{N_2}^\T /\sqrt{d}\rVert_{\op}\nonumber\\
&=O_{\PP}(\exp{(C\sqrt{\Log d})}),\label{eq:decompBJ}
\end{align}
where the last equality in \eqref{eq:decompBJ} follows by Lemma C.5 in \citep{mei2022generalization}. Moreover, for any fixed $\lambda>0$, it also determinstically holds that
\begin{align*}
\rVert(\Jb^\T \Jb+\lambda\Ib_N)^{-1}\Jb^\T \rVert_{\op}\leq2/\sqrt{\lambda},\quad \rVert(\Jb^\T \Jb+\lambda\Ib_N)^{-1} \rVert_{\op} \leq 1/\lambda.
\end{align*}
Now recall that
$$\begin{aligned}
K_{11}=&\psi_3\1_N^\T \bLambda_{d,0}(\Jb^\T \Jb+\lambda\Ib_N)^{-1}\bLambda_{d,0}\1_N,\\
K_{12}=&\1_N^\T \bLambda_{d,0}(\Jb^\T \Jb+\lambda\Ib_N)^{-1}\Jb^\T \1_n/\sqrt{d},\\
K_{22}=&\1_n^\T \Jb(\Jb^\T \Jb+\lambda\Ib_N)^{-1}\Jb^\T \1_n/n,\\
G_{11}=&\psi_3\1_N^\T \bLambda_{d,0}(\Jb^\T \Jb+\lambda\Ib_N)^{-1}{\bGamma}(\Jb^\T \Jb+\lambda\Ib_N)^{-1}\bLambda_{d,0}\1_N,\\
G_{12}=&\1_N^\T \bLambda_{d,0}(\Jb^\T \Jb+\lambda\Ib_N)^{-1}{\bGamma}(\Jb^\T \Jb+\lambda\Ib_N)^{-1}\Jb^\T \1_n/\sqrt{d},\\
G_{22}=&\1_n^\T \Jb(\Jb^\T \Jb+\lambda\Ib_N)^{-1}{\bGamma}(\Jb^\T \Jb+\lambda\Ib_N)^{-1}\Jb^\T \1_n/n.
\end{aligned}$$
Therefore we deterministically have
\begin{align}
|K_{12}|&\leq \rVert(\Jb^\T \Jb+\lambda\Ib_N)^{-1}\Jb^\T \rVert_{\op}\lVert\1_n\1_N^\T \bLambda_{d,0}/\sqrt{d}\rVert_{\op}=O_d(\sqrt{d/\lambda}).\label{eq:boundK12}
\end{align}
For $K_{22}$, by its definition, it is clear that $K_{22} > 0$. Moreover, we have
\begin{align*}
K_{22}&\leq \lambda_{\max}(\Jb(\Jb^\T \Jb+\lambda\Ib_N)^{-1}\Jb)\tr(\1_n\1_n^\T /n)\nonumber\\
&=\lambda_{\max}(\Ib_N-\lambda(\Jb^\T \Jb+\lambda\Ib_N)^{-1})
=1-\frac{\lambda}{\lVert\Jb^\T \Jb\rVert_{\op}+\lambda}.
\end{align*}
Therefore we have
\begin{align}\label{eq:boundK22}
    0< K_{22} \leq 1-\frac{\lambda}{\lVert\Jb^\T \Jb\rVert_{\op}+\lambda}.
\end{align}
For $K_{11}$, the condition $\mu_{1,0}^2+\mu_{2,0}^2>0$   
ensures that there exists $\sfc\in\{1,2\}$ such that $\mu^2_{\sfc,0} > 0$. 
By Lemma~\ref{lemma:connectioninGegenandHermite} (note that $B(d,0) = 1$), we have $\lambda_{d,0}^2(\sigma_{\sfc})\rightarrow\mu_{\sfc,0}^2$ as $d\rightarrow+\infty$.
Therefore for large enough $d$, we have $\lambda_{d,0}(\sigma_{\sfc})> \mu_{\sfc,0} / 2 >0$, and 
\begin{align}
K_{11}&\geq\psi_3\1_N^\T\bLambda_{d,0}^2\1_N\lambda_{\min}((\Jb^\T \Jb+\lambda\Ib_N)^{-1})\nonumber\\
&\geq \psi_3 \cdot (\mu_{\sfc,0}^2 / 4) \cdot N_\sfc \cdot \lambda_{\min}((\Jb^\T \Jb+\lambda\Ib_N)^{-1})\nonumber\\
&=\frac{\Omega_d(d)}{\lVert\Jb^\T \Jb\rVert_{\op}+\lambda}.\label{eq:boundK11}
\end{align}
Plugging \eqref{eq:boundK12}, \eqref{eq:boundK22}, \eqref{eq:boundK11} into \eqref{eq:Bcalcu} then gives
\begin{align*}
    |\cB | &=  \frac{\big| G_{22}(1+K_{12})^2+G_{11}(1-K_{22})^2+2G_{12}(1+K_{12})(1-K_{22}) \big| }{\big[(1+K_{12})^2+ K_{11}\cdot (1 -K_{22})\big]^2} \\
    &\leq \frac{\big| G_{22}(1+K_{12})^2+G_{11}(1-K_{22})^2+2G_{12}(1+K_{12})(1-K_{22}) \big| }{\big[ K_{11}\cdot (1 -K_{22})\big]^2} \\
    &\leq O_d(1)\cdot \frac{|G_{22}|\cdot d + |G_{12}| \cdot\sqrt{d}+ |G_{11}| }{d^2/(\lambda+\lVert\Jb\Jb^\T \rVert_{\op})^4}, 
\end{align*}
where we utilize the upper and lower bounds in \eqref{eq:boundK12}, \eqref{eq:boundK22}, \eqref{eq:boundK11} to obtain the last inequality.
For $G_{11},G_{12}$ and $G_{22}$, we have
\small
$$\begin{aligned}
\EE \big[|G_{11}|^{k}\big]^{1/k}&\leq\psi_3\big\lVert(\Jb^\T \Jb+\lambda\Ib_N)^{-1}\big\rVert_{\op}\big[\EE\lVert{{\bGamma} }\rVert_{\op}^{k}\big]^{1/k}\big\lVert(\Jb^\T \Jb+\lambda\Ib_N)^{-1}\big\rVert_{\op}\big\lVert\bLambda_{d,0}\1_N\1_N^\T \bLambda_{d,0}\big\rVert_{\op}=O_d(d),\\
\EE \big[|G_{12}|^{k}\big]^{1/k}&\leq\big\lVert(\Jb^\T \Jb+\lambda\Ib_N)^{-1}\big\rVert_{\op}\big[\EE\lVert{\bGamma }\rVert_{\op}^{k}\big]^{1/k}\big\lVert(\Jb^\T \Jb+\lambda\Ib_N)^{-1}\Jb^\T \big\rVert_{\op}\big\lVert\1_n\1_N^\T \bLambda_{d,0}/\sqrt{d}\big\rVert_{\op}=O_d(\sqrt{d}),\\
\EE \big[|G_{22}|^{k}\big]^{1/k}&\leq\big\lVert(\Jb^\T \Jb+\lambda\Ib_N)^{-1}\big\rVert_{\op}\big[\EE \|\bGamma\|_{\op}^{k}\big]^{1/k}\big\lVert(\Jb^\T \Jb+\lambda\Ib_N)^{-1}\Jb^\T \Jb\big\rVert_{\op}\tr\big(\1_n\1_n^\T /n\big)=O_d(1).
\end{aligned}$$
\normalsize
Thus by the  bounds above and the triangle inequality of the $L_k$-norm $\EE[ |\cdot|^k ]^{1/k}$, we have
$$\begin{aligned}
\big(\EE |\cB |^k\big)^{1/k}&\leq O_d(1)\cdot \frac{\EE \big[|G_{22}|^{k}\big]^{1/k}\cdot d+\EE \big[|G_{12}|^{k}\big]^{1/k}\cdot\sqrt{d}+\EE \big[|G_{11}|^{k}\big]^{1/k}}{d^2/(\lambda+\lVert\Jb\Jb^\T \rVert_{\op})^4}\\
& = O_d(1)\cdot \frac{d}{d^2/(\lambda+\lVert\Jb\Jb^\T \rVert_{\op})^4} = O_d\bigg(\frac{\big(\lambda+\lVert\Jb\Jb^\T \rVert_{\op}\big)^4}{d}\bigg)=O_d\bigg(\frac{\exp(C\sqrt{\log d})}{d}\bigg)\\
&=o_d(1),
\end{aligned}$$
and
\begin{align*}
    \EE |\cC |=&O_d\bigg(\frac{\big(\lambda+\lVert\Jb\Jb^\T \rVert_{\op}\big)^2}{ d}\bigg)=O_d\bigg(\frac{\exp(C\sqrt{\log d})}{d}\bigg)=o_d(1).
\end{align*}
This completes the proof. 

\subsubsection{Estimation for $\cD$}
\label{subsubsection:D}
The proof is similar to the calculations for $\cB$ and $\cC$ in the previous section.
Also we use a set of similar notations as previously which may however have slightly different values.  Let 
$$
\begin{array}{l@{\quad}l@{\quad}l@{\quad}l}
\Jb= \Zb-\1_n\1_N^\T \bLambda_{d,0}/\sqrt{d}, &\Eb_0=\Jb\Jb^\T+\lambda\Ib_n,&\\
\Tb_1=\Jb\bLambda_{d,0}\1_N/\sqrt{d},
&\Tb_2=\1_n,&\\
K_{11}=\Tb_1^\T \Eb_0^{-1}\Tb_1,&K_{12}=\Tb_1^\T \Eb_0^{-1}\Tb_2,&K_{22}=\Tb_2^\T \Eb_0^{-1}\Tb_2,\\
G_{11}=\Tb_1^\T \Eb_0^{-1}\bGamma \Eb_0^{-1}\Tb_1,&G_{12}=\Tb_1^\T \Eb_0^{-1}\bGamma \Eb_0^{-1}\Tb_2,&G_{22}=\Tb_2^\T \Eb_0^{-1}\bGamma \Eb_0^{-1}\Tb_2,
\end{array}
$$
where   $\bGamma\in\RR^{n\times n}$ is a symmetric matrix.  
We express $\cD$ with the terms defined above. Recall that   $\bUpsilon=(\Zb^{\T}\Zb+\lambda\Ib_N)^{-1}$and  further define $\bXi=(\Zb\Zb^\T+\lambda\Ib_n)^{-1}$. Clearly,    $\Zb\bUpsilon=\bXi\Zb$. Therefore we have
\begin{equation}
\label{eq:decomDalpha}
    \begin{split}
        \cD=\frac{1}{d}\tr(\Zb\bUpsilon\bLambda_{d,0}\1_N\1_N^\T \bLambda_{d,0}\bUpsilon  \Zb ^\T \bGamma  )=\frac{1}{d}\tr( \bXi\Zb\bLambda_{d,0}\1_N\1_N^\T \bLambda_{d,0}\Zb^\T  \bXi\bGamma  ).
    \end{split}
\end{equation}
We proceed to calculate $ \bXi$ and $\Zb\bLambda_{d,0}\1_N\1_N^\T \bLambda_{d,0}\Zb^\T $, respectively. 
Define $c={\1_N^\T \bLambda_{d,0}^2\1_N}/{d} = \Theta(1)$, $\Fb_1=[\Tb_1,\Tb_2,\Tb_2]\in \RR^{n \times 3}$, $\Fb_2=[\Tb_2,\Tb_1,c\Tb_2]\in \RR^{n \times 3}$. Then   we have
\begin{align}
 \bXi=&\Big(\big(\Jb+\1_n\1_N^\T \bLambda_{d,0}/\sqrt{d}\big)\big(\Jb+\1_n\1_N^\T \bLambda_{d,0}/\sqrt{d}\big)^\T +\lambda \Ib_n \Big)^{-1}\nonumber\\
=&\big(\Eb_0+\Fb_1\Fb_2^\T \big)^{-1}=~\Eb_0^{-1}-\Eb_0^{-1}\Fb_1(\Ib_3+\Fb_2^\T \Eb_0^{-1}\Fb_1)^{-1}\Fb_2^{T}\Eb_0^{-1},\label{eq:decomBXi}
\end{align}
where the last equality follows from  the Sherman-Morrison-Woodbury formula.
Moreover, we have
\begin{align}
    \Zb\bLambda_{d,0}\1_N\1_N^\T \bLambda_{d,0}\Zb^\T =&\Big(\Jb+\1_n\1_N^\T \bLambda_{d,0}/\sqrt{d}\Big)\bLambda_{d,0}\1_N\1_N^\T \bLambda_{d,0}\Big(\Jb+\1_n\1_N^\T \bLambda_{d,0}/\sqrt{d}\Big)^\T \nonumber\\
=&d\big(\Tb_1\Tb_1^\T +c\big(\Tb_2\Tb_1^\T +\Tb_1\Tb_2^\T \big)+c^2\Tb_2\Tb_2^\T  \big)\nonumber\\
=&d\cdot \big(\Tb_1 +c\Tb_2\big) \big(\Tb_1 +c\Tb_2\big)^\T.\label{eq:decompZLambZ}
\end{align}
Plugging \eqref{eq:decomBXi} and \eqref{eq:decompZLambZ} into \eqref{eq:decomDalpha}, we obtain 
\begin{align*}
\cD=&\tr\Big(\big(\Tb_1+c\Tb_2\big)^\T \big(\Eb_0^{-1}-\Eb_0^{-1}\Fb_1(\Ib_3+\Fb_2^\T \Eb_0^{-1}\Fb_1)^{-1}\Fb_2^{T}\Eb_0^{-1}\big)\\ 
&\quad\cdot{\bGamma } \big(\Eb_0^{-1}-\Eb_0^{-1}\Fb_1(\Ib_3+\Fb_2^\T \Eb_0^{-1}\Fb_1)^{-1}\Fb_2^{T}\Eb_0^{-1}\big)\big(\Tb_1+c\Tb_2\big)\Big).
\end{align*}
With similar calculation as in the proof of Lemma~\ref{lemma:decomp2}, we obtain that
\begin{align}
    \cD=&\frac{G_{11}(1+K_{12})^2+G_{22}(c-K_{11})^2+2G_{12}(1+K_{12})(c-K_{11})}{\big(1+2K_{12}+K_{12}^2+cK_{22}-K_{11}K_{22}\big)^2}.\label{eq:decomB2}
\end{align}
We then estimate the order for $K_{11}$, $K_{12}$, $K_{22}$, $G_{11}$, $G_{12}$ and $G_{22}$, respectively. 
For $K_{11}$, apparently we have $K_{11} > 0$. Moreover, 
\begin{align*}
    c-K_{11}=&\frac{1}{d}\1_N^\T \bLambda_{d,0}\big(\Ib_N-\Jb^\T (\Jb\Jb^\T +\lambda\Ib_n)^{-1}\Jb\big)\bLambda_{d,0}\1_N\\
\geq& c\big(1-\lambda_{\max}(\Jb^\T (\Jb\Jb^\T +\lambda\Ib_n)^{-1}\Jb)\big)=\frac{c\lambda}{\lambda+\lVert\Jb
\Jb^\T \rVert_{\op}}>0.
\end{align*}
Therefore we have
\begin{align}\label{eq:decomBorderK1}
c\geq c-K_{11} \geq \frac{c\lambda}{\lambda+\lVert\Jb
\Jb^\T \rVert_{\op}}>0.
\end{align}
Similarly, for $K_{12}$ and $K_{22}$ we have
\begin{align}
&|K_{12}|\leq \rVert(\Jb\Jb^\T +\lambda\Ib_n)^{-1}\Jb^\T \rVert_{\op}\lVert\1_n\1_N^\T \bLambda_{d,0}/\sqrt{d}\rVert_{\op}=O_d(\sqrt{d/\lambda}),\label{eq:decomBorderK2}\\
&K_{22}\geq n\lambda_{\min}\big(\big(\Jb\Jb^\T +\lambda\Ib_n\big)^{-1}\big)=\Omega(d)/(\lVert\Jb\Jb^\T  \rVert_{\op}+\lambda).\label{eq:decomBorderK3}
\end{align}
Plugging \eqref{eq:decomBorderK1}, \eqref{eq:decomBorderK2}, \eqref{eq:decomBorderK3} into \eqref{eq:decomB2} then gives
\begin{align*}
    |\cD | &=  \frac{\big| G_{11}(1+K_{12})^2+G_{22}(c-K_{11})^2+2G_{12}(1+K_{12})(c-K_{11}) \big| }{\big[(1+K_{12})^2+ K_{22}\cdot (c -K_{11})\big]^2} \\
    &\leq \frac{\big| G_{11}(1+K_{12})^2+G_{22}(c-K_{11})^2+2G_{12}(1+K_{12})(c-K_{11}) \big| }{\big[ K_{22}\cdot (c -K_{11})\big]^2} \\
    &\leq O_d(1)\cdot \frac{|G_{11}|\cdot d + |G_{12}| \cdot\sqrt{d}+ |G_{22}| \cdot c^2}{d^2/(\lambda+\lVert\Jb\Jb^\T \rVert_{\op})^4}, 
\end{align*}
where we utilize the upper and lower bounds in \eqref{eq:decomBorderK1}, \eqref{eq:decomBorderK2}, \eqref{eq:decomBorderK3} to obtain the last inequality. Now recall that
\begin{equation*}
    \begin{split}
G_{11}=&\1_N^\T \bLambda_{d,0}\Jb^\T (\Jb\Jb^\T +\lambda\Ib_n)^{-1}{\bGamma}(\Jb\Jb^\T +\lambda\Ib_n)^{-1}\Jb\bLambda_{d,0}\1_N/d,\\
G_{12}=&\1_N^\T \bLambda_{d,0}\Jb^\T (\Jb\Jb^\T +\lambda\Ib_n)^{-1}{\bGamma }(\Jb\Jb^\T +\lambda\Ib_n)^{-1}\1_n/\sqrt{d},\\
G_{22}=&\1_n(\Jb\Jb^n+\lambda\Ib_n)^{-1}{\bGamma}(\Jb\Jb^n+\lambda\Ib_n)^{-1}\1_n.
    \end{split}
\end{equation*}
Therefore we have 
\begin{align}
    \EE \big[|G_{11}|^{k}\big]^{1/k}\leq& c\cdot O_d(\lambda^{-1})\cdot \big[\EE\lVert{\bGamma }\rVert_{\op}^{k}\big]^{1/k}=O_d(1), \label{eq:decomBorderG1} \\
\EE \big[|G_{12}|^{k}\big]^{1/k}\leq& O_d(\lambda^{-3/2})\cdot\big[\EE\lVert{\bGamma }\rVert_{\op}^{k}\big]^{1/k}\cdot \lVert\1_n\1_N^\T \bLambda_{d,0}/\sqrt{d}\rVert_{\op}=O_d(\sqrt{d}), \label{eq:decomBorderG2} \\
\EE \big[|G_{22}|^{k}\big]^{1/k}\leq& O_d(\lambda^{-2})\cdot\big[\EE\lVert{\bGamma }\rVert_{\op}^{k}\big]^{1/k}\cdot\lVert\1_n\1_n^\T \rVert_{\op}=O_d(d). \label{eq:decomBorderG3} 
\end{align}
By the triangle inequality of the $L_k$-norm $\EE[ |\cdot|^k ]^{1/k}$, we have
\begin{equation*}
    \begin{split}
\big(\EE |\cD |^k\big)^{1/k}\leq& O_d(1)\cdot \frac{\EE \big[|G_{11}|^{k}\big]^{1/k}\cdot d+\EE \big[|G_{12}|^{k}\big]^{1/k}\cdot\sqrt{d}+\EE \big[|G_{22}|^{k}\big]^{1/k}\cdot c^2}{d^2/(\lambda+\lVert\Jb\Jb^\T \rVert_{\op})^4}\\
=&O_d(1)\cdot \frac{d}{d^2/(\lambda+\lVert\Jb\Jb^\T \rVert_{\op})^4}\\
=& O_d\bigg(\frac{\big(\lambda+\lVert\Jb\Jb^\T \rVert_{\op}\big)^4}{d}\bigg)
=O_d\bigg(\frac{\exp(C\sqrt{\log d})}{d}\bigg)=o_d(1),
    \end{split}
\end{equation*}
where the first equality follows by \eqref{eq:decomBorderG1}, \eqref{eq:decomBorderG2} and \eqref{eq:decomBorderG3}. 
This completes the proof.

\subsection{Proof of Lemma~\ref{lemma:replace}}
\label{sec:appendixreplace}
The first result follows by the rotation invariance of the learning problem. For any $\bbeta_{d} = [F_0,\bbeta_{1,d}^\T]^\T$ and $\tilde{\bbeta}_d = [F_0,\tilde{\bbeta}_{1,d}^\T]^\T$ with $\bbeta_{1,d}, \tilde{\bbeta}_{1,d} \in F_{1,d}\cdot\SSS^{d-1}$, there exists an orthogonal  matrix $\Pb$ such that $\Pb  \bbeta_{1,d} = \tilde{\bbeta}_{1,d}$. Then by definition, we have
\begin{align*}
     R_d(\Xb\Pb,\bTheta\Pb,\lambda,\bbeta_d,\bvarepsilon) = R_d(\Xb,\bTheta,\lambda,\tilde{\bbeta}_d,\bvarepsilon).
\end{align*}
Moreover, it is easy to check that
\begin{align*}
    \overline{R}_d(\Xb\Pb,\bTheta\Pb,\lambda,F_{1,d},\tau) = \overline{R}_d(\Xb,\bTheta,\lambda,F_{1,d},\tau).
\end{align*}
Since $(\Xb\Pb, \bTheta\Pb) \overset{\mathrm{d}}{=} (\Xb, \bTheta)$, we  see that conditional to  $\bbeta_d,\tilde{\bbeta}_d$, we have
\begin{align*}
    R_d(\Xb,\bTheta,\lambda,\bbeta_d,\bvarepsilon) - \overline{R}_d(\Xb,\bTheta,\lambda,F_{1,d},\tau) \overset{\mathrm{d}}{=}  R_d(\Xb,\bTheta,\lambda,\tilde{\bbeta}_d,\bvarepsilon) - \overline{R}_d(\Xb,\bTheta,\lambda,F_{1,d},\tau).
\end{align*}
This implies the first result in Lemma~\ref{lemma:replace}. 

If we assume that $\tilde{\bbeta}_{1,d}\sim\rmN(\0,[F_{1,d}^2/d]\Ib_d)$, then $F_{1,d}\cdot\tilde{\bbeta}_{1,d}/\|\tilde{\bbeta}_{1,d}\|_2\sim F_{1,d}\cdot\Unif(\SSS^{d-1})$. The proof of the second result in Lemma~\ref{lemma:replace} from Gaussian $\tilde{\bbeta}_{1,d}$ to spherical $\tilde{\bbeta}_{1,d}$ differs by the factor $\|\tilde{\bbeta}_{1,d}\|_2/F_{1,d}$.
Note that in high dimensions, the norm of Gaussian $\tilde{\bbeta}_{1,d}$ ($\|\tilde{\bbeta}_{1,d}\|_2$) is tightly concentrated on $F_{1,d}$. Therefore, it is not hard to translate the proof from Gaussian version to spherical version.

Based on the analysis above, without loss of generality we could assume $\tilde{\bbeta}_{1,d}\sim\rmN(\0,[F_{1,d}^2/d]\Ib_d)$ in the following of the proof. The lemma below helps us further handle the quadratic form of the variance which appears later.
\begin{lemma}
\label{lemma:quadratic}
Let $\Ab\in\RR^{d\times d}$, and define the random vector $\hb\sim\rmN(0,(F_{1,d}^2/d)\Ib_d)$. Then we have
\begin{align*}
    \Var_{\hb}(\hb^\T\Ab\hb)=\frac{F_{1,d}^4}{d^2}\big(\|\Ab\|_{F}^2+\tr(\Ab^2)\big).
\end{align*}
\end{lemma}
The proof of Lemma~\ref{lemma:quadratic} is given
at the end of this section.
With this lemma, we are well-prepared to prove the second result in Lemma~\ref{lemma:replace}. Recall the definitions
$$\begin{aligned}
\bsigma(\xb)=&\big(\sigma_1(\xb^\T \bTheta_1^\T /\sqrt{d}),\sigma_2(\xb^\T \bTheta_2^\T /\sqrt{d})\big)^\T \in\RR^N,~\quad \bUpsilon=( \Zb ^\T  \Zb +\lambda\Ib_N)^{-1},
\end{aligned}$$
and
$\tilde{\Vb}=\EE_{\xb}\bsigma(\xb)(\xb^\T \tilde{\bbeta}_{1,d}+F_0),~\Ub=\EE_{\xb}\bsigma(\xb)\bsigma(\xb)^\T .$ By the definition of the risk $R_d(\Xb,\bTheta,\lambda,\tilde{\bbeta}_{1,d},\bvarepsilon)$, we have
\begin{align}   
    \nonumber
     R_d(\Xb,\bTheta,\lambda,\tilde{\bbeta}_{1,d},\bvarepsilon) & =  \EE_{\xb}\big(\xb^\T \tilde{\bbeta}_{d}+F_0-\hat{\ab}(\lambda)^\T \bsigma(\xb)\big)^2\\   
     &=  F_0^2+F_{1,d}^2-2\Gamma_1+\Gamma_2+\Gamma_3-2\Gamma_4+2\Gamma_5,
    \label{eq:Rdef_Gammas}
\end{align} 
where
$$
\begin{array}{llll}
&\bbf=\Xb\tilde{\bbeta}_{1,d}+\1_nF_0,&\quad \Gamma_1=\bbf^\T \Zb\bUpsilon\tilde{\Vb}/\sqrt{d},&\quad \Gamma_2=\bbf^\T \big[\Ub\big]_\Zb \bbf/d,\\
&\Gamma_2=\bvarepsilon^\T \big[\Ub\big]_\Zb \bvarepsilon/d,&\quad \Gamma_4=\bvarepsilon^\T \Zb\bUpsilon\tilde{\Vb}/\sqrt{d},&\quad \Gamma_5=\bvarepsilon^\T \big[\Ub\big]_\Zb \bbf/d.
\end{array}
$$
Note that the terms $F_0^2$ and  $F_{1,d}^2$ in \eqref{eq:Rdef_Gammas} are constants, and therefore do not contribute to the variance of $R_d(\Xb,\bTheta,\lambda,\tilde{\bbeta}_{1,d},\bvarepsilon)$. In the following, we aim to show that $\EE_{\Xb,\bTheta}\big[\Var_{\tilde{\bbeta}_{1,d},\bvarepsilon}(\Gamma_k) \big]=o_d(1)$ for $k\in[5]$.
Consider first  the variance  of $\Gamma_1$. From Lemma~\ref{lemma:decomp1}, we have
$\tilde{\Vb}=\bLambda_{d,1}{\bTheta}\tilde{\bbeta}_{1,d}+\bLambda_{d,0}\1_NF_0$. Then
\begin{align*}
    \Var_{\tilde{\bbeta}_{1,d}}(\Gamma_1) &=\Var_{\tilde{\bbeta}_{1,d}}\Big((\Xb\tilde{\bbeta}_{1,d}+\1_nF_0)^\T \Zb\bUpsilon\big(\bLambda_{d,1}\bTheta\tilde{\bbeta}_{1,d}+\bLambda_{d,0}\1_NF_0\big)/\sqrt{d}~\Big)\\
    &= \frac{1}{d}\Var_{\tilde{\bbeta}_{1,d}}\Big(\tilde{\bbeta}_{1,d}^\T\Xb^T\Zb\bUpsilon\bLambda_{d,1}\bTheta\tilde{\bbeta}_{1,d}+\tilde{\bbeta}_{1,d}^\T\Xb\Zb\bUpsilon\bLambda_{d,0}\1_NF_0\\
    &\quad+F_0\1_n^\T\Zb\bUpsilon\bLambda_{d,1}\bTheta\tilde{\bbeta}_{1,d}+F_0\1_n^T\Zb\bUpsilon\bLambda_{d,0}\1_NF_0\Big)\\
    &\leq \frac{4}{d}\Var_{\tilde{\bbeta}_{1,d}}\Big(\tilde{\bbeta}_{1,d}^\T\Xb^T\Zb\bUpsilon\bLambda_{d,1}\bTheta\tilde{\bbeta}_{1,d}\Big)+\frac{4}{d}\Var_{\tilde{\bbeta}_{1,d}}\Big(\tilde{\bbeta}_{1,d}^\T\Xb\Zb\bUpsilon\bLambda_{d,0}\1_NF_0\Big)
    \\
    &\quad+\frac{4}{d}\Var_{\tilde{\bbeta}_{1,d}}\Big(F_0\1_n^\T\Zb\bUpsilon\bLambda_{d,1}\bTheta\tilde{\bbeta}_{1,d}\Big)+0
    \\
    & \leq  4F_{1,d}^4\cdot\underbrace{\frac{1}{d^3}\Big(\|\Xb^T\Zb\bUpsilon\bLambda_{d,1}\bTheta\|_{F}^2+\tr\big(\Xb^T\Zb\bUpsilon\bLambda_{d,1}\bTheta\Xb^T\Zb\bUpsilon\bLambda_{d,1}\bTheta\big)\Big)}_{I_1} \\
    &\quad+4F_{1,d}^2F_0^2\cdot\underbrace{\bigg(\frac{1}{d}\tr\bigg( \frac{\Xb\Xb^\T }{d}\big[\bLambda_{d,0}\1_N\1_N^\T \bLambda_{d,0}\big]_\Zb\bigg)+\frac{1}{d}\tr\bigg( \1_n\1_n^\T \bigg[\bLambda_{d,1}\frac{\bTheta\bTheta^\T }{d}\bLambda_{d,1}\bigg]_\Zb\bigg)\bigg)}_{I_2}.
\end{align*}
The first inequality holds from $\Var(a+b)\leq2\Var(a)+2\Var(b)$, $I_1$ comes from Lemma~\ref{lemma:quadratic} and $I_2$ comes from $\Var(a)\leq\EE a^2$. Note that $\lVert\bLambda_{d,1}\rVert_{\op}=O_d(1/\sqrt{d})$  and $\|\Zb\bUpsilon\|_{\op}\leq1/(2\sqrt{\lambda})$, we conclude that
\begin{align*}
    \EE_{\Xb,\bTheta} |I_1|&\leq \Big|\frac{1}{d^3}\EE_{\Xb,\bTheta}\tr(\Xb^T\Zb\bUpsilon\bLambda_{d,1}\bTheta\bTheta^\T\bLambda_{d,1}\bUpsilon\Zb^\T\Xb) \Big| +\Big|\frac{1}{d^3}\EE_{\Xb,\bTheta}\tr\big(\Xb^T\Zb\bUpsilon\bLambda_{d,1}\bTheta\Xb^T\Zb\bUpsilon\bLambda_{d,1}\bTheta\big)\Big|\\
    &\leq\frac{1}{4\lambda}\EE_{\Xb,\bTheta}\Big\|\bLambda_{d,1}\frac{\bTheta\bTheta^\T}{d}\bLambda_{d,1} \Big\|_{\op}\cdot \Big\|\frac{\Xb\Xb^\T}{d}\Big\|_{\op}+\frac{1}{4\lambda} \EE_{\Xb,\bTheta}\Big\|\frac{\bLambda_{d,1}\bTheta\Xb^\T}{d}\Big\|_{\op}^2=o_d(1).
\end{align*}
Furthermore from Lemma~\ref{lemma:decomp2}, 
 we have
$$
\EE_{\Xb,\bTheta}|I_2|=\EE_{\Xb,\bTheta}\bigg|\frac{1}{d}\tr\bigg( \frac{\Xb\Xb^\T }{d}\big[\bLambda_{d,0}\1_N\1_N^\T \bLambda_{d,0}\big]_\Zb\bigg)+\frac{1}{d}\tr\bigg( \1_n\1_n^\T \bigg[\bLambda_{d,1}\frac{\bTheta\bTheta^\T }{d}\bLambda_{d,1}\bigg]_\Zb\bigg)\bigg|=o_d(1).
$$
Thus we obtain 
$\EE_{\Xb,\bTheta}\big(\Var_{\bbeta_d}(\Gamma_1) \big)=o_d(1)$.
Similarly, we have for $\Gamma_2$,

\begin{align*}
    \Var_{\tilde{\bbeta}_{1,d}}(\Gamma_2)=&\Var_{\tilde{\bbeta}_{1,d}}\Big((\Xb\tilde{\bbeta}_{1,d}+\1_nF_0)^\T \big[\Ub\big]_\Zb(\Xb\tilde{\bbeta}_{1,d}+\1_nF_0)/d~\Big)\\
    =& \frac{1}{d^2}\Var_{\tilde{\bbeta}_{1,d}}\Big(\tilde{\bbeta}_{1,d}^\T\Xb^T\big[\Ub\big]_\Zb\Xb\tilde{\bbeta}_{1,d}+\tilde{\bbeta}_{1,d}^\T\Xb^\T\big[\Ub\big]_\Zb\1_nF_0\\
    &\quad+F_0\1_n^\T\big[\Ub\big]_\Zb\Xb\tilde{\bbeta}_{1,d}+F_0\1_n^T\big[\Ub\big]_\Zb\1_nF_0\Big)\\
    \leq&\frac{4}{d^2}\Var_{\tilde{\bbeta}_{1,d}}\Big(\tilde{\bbeta}_{1,d}^\T\Xb^T\big[\Ub\big]_\Zb\Xb\tilde{\bbeta}_{1,d}\Big)+\frac{8}{d^2}\Var_{\tilde{\bbeta}_{1,d}}\Big(\tilde{\bbeta}_{1,d}^\T\Xb^\T\big[\Ub\big]_\Zb\1_nF_0\Big)
    \\
    \leq&8F_{1,d}^4\cdot\underbrace{\frac{1}{d^4}\tr\bigg(\Xb^T\big[\Ub\big]_\Zb\Xb\Xb^T\big[\Ub\big]_\Zb\Xb\bigg)}_{I_3}+8F_{1,d}^2F_0^2\cdot\underbrace{\frac{1}{d^2}\tr\bigg(\big[\Ub\big]_\Zb\1_n\1_n^\T\big[\Ub\big]_\Zb\frac{\Xb\Xb^\T}{d}\bigg)}_{I_4}.
\end{align*}
The first inequality holds from $\Var(a+b)\leq2\Var(a)+2\Var(b)$, $I_3$ comes from Lemma~\ref{lemma:quadratic} and the symmetric of $\Xb^T\big[\Ub\big]_\Zb\Xb$, and $I_4$ comes from $\Var(a)\leq\EE a^2$. Define $\bGamma_{\Ub}=\Ub-\bLambda_{d,0}\1_N\1_N^\T\bLambda_{d,0}$, from Lemma~\ref{lemma:decomp1}, $\EE\lVert\bGamma_{\Ub}\rVert_{\op}^2=O_d(1)$. By replacing $\Ub$ by $\bLambda_{d,0}\1_N\1_N^\T\bLambda_{d,0}+\bGamma_{\Ub}$ in the terms $I_3$ and $I_4$, we obtain the following equalities:
\begin{align*}
I_3=&\frac{1}{d^2}\tr\bigg(\big[\bLambda_{d,0}\1_N\1_N^\T\bLambda_{d,0}+\bGamma_{\Ub}\big]_\Zb\frac{\Xb\Xb^\T}{d}\big[\bLambda_{d,0}\1_N\1_N^\T\bLambda_{d,0}+\bGamma_{\Ub}\big]_\Zb\frac{\Xb\Xb^\T}{d}\bigg)\\
    =&\underbrace{\frac{1}{d^2}\tr\bigg(\big[\bLambda_{d,0}\1_N\1_N^\T\bLambda_{d,0}\big]_\Zb\frac{\Xb\Xb^\T}{d}\big[\bLambda_{d,0}\1_N\1_N^\T\bLambda_{d,0}\big]_\Zb\frac{\Xb\Xb^\T}{d}\bigg)}_{J_1}\\
    &+\underbrace{\frac{2}{d^2}\tr\bigg(\big[\bLambda_{d,0}\1_N\1_N^\T\bLambda_{d,0}\big]_\Zb\frac{\Xb\Xb^\T}{d}\big[\bGamma_{\Ub}\big]_\Zb\frac{\Xb\Xb^\T}{d}\bigg)}_{J_2}+\underbrace{\frac{1}{d^2}\tr\bigg(\big[\bGamma_{\Ub}\big]_\Zb\frac{\Xb\Xb^\T}{d}\big[\bGamma_{\Ub}\big]_\Zb\frac{\Xb\Xb^\T}{d}\bigg)}_{J_3},\\
    I_4=&\frac{1}{d^2}\tr\bigg(\big[\bLambda_{d,0}\1_N\1_N^\T\bLambda_{d,0}+\bGamma_{\Ub}\big]_\Zb\1_n\1_n^\T\big[\bLambda_{d,0}\1_N\1_N^\T\bLambda_{d,0}+\bGamma_{\Ub}\big]_\Zb\frac{\Xb\Xb^\T}{d}\bigg)\\
    =&\underbrace{\frac{1}{d^2}\tr\bigg(\big[\bLambda_{d,0}\1_N\1_N^\T\bLambda_{d,0}\big]_\Zb\1_n\1_n^\T\big[\bLambda_{d,0}\1_N\1_N^\T\bLambda_{d,0}\big]_\Zb\frac{\Xb\Xb^\T}{d}\bigg)}_{K_1}\\
    &+\underbrace{\frac{2}{d^2}\tr\bigg(\big[\bLambda_{d,0}\1_N\1_N^\T\bLambda_{d,0}\big]_\Zb\1_n\1_n^\T\big[\bGamma_{\Ub}\big]_\Zb\frac{\Xb\Xb^\T}{d}\bigg)}_{K_2}+\underbrace{\frac{1}{d^2}\tr\bigg(\big[\bGamma_{\Ub}\big]_\Zb\1_n\1_n^\T\big[\bGamma_{\Ub}\big]_\Zb\frac{\Xb\Xb^\T}{d}\bigg)}_{K_3}.
\end{align*}
We investigate the terms $K_i$, $i=1,2,3$. The investigation of terms $J_i$, $i=1,2,3$ are quite similar, we omit the proof for $J_i$ for brevity. Consider first the term $K_2$. Due to $\EE\lVert\bGamma_{\Ub}\rVert_{\op}^2=O_d(1)$, it is true that
\begin{align*}
    \bigg(\EE\bigg\lVert\frac{\1_n\1_n^\T}{d}\big[\bGamma_{\Ub}\big]_\Zb\frac{\Xb\Xb^\T}{d} \bigg\rVert_{\op}^2\bigg)^{1/2}=&O_d(1)\cdot\bigg(\EE\bigg\lVert\bGamma_{\Ub}\frac{\Xb\Xb^\T}{d}\bigg\rVert_{\op}^2\bigg)^{1/2}\\
    =&O_d(1)\cdot\bigg(\EE\Big\lVert\bGamma_{\Ub}\Big\rVert_{\op}^2\bigg)^{1/2}\cdot\bigg(\EE\bigg\lVert\frac{\Xb\Xb^\T}{d}\bigg\rVert_{\op}^2\bigg)^{1/2}=O_d(1).
\end{align*}
The second equality comes from the independence of $\bGamma_{\Ub}$ and $\Xb$. Note that for any rank $1$ matrix $\Ab$, $|\tr\Ab|=\lVert\Ab\rVert_{\op}$, the term $K_2$ has the property

\begin{align*}
 \EE_{\Xb,\bTheta}&|K_2|=\EE_{\Xb,\bTheta}\bigg|\frac{2}{d^2}\tr\bigg(\bigg[\bLambda_{d,0}\1_N\1_N^\T\bLambda_{d,0}\big]_\Zb\1_n\1_n^\T\big[\bGamma_{\Ub}\big]_\Zb\frac{\Xb\Xb^\T}{d}\bigg) \bigg|\\
 &\leq \EE_{\Xb,\bTheta}\bigg(\frac{2}{d}\bigg\lVert\big[\bLambda_{d,0}\1_N\1_N^\T\bLambda_{d,0}\big]_\Zb \bigg\rVert_{\op}\cdot\bigg\lVert \frac{\1_n\1_n^\T}{d}\big[\bGamma_{\Ub}\big]_\Zb\frac{\Xb\Xb^\T}{d} \bigg\rVert_{\op}\bigg)\\
 &\leq \frac{2}{d}\bigg(\EE_{\Xb,\bTheta}\bigg\lVert\big[\bLambda_{d,0}\1_N\1_N^\T\bLambda_{d,0}\big]_\Zb \bigg\rVert_{\op}^2\cdot\EE_{\Xb,\bTheta}\bigg\lVert \frac{\1_n\1_n^\T}{d}\big[\bGamma_{\Ub}\big]_\Zb\frac{\Xb\Xb^\T}{d} \bigg\rVert_{\op}^2\bigg)^{1/2}\\
 &=o_d(1)\cdot O_d(1)=o_d(1).
\end{align*}
The equality comes from the estimation of $\cD$ in Lemma~\ref{lemma:decomp2} . For the term $K_1$, it is true from the estimation of $\cD$ in Lemma~\ref{lemma:decomp2} 
that 
\begin{align*}
    & \frac{1}{d}\bigg(\EE\bigg\lVert\frac{\1_n\1_n^\T}{d}\big[\bLambda_{d,0}\1_N\1_N^\T\bLambda_{d,0}\big]_\Zb\frac{\Xb\Xb^\T}{d} \bigg\rVert_{\op}^2 \bigg)^{1/2}\\
    &\quad \leq \frac{1}{d}\bigg\lVert\frac{\1_n\1_n^\T}{d}\bigg\rVert_{\op}\bigg(\EE\bigg\lVert\big[\bLambda_{d,0}\1_N\1_N^\T\bLambda_{d,0}\big]_\Zb\frac{\Xb\Xb^\T}{d}\bigg\rVert_{op}^2 \bigg)^{1/2}
    = O_d(1)\cdot o_d(1)=o_d(1).
\end{align*}
By repeating the arguments used previously for the term $K_2$ but for the consideration of  $K_1$, we have 
\begin{align*}
    &\EE_{\Xb,\bTheta}|K_1|=\frac{1}{d^2}\EE_{\Xb,\bTheta}\bigg|\tr\bigg(\big[\bLambda_{d,0}\1_N\1_N^\T\bLambda_{d,0}\big]_\Zb\1_n\1_n^\T\big[\bLambda_{d,0}\1_N\1_N^\T\bLambda_{d,0}\big]_\Zb\frac{\Xb\Xb^\T}{d}\bigg)\bigg|=o_d(1).
\end{align*}
For the term $K_3$, similarly we have
\begin{align*}
    \EE_{\Xb,\bTheta}|K_3|&=\EE_{\Xb,\bTheta}\bigg|\frac{1}{d^2}\tr\bigg(\big[\bGamma_{\Ub}\big]_\Zb\1_n\1_n^\T\big[\bGamma_{\Ub}\big]_\Zb\frac{\Xb\Xb^\T}{d}\bigg)\bigg|\\
    &\leq \frac{1}{d^2}\EE_{\Xb,\bTheta}\bigg( \big\lVert \big[\bGamma_{\Ub}\big]_\Zb\1_n\1_n^\T\big\rVert_{\op}\cdot\bigg\lVert\big[\bGamma_{\Ub}\big]_\Zb\frac{\Xb\Xb^\T}{d} \bigg\rVert_{\op} \bigg)\\
    &\leq \frac{1}{d^2}\bigg(\EE_{\Xb,\bTheta} \big\lVert \big[\bGamma_{\Ub}\big]_\Zb\1_n\1_n^\T\big\rVert_{\op}^2\cdot \EE_{\Xb,\bTheta} \bigg\lVert \big[\bGamma_{\Ub}\big]_\Zb\frac{\Xb\Xb^\T}{d} \bigg\rVert_{\op}^2 \bigg)^{1/2}\\
    &\leq\frac{1}{d^2}\big(\EE_{\Xb,\bTheta} \big\lVert \big[\bGamma_{\Ub}\big]_\Zb\1_n\1_n^\T\big\rVert_{\op}^2\big)^{1/2}\cdot O_d(1)\\
    &=o_d(1)\cdot O_d(1)=o_d(1).
\end{align*}
Now we conclude that $K_1$, $K_2$ and $K_3$ are all small terms under the expectation over $\Xb$ and $\bTheta$, we immediately get that 
\begin{align*}
   \EE_{\Xb,\bTheta}|I_4|=o_d(1).
\end{align*}
Similarly we get $\EE_{\Xb,\bTheta}|I_3|=o_d(1)$, thus we conclude that $\EE_{\Xb,\bTheta}\Var_{\tilde{\bbeta}_{1,d}}(\Gamma_2)=o_d(1)$. We omit the other terms for brevity. 
The proof of Lemma~\ref{lemma:replace} is complete.

\subsection*{Proof of Lemma~\ref{lemma:quadratic}}\label{subsec:proof_quadratic}

We have
\begin{align*}
    \EE[\hb^\T\Ab\hb]=\EE\tr(\Ab\hb\hb^\T)=\frac{F_{1,d}^2}{d}\tr(\Ab).
\end{align*}
Hence we have 
\begin{align*}
    \Var(\hb^\T\Ab\hb)&=\sum\limits_{i_1,i_2,i_3,i_4}\EE\Big[\hb_{i_1}\Ab_{i_1,i_2}\hb_{i_2}\hb_{i_3}\Ab_{i_3,i_4}\hb_{i_4}\Big]-\frac{F_{1,d}^4}{d^2}\tr(\Ab)^2\\
    &= \Bigg\{\bigg(\sum\limits_{i_1=i_2,i_3=i_4\atop i_1= i_3}+\sum\limits_{i_1=i_2,i_3=i_4\atop i_1\neq i_3}+\sum\limits_{i_1=i_3,i_2=i_4\atop i_1\neq i_2}+\sum\limits_{i_1=i_4,i_2=i_3\atop i_1\neq i_2}\bigg)
    \EE\Big[\hb_{i_1}\Ab_{i_1,i_2}\hb_{i_2}\hb_{i_3}\Ab_{i_3,i_4}\hb_{i_4}\Big]\Bigg\}\\
    &\quad-\frac{F_{1,d}^4}{d^2}\tr(\Ab)^2\\
    &=\frac{F_{1,d}^4}{d^2}\bigg(\sum_{i=1}^d\Ab_{i,i}^2\cdot\frac{d^2}{F_{1,d}^4}(\EE h_i^4)+\sum\limits_{i\neq j}\Ab_{i,i}\Ab_{j,j}+\sum\limits_{i\neq j}(\Ab_{i,j}\Ab_{i,j}+\Ab_{i,j}\Ab_{j,i})-\tr(\Ab)^2\bigg)\\
    &=\frac{F_{1,d}^4}{d^2}\bigg(\sum_{i=1}^d3\Ab_{i,i}^2+\sum\limits_{i\neq j}(\Ab_{i,j}\Ab_{i,j}+\Ab_{i,j}\Ab_{j,i})+\sum\limits_{i\neq j}\Ab_{i,i}\Ab_{j,j}-\tr(\Ab)^2\bigg)\\
    &=\frac{F_{1,d}^4}{d^2}\bigg(\sum\limits_{i,j}(\Ab_{i,j}\Ab_{i,j}+\Ab_{i,j}\Ab_{j,i})+\sum\limits_{i,j}\Ab_{i,i}\Ab_{j,j}-\tr(\Ab)^2\bigg)=\frac{F_{1,d}^4}{d^2}\big(\|\Ab\|_{F}^2+\tr(\Ab^2)\big).
\end{align*}
This proves Lemma~\ref{lemma:quadratic}.

\section{Proof of Proposition~\ref{prop:connectGd}}\label{sec:appendixconnectedGd}
  We first recall  the following matrix differential rules in \citet{petersen2008matrix}:
\begin{align}
\frac{\partial\det (\Yb)}{\partial x}=\det (\Yb)\cdot\tr\Big(\Yb^{-1}\cdot\frac{\partial\Yb}{\partial x}\Big),\quad\frac{\partial\Yb^{-1}}{\partial x}=-\Yb^{-1}\frac{\partial\Yb}{\partial x}\Yb^{-1}.\label{eq:matrix_derivative_rule}
\end{align}
Let $q_i,q_j$ be the elements in the vector $\qb$. Now  the matrix $\Ab = \Ab(\qb,\bmu)$ (see  Definition~\ref{def:linear pencil}) is linear  in  $\qb$, thus $\frac{\partial^2 \Ab}{\partial q_i\partial q_j}=0$. 
Therefore by the definition  $G_d(\xi)=\frac{1}{d}\Log\det\big(\Ab-\xi\Ib\big)$ and the matrix derivative rules in \eqref{eq:matrix_derivative_rule}, we have
\begin{align}
        \frac{\partial G_d}{\partial q_i}=&\frac{1}{d}\tr\Big(\big(\Ab-\xi\Ib\big)^{-1}\frac{\partial\Ab}{\partial  q_i}\Big),\label{eq:one_order}\\
        \frac{\partial^2 G_d}{\partial q_i\partial q_j}=&\frac{1}{d}\tr\bigg(\frac{\partial\big(\Ab-\xi\Ib\big)^{-1}}{\partial q_j}\frac{\partial\Ab}{\partial  q_i}\bigg)=-\frac{1}{d}\tr\Big(\big(\Ab-\xi\Ib\big)^{-1}\frac{\partial\Ab}{\partial  q_j}\big(\Ab-\xi\Ib\big)^{-1}\frac{\partial\Ab}{\partial  q_i}\Big).\label{eq:two_order}
\end{align}
By the Schur complement formula, we further have
\begin{equation*}
    \begin{split}
        (\Ab(\0,\bmu)-\xi\Ib_{\PO})^{-1}&=\begin{bmatrix}
        -\xi\Ib_{N}&\Zb^\T\\
        \Zb&-\xi\Ib_n
        \end{bmatrix}^{-1}=\begin{bmatrix}
        *&\Zb^\T(\Zb\Zb^\T-\xi^2\Ib_n)^{-1}\\
        (\Zb\Zb^\T-\xi^2\Ib_n)^{-1}\Zb&*
        \end{bmatrix},
    \end{split}
\end{equation*}
where we use $*$ to hide the irrelevant blocks in the matrix inverse. 
Moreover, by definition, it holds that
\begin{align*}
    &
    \frac{\partial\Ab(\qb,\bmu)}{\partial q_1} =\begin{bmatrix}
        \0&\frac{1}{d}\Mb_1\bTheta\Xb^\T\\
       \frac{1}{d}\Xb\bTheta^\T\Mb_1&\0\end{bmatrix},
    \quad\frac{\partial\Ab(\qb,\bmu)}{\partial q_4 } =\begin{bmatrix}
        \Mb_1\frac{\bTheta\bTheta^\T}{d}\Mb_1&~~\0\\
        \0&~~\0
    \end{bmatrix},\\
    &
    \frac{\partial\Ab(\qb,\bmu)}{\partial q_2 } =\begin{bmatrix}
     \Mb_2^2&~~\0\\
     \0&~~\0
    \end{bmatrix},
    \quad 
    \frac{\partial\Ab(\qb,\bmu)}{\partial q_3 } =\begin{bmatrix}
        \0&~~\0\\
        \0&~~\Ib_n
    \end{bmatrix},
    \quad
    \frac{\partial\Ab(\qb,\bmu)}{\partial q_5 } =\begin{bmatrix}
    \0&~~\0\\
    \0&~~\frac{\Xb\Xb^\T}{d}
    \end{bmatrix}.
\end{align*}
\if UT
\begin{align*}
    &\frac{\partial\Ab(\qb,\bmu)}{\partial q_2 }\Bigg|_{\qb=\0}=\begin{bmatrix}
    \Mb_2^2&~~\0\\
    \0&~~\0
    \end{bmatrix},\quad\frac{\partial\Ab(\qb,\bmu)}{\partial q_4 }\Bigg|_{\qb=\0}=\begin{bmatrix}
    \Mb_1\frac{\bTheta\bTheta^\T}{d}\Mb_1&~~\0\\
    \0&~~\0
    \end{bmatrix},\\
    &\frac{\partial\Ab(\qb,\bmu)}{\partial q_3 }\Bigg|_{\qb=\0}=\begin{bmatrix}
    \0&~~\0\\
    \0&~~\Ib_n
    \end{bmatrix},\quad\frac{\partial\Ab(\qb,\bmu)}{\partial q_5 }\Bigg|_{\qb=\0}=\begin{bmatrix}
    \0&~~\0\\
    \0&~~\frac{\Xb\Xb^\T}{d}
    \end{bmatrix},\quad\frac{\partial\Ab(\qb,\bmu)}{\partial q_1}\Bigg|_{\qb=\0}=\begin{bmatrix}
        \0&\frac{1}{d}\Mb_1\bTheta\Xb^\T\\
        \frac{1}{d}\Xb\bTheta^\T\Mb_1&\0\end{bmatrix}.
\end{align*}
\fi
By   plugging these derivatives   into \eqref{eq:one_order} and \eqref{eq:two_order}, 
we continue the calculation with $\xi=\xi^*$. 
Note that we have the identity $(\Zb\Zb^\T-(\xi^*)^2\Ib_n)^{-1}\Zb=\Zb(\Zb^\T\Zb+\lambda\Ib_N)^{-1}=\Zb\bUpsilon$. Then by \eqref{eq:one_order}, we have
\begin{align*}
\frac{\partial G_d(\xi^*;\qb,\bmu)}{\partial q_1}\Bigg|_{\qb=\0}=&\frac{1}{d}\tr\left(\begin{bmatrix}
        *&\Zb^\T(\Zb\Zb^\T-(\xi^*)^2\Ib_n)^{-1} \\
        (\Zb\Zb^\T-(\xi^*)^2\Ib_n)^{-1}\Zb&
        *
        \end{bmatrix}\begin{bmatrix}
        \0&\frac{1}{d}\Mb_1\bTheta\Xb^\T\\
        \frac{1}{d}\Xb\bTheta^\T\Mb_1&\0\end{bmatrix}\right)\\
        =&\frac{1}{d}\tr\left(\begin{bmatrix}
        *&\bUpsilon\Zb^\T \\
        \Zb\bUpsilon&
        *
        \end{bmatrix}\begin{bmatrix}
        \0&\frac{1}{d}\Mb_1\bTheta\Xb^\T\\
        \frac{1}{d}\Xb\bTheta^\T\Mb_1&\0\end{bmatrix}\right)=\frac{2}{d}\tr \Mb_1\frac{\bTheta\Xb^\T }{d}\Zb\bUpsilon,
\end{align*}
Similarly, by \eqref{eq:two_order}, we  have
\begin{align*}
-\frac{\partial^2G_d(\xi^*;\qb,\bmu)}{\partial q_4\partial q_5}\Bigg|_{\qb=\0} =&\tr\left(\begin{bmatrix}
        *&\bUpsilon\Zb^\T \\
        \Zb\bUpsilon&
        *
        \end{bmatrix} \begin{bmatrix}
\Mb_1\frac{\bTheta\bTheta^\T}{d}\Mb_1&\0\\
    \0&\0
\end{bmatrix}\begin{bmatrix}
        *&\bUpsilon\Zb^\T \\
        \Zb\bUpsilon&
        *
        \end{bmatrix}\begin{bmatrix}
    \0&\0\\
    \0&\frac{\Xb\Xb^\T}{d}
    \end{bmatrix}\right)\\
=&\frac{1}{d}\tr \Zb \bUpsilon\Mb_1\frac{\bTheta\bTheta^\T}{d}\Mb_1\bUpsilon  \Zb ^\T \frac{\Xb\Xb^\T }{d},\\
-\frac{\partial^2G_d(\xi^*;\qb,\bmu)}{\partial q_2\partial q_5}\Bigg|_{\qb=\0} =&\tr\left(\begin{bmatrix}
        *&\bUpsilon\Zb^\T \\
        \Zb\bUpsilon&
        *
        \end{bmatrix} \begin{bmatrix}
\Mb_2\Mb_2&\0\\
    \0&\0
\end{bmatrix}\begin{bmatrix}
        *&\bUpsilon\Zb^\T \\
        \Zb\bUpsilon&
        *
        \end{bmatrix}\begin{bmatrix}
    \0&\0\\
    \0&\frac{\Xb\Xb^\T}{d}
    \end{bmatrix}\right)\\
=&\frac{1}{d}\tr \Zb \bUpsilon\Mb_2\Mb_2\bUpsilon  \Zb ^\T \frac{\Xb\Xb^\T }{d},\\
-\frac{\partial^2G_d(\xi^*;\qb,\bmu)}{\partial q_3\partial q_4} \Bigg|_{\qb=\0}=&\tr\left(\begin{bmatrix}
        *&\bUpsilon\Zb^\T \\
        \Zb\bUpsilon&
        *
        \end{bmatrix} \begin{bmatrix}
\Mb_1\frac{\bTheta\bTheta^\T}{d}\Mb_1&\0\\
    \0&\0
\end{bmatrix}\begin{bmatrix}
        *&\bUpsilon\Zb^\T \\
        \Zb\bUpsilon&
        *
        \end{bmatrix}\begin{bmatrix}
    \0&\0\\
    \0&\Ib_n
    \end{bmatrix}\right)\\
=&\frac{1}{d}\tr \Zb \bUpsilon\Mb_1\frac{\bTheta\bTheta^\T}{d}\Mb_1\bUpsilon  \Zb ^\T ,\\
-\frac{\partial^2G_d(\xi^*;\qb,\bmu)}{\partial q_2\partial q_3} \Bigg|_{\qb=\0}=&\tr\left(\begin{bmatrix}
        *&\bUpsilon\Zb^\T \\
        \Zb\bUpsilon&
        *
        \end{bmatrix} \begin{bmatrix}
\Mb_2\Mb_2&\0\\
    \0&\0
\end{bmatrix}\begin{bmatrix}
        *&\bUpsilon\Zb^\T \\
        \Zb\bUpsilon&
        *
        \end{bmatrix}\begin{bmatrix}
    \0&\0\\
    \0&\Ib_n
    \end{bmatrix}\right)\\
=&\frac{1}{d}\tr \Zb \bUpsilon\Mb_2\Mb_2\bUpsilon  \Zb ^\T.
\end{align*}
The above equations complete the proof of Proposition~\ref{prop:connectGd}.

\section{Properties of the fixed point equation}\label{sec:appendixuniquesolution}
In this section, we justify the definition of $\mb(\xi;\qb,\bmu)$ below Definition~\ref{def:implicit1}, by 
proving that there exists a constant $\xi_0 > 0$, such that the fixed point equation \eqref{eq:m123} has a unique solution defined on $\{\xi:\Im(\xi) > \xi_0\}$ satisfying $| m_j(\xi) |\leq2\psi_j/\xi_0$
for $j=1,2,3$. 
The result is given in the following lemma.

\begin{lemma}
\label{lemma:uniquesolution}
Let $\sfFb(\mb;\xi,\qb,\bmu)$, $\qb\in\cQ$ be defined in Definition~\ref{def:implicit1}, and $\DD(r)=\{z:|z|<r\}$ be the disk of radius $r$ in the complex plane. There exists $\xi_0>0$ such that, for any $\xi\in\bbC_+ $ with $\Im(\xi)>\xi_0$, $\sfFb(\mb;\xi,\qb,\bmu)$ is $1/2$-Lipschitz continuous with respect to the $\ell_2$ norm, and the map $\mb\mapsto \sfFb(\mb;\xi,\qb,\bmu)$
admits a unique fixed point  in $\DD(2\psi_1/\xi_0)\times\DD(2\psi_2/\xi_0)\times\DD(2\psi_3/\xi_0)$.
\end{lemma}
Lemma \ref{lemma:uniquesolution} demonstrates that our definition of $\mb$ in Subsection~\ref{step3} as the unique fixed point of $\sfFb$ is valid.
 
\begin{proof}[Proof of Lemma~\ref{lemma:uniquesolution}]
We prove the existence and uniqueness of the solution by the Banach fixed point theorem when $\Im(\xi) \geq \xi_0$ for some sufficiently large $\xi_0$. To do so, we want to show that 
\begin{enumerate}[leftmargin = *]
    \item $\sfFb(\cdot;\qb,\bmu)$ maps domain $\DD(2\psi_1/\xi_0)\times\DD(2\psi_2/\xi_0)\times\DD(2\psi_3/\xi_0)$ into itself. 
    \item $\sfFb(\cdot;\qb,\bmu)$ is Lipschitz continuous with a Lipschitz   constant smaller than 1.
\end{enumerate}
For  $\sfF_1(\cdot;\qb,\bmu)$, 
by Definition~\ref{def:implicit1}, we have
\begin{align*}
    \sfF_1(\mb;\xi,\qb,\bmu)=&\frac{\psi_1}{-\xi+q_2\mu_{1,2}^2+H_1(\mb;\qb,\bmu)},
\end{align*}
where 
\begin{align}
    H_1(\mb;\qb,\bmu)=&-\mu_{1,2}^2m_3+\frac{1}{m_1+\frac{-\mu_{2,1}^2(1+q_1)^2m_2m_3+(1+\mu_{2,1}^2m_2q_4)(1+m_3q_5)}{\mu_{1,1}^2q_4(1+m_{3}q_5)-\mu_{1,1}^2(1+q_1)^2m_{3}}}.\label{eq:H(m)}
\end{align}
Note that $q_4,q_5\leq (1+q_1)/2$. Thus  for small enough $r_0$, we have for any $\mb\in\DD(r_0)^3$
\begin{align}\label{eq:fixedpoint_H1bound}
    |H_1(\mb;\qb,\bmu)|\leq 2+2|q_4|\mu_{1,1}^2.
\end{align}
Now as long as $\xi_0 \geq 4+4|q_4|\mu_{1,1}^2$, it is clear that for $\xi$ with $\Im(\xi) \geq \xi_0$ we have
\begin{align}\label{eq:fixedpoint_Imxi>H1}
    \Im(\xi) \geq \xi_0/ 2 + \xi_0/ 2 \geq \xi_0/ 2 +  2+2|q_4|\mu_{1,1}^2 \geq \xi_0/ 2 +  |H_1(\mb;\qb,\bmu)|.
\end{align} 
Therefore,  
\begin{align*}
    | \sfF_1(\mb;\xi,\qb,\bmu) |  
    \leq&\frac{\psi_1}{|\Im\big(\xi-q_2\mu_{1,2}^2-H_1(\mb;\qb,\bmu)\big)|}\\
    \leq & \frac{\psi_1}{\Im(\xi)-|H_1(\mb;\qb,\bmu)|}
    \leq \frac{2\psi_1}{\xi_0}
    ,
\end{align*}
where the last  inequality   follows from \eqref{eq:fixedpoint_Imxi>H1}.

Similarly, for $\sfF_2$ and $\sfF_3$ we show that   
$|\sfF_2(\mb;\xi,\qb,\bmu)|\leq2\psi_2/\xi_0$ provided $\Im(\xi)\geq\xi_0\geq4+4|q_4|\mu_{2,1}^2$, and $|\sfF_3(\mb;\xi,\qb,\bmu)|\leq2\psi_3/\xi_0$ provided $\Im(\xi)\geq\xi_0\geq4+4|q_5|$. 
Therefore if $\xi_0$ satisfies
$2\max\{\psi_1,\psi_2,\psi_3\}/\xi_0\leq r_0$ and $\xi_0 \geq 4+4 \max\{ |q_4|\mu_{1,1}^2, |q_4|\mu_{2,1}^2, |q_5| \}$, $\sfFb$ maps domain $\DD(2\psi_1/\xi_0)\times\DD(2\psi_2/\xi_0)\times\DD(2\psi_3/\xi_0)$ into itself.

As for the Lipschitz continuity of $\sfFb(\cdot;\qb,\bmu)$, note that
$$
\nabla_{\mb}\sfF_1(\mb;\xi,\qb,\bmu)= - \frac{\psi_1}{(-\xi+q_2\mu_{1,2}^2+H_1(\mb;\qb,\bmu))^2} \cdot \nabla_{\mb}H_1(\mb;\qb,\bmu).
$$
With the same calculation as above,  it is easy to see that when $\xi_0$ is sufficiently large, 
$\lVert\nabla_{\mb}H_1(\mb;\qb,\bmu)\rVert_2\leq C(\qb,\bmu)$ for all $\mb\in\DD(2\psi_1/\xi_0)\times\DD(2\psi_2/\xi_0)\times\DD(2\psi_3/\xi_0)$, where $C(\qb,\bmu)$ is a constant that only depends on $\qb$ and $\bmu$.
Thus for such  $\xi_0$ and $\xi$ with $\Im(\xi) \geq \xi_0$, 
$$
\big\lVert\nabla_{\mb}\sfF_1(\mb;\xi,\qb,\bmu)\big\rVert_2\leq \frac{C(\qb,\bmu) \cdot \psi_1}{\Im(\xi)-|H_1(\mb;\qb,\bmu)|} \leq \frac{4 C(\qb,\bmu) \cdot \psi_1 }{\xi_0} \leq \frac{1}{4},
$$
where we again utilize \eqref{eq:fixedpoint_Imxi>H1}. 
We can apply the same argument for $\sfF_2$ and $\sfF_3$, and conclude that $\sfFb$ is $1/2$-Lipschitz on $\mb\in\DD(2\psi_1/\xi_0)\times\DD(2\psi_2/\xi_0)\times\DD(2\psi_3/\xi_0)$. Therefore by Banach fixed point theorem, there exists a unique fixed point of $\sfFb$. Thus the solution of the implicit equations defined in Definition~\ref{def:implicit1} exists and is unique.
\end{proof}

\section{Proof of Proposition~\ref{prop:implicit1}}
\label{sec:appendixpropimplicit1}
The proof of Proposition~\ref{prop:implicit1} is split  into several sections. 
In Sections~\ref{sec:equivalenceinGaussianandSphere} and \ref{subsec:resolvent}, we give some useful preliminary results.
In Section~\ref{subsec:conclusion12}, we show that $\mb(\xi;\qb,\bmu)$ is analytic on $\{\xi:\Im(\xi) \geq \xi_0\}$, and then prove the first and second conclusions of Proposition~\ref{prop:implicit1}. In Section~\ref{subsec:proofpoly}, we  prove the point convergence of $ M_d(\xi;\qb,\bmu)$ to $\mb(\xi;\qb,\bmu)$ under the additional assumption that $\sigma_\sfc(x)$, $\sfc=1,2$ are polynomials. In Section~\ref{subsec:completeproof}, we extend this point convergence result to general activation functions satisfying Assumption~\ref{assump1}.
In Section~\ref{subsec:uniformconvergence}, we conclude the proof by showing the uniform convergence of $ M_d(\xi;\qb,\bmu)$ to $\mb(\xi;\qb,\bmu)$ on compact sets.

\subsection{Equivalence between Gaussian and spherical versions}
\label{sec:equivalenceinGaussianandSphere}
The first step in the proof of Proposition~\ref{prop:implicit1} is to relate the Stieltjes transform $M_d(\xi;\qb,\bmu)$  to the Stieltjes transform corresponding to Gaussian data and Gaussian random features. 

\begin{definition}\label{def:Gaussian_data_rf}
Let $(\overline{\btheta}_{a})_{a\in[N]}$ be i.i.d. standard Gaussian random vectors distributed as  $\rmN(\0,\Ib_d)$, and  $\overline{\bTheta}\in\RR^{N\times d}$ be the matrix whose $a^\text{th}$ row is given by $\overline{\btheta}_{a}$. Similarly, we denote $(\overline{\xb}_{i})_{i\in[n]}\sim_{\iid} \rmN(\0,\Ib_d)$, and let $\overline{\Xb}\in\RR^{n\times d}$ be the matrix whose $a^\text{th}$ row is $\overline{\xb}_{i}$.
\end{definition}
Given these definitions, our original data inputs and random feature parameters which are distributed uniformly on the sphere $\sqrt{d}\cdot\SSS^{d-1}$ can be represented as
\begin{align}\label{eq:x_theta_equiv_def}
    \xb_i=\sqrt{d}\cdot\frac{\overline{\xb}_{i}}{\lVert\overline{\xb}_{i}\rVert_2}\sim\Unif(\sqrt{d}\cdot\SSS^{d-1}), \text{ and } \btheta_a=\sqrt{d}\cdot\frac{\overline{\btheta}_{a}}{\lVert\overline{\btheta}_{a}\rVert_2}\sim\Unif(\sqrt{d}\cdot\SSS^{d-1})
\end{align}
for all $i \in [n]$ and $a \in [N]$. 
We can now consider the ``Gaussian version'' of the learning problem, where the data inputs are $(\overline{\xb}_{i})_{i\in[n]}$, and the double random feature model uses random parameters $(\overline{\btheta}_{a})_{a\in[N]}$ and activation functions
\begin{align}\label{eq:def_activation_Gaussian}
    \phi_{\sfc}(x)=\sigma_{\sfc}(x)-\EE_{G\sim\rmN(0,1)}[\sigma_{\sfc}(G)],\quad \sfc = 1,2.
\end{align}
For this version of the learning problem, we can similarly construct the linear pencil matrix $\overline{\Ab}(\qb,\bmu)$, which is the counterpart of the linear pencil matrix $\Ab(\qb,\bmu)$ defined in Definition~\ref{def:linear pencil}. 
\begin{definition}
\label{def:linear pencil Gaussian}
The linear pencil matrix $\overline\Ab(\qb,\bmu)\in\RR^{\PO\times \PO}$ ($\PO=N+n$) is defined as
\begin{equation*}
\overline{\Ab}(\qb,\bmu)=\begin{bmatrix}
q_2\mu_{1,2}^2\Ib_{N_1}+q_4\mu_{1,1}^2\frac{\overline{\bTheta}_1\overline{\bTheta}_1^\T }{d}&q_4\mu_{1,1}\mu_{2,1}\frac{\overline{\bTheta}_1\overline{\bTheta}_2^\T }{d}& \Jb_1^\T +q_1\tilde{\Jb}_1^\T \\
q_4\mu_{1,1}\mu_{2,1}\frac{\overline{\bTheta}_2\overline{\bTheta}_1^\T }{d}&q_2\mu_{2,2}^2\Ib_{N_2}+q_4\mu_{2,1}^2\frac{\overline{\bTheta}_2\overline{\bTheta}_2^\T }{d}& \Jb_2^\T +q_1\tilde{\Jb}_2^\T \\
 \Jb_1+q_1\tilde{\Jb}_1& \Jb_2+q_1\tilde{ \Jb}_2&q_3\Ib_n+q_5\frac{\overline{\Xb}~\overline{\Xb}^\T }{d}
\end{bmatrix},
\end{equation*}
where $\Jb_{\sfc}=\phi_{\sfc}\big(\overline{\Xb}~\overline{\bTheta}_{\sfc}^\T /\sqrt{d}\big)/\sqrt{d}$, $\tilde{\Jb}_{\sfc}=\frac{\mu_{\sfc,1}}{d}\overline{\Xb}~\overline{\bTheta}_{\sfc}^\T $, $\sfc=1,2$.
\end{definition}
We also define $ \overline{M}_d(\xi;\qb,\bmu)=\frac{1}{d}\tr\big[(\overline{\Ab}(\qb,\bmu)-\xi\Ib_{\PO})^{-1}\big]$ as the counterpart of the Stieltjes transform $ M_d(\xi;\qb,\bmu)$.
The following lemma establishes the equivalence between the two versions of the Stieltjes transforms. 
\begin{lemma}
\label{lemma:GaussianSphere}
Suppose that $\sigma_\sfc(x)$, $\sfc=1,2$, are polynomials. Then
for any fixed $\qb$ and $\xi\in\bbC_+$,  we have
\begin{equation*}
    \EE\big|\overline{M}_d(\xi;\qb,\bmu)-M_d(\xi;\qb,\bmu)\big|=o_d(1).
\end{equation*}
\end{lemma}


\begin{proof}[Proof of Lemma~\ref{lemma:GaussianSphere}]
Define $$\bDelta(\Ab,\overline{\Ab},\xi)=M_d(\xi;\qb,\bmu)-\overline{M}_d(\xi;\qb,\bmu),$$
and write $M_d(\xi;\qb,\bmu)$ and $\overline{M}_d(\xi;\qb,\bmu)$ as $M_d(\xi)$ and $\overline{M}_d(\xi)$ to simplify the notation. Then by definition we have
\begin{align}
    |\bDelta(\Ab,\overline{\Ab},\xi)|=&\big|\tr\big[ (\Ab-\xi\Ib)^{-1}(\Ab-\overline{\Ab})(\overline{\Ab}-\xi\Ib)^{-1} \big] \big|/d\nonumber\\
    \leq&\lVert(\Ab-\xi\Ib)^{-1}(\overline{\Ab}-\xi\Ib)^{-1}\rVert_{\op}\lVert \Ab-\overline{\Ab}\rVert_\star/d\nonumber\\
    \leq& \lVert \Ab-\overline{\Ab}\rVert_\star\cdot\frac{1}{d}\cdot\frac{1}{(\Im(\xi))^2},\label{eq:Delta_bound_AAbar}
\end{align}
where $\lVert\cdot\rVert_\star$ is the nuclear norm, the first inequality follows from the fact that $\tr(\Ub \Vb) \leq \|\Ub\|_{\op} \| \Vb \|_{\star}$ for all $\Ub \in \bbC^{N\times N}$ and Hermite $\Vb\in\bbC^{N\times N}$, and the second inequality follows from  the fact that $\Ab$ and   $\overline{\Ab}$ are real matrices. 
Because 
\begin{align*}
    |M_d(\xi)|&=\frac{1}{d} \Big| \tr \big(\Ab-\xi\Ib\big)^{-1} \Big| \leq \frac{\PO}{d} \big\|\Ab-\xi\Ib \big\|_{\op}\leq \PO/(d\cdot\Im(\xi)),\\
    |\overline{M}_d(\xi)|&=\frac{1}{d} \Big| \tr\big(\overline \Ab-\xi\Ib\big)^{-1}\Big| \leq \frac{\PO}{d} \big\|\overline \Ab-\xi\Ib \big\|_{\op}\leq \PO/(d\cdot\Im(\xi)),
\end{align*}
 $|\bDelta(\Ab,\overline{\Ab},\xi)|$ is deterministically upper bounded:
\begin{align}\label{eq:Delta_bound_deterministic}
    |\bDelta(\Ab,\overline{\Ab},\xi)|\leq|M_d(\xi)|+|\overline{M}_d(\xi)|\leq2\PO/(d\cdot\Im(\xi)).
\end{align}
Therefore, if we can prove $\lVert \Ab-\overline{\Ab}\rVert_\star/d=o_{\PP}(1)$, then according to \eqref{eq:Delta_bound_AAbar} and \eqref{eq:Delta_bound_deterministic}, we can conclude that $\EE|\bDelta(\Ab,\overline{\Ab},\xi)|=o_{d}(1)$ by the dominated convergence theorem. To this end,
we first recall the notations in Definitions \ref{def:notation1} and \ref{def:linear pencil} that for ${\sfc} = 1,2$,
\begin{align*}
    \Zb _{\sfc}=\sigma_{\sfc}\left(\Xb\bTheta_{\sfc}^\T /\sqrt d\right)/\sqrt d\in\RR^{n\times N_{\sfc}}, \qquad \tilde{\Zb}_\sfc=\frac{\mu_{\sfc,1}}{d}\Xb\bTheta_\sfc^\T.
\end{align*}
We also remind readers that $\Jb_{\sfc}=\phi_{\sfc}\big(\overline{\Xb}~\overline{\bTheta}_{\sfc}^\T /\sqrt{d}\big)/\sqrt{d}$,  $\tilde{\Jb}_{\sfc}=\frac{\mu_{\sfc,1}}{d}\overline{\Xb}~\overline{\bTheta}_{\sfc}^\T $ are the ``Gaussian version'' counterparts of $\Zb _{\sfc}$ and $\tilde\Zb _{\sfc}$ respectively. 
We further denote $\Zb_{\sfc,0}=\mu_{\sfc,0}\1_n\1_{N_{\sfc}}/\sqrt{d}$ and let 
$\Zb_{\sfc,\star}= \Zb _{\sfc} - \Zb_{\sfc,0}$ for $\sfc = 1,2$. Then by the definition of the functions $\phi_{1},\phi_{2}$, clearly we have $\Zb_{\sfc,\star}=\phi_{\sfc}(\Xb\bTheta^\T/\sqrt{d})/\sqrt{d}$ for $\sfc = 1,2$. 
With these notations, we can rewrite $\Ab-\overline{\Ab}$ as follows:
\begin{align*}
    \Ab-\overline{\Ab}=&q_5\begin{bmatrix}
\0&\0\\
\0&\frac{\Xb\Xb^\T-\overline{\Xb}~\overline{\Xb}^\T}{d}\end{bmatrix}
+q_4\begin{bmatrix}
\frac{\Mb_1\bTheta\bTheta^\T\Mb_1-\Mb_1\overline{\bTheta}~\overline{\bTheta}^\T\Mb_1}{d}&\0\\
\0&\0
\end{bmatrix}\\
&+q_1\begin{bmatrix}
\0&[\tilde{\Zb}_1,\tilde{\Zb}_2]^\T-[\tilde{\Jb}_1,\tilde{\Jb}_2]^\T\\
[\tilde{\Zb}_1,\tilde{\Zb}_2]^\T-[\tilde{\Jb}_1,\tilde{\Jb}_2]&\0
\end{bmatrix}+\begin{bmatrix}
\0&[\Zb_{1,0},\Zb_{2,0}]^\T\\
[\Zb_{1,0},\Zb_{2,0}]&\0
\end{bmatrix}\\
&+\begin{bmatrix}
\0&[\Zb_{1,\star},\Zb_{2,\star}]^\T-[\Jb_1,\Jb_2]^\T\\
[\Zb_{1,\star},\Zb_{2,\star}]-[\Jb_1,\Jb_2]&\0
\end{bmatrix}.
\end{align*}
Then by the triangle inequality and Cauchy-Schwarz inequality, we have
\begin{align*}
    \frac{\lVert\Ab-\overline{\Ab}\rVert_{\star}}{d}=&O_{\PP}( I_1 + I_2 + I_3 + I_4 + I_5 ),
\end{align*}
where
\begin{align*}
&I_1=\frac{1}{\sqrt d}\Bigg\lVert\begin{bmatrix}
\0&\0\\
\0&\frac{\Xb\Xb^\T-\overline{\Xb}~\overline{\Xb}^\T}{d}\end{bmatrix}\Bigg\rVert_{F}, \\
    &I_2= \frac{1}{\sqrt d}\Bigg\lVert\begin{bmatrix}
\frac{\Mb_1\bTheta\bTheta^\T\Mb_1-\Mb_1\overline{\bTheta}~\overline{\bTheta}^\T\Mb_1}{d}&\0\\
\0&\0 \end{bmatrix}\Bigg\rVert_{F}, \\
&I_3 = \frac{1}{\sqrt d}\Bigg\lVert\begin{bmatrix}
\0&[\tilde{\Zb}_1,\tilde{\Zb}_2]^\T-[\tilde{\Jb}_1,\tilde{\Jb}_2]^\T\\
[\tilde{\Zb}_1,\tilde{\Zb}_2]^\T-[\tilde{\Jb}_1,\tilde{\Jb}_2]&\0
\end{bmatrix}\Bigg\rVert_{F}, \\
&I_4 = \frac{1}{d}\Bigg\lVert\begin{bmatrix}
\0&[\Zb_{1,0},\Zb_{2,0}]^\T\\
[\Zb_{1,0},\Zb_{2,0}]&\0
\end{bmatrix}\Bigg\lVert_{\star},\\
&I_5 = \frac{1}{\sqrt d}\Bigg\lVert\begin{bmatrix}
\0&[\Zb_{1,\star},\Zb_{2,\star}]^\T-[\Jb_1,\Jb_2]^\T\\
[\Zb_{1,\star},\Zb_{2,\star}]-[\Jb_1,\Jb_2]&\0
\end{bmatrix}\Bigg\rVert_{F}.
\end{align*}
In the following, we bound the terms $I_1,\ldots, I_5$ separately.  For $I_1$,    let  $\Db_{\xb}=\diag(\sqrt{d}/\lVert\overline{\xb}_{1}\rVert_2,...,\sqrt{d}/\lVert\overline{\xb}_{n}\rVert_2)$. Then we have $\Xb=\Db_\xb\overline{\Xb}$ by \eqref{eq:x_theta_equiv_def}, and
\begin{align}
    I_1 =& \frac{1}{\sqrt{d}}\Big\lVert\frac{\Xb\Xb^\T-\overline{\Xb}~\overline{\Xb}^\T}{d}\Big\rVert_F\leq\Big\lVert\frac{\Xb\Xb^\T-\overline{\Xb}~\overline{\Xb}^\T}{d}\Big\rVert_{\op}=\Big\lVert\frac{\Db_\xb\overline{\Xb}~\overline{\Xb}^\T\Db_\xb-\overline{\Xb}~\overline{\Xb}^\T}{d}\Big\rVert_{\op}\nonumber\\
    =&\Big\lVert\frac{(\Db_\xb-\Ib_n)\overline{\Xb}~\overline{\Xb}^\T(\Db_\xb+\Ib_n)+\overline{\Xb}~\overline{\Xb}^\T\Db_\xb-\Db_\xb\overline{\Xb}~\overline{\Xb}^\T}{d}\Big\rVert_{\op}\nonumber\\
    \leq& \lVert\Db_\xb-\Ib_n\rVert_{\op} \cdot \Big\lVert\frac{\overline{\Xb}~\overline{\Xb}^\T}{d} \Big\rVert_{\op} \cdot \big(1 +\lVert\Db_\xb\rVert_{\op}\big)+\Big\lVert\frac{\overline{\Xb}~\overline{\Xb}^\T\Db_\xb-\Db_\xb\overline{\Xb}~\overline{\Xb}^\T}{d}\Big\rVert_{\op}\nonumber\\
    =& \Big\lVert\frac{\overline{\Xb}~\overline{\Xb}^\T\Db_\xb-\Db_\xb\overline{\Xb}~\overline{\Xb}^\T}{d}\Big\rVert_{\op}+o_{\PP}(1),\label{eq:I1_bound_eq1}
\end{align}
where the first inequality holds since the average of the $d$ squared eigenvalues of $(\Xb\Xb^\T-\overline{\Xb}~\overline{\Xb}^\T) / d$ is bounded by the largest one of them, and the last equality follows from  $\lVert\Db_\xb-\Ib_n\rVert_{\op}=O_{\PP}\Big(\sqrt{\frac{\log d}{d}}\Big)$ and $\lVert\Db_\xb\rVert_{\op}=O_{\PP}(1)$, which are direct consequences of   the definition of $\Db_\xb$. 
We further let $\widetilde{\Db}_\xb$ be the matrix whose elements $\big(\widetilde{\Db}_\xb\big)_{ij}$ satisfy $\big(\widetilde{\Db}_\xb\big)_{ij}=\big({\Db_\xb}\big)_{jj}-\big({\Db_\xb}\big)_{ii}$ for $i,j \in [n]$. Then we have $\lVert\widetilde{\Db}_\xb\rVert_{\max}=o_{\PP}(1)$, and
\begin{align}\label{eq:I1_bound_eq2}
    \Big\lVert\frac{\overline{\Xb}~\overline{\Xb}^\T\Db_\xb-\Db_\xb\overline{\Xb}~\overline{\Xb}^\T}{d}\Big\rVert_{\op}=\Big\lVert\widetilde{\Db}_\xb\odot\frac{\overline{\Xb}~\overline{\Xb}^\T}{d}\Big\rVert_{\op}
    \leq\lVert\widetilde{\Db}_\xb\rVert_{\max}\cdot\Big\lVert\frac{\overline{\Xb}~\overline{\Xb}^\T}{d}\Big\rVert_{\op}=o_{\PP}(1).
\end{align}
Plugging \eqref{eq:I1_bound_eq2} into \eqref{eq:I1_bound_eq1} completes the proof of $I_1=o_{\PP}(1)$. 
Similarly, it can be shown  that $I_2$ and $I_3$ are both $o_{\PP}(1)$. For $I_4$, by the definition that $\Zb_{\sfc,0}=\mu_{\sfc,0}\1_n\1_{N_{\sfc}}/\sqrt{d}$, $\sfc = 1,2$, it is clear that $\Zb_{\sfc,0}$ is rank-one and $\| \Zb_{\sfc,0}\|_{\op} = O_d(\sqrt{d})$. Therefore we have
$$
I_4=\frac{1}{d}\Bigg\lVert\begin{bmatrix}
\0&[\Zb_{1,0},\Zb_{2,0}]^\T\\
[\Zb_{1,0},\Zb_{2,0}]&\0
\end{bmatrix}\Bigg\lVert_{\star}=o_d(1).
$$
Finally, to prove $I_5 = o_d(1)$, it clearly suffices to show that 
\begin{align*}
    \frac{1}{\sqrt{d}}\lVert{\Zb}_{\sfc,\star}-\Jb_{\sfc}\rVert_F=o_{\PP}(1),\quad \sfc=1,2.
\end{align*}
Define $\overline{\Zb}_{\sfc,\star}= \phi_{\sfc}(\Xb~\overline{\bTheta}_{\sfc}^\T/\sqrt{d})/\sqrt{d} =\phi_{\sfc}(\Db_{\xb}\overline{\Xb}~\overline{\bTheta}_{\sfc}^\T/\sqrt{d})/\sqrt{d}$, $j\in[N]$ and $r_i=\sqrt{d}/\lVert\overline{\xb}_{i}\rVert_2$, $i\in[n]$. By the mean value theorem, for $\sfc=1,2$, $a\in [N_\sfc]$ and $i\in[n]$, there exists $\zeta_{ia}$ between $r_i$ and $1$, such that
$$\begin{aligned}
\overline{\Zb}_{\sfc,\star}-\Jb_{\sfc}=&\big[\phi_{\sfc}(r_i\la\overline{\xb}_{i},\overline{\btheta}_{a}\ra/\sqrt{d})/\sqrt{d}-\phi_{\sfc}(\la\overline{\xb}_{i},\overline{\btheta}_{a}\ra/\sqrt{d})/\sqrt{d}\big]_{i\in[n],a\in[N_{\sfc}]}\\
=&\big[(r_i-1)(\la\overline{\xb}_{i},\overline{\btheta}_{a}\ra/\sqrt{d})\phi_{\sfc}'(\zeta_{ij}\la\overline{\xb}_{i},\overline{\btheta}_{a}\ra/\sqrt{d})/\sqrt{d}\big]_{i\in[n],a\in[N_{\sfc}]}\\
=&(\Db_\xb-\Ib_n)\overline{\phi}_{\sfc}(\bzeta\odot(\overline{\Xb}~\overline{\bTheta}_{\sfc}^\T/\sqrt{d}))/\sqrt{d},
\end{aligned}$$
where $\bzeta=(\zeta_{ij})_{i\in[n],a\in[N_{\sfc}]}$ and $\overline{\phi}_{\sfc}(x)=x\phi_{\sfc}'(x)$. By Bernstein-type concentration inequalities \citep{vershynin2010introduction}, we have 
$$\begin{aligned}
&\lVert\Db_\xb-\Ib_n\rVert_{\op}=O_{\PP}\Big(\sqrt{\frac{\log d}{d}}\Big),~\lVert\bzeta\rVert_{\max}=O_{\PP}(1),~ \lVert\overline{\Xb}~\overline{\bTheta}_{\sfc}^\T/\sqrt{d}\rVert_{\max}=O_{\PP}(\sqrt{\log d}).
\end{aligned}$$
Moreover, note that we currently assume that the activation functions $\sigma_\sfc$, $\sfc=1,2$ are fixed polynomials, which implies that $\phi_{\sfc}$ are also fixed polynomials. Therefore, there exists a constant $M_0\in\NN$ such that
\begin{equation*}
    \begin{split}
        ||\overline{\Zb}_{\sfc,\star}-\Jb_{\sfc}||_F/\sqrt{d}\leq&\lVert\Db_\xb-\Ib_n\rVert_{\op}||\overline{\phi}_{\sfc}(\bzeta\odot(\overline{\Xb}~\overline{\bTheta}_{\sfc}^\T/\sqrt{d}))||_F/d
        =O_{\PP}((\log d)^{M_0}/\sqrt{d})=o_{\PP}(1).
    \end{split}
\end{equation*}
With exactly the same argument, it can be shown  that $||\overline{\Zb}_{\sfc,\star}-\Zb_{\sfc,\star}||_F/\sqrt{d}=o_{\PP}(1)$ (recall $\Zb_{\sfc,\star}=\phi_{\sfc}(\Xb\bTheta^\T/\sqrt{d})/\sqrt{d}$). Therefore we have
\begin{align*}
    \frac{1}{\sqrt{d}}\lVert{\Zb}_{\sfc,\star}-\Jb_{\sfc}\rVert_F\leq \frac{1}{\sqrt{d}}\lVert \overline{\Zb}_{\sfc,\star} - {\Zb}_{\sfc,\star}\rVert_F + \frac{1}{\sqrt{d}}\lVert \overline{\Zb}_{\sfc,\star} -\Jb_{\sfc}\rVert_F = 
    o_{\PP}(1)
\end{align*}
for $\sfc=1,2$. Finally  $I_5 = o_p(1)$ and  the proof of Lemma~\ref{lemma:GaussianSphere} is complete.
\end{proof}

\subsection{Calculation of the resolvent equations}
\label{subsec:resolvent}
Lemma~\ref{lemma:GaussianSphere} and its proof show the readers that the Stieltjes transforms of the empirical eigenvalue distributions of $\Ab$ and $\overline{\Ab}$ share the same asymptotics. Based on this result, we can equivalently consider the ``Gaussian version'' counterpart of the learning problem. 
Therefore, throughout Section~\ref{subsec:resolvent}, we directly consider the matrices $\overline\Xb$ and $\overline{\bTheta}$, whose elements are independently generated from standard normal $\rmN(0,1)$. In addition, the activation functions for the two types of random features are $\phi_{\sfc}(x)=\sigma_\sfc(x)-\mu_{\sfc,0}$, $\sfc = 1,2$, and the linear pencil matrix  is $\overline{\Ab}(\qb,\bmu)$ is given in Definition~\ref{def:linear pencil Gaussian}.
\if UT
$$
\overline{\Ab}(\qb,\bmu)=\begin{bmatrix}
q_2\mu_{1,2}^2\Ib_{N_1}+q_4\mu_{1,1}^2\frac{\overline{\bTheta}_1\overline{\bTheta}_1^\T }{d}&q_4\mu_{1,1}\mu_{2,1}\frac{\overline{\bTheta}_1\overline{\bTheta}_2^\T }{d}& \Jb_1^\T +q_1\tilde{\Jb}_1^\T \\
q_4\mu_{1,1}\mu_{2,1}\frac{\overline{\bTheta}_2\overline{\bTheta}_1^\T }{d}&q_2\mu_{2,2}^2\Ib_{N_2}+q_4\mu_{2,1}^2\frac{\overline{\bTheta}_2\overline{\bTheta}_2^\T }{d}& \Jb_2^\T +q_1\tilde{\Jb}_2^\T \\
 \Jb_1^\T +q_1\tilde{\Jb}_1^\T& \Jb_2^\T +q_1\tilde{\Jb}_2^\T&q_3\Ib_n+q_5\frac{\overline{\Xb}~\overline{\Xb}^\T }{d}
\end{bmatrix}.
$$
\fi

Moreover, for $\sfc=1,2$, let 
$\Phi_{\sfc}(x)=\phi_{\sfc}(x)+q_1\mu_{\sfc,1}x$, and it is easy to see $ \Jb_{\sfc}^\T +q_1\tilde{\Jb}_{\sfc}^\T=\Phi_{\sfc}\left(\overline{\Xb}~\overline{\bTheta}_{\sfc}^\T /\sqrt d\right)/\sqrt d$. We further denote $\phi_{\sfc,0}\triangleq\EE_{G\sim\rmN(0,1)}\{\Phi_{\sfc}(G)\}$, $\phi_{\sfc,1}\triangleq\EE_{G\sim\rmN(0,1)}\{G\Phi_{\sfc}(G)\}$, $
\phi_{\sfc,2}\triangleq\EE_{G\sim\rmN(0,1)}$
$\{\Phi_{\sfc}(G)^2\}-\phi_{\sfc,0}^2-\phi_{\sfc,1}^2$. By these definitions, it is easy to see that $\phi_{\sfc,0}=0$, $\phi_{\sfc,1}^2=\mu_{\sfc,1}^2(1+q_1)^2$, $\phi_{\sfc,2}^2=\mu_{\sfc,2}^2$. Importantly, the property that $\EE_{G\sim\rmN(0,1)}\{\Phi_{\sfc}(G)\} = \phi_{\sfc,0}=0$ enables the application of the following lemma, which is summarized from Section 4.3, Step 2 in \citet{cheng2013spectrum}.

\begin{lemma}\label{lemma:centralized_activation_odp}
    Suppose that $\Phi$ is a polynomial satisfying $\EE_{G\sim\rmN(0,1)}\{\Phi(G)\} = 0$, $\EE_{G\sim\rmN(0,1)}\{G\Phi(G)\} = 0$ and $\overline{\xb}_i, \overline{\btheta}_{a}\in\RR^d$, $i\in [n]$, $a\in [N]$ are standard Gaussian vectors. 
    Define matrix $\Eb \in \RR^{n\times N}$ elementwisely as
    \begin{align*}
        \big(\Eb_\sfc\big)_{i,a}=\frac{1}{\sqrt{d}}\bigg[\Phi\bigg(\frac{1}{\sqrt{d}}\langle \overline{\xb} _i, \overline{\btheta}_{a}\rangle \bigg)-\Phi\bigg(\frac{1}{\sqrt{d}}\langle(\overline{\xb}_i)_{[1:d-1]},(\overline{\btheta}_{a})_{[1:d-1]}\rangle\bigg)\bigg]
    \end{align*}
    for $i\in [n]$, $a\in [N]$.
    Then $ \| \Eb \|_{\op} = o_{\PP}(1) $. 
\end{lemma}
Lemma~\ref{lemma:centralized_activation_odp} formally shows the intuitive result that under the setting where $d,n$ grows proportionally, removing one entry in the random vectors does not change the asymptotic limit of polynomials. This enables us to apply the standard leave-one-out argument in random matrix theory. 

Our goal in this part of the proof is to calculate the resolvent equations of the Stieltjes transforms corresponding to the pencil  matrix  $\overline{\Ab}(\qb,\bmu)$. To do so, we define the following terms:
\begin{equation*}
    \begin{split}
        \overline{m}_{1,d}(\xi;\qb,\bmu)&=\EE\big[ \overline{M}_{1,d}(\xi;\qb,\bmu)\big],\quad \overline{M}_{1,d}(\xi;\qb,\bmu)=\frac{1}{d}\text{tr}_{[1:N_1]}\big[(\overline{\Ab}(\qb,\bmu)-\xi \Ib_{\PO})^{-1}\big],\\
        \overline{m}_{2,d}(\xi;\qb,\bmu)&=\EE\big[\overline{M}_{2,d}(\xi;\qb,\bmu)\big],\quad \overline{M}_{2,d}(\xi;\qb,\bmu)=\frac{1}{d}\text{tr}_{[N_1+1:N]}\big[ (\overline{\Ab}(\qb,\bmu)-\xi \Ib_{\PO})^{-1}\big],\\
        \overline{m}_{3,d}(\xi;\qb,\bmu)&=\EE\big[\overline{M}_{3,d}(\xi;\qb,\bmu)\big],\quad \overline{M}_{3,d}(\xi;\qb,\bmu)=\frac{1}{d}\text{tr}_{[N+1:\PO]}\big[ (\overline{\Ab}(\qb,\bmu)-\xi \Ib_{\PO})^{-1}\big].
    \end{split}
\end{equation*}
Standard argument in random matrix theory then gives us the concentration result
\begin{align*}
   \EE\big|\overline{M}_{i,d}(\xi;\qb,\bmu)- \overline{m}_{i,d}(\xi;\qb,\bmu)\big|=o_d(1)
\end{align*}
for any fixed $\xi\in\bbC_+$. Therefore, denoting $\overline{m}_d(\xi)=\sum\limits_{i=1}^3\overline{m}_{i,d}(\xi)$ and $ \overline{M}_d(\xi)=\sum\limits_{i=1}^3\overline{M}_{i,d}(\xi)$, (we drop the argument $\qb,\bmu$ for simplicity),
we have
\begin{align}\label{eq:concentration}
    \EE\big|\overline{M}_{d}(\xi)- \overline{m}_{d}(\xi)\big|=o_d(1)
\end{align}
for any fixed $\xi\in\bbC_+$.
A proof of this concentration can be found in \citet{hastie2022surprises,mei2022generalization}. 
Based on \eqref{eq:concentration}, to study $\overline{M}_{d}(\xi;\qb,\bmu)$, which is the Stieltjes transform of the empirical eigenvalue distribution of $\overline{A}(\qb,\bmu)$, it suffices to derive  the resolvent equations for $\overline{m}_{d}(\xi;\qb,\bmu)$. This is done in the following lemma.
\begin{lemma}\label{lemma:resolvent_calculation}
Let $\overline{\mb}_d(\xi)=[\overline{m}_{1,d}(\xi),\overline{m}_{2,d}(\xi),\overline{m}_{3,d}(\xi)]^\T$.  Then for any fixed $\xi\in\bbC_+$, the following  property holds:
\begin{equation*}
    \begin{split}
\| \overline{\mb}_d(\xi) - \sfFb(\overline{\mb}_d(\xi)) \|_2 = o_d(1).
    \end{split}
\end{equation*}
\end{lemma}

\begin{proof}[Proof of Lemma~\ref{lemma:resolvent_calculation}] 
Since $\overline{\mb}_d(\xi),\sfFb(\overline{\mb}_d(\xi)) \in \bbC^3 $,  Lemma~\ref{lemma:resolvent_calculation} essentially contains three results showing that the first, second, and third elements of $\overline{\mb}_d(\xi) - \sfFb(\overline{\mb}_d(\xi))$ are all asymptotically zero. 
Since the proofs of the three results are almost the same, we mainly focus on the proof of the first result. 
The proof consists of three  steps. The first step is to use the Schur complement formula to calculate $\overline{m}_{1,d}$. The second step is to simplify the formula of $\overline{m}_{1,d}$. The third step is to give the recursive equations of $\overline{m}_{1,d}$ based on the result of step 2.

\smallskip\noindent{\textbf{Step 1}.} 
We first use a leave-one-out argument to calculate $\overline{m}_{1,d}$.
Let $\overline{\Ab}_{\cdot,N_1} \in \RR^{\PO-1}$ be the $N_1^{\text{th}}$ column of $\overline{\Ab}$, with the $N_1^{\text{th}}$ entry removed. We further denote by $\overline{\Bb}\in\RR^{( \PO-1)\times( \PO-1)}$ the sub-matrix of $\overline{\Ab}$ obtained by removing the $N_1^{\text{th}}$ row and $N_1^{\text{th}}$ column in $\overline{\Ab}$. We can then treat $\overline{\Ab}$ as a $2\times 2$ block matrix formed by $\overline{\Ab}_{\cdot,N_1}$, $\overline{\Ab}_{\cdot,N_1}^\T$, $\overline{\Bb}$, and $\overline{\Ab}_{N_1,N_1} = q_2\mu_{1,2}^2+q_4\mu_{1,1}^2\lVert\overline{\btheta}_{N_1}\rVert_2^2/d$. Then by the Schur complement formula, we get
\begin{align}
    \overline{m}_{1,d}=\psi_1\EE\left(-\xi+q_2\mu_{1,2}^2+q_4\mu_{1,1}^2\lVert\overline{\btheta}_{N_1}\rVert_2^2/d-\overline{\Ab}_{\cdot,N_1}^\T (\overline{\Bb}-\xi\Ib_{ \PO-1})^{-1}\overline{\Ab}_{\cdot,N_1}\right)^{-1}.\label{eq:schurcomplement}
\end{align}
We decompose the vectors $\overline{\btheta}_{a}$, $a\in [N]$ and $\overline{\xb}_i$, $i\in [n]$ into components along the direction of $\overline{\btheta}_{N_1}$ and other orthogonal directions:
\begin{equation}
\label{eq:schurdecomp}
    \begin{split}
        &\overline{\btheta}_{a}=\eta_{a}\frac{\overline{\btheta}_{N_1}}{\lVert\overline{\btheta}_{N_1}\rVert}+\tilde{\btheta}_{a},~\langle\overline{\btheta}_{N_1},\tilde{\btheta}_{a}\rangle=0,~a\in[N]\backslash\{N_1 \},\\
&\overline{\xb}_i=u_i\frac{\overline{\btheta}_{N_1}}{\lVert\overline{\btheta}_{N_1}\rVert}+\tilde{\xb}_i,~\langle\overline{\btheta}_{N_1},\tilde{\xb}_i\rangle=0,~i\in[n].
    \end{split}
\end{equation}
Note that for any $a\in [N]\backslash\{N_1 \}$ and $i\in [n]$, $\eta_{a}$, $u_i$ are standard Gaussian and are independent of $\tilde{\btheta}_{a}$ and $\tilde{\xb}_i$. Moreover,
$\tilde{\btheta}_{a}$ and $\tilde{\xb}_i$ are conditionally independent on each other given $\overline{\btheta}_{N_1}$, with $\tilde{\btheta}_{a},\tilde{\xb}_i\sim N(0,P_{\bot})$, where $P_{\bot}$ is the projector orthogonal to $\overline{\btheta}_{N_1}$. We can then use the coefficients $\eta_{a}$, $a\in[N]\backslash\{N_1 \}$ and $u_i$, $i\in[n]$ to represent the entries of $\overline{\Ab}_{\cdot,N_1}$. We have  $\overline{\Ab}_{\cdot,N_1}=[\overline{\Ab}_{1,N_1},...,\overline{\Ab}_{ \PO-1,N_1}]^\T \in \RR^{ \PO-1}$ with
\begin{align}\label{eq:A_decom_final}
     \overline{\Ab}_{i,N_1}=\left\{\begin{aligned}
&\frac{q_4\mu_{1,1}^2\eta_i}{d}\lVert\overline{\btheta}_{N_1}\rVert_2, &&\text{if  } i\in[1,N_1-1],\\
&\frac{q_4\mu_{1,1}\mu_{2,1}\eta_{i+1}}{d}\lVert\overline{\btheta}_{N_1}\rVert_2, &&\text{if } i\in[N_1,N-1],\\
&\frac{1}{\sqrt{d}}\Phi_{1}\Big(\frac{1}{\sqrt{d}}u_{i-N+1}\lVert\overline{\btheta}_{N_1}\rVert_2\Big), &&\text{if } i\geq N.
\end{aligned}\right.
\end{align}
To calculate the resolvent equations, we need to further represent the matrix $\overline{\Bb}$ in \eqref{eq:schurcomplement} with $\eta_a$, $\tilde{\btheta}_a$, $u_i$, and $\tilde{\xb}_i$ for $a\in [N]\backslash\{N_1 \}$ and $i\in [n]$. 
Below we first list some additional notations for easier reference. We write 
$\bmeta_1=[\eta_1,...,\eta_{N_1-1}]\in\RR^{N_1-1}$, $\bmeta_2=[\eta_{N_1+1},...,\eta_{N}] \in\RR^{N_2}$,
$\bmeta=[\bmeta_1^\T ,\bmeta_2^\T ]^\T \in\RR^{N-1}$,
$\ub=[u_1,...,u_n]^\T\in\RR^n $,
$\tilde{\bTheta}_1=[\tilde{\btheta}_1,...,\tilde{\btheta}_{N_1-1}]^\T $,
$\tilde{\bTheta}_2=[\tilde{\btheta}_{N_1+1},...,\tilde{\btheta}_{N}]^\T $, $\tilde{\bTheta}=\begin{bmatrix}
\tilde{\bTheta}_1\\\tilde{\bTheta}_2
\end{bmatrix}\in\RR^{(N-1)\times d}$, $\tilde{\Mb}_1=\begin{bmatrix}
\mu_{1,1}\Ib_{N_1-1}&~\\
~&\mu_{2,1}\Ib_{N_2}
\end{bmatrix}$ and $\tilde{\Mb}_*=\begin{bmatrix}
\mu_{1,2}\Ib_{N_1-1}&~\\
~&\mu_{2,2}\Ib_{N_2}
\end{bmatrix}$. 
Now with \eqref{eq:schurdecomp} and the notations above, we can decompose $\overline{\Bb}_{[1:N-1],[1:N-1]}$ as follows:
\begin{align}
    \overline{\Bb}_{[1:N-1],[1:N-1]}=q_2\tilde{\Mb}_*\tilde{\Mb}_*+\frac{q_4}{d}\tilde{\Mb}_1\tilde{\bTheta}\tilde{\bTheta}^\T \tilde{\Mb}_1+\frac{q_4}{d}\tilde{\Mb}_1\bmeta\bmeta^\T \tilde{\Mb}_1\label{eq:B_decomp11}.
\end{align}
Moreover, for $i,j\in[n]$ we define
\begin{align*}
    \quad\big(\tilde{\Hb}\big)_{ij}=\frac{1}{d}\langle\tilde{\xb}_i,\tilde{\xb}_j\rangle.
\end{align*}
Then we can decompose $\overline{\Bb}_{[N:P-1],[N:P-1]}$ into
\begin{align}
    \overline{\Bb}_{[N:P-1],[N:P-1]}=q_3\Ib_n+q_5\tilde{\Hb}+\frac{q_5}{d}\ub\ub^\T.\label{eq:B_decomp22}
\end{align}
For $\overline{\Bb}_{[N:P-1],[1:N-1]}$, by definition we see that the elements in $\overline{\Bb}_{[N:P-1],[1:N-1]}$ are $(\Zb)_{i,a}$ for $a\in [N]\backslash\{N_1 \}$ and $i\in [n]$. Therefore, we have
\begin{align*}
    (\Zb)_{i,a}&=\frac{1}{\sqrt{d}}\Phi_{\sfc}\Big(\frac{1}{\sqrt{d}}\langle{\overline{\xb}}_i,{\overline{\btheta}}_{a}\rangle\Big)=\frac{1}{\sqrt{d}}\Phi_{\sfc}\Big(\frac{1}{\sqrt{d}}\langle\tilde{\xb}_i,\tilde{\btheta}_{a}\rangle+\frac{1}{d}u_i\eta_a\Big)\\
    &=\frac{1}{\sqrt{d}}\Phi_{\sfc}\Big(\frac{1}{\sqrt{d}}\langle\tilde{\xb}_i,\tilde{\btheta}_{a}\rangle\Big)+\frac{\phi_{\sfc,1}}{d}u_i\eta_a+\frac{1}{\sqrt{d}}\Big[\Phi_{\sfc,\bot}\Big(\frac{1}{\sqrt{d}}\langle\tilde{\xb}_i,\tilde{\btheta}_{a}\rangle+\frac{1}{\sqrt{d}}u_i\eta_{a}\Big)-\Phi_{\sfc,\bot}\Big(\frac{1}{\sqrt{d}}\langle\tilde{\xb}_i,\tilde{\btheta}_{a}\rangle\Big)\Big],
\end{align*}
where $\Phi_{\sfc,\bot}(x)=\Phi_{\sfc}(x)-\phi_{\sfc,1}x$, $\sfc=1$ when $a\leq N_1-1$ and $\sfc=2$ when $a\geq N_1+1$. 
By the symmetry of $\overline{\Bb}$, we can then decompose $\overline{\Bb}_{[N:P-1],[1:N-1]}$ and $\overline{\Bb}_{[1:N-1],[N:P-1]}^\T$ into
\begin{align}
    \overline{\Bb}_{[N:P-1],[1:N-1]}= \overline{\Bb}_{[1:N-1],[N:P-1]}^\T =\tilde{\Zb}+\frac1d\ub\bmeta \Mb_{\phi}+[\Eb_1,\Eb_2],\label{eq:B_decomp12}
\end{align}
where we define
\begin{align*}
    &\tilde{\Zb}=[\tilde{\Zb}_1,\tilde{\Zb}_2],\quad(\tilde{\Zb}_1)_{i,a}=\frac{1}{\sqrt{d}}\Phi_{1}\Big(\frac{1}{\sqrt{d}}\langle\tilde{\xb}_i,\tilde{\btheta}_{a}\rangle\Big), \quad(\tilde{\Zb}_2)_{i,a}=\frac{1}{\sqrt{d}}\Phi_{2}\Big(\frac{1}{\sqrt{d}}\langle\tilde{\xb}_i,\tilde{\btheta}_{a}\rangle\Big),\\
    &\Mb_\phi=\begin{bmatrix}
\phi_{1,1}\Ib_{N_1-1}&~\\
~&\phi_{2,1}\Ib_{N_2}
\end{bmatrix},\quad \big(\Eb_\sfc\big)_{i,a}=\frac{1}{\sqrt{d}}\Big[\Phi_{\sfc,\bot}\Big(\frac{1}{\sqrt{d}}\langle\tilde{\xb}_i,\tilde{\btheta}_{a}\rangle+\frac{1}{\sqrt{d}}u_i\eta_{a}\Big)-\Phi_{\sfc,\bot}\Big(\frac{1}{\sqrt{d}}\langle\tilde{\xb}_i,\tilde{\btheta}_{a}\rangle\Big)\Big]
\end{align*}
for $a\in [N]\backslash\{N_1 \}$, $i\in [n]$ and $\sfc \in \{1,2\}$.
Combining \eqref{eq:B_decomp11}, \eqref{eq:B_decomp22} and \eqref{eq:B_decomp12}, we decompose $\overline{\Bb}$ into
\begin{align}\label{eq:B_decom_final}
    \overline{\Bb}=\tilde{\Bb}+\bDelta+\Eb\in\RR^{(\PO-1)\times(\PO-1)},
\end{align}
where
\begin{align*}
\tilde{\Bb}=&\begin{bmatrix}
q_2\tilde{\Mb}_*\tilde{\Mb}_*+\frac{q_4}{d}\tilde{\Mb}_1\tilde{\bTheta}\tilde{\bTheta}^\T \tilde{\Mb}_1&\tilde{\Zb}^\T \\
\tilde{\Zb}&q_3\Ib_n+q_5\tilde{\Hb}
\end{bmatrix}\\
=&\begin{bmatrix}
q_2\mu_{1,2}^2\Ib_{N_1-1}+\frac{q_4\mu_{1,1}^2}{d}\tilde{\bTheta}_1\tilde{\bTheta}_1^\T &\frac{q_4\mu_{1,1}\mu_{2,1}}{d}\tilde{\bTheta}_1\tilde{\bTheta}_2^\T &\tilde{\Zb}_1^\T \\
\frac{q_4\mu_{1,1}\mu_{2,1}}{d}\tilde{\bTheta}_2\tilde{\bTheta}_1^\T &q_2\mu_{2,2}^2\Ib_{N_2}+\frac{q_4\mu_{2,1}^2}{d}\tilde{\bTheta}_2\tilde{\bTheta}_2^\T &\tilde{\Zb}_2^\T \\
\tilde{\Zb}_1&\tilde{\Zb}_2&q_3\Ib_n+q_5\tilde{\Hb}
\end{bmatrix},\\
\bDelta=&\begin{bmatrix}
\frac{q_4}{d}\tilde{\Mb}_1\bmeta\bmeta^\T \tilde{\Mb}_1&\frac1d\Mb_{\phi}\bmeta\ub^\T \\
\frac1d\ub\bmeta \Mb_{\phi}&\frac{q_5}{d}\ub\ub^\T 
\end{bmatrix}=\begin{bmatrix}
\frac{q_4\mu_{1,1}^2}{d}\bmeta_1\bmeta_1^\T &\frac{q_4\mu_{1,1}\mu_{2,1}}{d}\bmeta_1\bmeta_2^\T &\frac{\phi_{1,1}}{d}\bmeta_1\ub^\T \\\frac{q_4\mu_{1,1}\mu_{2,1}}{d}\bmeta_2\bmeta_1^\T &\frac{q_4\mu_{2,1}^2}{d}\bmeta_2\bmeta_2^\T &\frac{\phi_{2,1}}{d}\bmeta_2\ub^\T \\\frac{\phi_{1,1}}{d}\ub\bmeta_1^\T &\frac{\phi_{2,1}}{d}\ub\bmeta_2^\T &\frac{q_5}{d}\ub\ub^\T 
\end{bmatrix},~
\Eb=\begin{bmatrix}
\0&\0 & \Eb_1^\T \\
\0&\0 & \Eb_2^\T \\
\Eb_1&\Eb_2& \0
\end{bmatrix}.
\end{align*}
Clearly, by the definition of $\tilde{\Bb}$, 
the Stieltjes transform corresponding to $\tilde{\Bb}$ shares the same asymptotics as the Stieltjes transform corresponding to $\overline{\Ab}$.

\smallskip\noindent{\textbf{Step 2}.} 
According to our analysis in \textbf{Step 1}, we can then calculate $\overline{m}_{1,d}$ by \eqref{eq:schurcomplement}, in which the terms $\overline{\Ab}_{\cdot,N_1}$ and $\overline{\Bb}$ have the decompositions \eqref{eq:A_decom_final} and \eqref{eq:B_decom_final} respectively. In this step, we aim to further simplify the calculation by getting rid of the terms $\lVert\overline{\btheta}_{N_1}\rVert_2^2/d$ in \eqref{eq:schurcomplement} and $\Eb$ in \eqref{eq:B_decom_final}. 
Define
$$\begin{aligned}
w_0=&\left(-\xi+q_2\mu_{1,2}^2+q_4\mu_{1,1}^2-\overline{\Ab}_{\cdot,N_1}^\T (\overline{\Bb}-\xi\Ib_{ \PO-1})^{-1}\overline{\Ab}_{\cdot,N_1}\right)^{-1},\\
w_1=&\left(-\xi+q_2\mu_{1,2}^2+q_4\mu_{1,1}^2\lVert\overline{\btheta}_{N_1}\rVert_2^2/d-\overline{\Ab}_{\cdot,N_1}^\T (\overline{\Bb}-\xi\Ib_{ \PO-1})^{-1}\overline{\Ab}_{\cdot,N_1}\right)^{-1},\\
w_2=&\left(-\xi+q_2\mu_{1,2}^2+q_4\mu_{1,1}^2-\overline{\Ab}_{\cdot,N_1}^\T (\tilde{\Bb}+\bDelta-\xi\Ib_{ \PO-1})^{-1}\overline{\Ab}_{\cdot,N_1}\right)^{-1}.
\end{aligned}$$
Then by \eqref{eq:schurcomplement}, we have $\overline{m}_{1,d}=\psi_1\EE w_1$. We now give an upper bound of $|w_1-w_2|$. 
Recall that we consider a fixed $\xi\in\bbC_+$. Since $\overline{\Bb}$ is a real symmetric matrix, by diagonalizing $\overline{\Bb} $, it is easy to see that $\Im\big(\overline{\Ab}_{\cdot,N_1}^\T (\overline{\Bb}-\xi\Ib_{ \PO-1})^{-1}\overline{\Ab}_{\cdot,N_1}\big) \geq 0$. 
Therefore, we deterministically have 
\begin{align*}
    \Im(-w_1^{-1})= \Im(\xi)+\Im\big(\overline{\Ab}_{\cdot,N_1}^\T (\overline{\Bb}-\xi\Ib_{ \PO-1})^{-1}\overline{\Ab}_{\cdot,N_1}\big)\geq \Im(\xi).
\end{align*}
Thus we have $|w_1|\leq 1/\Im(\xi)$. Using a similar argument, we have $\max\{|w_0|,|w_1|,|w_2|\}\leq 1/\Im(\xi)$, which indicates that $|w_1-w_2|\leq 2/\Im(\xi)$. Moreover, we have 
\begin{align*}
|w_1-w_2|&\leq|w_1-w_0|+|w_0-w_2|\\
&\leq q_4\mu_{1,1}^2\big|w_1\big(\overline{\btheta}_{N_1}\rVert_2^2/d-1\big)w_0\big|+\big|w_1w_2\overline{\Ab}_{\cdot,N_1}^\T\big((\overline{\Bb}-\xi\Ib_{ \PO-1})^{-1}-(\tilde{\Bb}+\bDelta-\xi\Ib_{ \PO-1})^{-1}\big)\overline{\Ab}_{\cdot,N_1}\big|\\
&\leq q_4\mu_{1,1}^2|\lVert\overline{\btheta}_{N_1}\rVert_2^2/d-1|/\Im^2(\xi)+2\lVert\overline{\Ab}_{\cdot,N_1}\rVert_2^2\lVert \Eb\rVert_{\op}/\Im^4(\xi)
\end{align*}
By Lemma~\ref{lemma:centralized_activation_odp}, we have $\lVert\Eb_1\rVert_{\op}=o_{\PP}(1)$, $\lVert\Eb_2\rVert_{\op}=o_{\PP}(1)$. It is also easy to see that $\lVert\overline{\Ab}_{\cdot,N_1}\rVert_2^2=O_{\PP}(1)$ and $\| \overline{\btheta}_{N_1}\rVert_2^2/d-1=o_{\PP}(1)$. Therefore we have
$$
|w_1-w_2|=o_{\PP}(1).
$$
Combining with the fact that $|w_1-w_2|$ is deterministically bounded by $2/\Im(\xi)$, by the dominated convergence theorem, we have
$$
\EE|w_1- w_2|=o_d(1).
$$
Therefore $\overline{m}_{1,d}=\psi_1\EE w_2+o_{d}(1)$, and the derivation of the resolvent equations reduces to the calculation of $\EE w_2$.

\smallskip\noindent{\textbf{Step 3}. We calculate $\EE w_2$ to get the resolvent equations}.
For simplicity,  we give some notations which will be used later. Let  
$$
\vb=\overline{\Ab}_{\cdot,N_1},\quad\vb_i=\overline{\Ab}_{i,N_1}=\left\{\begin{aligned}
&\frac{q_4\mu_{1,1}^2\eta_i}{d}\lVert\overline{\btheta}_{N_1}\rVert_2, &&\text{if  } i\in[1,N_1-1],\\
&\frac{q_4\mu_{1,1}\mu_{2,1}\eta_{i+1}}{d}\lVert\overline{\btheta}_{N_1}\rVert_2, &&\text{if } i\in[N_1,N-1],\\
&\frac{1}{\sqrt{d}}\Phi_{1}\Big(\frac{1}{\sqrt{d}}u_{i-N+1}\lVert\overline{\btheta}_{N_1}\rVert_2\Big), &&\text{if } i\geq N,
\end{aligned}\right.
$$
and 
$$
\Ub=\frac{1}{\sqrt{d}}\begin{bmatrix}
\bmeta_1& & \\
&\bmeta_2&\\
& & \ub
\end{bmatrix}\in\RR^{(P-1)\times3},\quad \Mb=\begin{bmatrix}
q_4\mu_{1,1}^2&q_4\mu_{1,1}\mu_{2,1}&\phi_{1,1}\\
q_4\mu_{1,1}\mu_{2,1}&q_4\mu_{2,1}^2&\phi_{2,1}\\
\phi_{1,1}&\phi_{2,1}&q_5
\end{bmatrix}.
$$
By direct verification, we have 
$$
\bDelta=\Ub \Mb\Ub^\T .
$$
We now decompose $w_2$ into the terms related with $\tilde{\Bb}$, $\vb$ and $\Ub$.
By Schur complement formula, we have 
\begin{align}
    (\tilde{\Bb}+\Ub \Mb\Ub^\T -&\xi\Ib_{ \PO-1})^{-1}=(\tilde{\Bb}-\xi\Ib_{ \PO-1})^{-1} \nonumber\\
    & -(\tilde{\Bb}-\xi\Ib_{ \PO-1})^{-1}\Ub [\Mb^{-1}+\Ub^\T (\tilde{\Bb}-\xi\Ib_{ \PO-1})^{-1}\Ub ]^{-1}\Ub^\T (\tilde{\Bb}-\xi\Ib_{ \PO-1})^{-1}.\label{eq:Schur_for_Btilde}
\end{align}
Then $w_2$ can be rewritten as 
\begin{align}
        w_2=&\Big(-\xi+q_2\mu_{1,2}^2+q_4\mu_{1,1}^2-\vb^\T (\tilde{\Bb}+\Ub \Mb\Ub^\T -\xi\Ib_{ \PO-1})^{-1}\vb\Big)^{-1}\nonumber\\
=&\big[-\xi+q_2\mu_{1,2}^2+q_4\mu_{1,1}^2-\vb^\T (\tilde{\Bb}-\xi\Ib_{ \PO-1})^{-1}\vb\nonumber\\
&+\vb^\T (\tilde{\Bb}-\xi\Ib_{ \PO-1})^{-1}\Ub(\Mb^{-1}+\Ub^\T (\tilde{\Bb}-\xi\Ib_{ \PO-1})^{-1}\Ub)^{-1}\Ub^\T (\tilde{\Bb}-\xi\Ib_{ \PO-1})^{-1}\vb\big]^{-1},\label{eq:w2_woodberry}
\end{align}
where the first equation is the definition of $w_2$, and the second equation follows by \eqref{eq:Schur_for_Btilde}. 
To continue the calculation, we study the terms $\vb^\T (\tilde{\Bb}-\xi\Ib_{ \PO-1})^{-1}\vb$, $\vb^\T (\tilde{\Bb}-\xi\Ib_{ \PO-1})^{-1}\Ub$ and $\Ub^\T (\tilde{\Bb}-\xi\Ib_{\PO-1})^{-1}\Ub$ in the denominator of \eqref{eq:w2_woodberry}. To do so, we note that $\tilde{\Bb}$ is independent of $\vb$ and $\Ub$. Moreover, by the leave-one-out argument, the Stieltjes transform corresponding to $\tilde{\Bb}$ shares the same asymptotics as the Stieltjes transform corresponding to $\overline{\Ab}$.
Notice that $\eta_i$ is independent on $\tilde{\Bb}$ conditioned on $\overline{\btheta}_{N_1}$, and $\tilde{\Bb}$ is independent on $\overline{\btheta}_{N_1}$. We have
\begin{align*}
 &~ \EE\vb^\T (\tilde{\Bb}-\xi\Ib_{ \PO-1})^{-1}\vb=\EE\tr(\tilde{\Bb}-\xi\Ib_{ \PO-1})^{-1}\vb\vb^\T=\tr\Big(\EE(\tilde{\Bb}-\xi\Ib_{ \PO-1})^{-1}\EE\vb\vb^\T\Big)\\
 =&\tr\Bigg(\begin{bmatrix}
\frac{d\overline{m}_{1,d}}{N_1}\Ib_{N_1-1}&* &*\\
*&\frac{d\overline{m}_{2,d}}{N_2}\Ib_{N_2}&*\\
*&* &\frac{d\overline{m}_{3,d}}{n}\Ib_n
\end{bmatrix}\\&\cdot\frac{1}{d}\begin{bmatrix}
(q_4^2\mu_{1,1}^4+o_d(1))\Ib_{N_1-1}& &\\
&(q_4^2\mu_{1,1}^2\mu_{2,1}^2+o_d(1))\Ib_{N_2}&\\
& &(\phi_{1,1}^2+\phi_{1,2}^2+o_d(1))\Ib_n
\end{bmatrix}\Bigg)\\
=&
 q_4^2\mu_{1,1}^2(\mu_{1,1}^2\overline{m}_{1,d}+\mu_{2,1}^2\overline{m}_{2,d})+(\phi_{1,1}^2+\phi_{1,2}^2)\overline{m}_{3,d}+o_d(1),
\end{align*}
where the second equality follows from  the fact that 
$\EE\Phi_1^2\Big(\frac{1}{\sqrt{d}}u_{i-N+1}\lVert\overline{\btheta}_{N_1}\rVert_2\Big)=\phi_{1,1}^2+\phi_{1,2}^2+o_d(1)$, and we have denoted by `*' the blocks that are irrelevant to the calculation. By a concentration measure argument (see in \citet{tao2012topics} Section 2.4.3),  we have
\begin{align}
    \vb^\T (\tilde{\Bb}-\xi\Ib_{ \PO-1})^{-1}\vb=q_4^2\mu_{1,1}^2(\mu_{1,1}^2\overline{m}_{1,d}+\mu_{2,1}^2\overline{m}_{2,d})+(\phi_{1,1}^2+\phi_{1,2}^2)\overline{m}_{3,d}+o_{\PP}(1).\label{eq:term1}
\end{align}
After direct calculation with the same argument, we obtain that 
\begin{align}
  \vb^\T (\tilde{\Bb}-\xi\Ib_{ \PO-1})^{-1}\Ub=&\begin{bmatrix}
q_4\mu_{1,1}^2\overline{m}_{1,d}\\q_4\mu_{1,1}\mu_{2,1}\overline{m}_{2,d}\\ \phi_{1,1}\overline{m}_{3,d}
\end{bmatrix}^\T+o_{\PP}(1),\label{eq:term2}\\
\Ub^\T (\tilde{\Bb}-\xi\Ib_{ \PO-1})^{-1}\Ub=&\begin{bmatrix}
\overline{m}_{1,d}& &\\
&\overline{m}_{2,d}&\\
& &\overline{m}_{3,d}
\end{bmatrix}+o_{\PP}(1).\label{eq:term3}
\end{align}
Now since
$|w_2|\leq1/\xi_0$ is deterministically bounded, by dominated convergence theorem, we have the $L_1$ convergence of $w_2$ by plugging the main terms of \eqref{eq:term1}, \eqref{eq:term2} and \eqref{eq:term3} into \eqref{eq:w2_woodberry}. Further note that equation~\eqref{eq:w2_woodberry} has a part in the form of $\left(\Ab^{-1}+\Mb^{-1}\right)^{-1}$, where $$\Ab=\begin{bmatrix}
1/\overline{m}_{1,d}& &\\
&1/\overline{m}_{2,d}&\\
& &1/\overline{m}_{3,d}
\end{bmatrix},\quad\Mb=\begin{bmatrix}
q_4\mu_{1,1}^2&q_4\mu_{1,1}\mu_{2,1}&\phi_{1,1}\\
q_4\mu_{1,1}\mu_{2,1}&q_4\mu_{2,1}^2&\phi_{2,1}\\
\phi_{1,1}&\phi_{2,1}&q_5
\end{bmatrix}.$$ By the formula
$
\left(\Ab^{-1}+\Mb^{-1}\right)^{-1}=\Ab-\Ab\left(\Ab+\Mb\right)^{-1}\Ab
$, 
we have
\begin{align*}
\left(\Mb^{-1}+\Ab^{-1}\right)^{-1}=&\begin{bmatrix}
1/\overline{m}_{1,d}& &\\
&1/\overline{m}_{2,d}&\\
& &1/\overline{m}_{3,d}
\end{bmatrix}-\Ab\begin{bmatrix}
q_4\mu_{1,1}^2+1/\overline{m}_{1,d}&q_4\mu_{1,1}\mu_{2,1}&\phi_{1,1}\\
q_4\mu_{1,1}\mu_{2,1}&q_4\mu_{2,1}^2+1/\overline{m}_{2,d}&\phi_{2,1}\\
\phi_{1,1}&\phi_{2,1}&q_5+1/\overline{m}_{3,d}
\end{bmatrix}^{-1}\Ab.
\end{align*}
Denote $\bl=[q_4\mu_{1,1}^2\overline{m}_{1,d}\quad q_4\mu_{1,1}\mu_{2,1}\overline{m}_{2,d}\quad  \phi_{1,1}\overline{m}_{3,d}]^\T$. Then by plugging the equation above into \eqref{eq:w2_woodberry}, and combing it with \eqref{eq:term1}, \eqref{eq:term2} and \eqref{eq:term3}, we finally get
\begin{align*}
    \overline{m}_{1,d}=&\psi_1\EE w_2\\
    =&\psi_1\Big\{-\xi+q_2\mu_{1,2}^2+q_4\mu_{1,1}^2-q_4^2\mu_{1,1}^2(\mu_{1,1}^2\overline{m}_{1,d}+\mu_{2,1}^2\overline{m}_{2,d})\\
    &-(\phi_{1,1}^2+\phi_{1,2}^2)\overline{m}_{3,d}+\bl^\T\Ab\bl-\bl^\T\Ab(\Ab+\Mb)^{-1}\Ab\bl\Big\}^{-1}+o_d(1)\\
=&\psi_1\Bigg\{-\xi+q_2\mu_{1,2}^2+q_4\mu_{1,1}^2-\phi_{1,2}^2\overline{m}_{3,d}-\begin{bmatrix}
q_4\mu_{1,1}^2\\q_4\mu_{1,1}\mu_{2,1}\\ \phi_{1,1}
\end{bmatrix}^\T(\Ab+\Mb)^{-1}\begin{bmatrix}
q_4\mu_{1,1}^2\\q_4\mu_{1,1}\mu_{2,1}\\ \phi_{1,1}
\end{bmatrix}\Bigg\}^{-1}\\
&+o_d(1).
\end{align*}
Now note that $\phi_{\sfc,1}^2=\mu_{\sfc,1}^2(1+q_1)^2$, $\phi_{\sfc,2}^2=\mu_{\sfc,2}^2$, $\sfc=1,2$. 
Therefore with direct calculation, we have
\begin{equation}
\label{eq:implicit_d1}
    \begin{split}
      \overline{m}_{1,d}=&\psi_1\bigg\{-\xi+q_2\mu_{1,2}^2-\mu_{1,2}^2\overline{m}_{3,d}+\frac{H_{1,d}}{H_{D,d}}\bigg\}^{-1}+o_d(1),
    \end{split}
\end{equation}
where
$$\begin{aligned}
H_{1,d}=&\mu_{1,1}^2q_4(1+\overline{m}_{3,d}q_5)-\mu_{1,1}^2(1+q_1)^2\overline{m}_{3,d},\\
H_{D,d}=&(1+\mu_{1,1}^2\overline{m}_{1,d}q_4+\mu_{2,1}^2\overline{m}_{2,d}q_4)(1+\overline{m}_{3,d}q_5)-\mu_{2,1}^2(1+q_1)^2\overline{m}_{2,d}\overline{m}_{3,d}\\
&-\mu_{1,1}^2(1+q_1)^2\overline{m}_{1,d}\overline{m}_{3,d}.
\end{aligned}$$
The  equation above shows that the magnitude of the first element of $\overline{\mb}_d(\xi) - \sfFb(\overline{\mb}_d(\xi))$ is $o_d(1)$. With exactly the same proof, we also have 
\begin{equation}
\label{eq:implicit_d2}
    \begin{split}
      \overline{m}_{2,d}=&\psi_2\bigg\{-\xi+q_2\mu_{2,2}^2-\mu_{2,2}^2\overline{m}_{3,d}+\frac{H_{2,d}}{H_{D,d}}\bigg\}^{-1}+o_d(1),\\
\overline{m}_{3,d}=&\psi_3\bigg\{-\xi+q_3-\mu_{1,2}^2\overline{m}_{1,d}-\mu_{2,2}^2\overline{m}_{2,d}+\frac{H_{3,d}}{H_{D,d}}\bigg\}^{-1}+o_d(1),  
    \end{split}
\end{equation}
where
$$\begin{aligned}
H_{2,d}=&\mu_{2,1}^2q_4(1+\overline{m}_{3,d}q_5)-\mu_{2,1}^2(1+q_1)^2\overline{m}_{3,d},\\
H_{3,d}=&q_5(1+\mu_{1,1}^2\overline{m}_{1,d}q_4+\mu_{2,1}^2\overline{m}_{2,d}q_4)-\mu_{2,1}^2(1+q_1)^2\overline{m}_{2,d}-\mu_{1,1}^2(1+q_1)^2\overline{m}_{1,d}.
\end{aligned}$$
This completes the proof of Lemma~\ref{lemma:resolvent_calculation}.
\end{proof}

\subsection{Proof for conclusions 1 and 2 in Proposition~\ref{prop:implicit1}}
\label{subsec:conclusion12}
 We first introduce an important lemma about the property of Stieltjes transforms, which is given in \cite{hastie2022surprises}.
\begin{lemma}[Lemma 7 in \cite{hastie2022surprises}]
\label{lemma:Stieljes_property}
The functions $\xi\rightarrow \overline{m}_{i,d}(\xi)$, $i=1,2,3$, have the following properties:
\begin{enumerate}[leftmargin = *]
    \item $\overline{m}_{i,d}$, $i=1,2,3$ are analytical on $\bbC_+$, and map $\bbC_+$ into $\bbC_+$.
    \item Let $\Omega\subset\bbC_+$ be a set with an accumulation point. If $~\overline{m}_{i,d}\rightarrow m_i(\xi)$ for all $\xi\in\Omega$, then $m_i(\xi)$ has an unique analytic continuation to $\bbC_+$ and $\overline{m}_{i,d}\rightarrow m_i(\xi)$  for all $\xi\in\bbC_+$. Moreover, the convergence is uniform over compact sets $\Omega\subset \bbC_+$.
\end{enumerate}
\end{lemma}
We now give the proof of the conclusions in Proposition~\ref{prop:implicit1} that $\mb(\xi;\qb,\bmu)$ is analytic on $\{\xi:\Im(\xi)>\xi_0\}$ for some sufficiently large $\xi_0$, has unique analytic continuation to $\bbC_+$ and maps $\bbC_+$ to $\bbC_+^3$. Denote $\overline{\mb}_d=\overline{\mb}_d(\xi)=[ \overline{m}_{1,d}(\xi), \overline{m}_{2,d}(\xi), \overline{m}_{3,d}(\xi) ]^\T$. Then for any fixed $\xi\in\bbC_+$,  Lemma~\ref{lemma:resolvent_calculation} gives
\begin{equation}\label{eq:mbar_approximate_solution}
    || \overline{\mb}_d-\sfFb( \overline{\mb}_d)||_2=o_{d}(1).
\end{equation}
By Lemma~\ref{lemma:uniquesolution}, there exists a $\xi_0>0$, such that for all $\xi$ with $\Im(\xi)\geq\xi_0$, $\sfFb(\cdot)$ is $1/2$-Lipschitz with respect to $\ell_2$ norm. Moreover, for all $\xi$ with $\Im(\xi)\geq\xi_0$ we have
\begin{align*}
    \| \overline{\mb}_d - \mb \|_2 &= \| \overline{\mb}_d - \sfFb( \mb) \|_2 \\
    &\leq \| \overline{\mb}_d - \sfFb( \overline{\mb}_d) \|_2 + \| \sfFb( \overline{\mb}_d) - \sfFb( \mb) \|_2\\
    &\leq o_{d}(1) + \frac{1}{2}\cdot \|  \overline{\mb}_d - \mb  \|_2,
\end{align*}
where the equality is by the definition of $\mb$ as the unique fixed point of $\sfFb(\cdot)$, the first inequality is by triangle inequality, the second inequality is by \eqref{eq:mbar_approximate_solution} and the fact that $\sfFb(\cdot)$ is $1/2$-Lipschitz with respect to $\ell_2$ norm. Therefore we have $\| \overline{\mb}_d(\xi) - \mb(\xi) \|_2  =  o_{d}(1)$ for all $\xi$ with $\Im(\xi)\geq\xi_0$. The properties of Stieljes transforms (see Lemma~\ref{lemma:Stieljes_property}) then imply that $\mb(\xi;\qb,\bmu)$ is analytic in $\{\xi:\Im(\xi)>\xi_0\}$, and has a unique analytic continuation to $\bbC_+$. Moreover, the extended $\mb(\xi;\qb,\bmu)$ satisfies
\begin{align}\label{eq:transformcalculation_part1}
    || \overline{\mb}_d(\xi)-\mb(\xi)||_2=o_d(1)
\end{align}
for any fixed $\xi\in\bbC_+$. This implies  that $\mb(\xi;\qb,\bmu)$ is $\bbC_+\rightarrow\bbC_+^3$ by the definition of $\overline{\mb}_d$. The proof of Conclusion 1 is complete.

To prove Conclusion 2, we first prove that $\mb(\xi)$  is a continuity point of $\sfFb(\cdot)$ for any fixed $\xi\in\bbC_+$. For any fixed $\xi\in\bbC_+$, assume that $\mb(\xi)$  is not a continuity point of $\sfFb(\cdot)$ , by the definition of $\sfFb(\cdot)$ we have $\|\sfFb(\mb(\xi))\|_2=+\infty$.  Therefore, for any $M>0$, there exsits $\delta(\xi,M)>0$ ($\xi\in\bbC_+$ is fixed here), as long as  $\|\overline{\mb}_d(\xi)-\mb(\xi)\|<\delta(\xi,M)$, the inequality $\sfFb(\overline{\mb}_d(\xi))>M$ holds. 
Moreover, for the $\delta(\xi,M)$, there always exists $d_0$ such that $\|\overline{\mb}_{d}(\xi)-\mb(\xi)\|<\delta(\xi,M)$ for all $d>d_0$. That is: for any fixed $\xi\in\bbC_+$, and any large constant $M>0$, there always exists $d_0$ such that $\sfFb(\overline{\mb}_d(\xi))>M$ for $d>d_0$. Combined with \eqref{eq:mbar_approximate_solution}, there exists $d_1>0$, such that $|| \overline{\mb}_d(\xi)-\sfFb( \overline{\mb}_d(\xi))||_2<1$ for all $d>d_1$. Then for $d>\max\{d_0,d_1\}$, we have $\sfFb(\overline{\mb}_d(\xi))>M$ and $|| \overline{\mb}_d(\xi)-\sfFb( \overline{\mb}_d(\xi))||_2<1$. We have $\|\overline{\mb}_d(\xi)\|_2>M-1$ for $d>\max\{d_0,d_1\}$. On the other hand, $\|\overline{\mb}_d(\xi)\|_2\leq2(\psi_1+\psi_2+\psi_3)/\Im(\xi)$ from the definition of $\overline{\mb}_d(\xi)$. Note that $\xi\in\bbC_+$ is fixed here. Enlarging $M$ leads to a   contradiction. Therefore, $\mb(\xi)$  is the continuity point of $\sfFb(\cdot)$ for any fixed $\xi\in\bbC_+$.

For any fixed $\xi\in\bbC_+$, note that $\mb(\xi)$  is the continuity point of $\sfFb(\cdot)$. Let $d\rightarrow+\infty$,  \eqref{eq:mbar_approximate_solution} and \eqref{eq:transformcalculation_part1} give us 
\begin{align*}
    || {\mb}-\sfFb( {\mb})||_2=0.
\end{align*}
This means that $\sfFb(\mb(\xi;\qb,\bmu))=\mb(\xi;\qb,\bmu)$ for any fixed $\xi\in\bbC_+$. The proof of Conclusion 2 is complete.

\subsection{Point convergence for polynomial activation functions}
\label{subsec:proofpoly}
We now give the proof of point convergence under the additional assumption that the activation functions are  polynomials. We remind readers the ``Gaussian version'' of the problem defined   in Sections~\ref{sec:equivalenceinGaussianandSphere} and \ref{subsec:resolvent},  where the data inputs $\overline{\xb}_i$, $i\in[n]$ and $\overline{\btheta}_a$, $a\in [N]$ are defined in Definition~\ref{def:Gaussian_data_rf} and the activation functions $\phi_{1}(x), \phi_{2}(x)$ are given in \eqref{eq:def_activation_Gaussian}. We also remind readers that the ``Gaussian version'' and ``spherical version'' Stieltjes transforms of the empirical eigenvalue distributions of linear pencil matrices are denoted as $\overline{M}_d(\xi)$ and $M_d(\xi)$,  respectively. Importantly, the expectation of $\overline{M}_d(\xi)$ is denoted as $\overline{m}_d$, while $m(\xi;\qb,\bmu) = \sum_{i=1}^3 m_i(\xi) $, where $\mb=\mb(\xi) = (m_1(\xi),m_2(\xi),m_3(\xi))^\T$ is defined as the solution of \eqref{eq:m123} on $\{\xi: \Im(\xi) \geq \xi_0\}$ and then extended to $\bbC_+$ by analytic continuation.

By \eqref{eq:concentration}, for all fixed $\xi\in\bbC_+$, we have
\begin{align}
\EE\Bigg| \overline{M}_{d}(\xi)-\sum\limits_{i=1}^3\overline{m}_{i,d}(\xi)\Bigg|=o_d(1).\label{eq:secB32_result}
\end{align} 
In addition, by  Lemma~\ref{lemma:GaussianSphere}, when the activation functions are polynomials, we have 
\begin{align}\label{eq:transformcalculation_part3}
\EE\big|M_d(\xi)- \overline{M}_d(\xi)\big|=o_d(1).
\end{align}
Combining  \eqref{eq:transformcalculation_part1} \eqref{eq:secB32_result} and \eqref{eq:transformcalculation_part3} gives
\begin{align*}
    \EE\big| M_d(\xi)-m(\xi)\big|=o_d(1)
\end{align*}
for any fixed $\xi\in\bbC_+$,
which completes the proof of the point convergence for polynomial activation functions. 

\subsection{Point convergence for general activation functions satisfying Assumption~\ref{assump1}}
\label{subsec:completeproof} 
We now extend the result for polynomial activation functions to  general activation functions satisfying Assumption~\ref{assump1}.
Let $\tau_d$ be the marginal distribution of $\la\xb,\btheta\ra/\sqrt{d}$ for $\xb,\btheta\sim_{\iid}\Unif(\sqrt{d}\cdot\SSS^{d-1})$, and ${\overline{\tau}}_d$  the marginal distribution of $\la\overline\xb,\overline\btheta\ra/\sqrt{d}$ for $\overline\xb,\overline\btheta\sim_{\iid}\rmN(0,\Ib_d)$. For $j=1,2$, suppose that $\sigma_j$ are activation functions satisfying Assumption~\ref{assump1}. The idea here    is  to construct polynomial activation functions $\tilde{\sigma}_j$ to approximate $\sigma_j$. To do so, we  recall that $\mb= \mb(\xi;\qb,\bmu)$ solves  the implicit equations
\begin{align*}
    \mb = \sfFb(\mb;\xi,\qb,\bmu),
\end{align*}
where $\sfFb(\cdot;\xi,\qb,\bmu)$ is defined in Definition~\ref{def:implicit1}. When $\Im(\xi)>\xi_0$ for some large enough $\xi_0$, by the continuity of the solution of the fixed point equation with respect to $\bmu$, we have 
\begin{align*}
    \lim_{\tilde\bmu \rightarrow \bmu} \mb(\xi;\qb,\tilde\bmu) = \mb(\xi;\qb,\bmu).
\end{align*}
According to our proof in Section~\ref{subsec:conclusion12}, we can extend the definition of $\mb$ to $\xi\in \bbC_+$ with analytic continuation. Then with the same proof as in \citet{mei2022generalization} (see equation (10.56) in \citet{mei2022generalization}), for any fixed $\xi\in \bbC_+$ and any $\varepsilon > 0$, there exists $\delta = \delta(\varepsilon, \xi, \qb,\bmu)>0$ such that
%
\begin{align}\label{eq:m_continuity_mu}
    \| \mb(\xi;\qb,\tilde\bmu) - \mb(\xi;\qb,\bmu) \|_2 \leq \varepsilon
\end{align}
for all $\tilde\bmu$ with $\lVert\bmu-\tilde\bmu\rVert_2\leq\delta$. 
Now by Assumption~\ref{assump1}, for any fixed $\varepsilon>0$, we can choose a sufficiently large integer $\bar{k}$ and construct 
$$
\tilde{\sigma}_{\sfc}(x)=\sum\limits_{k=0}^{\bar k}\frac{\mu_{\sfc,k}}{k!}H_k(x),
$$
such that for $  {G\sim N(0,1)}$,
\begin{align}
    &\EE  [ \sigma_{\sfc}(G)-\tilde{\sigma}_{\sfc}(G) ]^2 \leq \varepsilon^2, \label{eq:sigma_difference_Gaussian1}\\
    &
    |\EE  \{\tilde\sigma_j(G)^2\}- \EE  \{\tilde\sigma_j(G)^2\}| \leq \delta^2 / 2.\label{eq:sigma_difference_Gaussian2}
\end{align}
Here, $\{H_k(x)\}$ are the  family of Hermite polynomials. 
Then by \eqref{eq:sigma_difference_Gaussian1} and Lemma~5 in \citet{ghorbani2021linearized}, we have
\begin{align}\label{eq:sigma_difference}
    \lVert\sigma_{\sfc}-\tilde{\sigma}_{\sfc} \rVert_{L^2(\tau_d)}\leq  \varepsilon
\end{align}
for $\sfc = 1,2$ and sufficiently large $d$, where 
we denote $\lVert\sigma_{\sfc}-\tilde{\sigma}_{\sfc} \rVert_{L^2(\nu)}=\int(\sigma_{\sfc}(x)-\tilde{\sigma}_{\sfc}(x))^2{\nu}(dx)$.

Following Definition~\ref{def:someconstant}, we can also define the parameters $\tilde\mu_{\sfc,0}, \tilde\mu_{\sfc,1}, \tilde\mu_{\sfc,2}$ corresponding to the polynomial activation functions $\tilde\sigma_\sfc$ by
\begin{align*}
    \tilde\mu_{j,0}\triangleq \EE \{\tilde\sigma_j(G)\},~ \tilde\mu_{j,1}\triangleq\EE \{G\tilde\sigma_j(G)\},~
    \tilde\mu_{j,2}^2\triangleq{\EE \{\tilde\sigma_j(G)^2\}-\tilde\mu_{j,0}^2-\tilde\mu_{j,1}^2}.
\end{align*}
Then by the definition of $\tilde\sigma_{\sfc}$, we have $\mu_{\sfc,0} = \tilde\mu_{\sfc,0}$, $\mu_{\sfc,1}=\tilde\mu_{\sfc,1}$ for $\sfc = 1,2$.  Moreover, we also have
\begin{align*}
    |\mu_{\sfc,2}- \tilde\mu_{\sfc,2}| \leq \sqrt{ | \mu_{\sfc,2}^2 - \tilde\mu_{\sfc,2}^2 |} = \sqrt{| \EE \{\sigma_j(G)^2 - \tilde\sigma_j(G)^2\}|} \leq \delta  / \sqrt{2}
\end{align*}
for $\sfc = 1,2$, where the first inequality follows from  $|a - b| \leq \sqrt{|a^2 - b^2|}$ for all $a,b > 0$, the equality follows by $\mu_{\sfc,0} = \tilde\mu_{\sfc,0}$, $\mu_{\sfc,1}=\tilde\mu_{\sfc,1}$ for $\sfc = 1,2$, and the last  inequality follows by \eqref{eq:sigma_difference_Gaussian2}. Therefore we have $\| \tilde\bmu - \bmu \|_2 \leq \delta$. 

Let $\tilde\mb(\xi) = [\tilde m_1(\xi), \tilde m_2(\xi), \tilde m_3(\xi)]^\T$ be the solution of the implicit equations
\begin{align*}
    \tilde\mb = \sfFb(\tilde\mb;\xi,\qb,\tilde\bmu),
\end{align*}
and let $\tilde{m}(\xi) = \tilde{m}_1(\xi) + \tilde{m}_2(\xi) + \tilde{m}_3(\xi)$, where we drop the arguments $\qb,\tilde\bmu$ in $\tilde\mb(\xi;\qb,\tilde\bmu)$ for notation simplification. 
Then by \eqref{eq:m_continuity_mu}, we have
\begin{align}\label{eq:mandgeneralm}
   \big|\tilde{m}(\xi)- m(\xi)\big|\leq 3\varepsilon.
\end{align}
Let $\tilde{\Ab}$ be the linear pencil matrix corresponding to $\tilde{\sigma}$ in Definition~\ref{def:linear pencil}, and define
$\tilde{M}_d(\xi)=(1/d)\cdot \tr[(\tilde{\Ab}-\xi\Ib)^{-1}]$. 
Then we have
\begin{align}
    \EE\big[|M_d(\xi)-\tilde{M}_d(\xi)|\big] 
    =& \frac{1}{d}\EE\big[\big|\tr[(\Ab-\xi\Ib)^{-1}(\tilde{\Ab}-\Ab)(\tilde{\Ab}-\xi\Ib)^{-1}]\big|\big]\nonumber
    \\
    \leq&\frac{1}{d}\EE\big[\lVert(\Ab-\xi\Ib)^{-1}(\tilde{\Ab}-\xi\Ib)^{-1}\rVert_{\op}\lVert\tilde{\Ab}-\Ab\rVert_{\star}\big]\nonumber\\
    \leq&\big[1/(\Im(\xi)^2\big]\cdot \PO^{-1/2} \cdot \EE\{\lVert\tilde{\Ab}-\Ab\rVert_F^2\}^{1/2}\nonumber\\
    \leq&C'(\xi,\bpsi)\cdot \big[1/(\Im(\xi)^2\big]\cdot d^{-1/2} \cdot \EE\{\lVert\tilde{\Ab}-\Ab\rVert_F^2\}^{1/2}\nonumber\\
    \leq& C''(\xi,\qb)\cdot\big(\lVert\sigma_1-\tilde{\sigma}_1\rVert_{L^2(\tau_d)}+\lVert\sigma_2-\tilde{\sigma}_2\rVert_{L^2(\tau_d)}\big),\label{eq:M_difference}
\end{align}
where $C'(\xi,\bpsi)> 0$ is a constant only depending on $\xi$ and $\bpsi$, and $C''(\xi,\qb) > 0$ only depends on $\xi$, $\qb$ and $\bpsi$. Here the second inequality above follows by Cauchy-Schwarz inequality, the third inequality follows by $P = N_1 + N_2 + n$ and the assumption that $N_1,N_2,n,d$ goes to infinity proportionally, and the last inequality follows by the definitions of $\tilde{\Ab}$ and $\Ab$. 
Therefore, by \eqref{eq:sigma_difference} and \eqref{eq:M_difference}, we have
\begin{equation}\label{eq:M_difference_bound}
\begin{split}
\EE|M_d(\xi)-\tilde{M}_d(\xi)|\leq 2C''(\xi,\qb)\cdot\varepsilon
\end{split}
\end{equation}
for sufficiently large $d$. Moreover, since $\tilde\sigma_\sfc$, $\sfc = 1,2$ are polynomial activation functions, by the results in Appendix~\ref{subsec:proofpoly}, we have
\begin{align}\label{eq:polynomial_stieltjes_forgeneralproof}
    \EE\big| \tilde{M}_d(\xi)-\tilde{m}(\xi)\big|=o_d(1).
\end{align}
Combining \eqref{eq:mandgeneralm}, \eqref{eq:M_difference_bound} and \eqref{eq:polynomial_stieltjes_forgeneralproof} and taking $d \rightarrow \infty$, we have
\begin{align*}
\limsup_{d\rightarrow+\infty}\EE\big|M_d(\xi)- m(\xi)\big|\leq (2C''(\xi,\qb)+ 3)\cdot\varepsilon
\end{align*}
for all fixed $\xi\in\bbC_+$.  Taking $\varepsilon\rightarrow0^+$, we conclude that
$\lim\limits_{d\rightarrow\infty}\EE\big|\tilde{M}_d(\xi)-\tilde m(\xi)\big|=0$, which proves the point convergence for general activation functions.

\subsection{Uniform convergence on compact sets}
\label{subsec:uniformconvergence}
In this section, we aim to prove that on compact sets the point convergence established above could be extended to uniform convergence.  Consider a  compact set $\Omega\subset \bbC_+$. From the proof above we have \begin{align*}\lim\limits_{d\rightarrow+\infty}|\EE M_d(\xi;\qb,\bmu)-m(\xi;\qb,\bmu)|=0.\end{align*} 
Then from Lemma~\ref{lemma:Stieljes_property}, we have
\begin{align}\label{eq:eq_compact1}
\lim\limits_{d\rightarrow+\infty}\sup\limits_{\xi\in\Omega}|\EE M_d(\xi;\qb,\bmu)-m(\xi;\qb,\bmu)|=0.
\end{align}
Moreover, by the definition of $M_d(\xi;\qb,\bmu)$, we have 
\begin{align*}
    |M_d(\xi_1;\qb,\bmu)-M_d(\xi_2;\qb,\bmu)|&=\frac{1}{d}\big|\tr\big((\Ab-\xi_1\Ib)^{-1}(\xi_1-\xi_2)(\Ab-\xi_2\Ib)^{-1}\big)\big|\\
    &\leq \frac{\PO}{d\cdot \Im(\xi_1)\Im(\xi_2)}\cdot |\xi_1-\xi_2|.
\end{align*}
Since $P$ is proportional to $d$, there exists a  constant $L_0$ that only depends on $\psi_1,\psi_2,\psi_3$ and $\Omega$, such that $M_d(\xi;\qb,\bmu)$ is $L_0$-Lipschitz for all $d\in \NN$. 
Then by the compactness of $\Omega$, for any $\varepsilon>0$, there exists a finite set $\cN_{\varepsilon}(\Omega)\subset\bbC_+$, that is an $\varepsilon/L_0$ covering of the compact set $\Omega$. Specifically, for any $\xi\in\Omega$, there exists a $\xi_*\in\cN_{\varepsilon}(\Omega)$ such that $|\xi-\xi_*|<\varepsilon/L_0$. Therefore 
\begin{equation}\label{eq:eq_compact2}
    \begin{split}
            &\sup\limits_{\xi\in\Omega}\inf\limits_{\xi_*\in\cN_{\varepsilon}(\Omega)}|M_d(\xi;\qb,\bmu)-M_d(\xi_*;\qb,\bmu)|\leq\varepsilon,\\
    &\sup\limits_{\xi\in\Omega}\inf\limits_{\xi_*\in\cN_{\varepsilon}(\Omega)}|\EE M_d(\xi;\qb,\bmu)-\EE M_d(\xi_*;\qb,\bmu)|\leq\varepsilon
    \end{split}
\end{equation}
for all $d\in \NN$. 
Moreover,  since $\cN_{\varepsilon}(\Omega)$ is  finite, the number of $\xi_*$ is finite. Similar to the proof of \eqref{eq:concentration}, we have
\begin{align*}
    \sup\limits_{\xi_*\in\cN_{\varepsilon}(\Omega)}|M_d(\xi_*;\qb,\bmu)-\EE M_d(\xi_*;\qb,\bmu)|=o_{\PP}(1).
\end{align*}
Now since 
$| M_d(\xi_*;\qb,\bmu)|\leq P/(d\cdot \Im(\xi_*))\leq \PO/(d\cdot\inf\limits_{\xi\in\Omega}\Im(\xi))$, $| M_d(\xi_*;\qb,\bmu)|$ is bounded by some constant. By the dominated convergence theorem, we have
\begin{align}
    \label{eq:eq_compact3}
    \EE\sup\limits_{\xi_*\in\cN_{\varepsilon}(\Omega)}|M_d(\xi_*;\qb,\bmu)-\EE M_d(\xi_*;\qb,\bmu)|=o_{d}(1).
\end{align}
Combining \eqref{eq:eq_compact1}, \eqref{eq:eq_compact2} and \eqref{eq:eq_compact3}, we obtain
\begin{align*}
    &\EE\big[\sup_{\xi\in\Omega}\big|M_d(\xi;\qb,\bmu)-m(\xi;\qb,\bmu)\big|\big]\\
    &\quad=\EE\bigg\{\sup_{\xi\in\Omega}\inf\limits_{\xi_*\in\cN_{\varepsilon}(\Omega)}\big|M_d(\xi;\qb,\bmu)+ M_d(\xi_*;\qb,\bmu)- M_d(\xi_*;\qb,\bmu)+\EE M_d(\xi_*;\qb,\bmu)\\
    &\qquad-\EE M_d(\xi_*;\qb,\bmu)+\EE M_d(\xi;\qb,\bmu)-\EE M_d(\xi;\qb,\bmu)-m(\xi;\qb,\bmu)\big|\bigg\} \\
    &\quad\leq  \EE\bigg\{\sup\limits_{\xi\in\Omega}\inf\limits_{\xi_*\in\cN_{\varepsilon}(\Omega)}\big[|M_d(\xi;\qb,\bmu)-M_d(\xi_*;\qb,\bmu)|+|\EE M_d(\xi;\qb,\bmu)-\EE M_d(\xi_*;\qb,\bmu)|\\
    &\qquad+|M_d(\xi_*;\qb,\bmu)-\EE M_d(\xi_*;\qb,\bmu)|+|\EE M_d(\xi;\qb,\bmu)-m(\xi;\qb,\bmu)|\big] \bigg\}\\
    &\quad\leq2\varepsilon+o_d(1).
\end{align*}
Taking $d\rightarrow+\infty$, we have
\begin{align*}
    \lim\limits_{d\rightarrow+\infty}\EE\bigg[\sup_{\xi\in\Omega}\big|M_d(\xi;\qb,\bmu)-m(\xi;\qb,\bmu)\big|\bigg]\leq 2\varepsilon.
\end{align*}
Therefore, taking $\varepsilon\rightarrow 0^+$ proves Conclusion 3 in Proposition~\ref{prop:implicit1}. The proof of Proposition~\ref{prop:implicit1} is complete.

\section{Proof of Proposition~\ref{prop:substitutiongd}}
\label{sec:appendixsubstitution}
We first present some lemmas in Section~\ref{subsec:somelemmas}, and then complete the proof in Section~\ref{subsec:appendixsubstitutiongd}. Recall that we assume $\qb\in\mathcal{Q}$ (see Definition~\ref{def:linear pencil}).

\subsection{Preliminary lemmas}\label{subsec:somelemmas} 
The lemma below presents some additional properties of the function $\mb(\xi) = [m_1(\xi), m_2(\xi), m_3(\xi)]^\T$ defined in Proposition~\ref{prop:implicit1}.
\begin{lemma}
\label{lemma:Log1}
Let $\mb(\xi) = [m_1(\xi), m_2(\xi), m_3(\xi)]^\T$ defined on $\xi\in\bbC_+$ be the analytic continuation of the solution of the implicit equations $\mb = \sfFb(\mb;\xi,\qb,\bmu)$ defined in Proposition~\ref{prop:implicit1}. Then for any fixed $\xi_r\in\RR$ and $j = 1,2,3$, we have
$$
\lim_{u\rightarrow+\infty}|m_j(\xi_r +\rmi u)\cdot (\xi_r +\rmi u)+\psi_j|=0,
$$
\end{lemma}
\begin{proof}[Proof of Lemma~\ref{lemma:Log1}]
We denote $\xi_u = \xi_r +\rmi u$ for $u > 0$, and use the same definition of $H_1(\mb;\qb,\bmu)$ as in \eqref{eq:H(m)} that 
\begin{align*}
    H_1(\mb;\qb,\bmu)=&-\mu_{1,2}^2m_3+\frac{1}{m_1+\frac{-\mu_{2,1}^2(1+q_1)^2m_2m_3+(1+\mu_{2,1}^2m_2q_4)(1+m_3q_5)}{\mu_{1,1}^2q_4(1+m_{3}q_5)-\mu_{1,1}^2(1+q_1)^2m_{3}}}.
\end{align*}
Then we have 
\begin{align}\label{eq:H1toF1}
    \sfF_1(\mb;\xi,\qb,\bmu)=&\frac{\psi_1}{-\xi+q_2\mu_{1,2}^2+H_1(\mb;\qb,\bmu)}.
\end{align}
Moreover, define $\check{m}_1(\xi)=-\psi_1/\xi$, $\check{m}_2(\xi)=-\psi_2/\xi$, and $\check{m}_3=-\psi_3/\xi$, and denote $\check{\mb}(\xi)=[\check{m}_1(\xi),\check{m}_2(\xi),\check{m}_3(\xi)]^\T$. Then clearly we have $\lim\limits_{u\rightarrow +\infty} \check{\mb}(\xi_u) = \mathbf{0}$. By the definition of $H_1(\mb;\qb,\bmu)$, with simple calculations, we can see that $\lim\limits_{u\rightarrow+\infty}H_1(\check{\mb}(\xi_u);\qb,\bmu)=q_4\mu_{1,1}^2$. 
Thus by \eqref{eq:H1toF1}, we have 
\begin{align*}
    | \xi_u \cdot [\check{m}_1-\sfF_1(\check{\mb}(\xi_u);\xi_u,\qb,\bmu)] | &=\psi_1\cdot \bigg| \frac{q_2\mu_{1,2}^2+H_1(\check{\mb}(\xi_u);\qb,\bmu)}{\xi_u-q_2\mu_{1,2}^2-H_1(\check{\mb}(\xi_u);\qb,\bmu)} \bigg| = O_u\bigg(\frac{1 }{u}\bigg).
\end{align*}
Similarly, we can show that $| \xi_u \cdot [\check{m}_j-\sfF_1(\check{\mb}(\xi_u);\xi_u,\qb,\bmu)] | = O_u(1/u)$, $j = 2,3$. Therefore we have
\begin{align}\label{eq:proof_Log1_eq1}
    \xi_u\cdot ||\check{\mb}(\xi_u)-\sfFb(\check{\mb}(\xi_u);\xi_u,\qb,\bmu)||_2 = O_u(u^{-1}). 
\end{align}
Moreover, by Lemma~\ref{lemma:uniquesolution}, there exists a sufficiently large $\xi_0$ such that for any $K\geq\xi_0$, 
$\Fb(\cdot ;\xi_u,\qb,\bmu)$ is $1/2$-Lipschitz on the domain $\DD(2\psi_1/\xi_0)\times\DD(2\psi_2/\xi_0)\times\DD(2\psi_3/\xi_0)$. Therefore for sufficiently large $K$, 
\begin{align*}
&\lVert\check{\mb}(\xi_u)-\mb(\xi_u)\rVert_2\\
&~=\lVert \sfFb(\check{\mb}(\xi_u);\xi_u,\qb,\bmu)-\sfFb(\mb(\xi_u);\xi_u,\qb,\bmu)+\check{\mb}(\xi_u)-\sfFb(\check{\mb}(\xi_u);\xi_u,\qb,\bmu)\rVert_2\\
&~\leq\lVert\sfFb(\check{\mb}(\xi_u);\xi_u,\qb,\bmu)-\sfFb(\mb(\xi_u);\xi_u,\qb,\bmu)\rVert_2+\lVert\check{\mb}(\xi_u)-\sfFb(\check{\mb}(\xi_u);\xi_u,\qb,\bmu)\rVert_2\\
&~\leq\lVert\check{\mb}(\xi_u)-\mb(\xi_u)\rVert_2/2+\lVert\check{\mb}(\xi_u)-\sfFb(\check{\mb}(\xi_u);\xi_u,\qb,\bmu)\rVert_2,
\end{align*}
where the first equality follows by the definition of $\mb(\xi_u)$ as the fixed point of $\sfFb(\cdot;\xi_u,\qb,\bmu)$, the first inequality follows by triangle inequality, and the second inequality follows by the $1/2$-Lipschitz continuity of $\sfFb(\cdot ;\xi_u,\qb,\bmu)$ on the domain $\DD(2\psi_1/\xi_0)\times\DD(2\psi_2/\xi_0)\times\DD(2\psi_3/\xi_0)$ (note that $\mb(\xi_u)$ is automatically in this domain according to Lemma~\ref{lemma:uniquesolution}, and $\check{\mb}(\xi_u)$ is also in this domain by its definition).
Rearranging terms then gives
\begin{align}\label{eq:proof_Log1_eq2}
    \lVert\check{\mb}(\xi_u)-\mb(\xi_u)\rVert_2 \leq 2 \lVert\check{\mb}(\xi_u)-\sfFb(\check{\mb}(\xi_u);\xi_u,\qb,\bmu)\rVert_2.
\end{align}
Thus for $j = 1,2,3$, we have
\begin{align*}
    |m_j(\xi_r +\rmi u)\cdot (\xi_r +\rmi u)+\psi_j| &= \xi_u \cdot |m_j(\xi_u) - \check{m}_j(\xi_u) | \\
    & \leq  2\xi_u \cdot \lVert\check{\mb}(\xi_u)-\sfFb(\check{\mb}(\xi_u);\xi_u,\qb,\bmu)\rVert_2\\
    &= O_u(u^{-1}),
\end{align*}
where the first inequality follows by \eqref{eq:proof_Log1_eq2}, and the second equality follows by \eqref{eq:proof_Log1_eq1}. This completes the proof.
\end{proof}
The following lemma shows the asymptotics of the functions $G_d(\rmi u;\qb,\bmu)$ and $g(\rmi u;\qb,\bmu)$ (defined in Definition~\ref{def:linear pencil} and \eqref{def:gd} respectively) as $u$ goes to infinity.
\begin{lemma}
\label{lemma:Log2}
Let $G_d(\xi;\qb,\bmu)$ be defined in Definition~\ref{def:linear pencil} and  $g(\xi;\qb,\bmu)$  defined in \eqref{def:gd}. The following limits hold:
$$\begin{aligned}
&\lim_{u\rightarrow+\infty}\sup_{d\geq1}\EE|G_d(\rmi u ;\qb,\bmu)-(\psi_1+\psi_2+\psi_3)\Log(-\rmi u )|=0,\\
&\lim_{u\rightarrow+\infty}|g(\rmi u ;\qb,\bmu)-(\psi_1+\psi_2+\psi_3)\Log(-\rmi u )|=0.
\end{aligned}$$
\end{lemma}

\begin{proof}[Proof of Lemma~\ref{lemma:Log2}]
The real and imaginary parts of $G_d(\rmi u ;\qb,\bmu)-(\psi_1+\psi_2+\psi_3)\Log(-\rmi u )$ are
$$\begin{aligned}
&\Big|\Re\Big[\frac{1}{\PO}\sum\big(\Log(\lambda_i(\Ab)-\rmi u )-\Log(-\rmi u )\big)\Big]\Big|=\frac{1}{2\PO}\sum_{i=1}^\PO \Log(1+\lambda_i(A)^2/u^2)\leq\frac{\lVert\Ab\rVert_F^2}{2\PO u^2},\\
&\Big|\Im\Big[\frac{1}{\PO}\sum\big(\Log(\lambda_i(\Ab)-\rmi u )-\Log(-\rmi u )\big)\Big]\Big|=\frac{1}{\PO}\sum_{i=1}^\PO \arctan(\lambda_i(\Ab)/u)\leq\frac{\lVert A\rVert_F}{\PO^{1/2}u}.
\end{aligned}$$
By the definition of the linear pencil matrix $\Ab$, it is 
easy to see that $\frac{1}{\PO}\EE[\lVert \Ab\rVert_F^2]=O_d(1)$, thus 
$$
\lim_{u\rightarrow+\infty}\sup_{d\geq1}\EE|G_d(\rmi u ;\qb,\bmu)-(\psi_1+\psi_2+\psi_3)\Log(-\rmi u )|=0.
$$
For the asymptotics of $g(\rmi u ;\qb,\bmu)$, note that
\begin{align*}
    L(\xi,z_1,z_2,z_3;\qb,\bmu)=&L_1(z_1,z_2,z_3;\qb,\bmu)+L_2(\xi,z_1,z_2,z_3;\qb,\bmu),
\end{align*}
where
$$\begin{aligned}
&L_1(z_1,z_2,z_3;\qb,\bmu)\\
&\quad=\Log\left[(1+\mu_{1,1}^2z_1q_4+\mu_{2,1}^2z_2q_4)(1+z_3q_5)-\mu_{1,1}^2(1+q_1)^2z_1z_3-\mu_{2,1}^2(1+q_1)^2z_2z_3\right]\\
&\quad\quad-\mu_{1,2}^2z_1z_3-\mu_{2,2}^2z_2z_3+q_2\mu_{1,2}^2z_1+q_2\mu_{1,2}^2z_2+q_3z_3,\\
&L_2(\xi,z_1,z_2,z_3;\qb,\bmu)\\
&\quad=-\psi_1\Log(z_1/\psi_1)-\psi_2\Log(z_2/\psi_2)-\psi_3\Log(z_3/\psi_3)-\xi(z_1+z_2+z_3)-\psi_1-\psi_2-\psi_3.
\end{aligned}$$
We now calculate the limits of $L_1(m_1(\rmi u ),m_2(\rmi u ),m_3(\rmi u );\qb,\bmu)$ and $L_2(\rmi u ,m_1(\rmi u ),m_2(\rmi u ),m_3(\rmi u );\qb,\bmu)$ separately. For $L_1$, 
by Lemma~\ref{lemma:Log1}, we have
$$
\lim_{u\rightarrow+\infty}m_1(\rmi u )=0,~\lim_{u\rightarrow+\infty}m_2(\rmi u )=0,~\lim_{u\rightarrow+\infty}m_3(\rmi u )=0,
$$
which immediately implies that 
$$
\lim_{u\rightarrow+\infty}L_1(m_1(\rmi u ),m_2(\rmi u ),m_3(\rmi u );\qb,\bmu)=0.
$$
For $L_2$, note that by Lemma~\ref{lemma:Log1} we also have
$$
\lim_{u\rightarrow+\infty}|m_1(\rmi u )\rmi u +\psi_1|=0,~\lim_{u\rightarrow+\infty}|m_2(\rmi u )\rmi u +\psi_2|=0,~\lim_{u\rightarrow+\infty}|m_3(\rmi u )\rmi u +\psi_3|=0.
$$
Therefore,
$$\begin{aligned}
|L_2(&\rmi u ,m_1(\rmi u ),m_2(\rmi u ),m_3(\rmi u );\qb,\bmu)-(\psi_1+\psi_2+\psi_3)\Log(-\rmi u )|\\
\leq &\psi_1|\Log(-\rmi u m_1(\rmi u )/\psi_1)|+\psi_2|\Log(-\rmi u m_2(\rmi u )/\psi_2)|+\psi_3|\Log(-\rmi u m_3(\rmi u )/\psi_3)|\\
&+|\psi_1+\rmi u m_1(\rmi u )|+|\psi_2+\rmi u m_2(\rmi u )|+|\psi_3+\rmi u m_3(\rmi u )|\rightarrow0,
\end{aligned}$$
which completes the proof.
\end{proof}

The following lemma gives an important identity between $g(\xi;\qb,\bmu)$ and $m(\xi;\qb,\bmu)$.
\begin{lemma}\label{lemma:q_derivative_m}
For all $\xi\in \bbC_+$, it holds that 
$$\frac{\partial g}{\partial\xi} (\xi;\qb,\bmu)
=-(m_1 + m_2 + m_3) (\xi;\qb,\bmu)= -m(\xi;\qb,\bmu).$$
\end{lemma}
\begin{proof}[Proof of Lemma~\ref{lemma:q_derivative_m}]
By the definition of $L(\xi,z_1,z_2,z_3;\qb,\bmu)$, it is easy to see that 
$$\begin{aligned}
&\partial_{z_1}L(\xi,z_1,z_2,z_3;\qb,\bmu)\\
&\qquad=
-\mu_{1,2}^2z_3+q_2\mu_{1,2}^2-\psi_1/z_1-\xi\\
&\qquad\quad+\frac{\mu_{1,1}^2q_4(1+z_{3}q_5)-\mu_{1,1}^2(1+q_1)^2z_{3}}{(1+\mu_{1,1}^2z_{1}q_4+\mu_{2,1}^2z_{2}q_4)(1+z_{3}q_5)-\mu_{2,1}^2(1+q_1)^2z_{2}z_{3}-\mu_{1,1}^2(1+q_1)^2z_{1}z_{3}}\\
&\qquad=\psi_1\bigg(\frac{1}{\sfF_1(\zb)}-\frac{1}{z_1}\bigg), \\
\end{aligned}$$
$$\begin{aligned}
&\partial_{z_2}L(\xi,z_1,z_2,z_3;\qb,\bmu)\\
&\qquad=
-\mu_{2,2}^2z_3+q_2\mu_{2,2}^2-\psi_2/z_2-\xi\\
&\qquad\quad+\frac{\mu_{2,1}^2q_4(1+z_{3}q_5)-\mu_{2,1}^2(1+q_1)^2z_{3}}{(1+\mu_{2,1}^2z_{2}q_4+\mu_{1,1}^2z_{1}q_4)(1+z_{3}q_5)-\mu_{1,1}^2(1+q_1)^2z_{1}z_{3}-\mu_{2,1}^2(1+q_1)^2z_{2}z_{3}}\\
&\qquad=\psi_2\bigg(\frac{1}{\sfF_2(\zb)}-\frac{1}{z_2}\bigg), \\
\end{aligned}$$
$$\begin{aligned}
&\partial_{z_3}L(\xi,z_1,z_2,z_3;\qb,\bmu)\\
&\qquad=
-\mu_{1,2}^2z_1-\mu_{2,2}^2z_2+q_3-\psi_3/z_3-\xi\\
&\qquad\quad+\frac{q_5(1+\mu_{1,1}^2z_{1}q_4+\mu_{2,1}^2z_{2}q_4)-\mu_{2,1}^2(1+q_1)^2z_{2}-\mu_{1,1}^2(1+q_1)^2z_{1}}{(1+\mu_{2,1}^2z_{2}q_4+\mu_{1,1}^2z_{1}q_4)(1+z_{3}q_5)-\mu_{1,1}^2(1+q_1)^2z_{1}z_{3}-\mu_{2,1}^2(1+q_1)^2z_{2}z_{3}}\\
&\qquad=\psi_3\bigg(\frac{1}{\sfF_3(\zb)}-\frac{1}{z_3}\bigg),
\end{aligned}$$
where we utilize the definition of $\sfFb$ in Definition~\ref{def:implicit1} and write $\zb=[z_1,z_2,z_3]$. Then by Proposition~\ref{prop:implicit1}, we have
$$
\nabla_{\zb}L(\xi,\zb;\qb,\bmu)|_{\zb=\mb}\equiv\0
$$
for all $\xi\in \bbC_+$. 
By the formula of implicit differentiation, we have 
\begin{align*}
   \frac{\partial g(\xi;\qb,\bmu)}{\partial\xi} 
    & = \big[ \big\la  \nabla_{\zb} L(\xi,\zb;\qb,\bmu) |_{ \zb = \mb} , \partial_{\xi} \mb \big\ra + \partial_{\xi} L(\xi,\zb;\qb,\bmu) |_{\zb = \mb} \big] \\
    & = 0 +\frac{\rmd  L(\xi,\zb;\qb,\bmu)}{\rmd\xi} \Big|_{\zb = \mb}=-m(\xi;\qb,\bmu).
\end{align*}
This completes the proof of Lemma~\ref{lemma:q_derivative_m}.
\end{proof}

The following lemma further shows that the derivatives of $G_d$ and $g_d$ are asymptotically bounded.
\begin{lemma}
\label{lemma:Log3}
For fixed $\xi\in\bbC_+$, the following limits hold:
$$\begin{aligned}
&\limsup_{d\rightarrow+\infty}\big\{\EE\sup_{\qb\in\cQ}\lVert \nabla_\qb G_d(\xi;\qb,\bmu)\rVert_2\big\}+\sup_{\qb\in\cQ}\lVert\nabla_\qb g(\xi;\qb,\bmu)\rVert_2<+\infty,\\
&\limsup_{d\rightarrow+\infty}\big\{\EE\sup_{\qb\in\cQ}\lVert \nabla_\qb^2 G_d(\xi;\qb,\bmu)\rVert_{\op}\big\}+\sup_{\qb\in\cQ}\lVert\nabla_\qb^2 g(\xi;\qb,\bmu)\rVert_{\op}<+\infty,\\
&\limsup_{d\rightarrow+\infty}\big\{\EE\sup_{\qb\in\cQ}\lVert \nabla_\qb^3 G_d(\xi;\qb,\bmu)\rVert_{\op}\big\}+\sup_{\qb\in\cQ}\lVert\nabla_\qb^3 g(\xi;\qb,\bmu)\rVert_{\op}<+\infty.
\end{aligned}$$
\end{lemma}

\begin{proof}[Proof of Lemma~\ref{lemma:Log3}]
Let $\xi=\xi_r+\rmi u$, where $\xi_r\in\RR$ and $u\in\RR_+$ are both fixed. We also denote
\begin{align*}
    \Sbb_1&=\begin{bmatrix}
\0&\0&\tilde{\Zb}_1^\T \\
\0&\0&\tilde{\Zb}_2^\T \\
\tilde{\Zb}_1&\tilde{\Zb}_2&\0
\end{bmatrix}=\begin{bmatrix}
\0&\tilde{\Zb}\\
\tilde{\Zb}&\0
\end{bmatrix},\quad \Sbb_2=\begin{bmatrix}
        \Mb_2\Mb_2&\0\\
        \0&\0
        \end{bmatrix},\\
\Sbb_3&=\begin{bmatrix}
\0&\0\\
\0&\Ib_n
\end{bmatrix},\quad \Sbb_4=\begin{bmatrix}
\Mb_1\frac{\bTheta\bTheta^\T}{d}\Mb_1&\0\\
\0&\0
\end{bmatrix},\quad \Sbb_5=\begin{bmatrix}
\0&\0\\
\0&\frac{\Xb\Xb^\T}{d}
\end{bmatrix}.
\end{align*}
Then $\Sbb_1,\ldots, \Sbb_5$ are not related to $\qb$, and it is easy to see that
\begin{align*}
\limsup_{d\rightarrow+\infty}\EE\sup_{\qb\in\cQ}\lVert \Sbb_i\rVert_{\op}^{2k}<+\infty
\end{align*}
for any fixed $k\in\NN$. Moreover, define $\Rb=\Rb(\qb)=(\Ab(\qb)-\xi_r\Ib_{\PO}-\rmi u\Ib_{\PO})^{-1}$. Since $\Ab(\qb)$ is a real symmetric matrix, the imaginary parts in the eigenvalues of $\Ab(\qb)-\xi_r\Ib_{\PO}-\rmi u\Ib_{\PO}$ are all $-\rmi u$, and hence we deterministically have
\begin{align}\label{eq:R_op_bound}
    \sup_{\qb}\lVert\Rb\rVert_{\op}\leq1/u.
\end{align}
Therefore, by \eqref{eq:one_order}, \eqref{eq:two_order} and the definition of the linear pencil matrix $\Ab$, we have
\begin{equation*}
    \begin{split}
        &\EE\sup_{\qb\in\cQ}|\partial_{q_i}G_d(\xi;\qb)|=\EE\sup_{\qb\in\cQ}\frac{1}{d}|\tr(\Rb\Sbb_i)|\leq\EE\sup_{\qb\in\cQ}\frac{1}{u}[\lVert\Sbb_i\rVert_{\op}]=O_d(1),\\   
        &\EE\sup_{\qb\in\cQ}|\partial_{q_i,q_j}^2G_d(\xi;\qb)|=\EE\sup_{\qb\in\cQ}\frac{1}{d}|\tr(\Rb\Sbb_i\Rb\Sbb_j)|
        \leq \bigg(\EE\sup_{\qb\in\cQ}\frac{1}{u^2}[\lVert\Sbb_i\rVert_{\op}^2\lVert\Sbb_j\rVert_{\op}^2]\bigg)^{\frac12}=O_d(1).
        \end{split}
\end{equation*}
Similarly, for the third order derivatives, we also have
\begin{equation*}
    \begin{split}
        \EE\sup_{\qb\in\cQ}|\partial_{q_i,q_j,q_l}^3G_d(\xi;\qb)|=&\EE\Big\{\sup_{\qb\in\cQ}\frac{1}{d}|\tr(\Rb\Sbb_i\Rb\Sbb_j\Rb\Sbb_l)+|\tr(\Rb\Sbb_i\Rb\Sbb_l\Rb\Sbb_j)|\Big\}\\
        \leq&\frac{2}{u^3}\Big(\EE\sup_{\qb\in\cQ}[\lVert\Sbb_i\rVert_{\op}^4\lVert\Sbb_j\rVert_{\op}^4\lVert\Sbb_j\rVert_{\op}^4]\Big)^{\frac14}=O_d(1).
        \end{split}
\end{equation*}
This completes the Proof for $G_d(\xi;\qb,\bmu)$.
As for $g(\xi;\qb,\bmu)$, we first show that if $\qb_1\neq\qb_2$, the following property holds: 
\begin{align*}
    \frac{|m(\xi;\qb_1,\bmu)-m(\xi;\qb_2,\bmu)|}{\lVert\qb_1-\qb_2\rVert_2}&=\frac{\Big|\lim\limits_{d\rightarrow\infty}\EE\big(M_d(\xi;\qb_1,\bmu)-M_d(\xi;\qb_2,\bmu)\big)\Big|}{\lVert\qb_1-\qb_2\rVert_2}\\
    &=\lim\limits_{d\rightarrow\infty}\frac{\Big|\EE\big[\tr\big(\Rb(\qb_1)-\Rb(\qb_2)\big)\big]\Big|}{d\lVert\qb_1-\qb_2\rVert_2} \\
    &= \lim\limits_{d\rightarrow\infty}\frac{\Big|\EE\big[\tr\big(\Rb(\qb_1)\big(\Ab(\qb_1)-\Ab(\qb_2)\big)\Rb(\qb_2)\big)\big]\Big|}{d\lVert\qb_1-\qb_2\rVert_2} \\
    &\leq  \lim\limits_{d\rightarrow\infty}\frac{\PO}{d}\cdot\EE\frac{\lVert\Ab(\qb_1)-\Ab(\qb_2)\rVert_{\op}}{u^2\cdot\lVert\qb_1-\qb_2\rVert_2} <+\infty,
\end{align*}
where the first equality follows by Proposition~\ref{prop:implicit1}, the third equality follows by the identity $\Ab^{-1} - \Bb^{-1} =\Ab^{-1}(\Bb-\Ab)\Bb^{-1} $ for any invertible matrices $\Ab,\Bb$, the first inequality follows by $|\tr(\Ab\Bb)| \leq \PO\cdot \| \Ab \|_{\op}\| \Bb \|_{\op}$ for all $\Ab,\Bb\in \bbC^{\PO\times \PO}$ and \eqref{eq:R_op_bound}, and the last inequality follows by the linearity of $\Ab(\qb)$ in $\qb$. Therefore we have $\sup_{\qb\in\cQ}\lVert \nabla_\qb m(\xi;\qb,\bmu)\rVert_2<+\infty$. Similarly, we can also show that $\sup_{\qb\in\cQ}\lVert \nabla_\qb^j m(\xi;\qb,\bmu)\rVert_{\op} <+\infty$ for any fixed $\xi\in\bbC_+$ and $j=2,3$. Moreover, by Lemma~\ref{lemma:q_derivative_m}, we have
\begin{align*}
    \frac{d}{d\xi}g(\xi;\qb,\bmu)=-m(\xi;\qb,\bmu).
\end{align*}
Then $\sup_{\qb\in\cQ}\lVert \nabla_\qb^j m(\xi;\qb,\bmu)\rVert_{\op}<+\infty$ indicates that $\sup_{\qb\in\cQ}\lVert \nabla_\qb^j g(\xi;\qb,\bmu)\rVert_{\op}<+\infty$. 
This completes the proof of Lemma~\ref{lemma:Log3}.
\end{proof}

Finally, we present a classic result which shows that the derivatives of a function in a compact region can be upper bounded by the function value and the second derivatives of the function in the region.
\begin{lemma}[lemma 11.4 in \cite{mei2022generalization}]
\label{lemma:Log4}
Let $f\in C^2([a,b])$. Then we have
$$
\sup_{x\in[a,b]}|f'(x)|\leq\Big|\frac{f(a)-f(b)}{a-b}\Big|+\frac{1}{2}\sup_{x\in[a,b]}|f''(x)|\cdot|a-b|.
$$ Moreover, letting, $f\in C^2(\bB(\xb_0,2r))$ where $\bB(\xb_0,2r)=\{\xb\in\RR^d:\lVert\xb-\xb_0\rVert_2\leq r\}$ with a point $\xb_0$, we have
$$
\sup_{\xb\in\bB(\xb_0,2r)}\lVert\nabla f(\xb)\rVert_2\leq r^{-1}\sup_{\xb\in\bB(\xb_0,2r)}|f(\xb)|+2r\sup_{\xb\in\bB(\xb_0,2r)}\lVert\nabla^2f(\xb)\rVert_{\op}.
$$
\end{lemma}

\subsection{Completion of the proof}
\label{subsec:appendixsubstitutiongd}
By Lemma~\ref{lemma:q_derivative_m}, we have 
$\frac{\partial g}{\partial \xi}(\xi;\qb,\bmu)=-m(\xi;\qb,\bmu)$.
Hence, for $\xi\in\bbC_+$, $u \in\RR_+$, and any compact continuous path $c(\xi,\rmi u )$ connecting $\xi$ and $\rmi u$, we have
$$
g(\xi;\qb,\bmu)-g(\rmi u ;\qb,\bmu)=\int_{c(\xi,\rmi u )}m(x;\qb,\bmu)dx.
$$
Moreover, from Definition~\ref{def:linear pencil}, we also have $\frac{\rmd G_d(\xi;\qb,\bmu)}{\rmd\xi}=-M_d(\xi;\qb,\bmu)$, and
$$
G_d(\xi;\qb,\bmu)-G_d(\rmi u ;\qb,\bmu)=\int_{c(\xi,\rmi u )}M_d(x;\qb,\bmu)dx.
$$
The two equations imply that 
\begin{align}
\EE[|&G_d(\xi;\qb,\bmu)-g(\xi;\qb,\bmu)|]\nonumber\\
&\leq\EE|G_d(\rmi u ;\qb,\bmu)-g(\rmi u ;\qb,\bmu)|+\int_{c(\xi,\rmi u )}\EE|M_d(x;\qb,\bmu)-m(x;\qb,\bmu)|dx.\label{eq:G_g_eq1}
\end{align}
Therefore, by taking supremum limit on both sides above and  using Proposition~\ref{prop:implicit1},
we have
\begin{align*}
&\limsup_{d\rightarrow +\infty} ~ \EE[|G_d(\xi;\qb,\bmu)-g(\xi;\qb,\bmu)|]\\
&\qquad\leq \limsup_{d\rightarrow +\infty}\EE|G_d(\rmi u ;\qb,\bmu)-g(\rmi u ;\qb,\bmu)|+\limsup_{d\rightarrow +\infty}\int_{c(\xi,\rmi u )}\EE|M_d(x;\qb,\bmu)-m(x;\qb,\bmu)|dx\\
&\qquad = \limsup_{d\rightarrow +\infty}\EE|G_d(\rmi u ;\qb,\bmu)-g(\rmi u ;\qb,\bmu)|.
\end{align*}
Now  the left hand side above does not depend on $u$. Moreover, by Lemma~\ref{lemma:Log2}, we have
$$
\lim_{u\rightarrow+\infty}\limsup_{d\rightarrow +\infty}\EE|G_d(\rmi u ;\qb,\bmu)-g(\rmi u ;\qb,\bmu)|=0.
$$
Therefore 
\begin{align}
    \lim_{d\rightarrow+\infty}\EE[|G_d(\xi;\qb,\bmu)-g(\xi;\qb,\bmu)|]=0,\label{eq:subtitude1}
\end{align}
which proves the first equality in Proposition~\ref{prop:substitutiongd}. 

Next, we will prove the second and third equalities in Proposition~\ref{prop:substitutiongd}.
Define $V_d(\qb)=G_d(\xi;\qb,\bmu)-g(\xi;\qb,\bmu)$. Then by Lemma~\ref{lemma:Log4}, we have
\begin{align}
\sup_{\tilde{\qb}\in\Bb(\0,\varepsilon)}\lVert\nabla V_d(\tilde{\qb})\rVert_2\leq \frac{\sup_{\tilde{\qb}\in\Bb(\0,\varepsilon)}|V_d(\tilde{\qb})|}{\varepsilon}+2\varepsilon\sup_{\tilde{\qb}\in\Bb(\0,\varepsilon)}\lVert\nabla^2 V_d(\tilde{\qb})\rVert_{\op}.\label{eq:substitude_key}
\end{align}
By equation \eqref{eq:subtitude1}, Lemma~\ref{lemma:Log3} and the covering number argument (similar to Section~\ref{subsec:uniformconvergence} and the proof in Section 11.2 in \citet{mei2022generalization}), we get that $\lim\limits_{d\rightarrow+\infty}\EE\sup\limits_{\tilde{\qb}\in\cQ_{\star}}|V_d(\tilde{\qb})|=0$.
Again from Lemma~\ref{lemma:Log3} and its proof, we already have \begin{align*}\lim_{d\rightarrow+\infty}\EE\sup_{\tilde{\qb}\in\Bb(\0,\varepsilon)}|\nabla^2V_d(\tilde{\qb})|<C,
\end{align*} 
for some absolute value $C$. Therefore, by \eqref{eq:substitude_key}, we have
$$
\lim_{d\rightarrow+\infty}\EE[\lVert \partial_{\qb}G_d(\rmi u;\qb,\bmu)|_{\qb=\0}-\partial_{\qb}g(\rmi u;\qb,\bmu)|_{\qb=\0}\rVert_2]\leq C\varepsilon.
$$
Taking $\varepsilon\rightarrow 0^+$, we have 
\begin{align*}
    \lim_{d\rightarrow+\infty}\EE[\lVert \partial_{\qb}G_d(\rmi u;\qb,\bmu)|_{\qb=\0}-\partial_{\qb}g(\rmi u;\qb,\bmu)|_{\qb=\0}\rVert_2]=0.
\end{align*}
This completes the proof of the second equality in Proposition~\ref{prop:substitutiongd}. The proof of the third equation in Proposition~\ref{prop:substitutiongd} follows by a similar argument.

\section{Proof of Proposition~\ref{prop:existence_uniqueness_nu}}
\label{sec:appendix_nu}
The existence result is obtained by directly checking that $\mb(\xi;\0,\bmu)$ satisfies the two properties stated in Proposition~\ref{prop:existence_uniqueness_nu} as follows. The second property follows by the original definition of $\mb(\xi;\qb,\bmu)$ on $\{\xi: \Im(\xi) \geq \xi_0\}$. For the first property, by Proposition~\ref{prop:implicit1},  the analytic continuation of $\mb$ satisfies $\mb(\xi,\0,\bmu)\equiv\sfFb[\mb(\xi,\0,\bmu);\xi,\0,\bmu]$ for all $\xi\in \bbC_+$. This directly implies that $\mb(\xi,\0,\bmu)$ solves the system \eqref{eq:implicit2} for all $\xi\in \bbC_+$, which verifies the first property in Proposition~\ref{prop:existence_uniqueness_nu}. Moreover, by Proposition~\ref{prop:implicit1}, $\mb(\xi;\0,\bmu)$ is analytic, and $\mb(\xi;\0,\bmu) \in \bbC_+^3$ for all $\xi\in \bbC_+$. Therefore $\mb(\xi;\0,\bmu)$ is indeed an analytic function from $\bbC_+$ to $\bbC_+^3$. This completes the proof of existence.

To show the uniqueness, suppose that an analytic function $\bnu:$ $\bbC_+\rightarrow\bbC_+^3$ satisfies the properties satisfied in Proposition~\ref{prop:implicit1}. It then suffices to show that $\bnu(\xi;\bmu) \equiv\mb(\xi;\0,\bmu)$. By Lemma~\ref{lemma:uniquesolution}, we clearly  have $\bnu(\xi;\bmu) =  \mb(\xi;\0,\bmu)$ for all $\xi$ with $\Im(\xi) > \xi_0$. Now since both $\bnu(\xi;\bmu)$ and $\mb(\xi;\0,\bmu)$ are analytic on $\bbC_+$, the result follows by the uniqueness of analytic continuation.

We denote by $\bnu^*=\bnu(\sqrt{\lambda}\cdot\rmi;\bmu)$. From  the definition of $\overline{\mb}_d(\xi)$ in Lemma~\ref{lemma:resolvent_calculation}, we easily get that the elements in  $\overline{\mb}_d(\xi)$ are purely imaginary in the upper half-plane of $\bbC$ when $\qb=0$ and $\xi=\sqrt{\lambda}\cdot\rmi$. \eqref{eq:transformcalculation_part1} further indicates that elements in $\bnu^*$ are purely imaginary. Based on the proof  above, we have $\nu_j^*/\rmi\in\RR_+$.

\section{Proofs of Lemmas and Propositions in Appendix~\ref{sec:appendixthm-m}}\label{sec:addionalproofs-m}

\subsection{Proof of Lemma~\ref{lemma:uniquesolution-m}}
\label{sec:appendixuniquesolution-m}

When $\Im(\xi) \geq \xi_0$ for some sufficiently large $\xi_0$, we prove the existence and uniqueness of the solution by the Banach fixed point theorem. To do so, we want to show that 
\begin{enumerate}[leftmargin = *]
    \item $\sfFb(\cdot;\qb,\bmu)$ maps domain $\DD(2\psi_1/\xi_0)\times\cdots\times\DD(2\psi_K/\xi_0)\times\DD(2\psi_{K+1}/\xi_0)$ into itself. 
    \item $\sfFb(\cdot;\qb,\bmu)$ is Lipschitz continuous with a Liptichz constant smaller than 1.
\end{enumerate}
For  $\sfF_1(\cdot;\qb,\bmu)$, 
by Definition~\ref{def:implicit1-m}, we have
\begin{align*}
    \sfF_1(\mb;\xi,\qb,\bmu)=&\frac{\psi_1}{-\xi+q_2\mu_{1,2}^2+H_1(\mb;\qb,\bmu)},
\end{align*}
where 
\begin{align}
    H_1(\mb;\qb,\bmu)=&-\mu_{1,2}^2m_3+\frac{1}{m_1+\frac{-\sum_{\sfc=2}^K\mu_{\sfc,1}^2(1+q_1)^2m_\sfc m_{K+1}+(1+\sum_{\sfc=2}^K\mu_{\sfc,1}^2m_\sfc q_4)(1+m_{K+1}q_5)}{\mu_{1,1}^2q_4(1+m_{3}q_5)-\mu_{1,1}^2(1+q_1)^2m_{3}}}
    \label{eq:H(m)-m}
\end{align}
Note that $q_4,q_5\leq(1+q_1)/2$, it is easy to see that for $r_0$ small enough and  $\mb\in\DD(r_0)\times\DD(r_0)\times\DD(r_0)$, we have
\begin{align}\label{eq:fixedpoint_H1bound-m}
    |H_1(\mb;\qb,\bmu)|\leq 2+2|q_4|\mu_{1,1}^2.
\end{align}
Now as long as $\xi_0 \geq 4+4|q_4|\mu_{1,1}^2$, it is clear that for $\xi$ with $\Im(\xi) \geq \xi_0$ we have
\begin{align}\label{eq:fixedpoint_Imxi>H1-m}
    \Im(\xi) \geq \xi_0/ 2 + \xi_0/ 2 \geq \xi_0/ 2 +  2+2|q_4|\mu_{1,1}^2 \geq \xi_0/ 2 +  |H_1(\mb;\qb,\bmu)|,
\end{align}
where the last inequality follows by \eqref{eq:fixedpoint_H1bound-m}. 
Therefore we have
\begin{align*}
    | \sfF_1(\mb;\xi,\qb,\bmu) |  
    \leq~&\frac{\psi_1}{|\Im\big(\xi-q_2\mu_{1,2}^2-H_1(\mb;\qb,\bmu)\big)|}\\
    \leq~ & \frac{\psi_1}{\Im(\xi)-|H_1(\mb;\qb,\bmu)|}
    \leq \frac{2\psi_1}{\xi_0}
    ,
\end{align*}
where the  inequalities  follow from \eqref{eq:fixedpoint_Imxi>H1-m}.


Similarly, for $\sfF_\sfc$, $\sfc=2,\cdots,K+1$, we also have
$|\sfF_\sfc(\mb;\xi,\qb,\bmu)|\leq2\psi_\sfc/\xi_0$ provided $\xi_0 \geq 4+4 \max_\sfc\{ |q_4|\mu_{\sfc,1}^2, |q_5| \}$. 
Therefore if $\xi_0$ satisfies
$2\max\{\psi_1,\cdots,\psi_{K+1}\}/\xi_0\leq r_0$ and $\xi_0 \geq 4+4 \max_\sfc\{ |q_4|\mu_{\sfc,1}^2, |q_5| \}$,  it is clear that $\sfFb$ maps domain $\DD(2\psi_1/\xi_0)\times\cdots\times\DD(2\psi_{K+1}/\xi_0)$ into itself.

As for the Lipschitz continuity of $\sfFb(\cdot;\qb,\bmu)$, note that
$$
\nabla_{\mb}\sfF_1(\mb;\xi,\qb,\bmu)= - \frac{\psi_1}{(-\xi+q_2\mu_{1,2}^2+H_1(\mb;\qb,\bmu))^2} \cdot \nabla_{\mb}H_1(\mb;\qb,\bmu).
$$
It is easy to see that when $\xi_0$ is sufficiently large, 
$\lVert\nabla_{\mb}H_1(\mb;\qb,\bmu)\rVert_2\leq C(\qb,\bmu)$ for all $\mb\in\DD(2\psi_1/\xi_0)\times\cdots\times\DD(2\psi_{K+1}/\xi_0)$, where $C(\qb,\bmu)$ is a constant that only depends on $\qb$ and $\bmu$.
Thus when $\xi_0$ is sufficiently large, for $\xi$ with $\Im(\xi) \geq \xi_0$, 
$$
\big\lVert\nabla_{\mb}\sfF_1(\mb;\xi,\qb,\bmu)\big\rVert_2\leq \frac{C(\qb,\bmu) \cdot \psi_1}{\Im(\xi)-|H_1(\mb;\qb,\bmu)|} \leq \frac{4 C(\qb,\bmu) \cdot \psi_1 }{\xi_0} \leq \frac{1}{4K},
$$
where we again utilize \eqref{eq:fixedpoint_Imxi>H1-m}. 
We can apply the same argument for $\sfF_2,\ldots,\sfF_{K+1}$, and conclude that $\sfFb$ is $\frac12$-Lipschitz on $\mb\in\DD(2\psi_1/\xi_0)\times\cdots\times\DD(2\psi_{K+1}/\xi_0)$. Therefore by Banach fixed point theorem, there exists a unique fixed point of $\sfFb$. Thus the fixed point of the functions defined in Definition~\ref{def:implicit1-m} exists and is unique.

\subsection{Proof of Proposition~\ref{prop:implicit1-m}}
\label{subsec:appendixpropimplicit1-m}
Following  the same argument as in Lemma~\ref{lemma:GaussianSphere}, we may assume that all the elements in $\overline{\Xb}$ and $\overline{\bTheta}_\sfc$ are independently generated from standard normal $\rmN(0,1)$, and the activation functions are polynomials and centralized as $\phi_\sfc(x)=\sigma_\sfc(x)-\mu_{\sfc,0}$. 
The linear  pencil matrix of this Gaussian version is defined as 
\begin{align*}
    \overline{\Ab}(\qb,\bmu)=\begin{bmatrix}
q_2\mu_{1,2}^2\Ib_{N_1}+q_4\mu_{1,1}^2\frac{\overline{\bTheta}_1\overline{\bTheta}_1^\T }{d}&\cdots&q_4\mu_{1,1}\mu_{K,1}\frac{\overline{\bTheta}_1\overline{\bTheta}_K^\T }{d}& \overline{\Zb}_1^\T \\\vdots&\ddots&\vdots&\vdots\\
q_4\mu_{K,1}\mu_{1,1}\frac{\overline{\bTheta}_K\overline{\bTheta}_1^\T }{d}&\cdots&q_2\mu_{K,2}^2\Ib_{N_K}+q_4\mu_{K,1}^2\frac{\overline{\bTheta}_K\overline{\bTheta}_K^\T }{d}& \overline{\Zb}_K^\T  \\
 \overline{\Zb}_1&\cdots& \overline{\Zb}_K&q_3\Ib_n+q_5\frac{\overline{\Xb}~\overline{\Xb}^\T }{d}
\end{bmatrix}.
\end{align*}
Here $ \overline{\Zb} _\sfc=\Phi_\sfc\left(\overline{\Xb}~\overline{\bTheta}_\sfc^\T /\sqrt d\right)/\sqrt d\in\RR^{n\times N_\sfc}$, and $\Phi_\sfc(x)$ is defined as  $\Phi_\sfc(x)=\phi_\sfc(x)+q_1\mu_{\sfc,1}x$. Moreover, for $\sfc = 1,\ldots,K$ and with $ G\sim\rmN(0,1)$, we denote $\phi_{\sfc,0}\triangleq\EE \{\Phi_{\sfc}(G)\}$, $\phi_{\sfc,1}\triangleq\EE \{G\Phi_{\sfc}(G)\}$, $
\phi_{\sfc,2}\triangleq\EE \{\Phi_{\sfc}(G)^2\}-\phi_{\sfc,0}^2-\phi_{\sfc,1}^2$.
It is easy to see $\phi_{\sfc,0}=0$, $\phi_{\sfc,1}^2=\mu_{\sfc,1}^2(1+q_1)^2$, $\phi_{\sfc,2}^2=\mu_{\sfc,2}^2$.

We remind readers  that
$\cN_{\sfc}$ is the index set of units that use  the $\sfc$-th activation function $\sigma_{\sfc}$. Define the following terms:
\begin{equation*}
    \begin{split}
        \overline{m}_{\sfc,d}(\xi;\qb,\bmu)&=\EE\big[\overline{M}_{\sfc,d}(\xi;\qb,\bmu)\big],\quad \overline{M}_{\sfc,d}(\xi;\qb,\bmu)=\frac{1}{d}\text{tr}_{\cN_\sfc}\big[\overline{\Ab}(\qb,\bmu)-\xi \Ib_{\PO})^{-1}\big],\quad \sfc=1,\ldots,K\\
        \overline{m}_{K+1,d}(\xi;\qb,\bmu)&=\EE\big[\overline{M}_{K+1,d}(\xi;\qb,\bmu)\big],\quad \overline{M}_{K+1,d}(\xi;\qb,\bmu)=\frac{1}{d}\text{tr}_{[N+1:\PO]}\big[\overline{\Ab}(\qb,\bmu)-\xi \Ib_{\PO})^{-1}\big].    
    \end{split}
\end{equation*}
With the  same argument as in Lemma~\ref{lemma:GaussianSphere}, we obtain   
\begin{align*}
    \EE\Big|\sum\limits_{\sfc=1}^{K+1}\overline{M}_{\sfc,d}(\xi;\qb,\bmu)-{M}_d(\xi;\qb,\bmu) \Big|=o_d(1),\quad\text{for any fixed }\xi\in\bbC_+.
\end{align*}
Next,  by contraction properties we have   
\begin{align*}
   \EE\Big|\overline{M}_{\sfc,d}(\xi;\qb,\bmu)- \overline{m}_{\sfc,d}(\xi;\qb,\bmu)\Big|=o_d(1), \quad\text{for any fixed }\xi\in\bbC_+.
\end{align*}

To study $\overline{M}_{d}(\xi;\qb,\bmu)$, which is the Stieltjes transform of the empirical eigenvalue distribution of $\overline{A}(\qb,\bmu)$, it suffices to derive  the resolvent equations for $\overline{m}_{d}(\xi;\qb,\bmu)$ here. This is done by the following lemma.
\begin{lemma}\label{lemma:resolvent_calculation-m}
Let $\overline{\mb}_d(\xi)=[\overline{m}_{1,d}(\xi),\ldots,\overline{m}_{K+1,d}(\xi)]^\T$.  Then for any fixed $\xi\in\bbC_+$, the following  property holds:
\begin{equation*}
    \begin{split}
\| \overline{\mb}_d(\xi) - \sfFb(\overline{\mb}_d(\xi)) \|_2 = o_d(1).
    \end{split}
\end{equation*}
\end{lemma}
\begin{proof}[Proof of Lemma~\ref{lemma:resolvent_calculation-m}] 
Since $\overline{\mb}_d(\xi),\sfFb(\overline{\mb}_d(\xi)) \in \bbC^{K+1} $,  Lemma~\ref{lemma:resolvent_calculation-m} essentially contains  results showing that each element of $\overline{\mb}_d(\xi) - \sfFb(\overline{\mb}_d(\xi))$ is asymptotically zero. 
Since the proofs of the results are almost the same, we mainly focus on the proof of the first element $\overline{m}_{1,d}$. 
The proof still consists of three main steps similar to the proof of Lemma~\ref{lemma:resolvent_calculation}.

\smallskip\noindent{\textbf{Step 1}.}  We first use a leave-one-out argument to calculate $\overline{m}_{1,d}$. Let $\overline{\Ab}_{\cdot,N_1}$ be the $N_1^{\text{th}}$ column of $\overline{\Ab}$, with the $N_1^{\text{th}}$ entry removed. We further denote $\overline{\Bb}\in\RR^{( \PO-1)\times( \PO-1)}$ the matrix from $\overline{\Ab}$ by removing the $N_1^{\text{th}}$ row and $N_1^{\text{th}}$ column. From the Schur complement formula, we get
\begin{align}
    \overline{m}_{1,d}=\psi_1\EE\left(-\xi+q_2\mu_{1,2}^2+q_4\mu_{1,1}^2\lVert\overline{\btheta}_{N_1}\rVert_2^2/d-\overline{\Ab}_{\cdot,N_1}^\T (\overline{\Bb}-\xi\Ib_{ \PO-1})^{-1}\overline{\Ab}_{\cdot,N_1}\right)^{-1}.\label{eq:schurcomplement-m}
\end{align}
We decompose the vectors $\overline{\btheta}_{a}$, $a\in [N]$ and $\overline{\xb}_i$, $i\in [n]$ into components along the direction of $\overline{\btheta}_{N_1}$ and the other orthogonal directions:
\begin{equation}
\label{eq:schurdecomp-m}
    \begin{split}
        &\overline{\btheta}_{a}=\eta_{a}\frac{\overline{\btheta}_{N_1}}{\lVert\overline{\btheta}_{N_1}\rVert}+\tilde{\btheta}_{a},~\langle\overline{\btheta}_{N_1},\tilde{\btheta}_{a}\rangle=0,a\in[N]\backslash\{N_1 \},\\
&\overline{\xb}_i=u_i\frac{\overline{\btheta}_{N_1}}{\lVert\overline{\btheta}_{N_1}\rVert}+\tilde{\xb}_i,~\langle\overline{\btheta}_{N_1},\tilde{\xb}_i\rangle=0,i\in[n].
    \end{split}
\end{equation}

Note that for any $a\in [N]\backslash\{N_1 \}$ and $i\in [n]$, $\eta_{a}$, $u_i$ are standard Gaussian and are independent of $\tilde{\btheta}_{a}$ and $\tilde{\xb}_i$. Moreover,
$\tilde{\btheta}_{a}$ and $\tilde{\xb}_i$ are conditionally independent on each other given $\overline{\btheta}_{N_1}$, with $\tilde{\btheta}_{a},\tilde{\xb}_i\sim N(0,P_{\bot})$, where $P_{\bot}$ is the projector orthogonal to $\overline{\btheta}_{N_1}$. We then have $\overline{\Ab}_{\cdot,N_1}=(\overline{\Ab}_{1,N_1},...,\overline{\Ab}_{ \PO-1,N_1})^\T \in \RR^{ \PO-1}$ with
\begin{align*}
     \overline{\Ab}_{i,N_1}=\left\{\begin{aligned}
&\frac{q_4\mu_{1,1}^2\eta_i}{d}\lVert\overline{\btheta}_{N_1}\rVert_2, &&\text{if  } i\in[1,N_1-1],\\
&\frac{q_4\mu_{1,1}\mu_{\sfc,1}\eta_{i+1}}{d}\lVert\overline{\btheta}_{N_1}\rVert_2, &&\text{if } i+1\in\cN_\sfc,~\sfc\geq2,\\
&\frac{1}{\sqrt{d}}\Phi_{1}\Big(\frac{1}{\sqrt{d}}u_{i-N+1}\lVert\overline{\btheta}_{N_1}\rVert_2\Big), &&\text{if } i\geq N.
\end{aligned}\right.
\end{align*}
To calculate the resolvent equations, we need to further represent the matrix $\overline{\Bb}$ in \eqref{eq:schurcomplement-m} with $\eta_a$, $\tilde{\btheta}_a$, $u_i$, and $\tilde{\xb}_i$ for $a\in [N]\backslash\{N_1 \}$ and $i\in [n]$. 
Below we first list some additional notations for easier reference. Write
$\bmeta_1=[\eta_1,...,\eta_{N_1-1}]\in\RR^{N_1-1}$, $\bmeta_\sfc=(\eta_{\cN_\sfc})\in\RR^{N_\sfc}$, $\sfc=2,\ldots,K$
$\bmeta=[\bmeta_1^\T ,\ldots,\bmeta_K^\T ]^\T \in\RR^{N-1}$,
$\ub=(u_1,...,u_n)^\T\in\RR^n $,
$\tilde{\bTheta}_1=[\tilde{\btheta}_1,...,\tilde{\btheta}_{N_1-1}]^\T $,
$\tilde{\bTheta}_{\sfc}=[\tilde{\btheta}_{\cN_\sfc}]^\T $, 
{\small
\[ 
\tilde{\bTheta}=\begin{bmatrix}
\tilde{\bTheta}_1\\\vdots\\\tilde{\bTheta}_K
\end{bmatrix}\in\RR^{(N-1)\times d}, \quad \tilde{\Mb}_1=\begin{bmatrix}
\mu_{1,1}\Ib_{N_1-1}& &\\
&\ddots &\\
& & \mu_{K,1}\Ib_{N_K}
\end{bmatrix},\quad  
\tilde{\Mb}_*=\begin{bmatrix}
\mu_{1,2}\Ib_{N_1-1}& &\\
&\ddots &\\
& & \mu_{K,2}\Ib_{N_K}
\end{bmatrix}. 
\]}
With $\overline{\Bb}$ defined previously and \eqref{eq:schurdecomp-m},   $\overline{\Bb}_{[1:N-1],[1:N-1]}$ is decomposed into
\begin{align}
    \overline{\Bb}_{[1:N-1],[1:N-1]}=q_2\tilde{\Mb}_*\tilde{\Mb}_*+\frac{q_4}{d}\tilde{\Mb}_1\tilde{\bTheta}\tilde{\bTheta}^\T \tilde{\Mb}_1+\frac{q_4}{d}\tilde{\Mb}_1\bmeta\bmeta^\T \tilde{\Mb}_1\label{eq:B_decomp11-m}.
\end{align}
Moreover, for $i,j\in[n]$ and $a\in[N]\backslash\{N_1 \}$, we define
\begin{align*}
    \quad\big(\tilde{\Hb}\big)_{ij}=\frac{1}{d}\langle\tilde{\xb}_i,\tilde{\xb}_j\rangle.
\end{align*}
Then we could decompose $\overline{\Bb}_{[N:P-1],[N:P-1]}$ into
\begin{align}
    \overline{\Bb}_{[N:P-1],[N:P-1]}=q_3\Ib_n+q_5\tilde{\Hb}+\frac{q_5}{d}\ub\ub^\T.\label{eq:B_decomp22-m}
\end{align}
$\overline{\Bb}_{[N:P-1],[1:N-1]}=\overline{\Bb}_{[1:N-1],[N:P-1]}^\T$ holds due to the symmetry of $\overline{\Bb}$ . For $i,j\in[n]$ and $a\in\cN_\sfc\backslash\{N_1 \}$, elementally we have
\begin{align*}
    (\overline{\Zb})_{i,a}&=\frac{1}{\sqrt{d}}\Phi_{\sfc}\Big(\frac{1}{\sqrt{d}}\langle{\overline{\xb}}_i,{\overline{\btheta}}_{a}\rangle\Big)=\frac{1}{\sqrt{d}}\Phi_{\sfc}\Big(\frac{1}{\sqrt{d}}\langle\tilde{\xb}_i,\tilde{\btheta}_{a}\rangle+\frac{1}{d}u_i\eta_a\Big)\\
    &=\frac{1}{\sqrt{d}}\Phi_{\sfc}\Big(\frac{1}{\sqrt{d}}\langle\tilde{\xb}_i,\tilde{\btheta}_{a}\rangle\Big)+\frac{\phi_{\sfc,1}}{d}u_i\eta_a+\frac{1}{\sqrt{d}}\Big[\Phi_{\sfc,\bot}\Big(\frac{1}{\sqrt{d}}\langle\tilde{\xb}_i,\tilde{\btheta}_{a}\rangle+\frac{1}{\sqrt{d}}u_i\eta_{a}\Big)-\Phi_{\sfc,\bot}\Big(\frac{1}{\sqrt{d}}\langle\tilde{\xb}_i,\tilde{\btheta}_{a}\rangle\Big)\Big],
\end{align*}
where $\Phi_{\sfc,\bot}(x)=\Phi_{\sfc}(x)-\phi_{\sfc,1}x$. By the symmetry of $\overline{\Bb}$, we can then decompose $\overline{\Bb}_{[N:P-1],[1:N-1]}$ into
\begin{align}
    \overline{\Bb}_{[N:P-1],[1:N-1]}=\tilde{\Zb}+\frac1d\ub\bmeta \Mb_{\phi}+[\Eb_1,\Eb_2].\label{eq:B_decomp12-m}
\end{align}
Here, we define
\begin{align*}
    &\tilde{\Zb}=[\tilde{\Zb}_1,\ldots,\tilde{\Zb}_K],\quad(\tilde{\Zb}_\sfc)_{i,a}=\frac{1}{\sqrt{d}}\Phi_{\sfc}\Big(\frac{1}{\sqrt{d}}\langle\tilde{\xb}_i,\tilde{\btheta}_{a}\rangle\Big),\quad\Mb_\phi=\begin{bmatrix}
\phi_{1,1}\Ib_{N_1-1}& &\\
&\ddots &\\
& & \phi_{K,1}\Ib_{N_K}
\end{bmatrix} \\
    & \big(\Eb_\sfc\big)_{i,a}=\frac{1}{\sqrt{d}}\Big[\Phi_{\sfc,\bot}\Big(\frac{1}{\sqrt{d}}\langle\tilde{\xb}_i,\tilde{\btheta}_{a}\rangle+\frac{1}{\sqrt{d}}u_i\eta_{a}\Big)-\Phi_{\sfc,\bot}\Big(\frac{1}{\sqrt{d}}\langle\tilde{\xb}_i,\tilde{\btheta}_{a}\rangle\Big)\Big].
\end{align*}
Combined \eqref{eq:B_decomp11-m}, \eqref{eq:B_decomp22-m} and \eqref{eq:B_decomp12-m}, we decompose $\overline{\Bb}$ into
\begin{align*}
\overline{\Bb}=\tilde{\Bb}+\bDelta+\Eb\in\RR^{(\PO-1)\times(\PO-1)},
\end{align*}
where
\begin{align*}
\tilde{\Bb}=&\begin{bmatrix}
q_2\tilde{\Mb}_*\tilde{\Mb}_*+\frac{q_4}{d}\tilde{\Mb}_1\tilde{\bTheta}\tilde{\bTheta}^\T \tilde{\Mb}_1&\tilde{\Zb}^\T \\
\tilde{\Zb}&q_3\Ib_n+q_5\tilde{\Hb}
\end{bmatrix}\\
=&\begin{bmatrix}
q_2\mu_{1,2}^2\Ib_{N_1}+q_4\mu_{1,1}^2\frac{\tilde{\bTheta}_1\tilde{\bTheta}_1^\T }{d}&\cdots&q_4\mu_{1,1}\mu_{K,1}\frac{\tilde{\bTheta}_1\tilde{\bTheta}_K^\T }{d}& \tilde{\Zb}_1^\T \\\vdots&\ddots&\vdots&\vdots\\
q_4\mu_{K,1}\mu_{1,1}\frac{\tilde{\bTheta}_K\tilde{\bTheta}_1^\T }{d}&\cdots&q_2\mu_{K,2}^2\Ib_{N_K}+q_4\mu_{K,1}^2\frac{\tilde{\bTheta}_K\tilde{\bTheta}_K^\T }{d}& \tilde{\Zb}_K^\T  \\
\tilde{\Zb}_1&\cdots& \tilde{\Zb}_K&q_3\Ib_n+q_5\tilde{\Hb}
\end{bmatrix},\\
\bDelta=&\begin{bmatrix}
\frac{q_4}{d}\tilde{\Mb}_1\bmeta\bmeta^\T \tilde{\Mb}_1&\frac1d\Mb_{\phi}\bmeta\ub^\T \\
\frac1d\ub\bmeta \Mb_{\phi}&\frac{q_5}{d}\ub\ub^\T 
\end{bmatrix}\\
=&\begin{bmatrix}
\frac{q_4\mu_{1,1}^2}{d}\bmeta_1\bmeta_1^\T&\cdots&\frac{q_4\mu_{1,1}\mu_{K,1}}{d}\bmeta_1\bmeta_K^\T & \frac{\phi_{1,1}}{d}\bmeta_1\ub^\T \\\vdots&\ddots&\vdots&\vdots\\
\frac{q_4\mu_{K,1}\mu_{1,1}}{d}\bmeta_K\bmeta_1^\T &\cdots&\frac{q_4\mu_{K,1}^2}{d}\bmeta_{K}\bmeta_{K}^\T& \frac{\phi_{K,1}}{d}\bmeta_K\ub^\T  \\
\frac{\phi_{1,1}}{d}\ub\bmeta_1^\T&\cdots& \frac{\phi_{K,1}}{d}\ub\bmeta_K^\T&\frac{q_5}{d}\ub\ub^\T
\end{bmatrix},
\quad
\Eb=\begin{bmatrix}
\0& \cdots &\0 &\Eb_1^\T \\
\vdots&\vdots &\vdots &\vdots \\
\0& \cdots &\0 &\Eb_K^\T\\
\Eb_1&\cdots&\Eb_K&\0
\end{bmatrix}.
\end{align*}
Clearly, by the definition of $\tilde{\Bb}$, 
the Stieltjes transform corresponding to $\tilde{\Bb}$ shares the same asymptotics as the Stieltjes transform corresponding to $\overline{\Ab}$.

\smallskip\noindent\textbf{Step 2.} Define $w_2=\left(-\xi+q_2\mu_{1,2}^2+q_4\mu_{1,1}^2-\overline{\Ab}_{\cdot,N_1}^\T (\tilde{\Bb}+\bDelta-\xi\Ib_{ \PO-1})^{-1}\overline{\Ab}_{\cdot,N_1}\right)^{-1}$. Similar to the argument in Section~\ref{subsec:resolvent}, we have $\overline{m}_{1,d}=\psi_1\EE w_2+o_{d}(1)$.


\smallskip\noindent\textbf{Step 3.}  We calculate $\EE w_2$ by mathematical induction. Similar to Section~\ref{subsec:resolvent},  we give some notations which will be used in the following calculation on $\EE w_2$. Let
$$
\vb=\overline{\Ab}_{\cdot,N_1},\quad\vb_i=\overline{\Ab}_{i,N_1}=\left\{\begin{aligned}
&\frac{q_4\mu_{1,1}^2\eta_i}{d}\lVert\overline{\btheta}_{N_1}\rVert_2, &&\text{if  } i\in[1,N_1-1],\\
&\frac{q_4\mu_{1,1}\mu_{\sfc,1}\eta_{i+1}}{d}\lVert\overline{\btheta}_{N_1}\rVert_2, &&\text{if } i\in\cN_\sfc-1,~\sfc\geq2,\\
&\frac{1}{\sqrt{d}}\Phi_{1}\Big(\frac{1}{\sqrt{d}}u_{i-N+1}\lVert\overline{\btheta}_{N_1}\rVert_2\Big), &&\text{if } i\geq N,
\end{aligned}\right.
$$
and 
$$
\Ub=\frac{1}{\sqrt{d}}\begin{bmatrix}
\bmeta_1& & & &\\
&\bmeta_2& & &\\
& & \ddots& & \\
& & &\bmeta_K &\\
& & & &\ub
\end{bmatrix}\in\RR^{(P-1)\times (K+1)},\quad \Mb=\begin{bmatrix}
q_4\mu_{1,1}^2&\cdots &q_4\mu_{1,1}\mu_{K,1}&\phi_{1,1}\\
\vdots&\ddots & \vdots&\vdots \\
q_4\mu_{1,1}\mu_{K,1}&\cdots&q_4\mu_{K,1}^2&\phi_{K,1}\\
\phi_{1,1}&\cdots&\phi_{K,1}&q_5
\end{bmatrix},
$$
respectively. Then after direct calculation, we have the decomposition of $\bDelta$ as
$$
\bDelta=\Ub \Mb\Ub^\T .
$$
Similar to \eqref{eq:w2_woodberry}, we again get that
\begin{equation}
\label{eq:w2_woodberry-m}
    \begin{split}
        w_2=&\Big(-\xi+q_2\mu_{1,2}^2+q_4\mu_{1,1}^2-\vb^\T (\tilde{\Bb}-\xi\Ib_{ \PO-1})^{-1}\vb\\
&+\vb^\T (\tilde{\Bb}-\xi\Ib_{ \PO-1})^{-1}\Ub(\Mb^{-1}+\Ub^\T (\tilde{\Bb}-\xi\Ib_{ \PO-1})^{-1}\Ub)^{-1}\Ub^\T (\tilde{\Bb}-\xi\Ib_{ \PO-1})^{-1}\vb\Big)^{-1}.
    \end{split}
\end{equation}
To continue the calculation, we still require to study the terms $\vb^\T (\tilde{\Bb}-\xi\Ib_{ \PO-1})^{-1}\vb$, $\vb^\T (\tilde{\Bb}-\xi\Ib_{ \PO-1})^{-1}\Ub$ and $\Ub^\T (\tilde{\Bb}-\xi\Ib_{\PO-1})^{-1}\Ub$ in the denominator of \eqref{eq:w2_woodberry-m}. To do so, we note that $\tilde{\Bb}$ is independent on $\vb$ and $\Ub$. Moreover, by the leave-one-out argument, the Stieltjes transform corresponding to $\tilde{\Bb}$ shares the same asymptotics as the Stieltjes transform corresponding to $\overline{\Ab}$.
Notice that $\eta_i$ is independent on $\tilde{\Bb}$ conditioned on $\overline{\btheta}_{N_1}$, and $\tilde{\Bb}$ is independent on $\overline{\btheta}_{N_1}$.
Similar to \eqref{eq:term1}-\eqref{eq:term3}, we have
\begin{align}
\vb^\T (\tilde{\Bb}-\xi\Ib_{ \PO-1})^{-1}\vb=&q_4^2\mu_{1,1}^2\Big(\sum\limits_{\sfc=1}^K\mu_{\sfc,1}^2\overline{m}_{\sfc,d}\Big)+(\phi_{1,1}^2+\phi_{1,2}^2)\overline{m}_{K+1,d}+o_{\PP}(1),\label{eq:term1-m}\\
  \vb^\T (\tilde{\Bb}-\xi\Ib_{ \PO-1})^{-1}\Ub=&\begin{bmatrix}
q_4\mu_{1,1}^2\overline{m}_{1,d}&~\cdots&~q_4\mu_{1,1}\mu_{K,1}\overline{m}_{K,d}&~ ~\phi_{1,1}\overline{m}_{K+1,d}
\end{bmatrix}+o_{\PP}(1),\label{eq:term2-m}\\
\Ub^\T (\tilde{\Bb}-\xi\Ib_{ \PO-1})^{-1}\Ub=&\begin{bmatrix}
\overline{m}_{1,d}& &\\
&\ddots&\\
& &\overline{m}_{K+1,d}
\end{bmatrix}+o_{\PP}(1).\label{eq:term3-m}
\end{align}
Since $|w_2|\leq\Im(\xi)$ is deterministically bounded, by dominated convergence theorem, we have the $L_1$ convergence of $w_2$ by  plugging \eqref{eq:term1-m}-\eqref{eq:term3-m} into \eqref{eq:w2_woodberry-m}.  We have
\begin{align}\label{eq:m1d_calcu1-m}
    \overline{m}_{1,d}=\psi_1\big\{-\xi+q_2\mu_{1,2}^2+q_4\mu_{1,1}^2-\phi_{1,2}^2\overline{m}_{K+1,d}-\bl_K^\T\Mb_K^{-1}\bl_K\big\}^{-1}+o_d(1).
\end{align}
Here we define $\bl_K=\begin{bmatrix}
q_4\mu_{1,1}^2&\quad\cdots&\quad q_4\mu_{1,1}\mu_{K,1}&\quad \phi_{1,1}
\end{bmatrix}^\T\in\RR^{(K+1)\times1}$, and $$\Mb_K=\begin{bmatrix}
q_4\mu_{1,1}^2+\frac{1}{\overline{m}_{1,d}}&\cdots&q_4\mu_{1,1}\mu_{K,1}&\phi_{1,1}\\
\vdots&\ddots&\vdots&\vdots\\
q_4\mu_{1,1}\mu_{K,d}&\cdots&q_4\mu_{K,1}^2+\frac{1}{\overline{m}_{K,d}}&\phi_{K,1}\\
\phi_{1,1}&\cdots&\phi_{K,1}&q_5+\frac{1}{\overline{m}_{K+1,d}}
\end{bmatrix}.$$
Note that $\phi_{\sfc,1}=\mu_{\sfc,1}(1+q_1)$, $\phi_{\sfc,2}=\mu_{\sfc,2}$, we aim to prove the following equality:
\begin{align}\label{eq:target}
    q_4\mu_{1,1}^2-\bl_K^\T\Mb_K^{-1}\bl_K=\frac{\mu_{1,1}^2q_4(1+q_5\overline{m}_{K+1,d})-\mu_{1,1}^2(1+q_1)^2\overline{m}_{K+1,d}}{\Big(1+q_4\sum\limits_{\sfc=1}^K\mu_{\sfc,1}^2\overline{m}_{\sfc,d}\Big)(1+q_5\overline{m}_{K+1,d})-(1+q_1)^2\sum\limits_{\sfc=1}^K\mu_{\sfc,1}^2\overline{m}_{\sfc,d}\overline{m}_{K+1,d}}.
\end{align}
We prove \eqref{eq:target} by mathematical induction. For $K=2$, \eqref{eq:target} holds  from Section~\ref{subsec:resolvent}.  We assume that
\begin{align}\label{eq:K-1}
    q_4\mu_{1,1}^2-\bl_{K-1}^\T\Mb_{K-1}^{-1}\bl_{K-1}=\frac{\mu_{1,1}^2q_4(1+q_5\overline{m}_{K,d})-\mu_{1,1}^2(1+q_1)^2\overline{m}_{K,d}}{\Big(1+q_4\sum\limits_{\sfc=1}^{K-1}\mu_{\sfc,1}^2\overline{m}_{\sfc,d}\Big)(1+q_5\overline{m}_{K,d})-(1+q_1)^2\sum\limits_{\sfc=1}^{K-1}\mu_{\sfc,1}^2\overline{m}_{\sfc,d}\overline{m}_{K,d}}
\end{align}
holds under the case $K-1$. We aim to prove \eqref{eq:target} for general  $K$ under the assumption that \eqref{eq:K-1} holds.
To prove so, define $\bmu_K=[\mu_{1,1},\ldots,\mu_{K,1}]^\T$. The vector $\bl_K$ could be separated into $\bl_K=[q_4\mu_{1,1}\cdot\bmu_K^\T\quad(1+q_1)\mu_{1,1}]^\T$. If we further define
$$
\Vb_0=\begin{bmatrix}
q_4\mu_{1,1}^2+\frac{1}{\overline{m}_{1,d}}&\cdots&q_4\mu_{1,1}\mu_{K,1}\\
\vdots&\ddots&\vdots\\
q_4\mu_{1,1}\mu_{K,1}&\cdots&q_4\mu_{K,1}^2+\frac{1}{\overline{m}_{K,d}}
\end{bmatrix},
$$
the target equation~\eqref{eq:target} could be rewritten as 
\begin{align}\label{eq:Induction1}
q_4\mu_{1,1}^2-\bl_K^\T\Mb_K^{-1}\bl_K=q_4\mu_{1,1}^2-\begin{bmatrix}q_4\mu_{1,1}\cdot\bmu_K\\(1+q_1)\mu_{1,1}\end{bmatrix}^\T\begin{bmatrix}
\Vb_0 &(1+q_1)\bmu_K^\T\\
(1+q_1)\bmu_K&q_5+\frac{1}{\overline{m}_{K+1,d}}
\end{bmatrix}^{-1}\begin{bmatrix}q_4\mu_{1,1}\cdot\bmu_K\\(1+q_1)\mu_{1,1}\end{bmatrix}.
\end{align}
Clearly, the formula in \eqref{eq:Induction1} requires us to investigate $\bmu_K^\T\Vb_0^{-1}\bmu_K$ first. Under the case $K-1$, \eqref{eq:K-1} holds from the induction hypothesis.  Thus if we set $(1+q_1)=q_4\mu_{K,1},~ q_5=q_4\mu_{K,1}^2$, we have $\bl_{K-1}=q_4\mu_{1,1}\cdot\bmu_K$. Plugging $\bl_{K-1}=q_4\mu_{1,1}\cdot\bmu_K$ into \eqref{eq:K-1} we obtain that

\begin{align}
  & ~~\bmu_K^\T\Vb_0^{-1}\bmu_K 
  \label{eq:Induction2}
  \\ 
  =&~~ \frac{1}{q_4^2\mu_{1,1}^2}\Bigg(q_4\mu_{1,1}^2-\frac{\mu_{1,1}^2q_4(1+\mu_{K,1}^2q_4\overline{m}_{K,d})-\mu_{1,1}^2q_4^2\mu_{K,1}^2\overline{m}_{K,d}}{\Big(1+q_4\sum\limits_{\sfc=1}^{K-1}\mu_{\sfc,1}^2\overline{m}_{\sfc,d}\Big)(1+\mu_{K,1}^2q_4\overline{m}_{K,d})-q_4^2\mu_{K,1}^2\sum\limits_{\sfc=1}^{K-1}\mu_{\sfc,1}^2\overline{m}_{\sfc,d}\overline{m}_{K,d}}\Bigg)\nonumber\\
=&~~\frac{1}{q_4\mu_{1,1}^2}\Bigg(\mu_{1,1}^2-\frac{\mu_{1,1}^2}{1+q_4\sum\limits_{\sfc=1}^{K}\mu_{\sfc,1}^2\overline{m}_{\sfc,d}}\Bigg)=\frac{\sum\limits_{\sfc=1}^{K}\mu_{\sfc,1}^2\overline{m}_{\sfc,d}}{1+q_4\sum\limits_{\sfc=1}^{K}\mu_{\sfc,1}^2\overline{m}_{\sfc,d}}. \nonumber 
\end{align}
Therefore, for the case $K$, we have
\begin{align*}
    q_4\mu_{1,1}^2-\bl_K^\T\Mb_K^{-1}\bl_K
    =&q_4\mu_{1,1}^2-\begin{bmatrix}q_4\mu_{1,1}\cdot\bmu_K\\(1+q_1)\mu_{1,1}\end{bmatrix}^\T\begin{bmatrix}
\Vb_0 &(1+q_1)\bmu_K^\T\\
(1+q_1)\bmu_K&q_5+\frac{1}{\overline{m}_{K+1,d}}
\end{bmatrix}^{-1}\begin{bmatrix}q_4\mu_{1,1}\cdot\bmu_K\\(1+q_1)\mu_{1,1}\end{bmatrix}\\
=&\frac{\mu_{1,1}^2q_4(1+q_5\overline{m}_{K,d})-\mu_{1,1}^2(1+q_1)^2\overline{m}_{K,d}}{\Big(1+q_4\sum\limits_{\sfc=1}^{K-1}\mu_{\sfc,1}^2\overline{m}_{\sfc,d}\Big)(1+q_5\overline{m}_{K,d})-(1+q_1)^2\sum\limits_{\sfc=1}^{K-1}\mu_{\sfc,1}^2\overline{m}_{\sfc,d}\overline{m}_{K,d}}.
\end{align*}
Here, the first equality directly comes from \eqref{eq:Induction1}, and the second equality comes from  Schur complement and \eqref{eq:Induction2} after direct calculation. We completed the mathematical induction for the general case $K$ and equation~\eqref{eq:target} is proved. Then we have
\begin{align*}
      &\overline{m}_{1,d}=\psi_1\bigg\{-\xi+q_2\mu_{1,2}^2-\mu_{1,2}^2\overline{m}_{K+1,d}+\frac{H_{1,d}}{H_{D,d}}\bigg\}^{-1}+o_d(1),
\end{align*}
where
$$\begin{aligned}
H_{1,d}=&\mu_{\sfc,1}^2q_4(1+q_5\overline{m}_{K+1,d})-\mu_{\sfc,1}^2(1+q_1)^2\overline{m}_{K+1,d},\quad\sfc=1,\ldots,K,\\
H_{D,d}=&\Big(1+q_4\sum\limits_{\sfc=1}^K\mu_{\sfc,1}^2\overline{m}_{\sfc,d}\Big)(1+q_5\overline{m}_{K+1,d})-(1+q_1)^2\sum\limits_{\sfc=1}^K\mu_{\sfc,1}^2\overline{m}_{\sfc,d}\overline{m}_{K+1,d}.
\end{aligned}$$
After similar argument, we conclude that 
\begin{align*}
      &\overline{m}_{\sfc,d}=\psi_\sfc\bigg\{-\xi+s_\sfc\mu_{\sfc,2}^2-\mu_{\sfc,2}^2\overline{m}_{K+1,d}+\frac{H_{\sfc,d}}{H_{D,d}}\bigg\}^{-1}+o_d(1),\quad\sfc=1,\ldots,{K+1},\\
&\overline{m}_{K+1,d}=\psi_{K+1}\bigg\{-\xi+q_3-\sum\limits_{\sfc=1}^K\mu_{\sfc,2}^2\overline{m}_{\sfc,d}+\frac{H_{{K+1},d}}{H_{D,d}}\bigg\}^{-1}+o_d(1),
\end{align*}
where
$$\begin{aligned}
H_{\sfc,d}=&\mu_{\sfc,1}^2q_4(1+q_5\overline{m}_{K+1,d})-\mu_{\sfc,1}^2(1+q_1)^2\overline{m}_{K+1,d},\quad\sfc=1,\ldots,K,\\
H_{{K+1},d}=&q_5\Big(1+q_4\sum\limits_{\sfc=1}^K\mu_{\sfc,1}^2\overline{m}_{\sfc,d}\Big)-(1+q_1)^2\sum\limits_{\sfc=1}^K\mu_{\sfc,1}^2\overline{m}_{\sfc,d}.
\end{aligned}$$
We get that each element of $\overline{\mb}_{d}(\xi)-\sfFb(\overline{\mb}_{d}(\xi))$ is asymptotically zero. Therefore $\|\overline{\mb}_{d}(\xi)-\sfFb(\overline{\mb}_{d}(\xi))\|_2=o_d(1)$. The remaining arguments are similar to those in Section~\ref{sec:appendixpropimplicit1} and the details are skipped.  Wrapping all together, we complete the proof of Proposition~\ref{prop:implicit1-m}.
\end{proof}
\end{document}